\documentclass{article} 
\pdfoutput=1
\usepackage{nips13submit_e,times}
\usepackage{url}
\usepackage{graphicx}
\usepackage{amssymb,amsmath,mathbbol}
\usepackage[margin=2.5cm]{geometry}
\usepackage{rotating}
\usepackage{multirow}
\usepackage{bbding}
\usepackage{color}

\newtheorem{theorem}{Theorem}[section]
\newcommand{\ignore}[1]{}

\definecolor{myblue}{rgb}{0.2,0,0.5} 
\definecolor{mypurple}{rgb}{1.0,0,0.4} 

\title{The HIM glocal metric and kernel for \\ network comparison and classification}

\author{
G. Jurman\thanks{Corresponding author} ${\;}^1$, 
R. Visintainer${}^1$,
M. Filosi${}^{1,2}$,
S. Riccadonna${}^3$,
C. Furlanello${}^1$
\\
\begin{tabular}{c}
${}^1$ Fondazione Bruno Kessler, Trento, Italy\\
${}^2$ CIBIO, University of Trento, Italy\\
${}^3$ Research and Innovation Centre, Fondazione Edmund Mach, San Michele all'Adige, Italy 
\end{tabular}
\\
\begin{tabular}{c}
\texttt{\{jurman,visintainer,filosi,furlan\}@fbk.eu}\\
\texttt{samantha.riccadonna@fmach.it}
\end{tabular}
}

\nipsfinalcopy 

\begin{document}
\maketitle
\begin{abstract}
Due to the ever rising importance of the network paradigm across several areas of science, comparing and classifying graphs represent essential steps in the networks analysis of complex systems.
Both tasks have been recently tackled via quite different strategies, even tailored \textit{ad-hoc} for the investigated problem.
Here we deal with both operations by introducing the Hamming-Ipsen-Mikhailov (HIM) distance, a novel metric to quantitatively measure the difference between two graphs sharing the same vertices.
The new measure combines the local Hamming distance and the global spectral Ipsen-Mikhailov distance so to overcome the drawbacks affecting the two components separately.
Building then the HIM kernel function derived from the HIM distance it is possible to move from network comparison to network classification via the Support Vector Machine (SVM) algorithm.
Applications of HIM distance and HIM kernel in computational biology and social networks science demonstrate the effectiveness of the proposed functions as a general purpose solution.
\end{abstract}
\section{Introduction}
\label{sec:intro}
The arising prevalence of the network paradigm \cite{barabasi12network} as the elective model for complex systems analysis in different workfields has strongly contributed in stimulating graph theoretical techniques in the recent scientific literature.
Methods based on graph properties have spread through the static and dynamic analysis of different economical, chemical and biological system, computer networking, social networks and neuroscience.
As a relevant example, it is worthwhile mentioning the rapid diffusion, in computational biology, of the differential network analysis \cite{sharan06modeling,ideker12differential,yoon12comparative,csermely13structure,chuang07network,yang13network,pavlopoulos11using,barla12machine,barla13machine}.
In particular, two key tasks constitute the backbone of most of the aforementioned analysis techinques, namely network comparison and network classification, and they both rely on the basic idea of measuring the similarity between two graphs.

Network comparison consists in the quantification of the difference between two homogeneous objects in some network space, while the aim of network classification is to predictively discriminate graphs belonging to different classes, for instance by means of machine learning algorithms.
Network comparison has its roots in the quantitative description of main properties of a graph (\textit{e.g.}, degree distribution), which can be encoded into a feature vector \cite{xiao08structure}, thus providing a convenient representation for classification tasks (see for instance \cite{dehmer13discrimination} for a very recent approach).
As a major alternative strategy, one can adopt a direct comparison method stemming from the graph isomorphism problem, by defining a suitable similarity measure on the topology of the underlying (possibly directed and/or weighted) graphs.
This line of study dates back to the 70’s with the theory of graph distances, regarding both inter- and intra-graphs metrics \cite{entringer76distance}.
Since then, a wide range of similarity measures has been defined, based on very different graph indicators. 
To mention some of the most important metrics, we list the family of edit distances, evaluating the minimum cost of transformation of one graph into another by means of the usual edit operations (insertion and deletion of links), the family of common network subgraphs, looking for shared structures between the graphs and the family of spectral measures, relying on functions of the eigenvalues of one of the graph connectivity matrices.
Similarly, graph classification can be tackled by a number of different techniques, for instance nearest neighbours on Euclidean distance of the features' vectors of the graphs \cite{zhu11classifying,aliakbary13learning,chen12discovery}, or Support Vector Machine with the graph Laplacian as a regularization term \cite{chen11identifying}, or via different subgraph-based lerning algorithms \cite{thorat13survey}. 
However, in general the most efficient techniques use a kernel machine, where the kernel itself corresponds to a scalar product (and hence a distance) in a suitable Hilbert space \cite{mahe04extension,gaertner06short,gaertner07kernel,borgwardt07graph,ketkar09empirical,vishwanathan10graph,tsuda10graph,vert05supervised,vert03graph}.
For more recent advances, we cite the Weisfeiler-Lehman graph kernel \cite{shervashidze11weisfeiler}, and its use in neuroimaging classification for discriminating mild cognitive impairment from Alzheimer's disease \cite{jie13integration}.
This last citation stands as an example of the increasing interest for these techniques recently appearing in neurosciences \cite{richiardi13machine,su13discriminative}.
 
In the present work we propose a novel solution to both the comparison and the classification tasks by introducing the novel HIM metric for comparing graphs (even directed and weighted) and a graph kernel induced by the HIM measure.
The HIM distance is defined as the one-parameter family of product metrics linearly combining -- by a non-negative real factor $\xi$ -- the normalized Hamming distance H \cite{dougherty10validation,tun06metabolic,iwayama12characterizing,morris08specification} and the normalized Ipsen-Mikhailov distance IM \cite{ipsen02evolutionary}; the product metric is normalized by the factor $\sqrt{1+\xi}$ to set its upper bound to 1.
In absence of a gold standard driving the search for the optimal weight ratio, we decided for an equal contribution from the two components $\xi=1$ as the most natural choice.
The Hamming distance is the simplest member of the family of edit distances, evaluating the occurrence of matching links in the compared networks: by definition, it is a local measure of dissimilarity between graphs, because it only focusses on the links as independent entities, disregarding the overall structure.
On the other hand, the spectral distances are global measures, evaluating the differences between the whole network structures: however, they cannot discriminate between isospectral non-identical graphs: for a recent spectral approach, see \cite{rajendran13analysis}.
In the comparative review \cite{jurman11introduction}, the properties of the existing graph spectral distances were studied, and the Ipsen-Mikhailov metric emerged as the more reliable and stable.
The combination of the two components within a single metric allows overcoming their drawbacks and obtaining a measure which is simultaneously global and local.
Moreover, the imposed normalization limits the values of the HIM distance between zero (reached only by comparing identical networks) and one (attained when comparing a clique and the empty graph), regardless of the number of vertices.
Finally, the HIM distance can also be applied to multilayer networks \cite{kivela13multilayer,dedomenico13mathematical}, since a rigorous definition of their Laplacian has just been proposed \cite{sole-ribalta13spectral,sanchezgarcia13dimensionality}. 
By a Gaussian-like map \cite{cortes03positive}, the HIM distance generates the HIM kernel. 
Plugging the HIM kernel \cite{shawe-taylor04kernel} into a Support Vector Machine gives us a classification algorithm based on the HIM distance, to be used as is or together with other graph kernels in a Multi-Kernel Learning framework to increase the classification performance and to enhance the interpretability of the results \cite{kloft11lp}.
Note that, although positive definiteness does not hold globally for the HIM kernel, this property can be guaranteed on the given training data, thus leading to positive definite matrices suitable for the convergence of the SVM optimizer.
 
To conclude with, we present some applications of the HIM distance and the HIM kernel to some real datasets belonging to different areas of science. 
These examples support the positive impact of the HIM suite as general analysis tool whenever it is required to extract information from the quantitative evaluation of the difference among diverse instances of a complex system.

We also provide for analysis the R \cite{R2013} package \textit{nettools} including functions to compute the HIM distance.
The package is provided as a working beta version and it is accessible on GitHub at \url{https://github.com/filosi/nettools.git}.
To reduce computing time, the software can be used on multicore workstations and on high performance computing (HPC) clusters.  
\section{The HIM family of distances}
\label{sec:him}
\subsection{Notations}
\label{ssec:notations}
Let $\mathcal{N}_1$ and $\mathcal{N}_2$ be two simple networks on $N$ nodes, described by the corresponding adjacency matrices $A^{(1)}$ and $A^{(2)}$, with $a^{(1)}_{ij}, a^{(2)}_{ij}\in\mathcal{F}$, where $\mathcal{F}=\mathbb{F}_2=\{0,1\}$ for unweighted graphs and $\mathcal{F}=[0,1]$ for weighted networks. 
Let then $\mathbb{I}_N$ be the $N\times N$ identity matrix $\mathbb{I}_N = \left( \begin{smallmatrix} 1&0&\cdots & 0 \\ 0&1&\cdots&0 \\ &\cdots \\ 0&0&\cdots &1  \end{smallmatrix} \right)$, let $\mathbb{1}_N$ be the $N\times N$ unitary matrix with all entries equal to one and let $\mathbb{0}_N$ be the $N\times N$  null matrix with all entries equal to zero. 
Denote then by $\mathcal{E}_N$ the empty network with $N$ nodes and no links (with adjacency matrix $\mathbb{0}_N$) and by $\mathcal{F}_N$ the clique (undirected full network) with $N$ nodes and all possible $N(N-1)$ links, whose adjacency matrix is $\mathbb{1}_N-\mathbb{I}_N$.
For an undirected network, its adjacency matrix is symmetric.
For a directed network $\mathcal{N}^\uparrow$, following the convention in \cite{liu11controllability}, a link ${i}\rightarrow{j}$ is represented by setting $a_{ji}=1$ in the corresponding adjacency matrix $A_{\mathcal{N}^\uparrow}$, which thus is, in general, not symmetric.

For instance, the matrix $A_{\mathcal{N}^\uparrow}=\mathbb{1}_N-\mathbb{I}_N$ represents the full directed network $\mathcal{F}^\uparrow_N$, with all possible $N^2-N$ directed links ${i}\rightarrow{j}$.
\subsection{The Hamming distance}
\label{ssec:hamming}
The Hamming distance is one of the most common dissimilarity measures in coding and string theory, recently used also for (biological) network comparison \cite{dougherty10validation,tun06metabolic,morris08specification,iwayama12characterizing}.
Since the Hamming measure basically evaluates the presence/absence of matching links on the two networks being compared, it has a simple expression in terms of the neworks' adjacency matrices.
This is not the case for many other members of the edit distance family, whose computation is known to be a NP-hard task.
The definition of the normalized Hamming distance H is in fact the following:
\begin{equation}
\label{eq:hamming}
\textrm{H}(\mathcal{N}_1,\mathcal{N}_2) = 
\frac{\textrm{Hamming}(\mathcal{N}_1,\mathcal{N}_2)}{\textrm{Hamming}(\mathcal{E}_N,\mathcal{F}_N)} = 
\frac{\textrm{Hamming}(\mathcal{N}_1,\mathcal{N}_2)}{N(N-1)} = 
\frac{1}{N(N-1)}\sum_{1\leq i\not = j\leq N} \vert A^{(1)}_{ij} - A^{(2)}_{ij} \vert\ ,
\end{equation}
where the normalization factor $N(N-1)$ bound the range of the function H in the interval $[0,1]$. 
The lower bound $0$ is attained only for identical networks $A^{(1)}=A^{(2)}$, while the upper bound $1$ is reached whenever the two networks are complementary $A^{(1)}+A^{(2)}=\mathbb{1}_N-\mathbb{I}_N=\left( \begin{smallmatrix} 0&1&\cdots & 1 \\ 1&0&\cdots&1 \\ &\cdots \\ 1&1&\cdots &0  \end{smallmatrix} \right)$.
When $\mathcal{N}_1$ and $\mathcal{N}_2$ are unweighted networks, $\textrm{H}(\mathcal{N}_1,\mathcal{N}_2)$ is just the fraction of different matching links over the total number $N(N-1)$ of possible links between the two graphs.
\subsection{The Ipsen-Mikhailov distance}
\label{ssec:ipsen}
\begin{figure}[!t]
\begin{center}
\begin{tabular}{ccccc}
\raisebox{-1.8cm}{\includegraphics[width=0.17\textwidth]{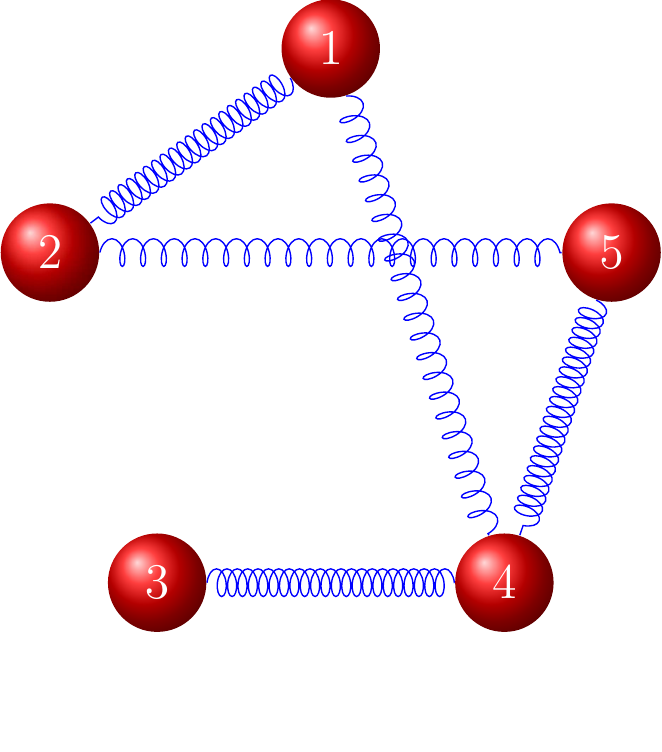}} &
\phantom{aaaaa} &
$
\begin{pmatrix}
0 & 1 & 0 & \frac{1}{2} & 0 \\
1 & 1 & 0 & 0 & \frac{1}{2}  \\
0 & 0 & 0 & 1  & 0 \\
\frac{1}{2} & 0 & 1 & 0 & 1 \\
0 & \frac{1}{2} & 0 & 1 & 0 
\end{pmatrix}
$ &
\phantom{aaaaa} &
\raisebox{-1.8cm}{\includegraphics[width=0.35\textwidth]{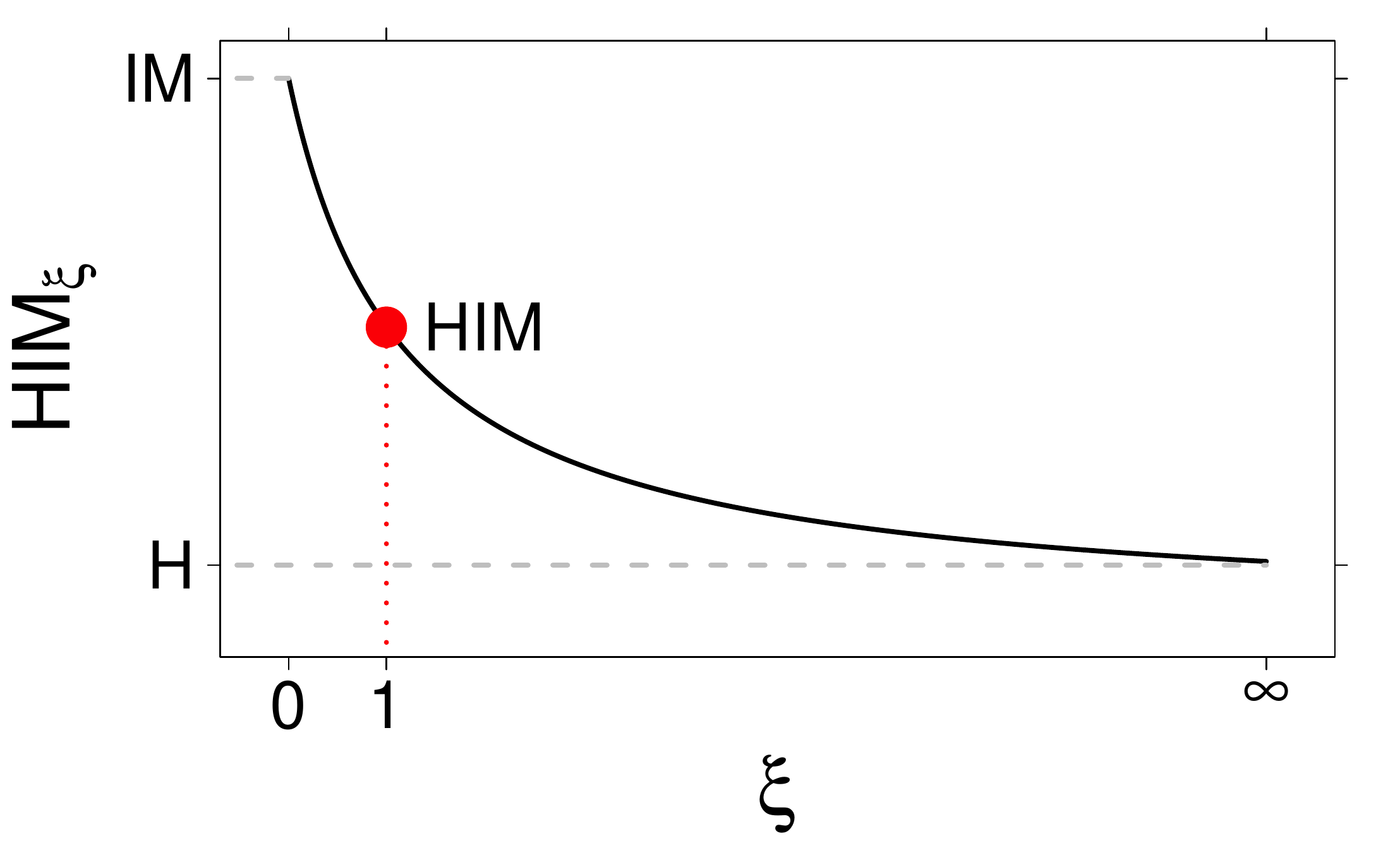}} \\
(a) & & (b) & & (c)
\end{tabular}
\end{center}
\caption{A five nodes network as a oscillatory system (a) and the corresponding adjacency matrix (b), with two different edge weights 1 and $\frac{1}{2}$, represented by different springs. In panel (c), the product metric $\textrm{HIM}_{\xi}=\frac{1}{\sqrt{1+\xi}} \sqrt{ \textrm{H}^2 + \xi \textrm{IM}^2}$ as a function of $\xi\in[0,\infty)$.}
\label{fig:springs}
\end{figure}

Originally introduced in \cite{ipsen02evolutionary} as a tool for network reconstruction from its Laplacian spectrum, the definition of the Ipsen-Mikhailov $\textrm{IM}$ metric follows the dynamical interpretation of an $N$ nodes network as an $N$ molecules system connected by identical elastic strings as in Fig.~\ref{fig:springs}(a-b), where the pattern of connections is defined by the adjacency matrix $A$ of the corresponding network.
The dynamical system is described by the set of $N$ differential equations 
\begin{equation}
\label{eq:ipsen_model}
\ddot{x}_i+\sum_{j=1}^N A_{ij}(x_i-x_j)=0\quad\textrm{for\ }i=0,\cdots,N-1\ .
\end{equation}
We recall that the Laplacian matrix $L$ of an undirected network is defined as the difference between the degree $D$ and the adjacency $A$ matrices $L=D-A$, where $D$ is the diagonal matrix with vertex degrees as entries.
$L$ is positive semidefinite and singular \cite{chung97spectral,spielman09spectral,tonjes09perturbation,atay06network}, so its eigenvalues are $0 = \lambda_0 \leq \lambda_1\leq \cdots\leq \lambda_{N-1}$.
The vibrational frequencies $\omega_i$ for the network model in Eq.~\ref{eq:ipsen_model} are given by the square root of the eigenvalues of the Laplacian matrix of the network: $\lambda_i = \omega^2_i$, with $\lambda_0=\omega_0=0$. 
In \cite{chung97spectral}, the Laplacian spectrum is called the vibrational spectrum.
Estimates (actual and asymptotic) of the eigenvalues distribution are available for complex networks \cite{rodgers05eigenvalue}, while the relations between the spectral properties and the structure and the dynamics of a network are discussed in \cite{jost02evolving, jost07dynamical, almendral07dynamical}.
The spectral density for a graph as the sum of Lorentz distributions is defined as 
\begin{displaymath}
\rho(\omega,\gamma)=K\sum_{i=1}^{N-1} \frac{\gamma}{(\omega-\omega_i)^2+\gamma^2}\ ,
\end{displaymath}
where $\gamma$ is the common width and $K$ is the normalization constant defined by the condition $\displaystyle{\int_0^\infty \rho(\omega,\gamma)\textrm{d}\omega =1}$, and thus
\begin{displaymath}
K = \frac{1}{\gamma\displaystyle{\sum_{i=1}^{N-1} \int_0^\infty \frac{\textrm{d}\omega}{(\omega-\omega_i)^2+\gamma^2} }}\ .
\end{displaymath}
The scale parameter $\gamma$ specifies the half-width at half-maximum, which is equal to half the interquartile range. 
An example of Lorentz distribution for two networks is shown In Fig.~\ref{fig:lorentz}.
Then the spectral distance $\epsilon_\gamma$ between two graphs $\mathcal{N}_1$ and $\mathcal{N}_2$ on $N$ nodes with densities $\rho_{\mathcal{N}_1}(\omega,\gamma)$ and $\rho_{\mathcal{N}_2}(\omega,\gamma)$ can then be defined as 
\begin{displaymath}
\epsilon_\gamma(\mathcal{N}_1,\mathcal{N}_2) = \sqrt{\int_0^\infty \left[\rho_{\mathcal{N}_1}(\omega,\gamma)-\rho_{\mathcal{N}_2}(\omega,\gamma)\right]^2 \textrm{d}\omega}\ .
\end{displaymath}
The highest value of $\epsilon_\gamma$ is reached, for each $N$, when evaluating the distance between $\mathcal{E}_N$ and $\mathcal{F}_N$.
Denoting by $\overline{\gamma}$ the unique solution of 
\begin{equation}
\label{eq:gamma_implicit}
\epsilon_\gamma(\mathcal{E}_N, \mathcal{F}_N) = 1\ , 
\end{equation}
the normalized Ipsen-Mikahilov distance between two undirected (possibly weighted) networks can be defined as
\begin{equation}
\label{eq:epsilon}
\textrm{IM}(\mathcal{N}_1,\mathcal{N}_2)=\epsilon_{\overline\gamma}(\mathcal{N}_1,\mathcal{N}_2) = \sqrt{\int_0^\infty \left[\rho_{\mathcal{N}_1}(\omega,\overline{\gamma})-\rho_{\mathcal{N}_2}(\omega,\overline{\gamma})\right]^2 \textrm{d}\omega}\ ,
\end{equation}
so that $\textrm{IM}$ is bounded between 0 and 1, with upper bound attained only for $\{\mathcal{N}_1,\mathcal{N}_2\}=\{\mathcal{E}_N,\mathcal{F}_N\}$.
A detailed description of the uniqueness of the solution of Eq.~\ref{eq:gamma_implicit} is described in Appendix \ref{sec:appendix}.
Isospectral networks (and thus also isomorphic networks) cannot be distinguished by this class of measures, so this is a distance between classes of isospectral graphs. 
Although the number of isospectral networks is negligible for large number of nodes \cite{haemers04enumeration}, their fraction is relevant for smaller networks.
The case of directed networks is discussed in a later paragraph.
\begin{figure}[!t]
\centering
\includegraphics[width=0.8\textwidth]{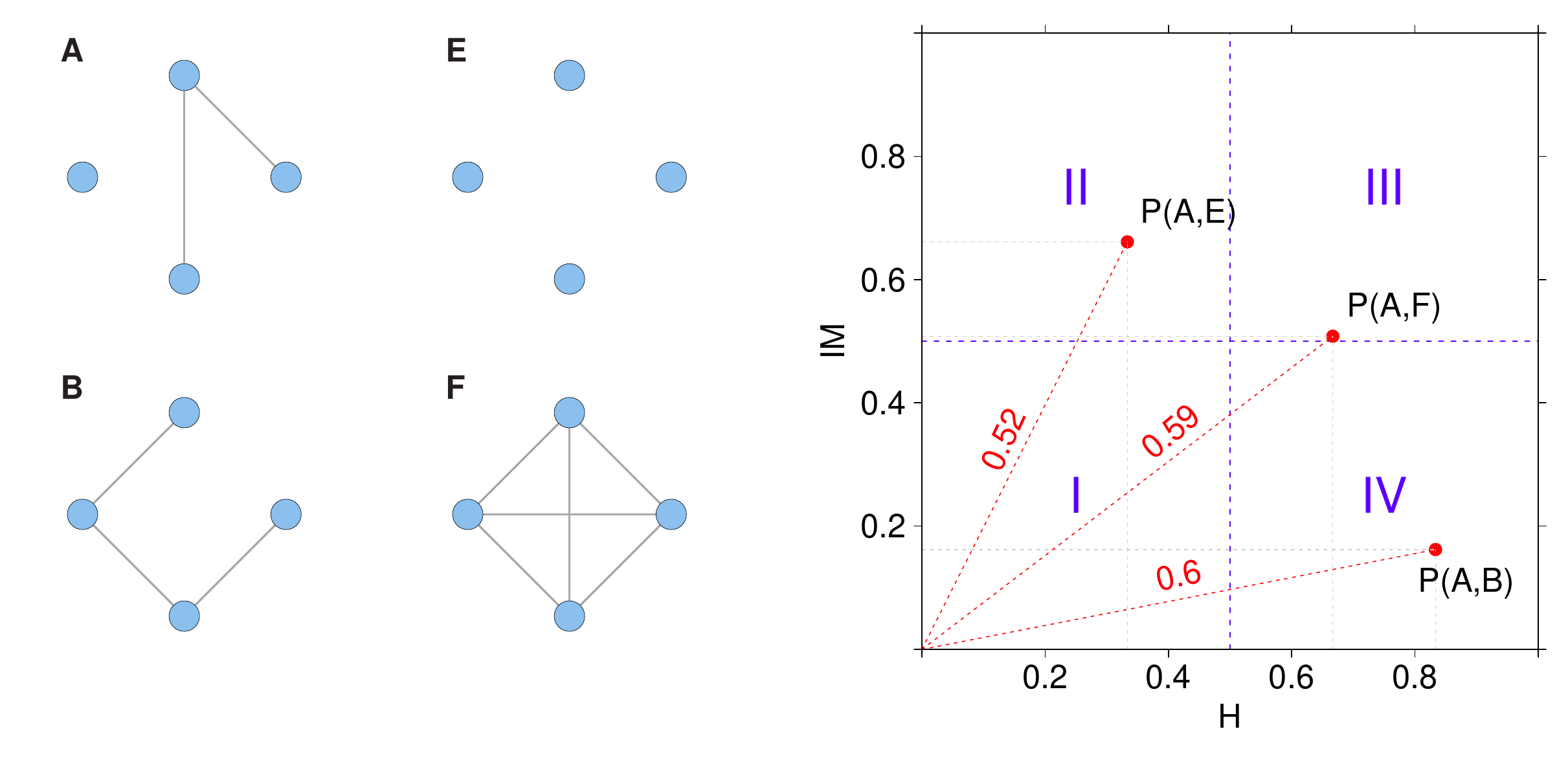}
\caption{Representation of the HIM distance in the Ipsen-Mikhailov (IM axis) and Hamming (H axis) distance space between networks A versus B, E and F, where E is the empty network and F is the clique.}
\label{fig:himspace}
\end{figure}
\subsection{The HIM distance}
\label{ssec:him}
Consider now two copies of the space $\pmb{N}(N)$ of all simple undirected networks on $N$ nodes, and endow the first copy with the Hamming metric H and the second copy with the Ipsen-Mikhailov distance IM. 
Then the two obtained pairs $(\pmb{N}(N),\textrm{H})$ and $(\pmb{N}(N),\textrm{IM})$ are metric spaces.
Define now on their Cartesian product the one-parameter HIM function as the $L_2$ (Euclidean) product metric \cite{deza09encyclopedia} combining H and $\sqrt{\xi}\cdot$ IM, normalized by the factor $\frac{1}{\sqrt{1+\xi}}$, for $\xi\in [0,+\infty)$.
Via the natural correspondence of the same network in the two spaces, the HIM function becomes a distance on $\pmb{N}(N)$:
\begin{equation}
\label{eq:glocal}
\textrm{HIM}_{\xi}(\mathcal{N}_1,\mathcal{N}_2) =  \frac{1}{\sqrt{1+\xi}}  || (\textrm{H}(\mathcal{N}_1,\mathcal{N}_2) , \sqrt{\xi}\cdot\textrm{IM}(\mathcal{N}_1,\mathcal{N}_2)) ||_{2} = \frac{1}{\sqrt{1+\xi}} \sqrt{ \textrm{H}^2(\mathcal{N}_1,\mathcal{N}_2) + \xi\cdot \textrm{IM}^2(\mathcal{N}_1,\mathcal{N}_2) } \ ,
\end{equation}
where in what follows we will omit the subscript $\xi$ when it is equal to one.
Obviously, $\textrm{HIM}_0 = \textrm{H}$ and $\displaystyle{\lim_{\xi\to +\infty} \textrm{HIM}_\xi= \textrm{IM}}$ (see Fig.~\ref{fig:springs}(c)); apart from values of $\xi$ close to the bounds $\{0, +\infty\}$ where the prevalence of one of the factors becomes dominant, the qualitative impact of $\xi$ is minimal in practice when using $\textrm{HIM}_\xi$ as a distance.
In what follows, when no \textit{a priori} hypothesis supports unbalancing the metric towards one of the two components, $\xi=1$ will be assumed. 
However, the impact of $\xi$ is definitely more relevant when $\textrm{HIM}_\xi$ is used to generate a kernel function to be used for classification purposes, as we will show in a later section.
The metric $\textrm{HIM}_\xi(\mathcal{N}_1,\mathcal{N}_2)$ is bounded in the interval $[0,1]$, with lower bound attained for every couple of identical networks, and upper bound attained only on the pair $(\mathcal{E}_N, \mathcal{F}_N)$.
Moreover, all distances $\textrm{HIM}_\xi$ will be nonzero for non-identical isomorphic/isospectral graphs.

Consider now the $[0,1]\times[0,1]$ Hamming/Ipsen-Mikhailov (H/IM) space, where a point $P$ has coordinates $(\textrm{H}(\mathcal{N}_1,\mathcal{N}_2),\textrm{IM}(\mathcal{N}_1,\mathcal{N}_2))$, and the distance of $P$ from the origin is $\sqrt{2}\cdot\textrm{HIM}(\mathcal{N}_1,\mathcal{N}_2)$.
If we (roughly) split the Hamming/Ipsen-Mikhailov space into four main zones I,II,III,IV as in Fig.~\ref{fig:himspace}, two networks whose distances correspond to a point in zone I are quite close both in terms of matching links and of structure, while those falling in the zone III are very different with respect to both measures. 
Networks corresponding to a point in zone II have many common links, but their structure is rather different (for instance, they have a different number of connected components), while a point in zone IV indicates two networks with few common links, but with similar structure (\textit{e.g.}, isospectral non-identical graphs).
In Fig.~\ref{fig:himspace} we show some examples of points in the Hamming/Ipsen-Mikhailov space.
\subsection{The directed network case}
\label{ssec:directed}
\begin{figure}[!t]
\begin{center}
\begin{tabular}{cccc}
\includegraphics[width=0.2\textwidth]{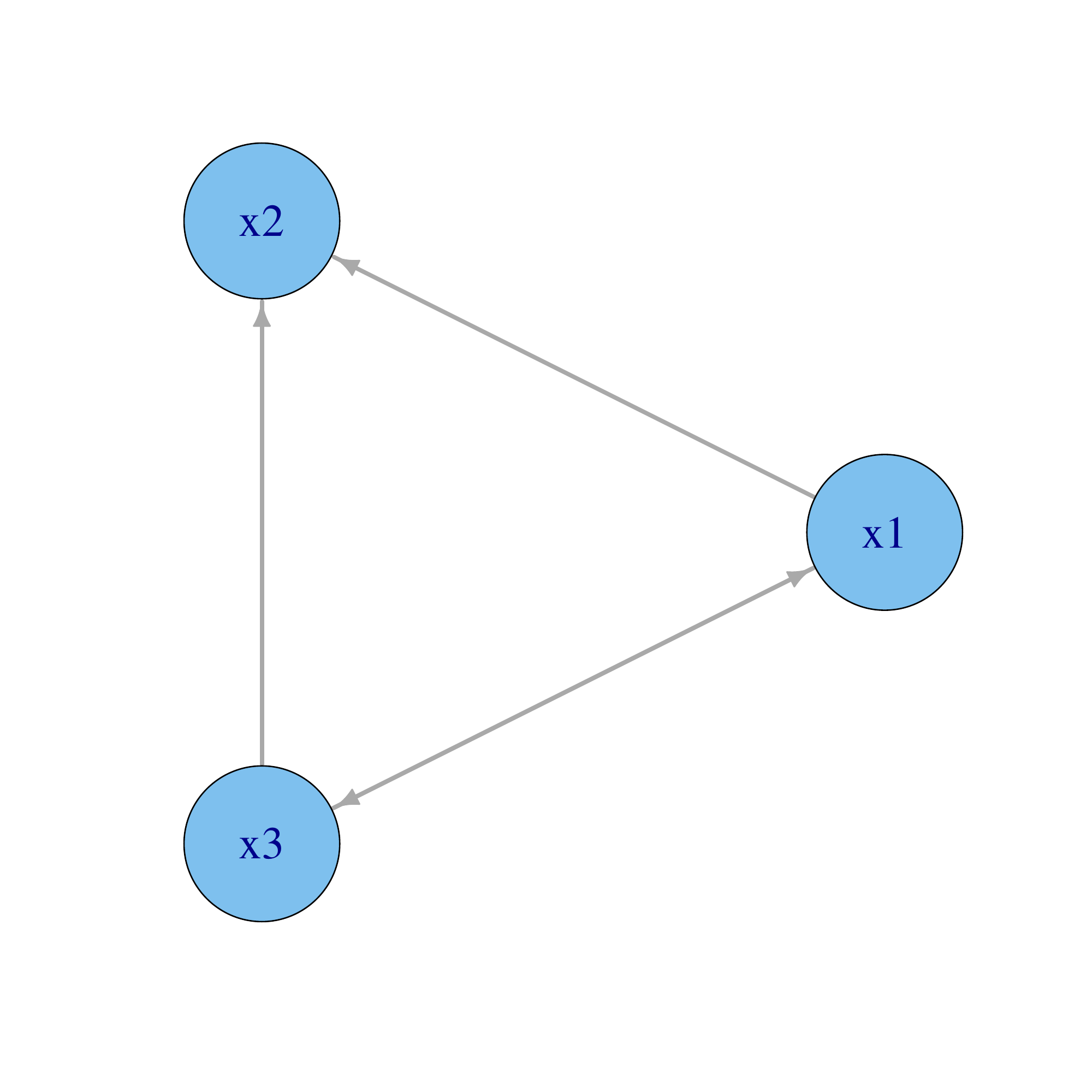} & 
\raisebox{1.5cm}{$\begin{pmatrix} 0&0&1\\ 1&0&1\\ 1&0&0 \end{pmatrix}$} &
\includegraphics[width=0.2\textwidth]{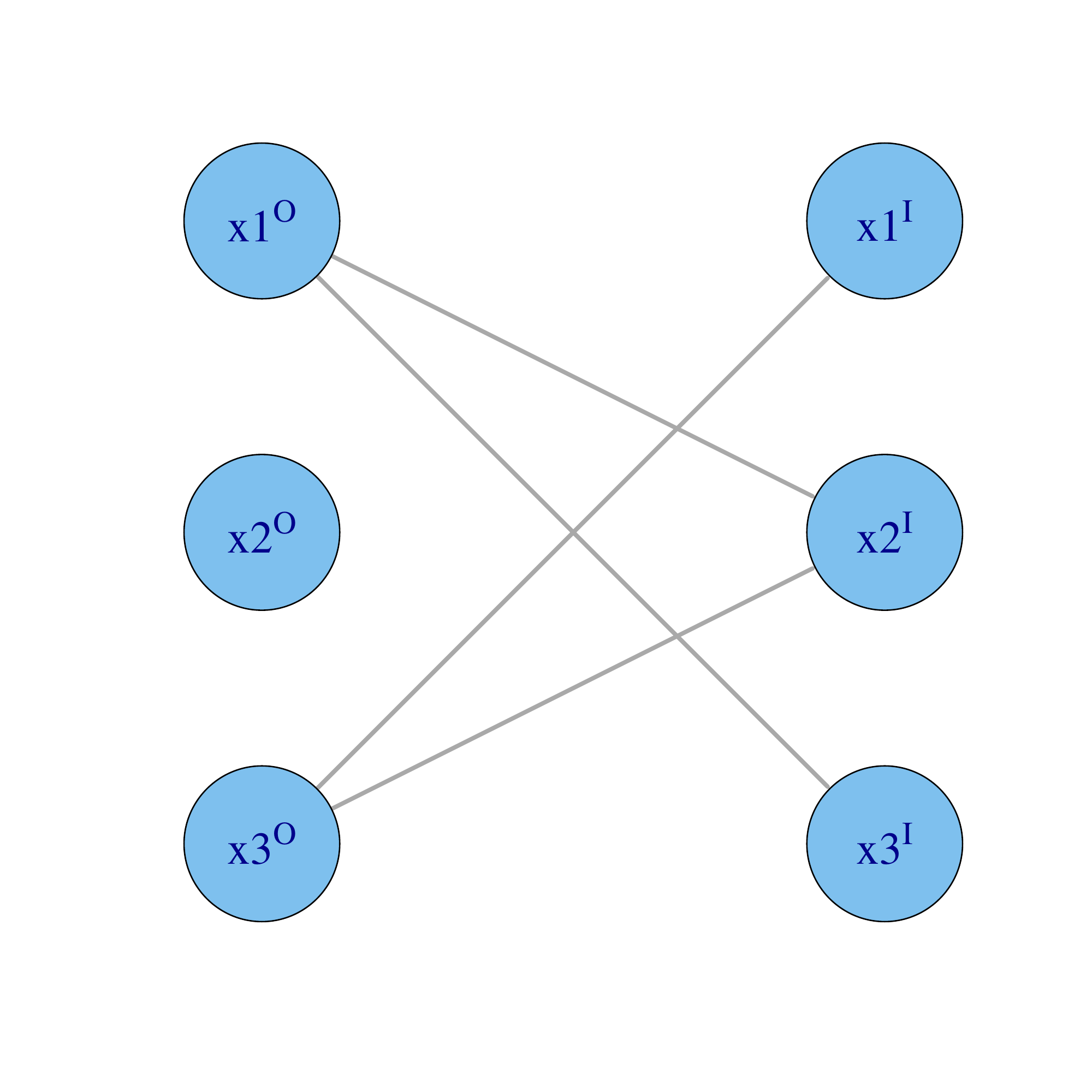} & 
\raisebox{1.5cm}{
$\begin{pmatrix} 
0&0&0&0&1&1\\
0&0&0&0&0&0\\
0&0&0&1&1&0\\
0&0&1&0&0&0\\
1&0&1&0&0&0\\
1&0&0&0&0&0
\end{pmatrix}$
}
\\
\\
\\
\multicolumn{2}{c}{$D^\uparrow$} & \multicolumn{2}{c}{$\hat{D}^\uparrow$} \\
\end{tabular}
\caption{A directed network $D^\uparrow$ on three nodes and the equivalent undirected network $\hat{D}^\uparrow$ on six nodes, together with their adjacency matrices.}
\label{fig:dir2undir}
\end{center}
\end{figure}

In this situation, the connectivity matrices are not symmetric, thus the Laplacian spectrum lies in $\mathbb{C}$. 
Hence, computing the Ipsen-Mikhailov distance would require extending the Lorentzian distribution to the complex plane.
A simpler solution can be obtained by transforming the directed network $D^\uparrow$ into an undirected (bipartite) one $\hat{D}^\uparrow$, as in \cite{liu11controllability}. 
For each node $x_i$ in $D^\uparrow$, the graph $\hat{D}^\uparrow$ has two nodes $x_i^I$ and $x_i^O$ (where I and O stand for In and Out respectively) and for each directed link $x_i\longrightarrow x_j$ in $D^\uparrow$ there is a link $x_i^O -  x_j^I$ in $\hat{D}^\uparrow$.
If the adjacency matrix for $D^\uparrow$ is $A_{D^\uparrow}$, the corresponding matrix for $\hat{D}^\uparrow$ is $A_{\hat{D}^\uparrow}=\left( \begin{smallmatrix}0 &  A^T_{D^\uparrow} \\ A_{D^\uparrow} & 0\end{smallmatrix}\right)$, with respect to the node ordering $x_1^O, x_2^O, \ldots x_n^O, x_1^I, \ldots, x_n^I$. 
An example of the above transformation is shown in Fig.~\ref{fig:dir2undir}.
Thus it is possible to define $\textrm{HIM}(\mathcal{N}^\uparrow_1,\mathcal{N}^\uparrow_2)$ as $\textrm{HIM}(\hat{\mathcal{N}}^\uparrow_1,\hat{\mathcal{N}}^\uparrow_2)$ after substituing the normalizing factors $\overline{\eta}$ and $\overline{\gamma}$ with the corresponding $\overline{\eta}^\uparrow$ and $\overline{\gamma}^\uparrow$ derived by imposing the conditions $\textrm{Hamming}(\hat{\mathcal{E}}_N,\hat{\mathcal{F}}_N)/\overline{\eta}^\uparrow=1$ and $\epsilon_{\overline{\gamma}^\uparrow}(\hat{\mathcal{E}}_N,\hat{\mathcal{F}}_N)=1$, so that $\textrm{HIM}(\hat{\mathcal{E}}_N,\hat{\mathcal{F}}_N)=1$ by using Eq.~(\ref{eq:glocal}).
It is immediate to compute $\bar{\eta}^\uparrow = 2N(N-1)$, while $\bar{\gamma}^\uparrow$ can be numerically computed as for $\bar{\gamma}$: details are given in Appendix \ref{sec:appendixb}.
\begin{figure}[!t]
\begin{center}
\begin{tabular}{cccc}
$
A^{I_1} = \left(\begin{smallmatrix}
0&1&0&0&1&0&0&1\\
1&0&0&0&0&0&1&1\\
0&0&0&0&0&1&1&0\\
0&0&0&0&1&1&0&0\\
1&0&0&1&0&0&0&0\\
0&0&1&1&0&0&0&0\\
0&1&1&0&0&0&0&1\\
1&1&0&0&0&0&1&0\\
\end{smallmatrix}
\right)
$ &
\raisebox{-1.5cm}{\includegraphics[width=0.2\textwidth]{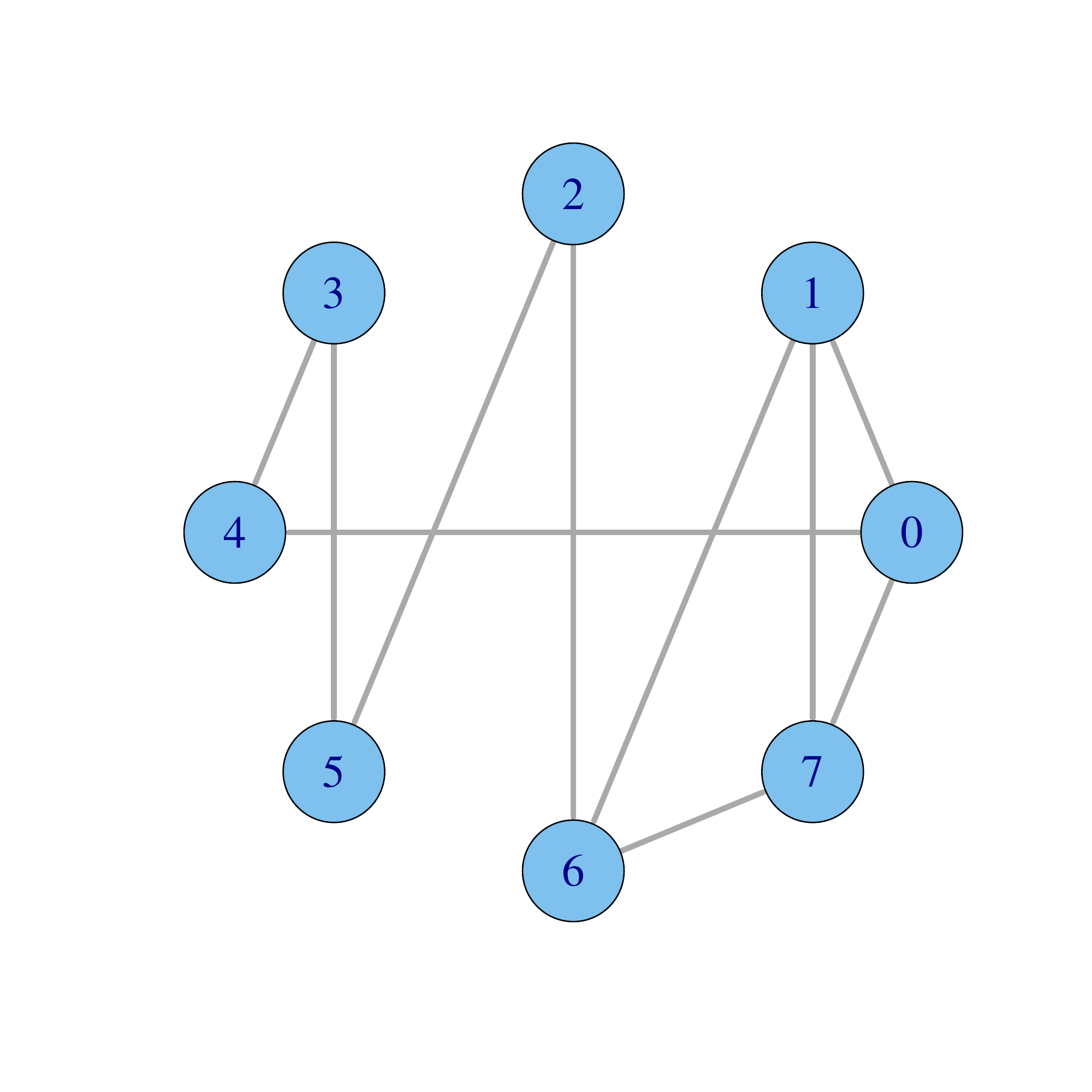}} &
$
A^{I_2} = \left(\begin{smallmatrix}
0&1&0&0&0&1&1&0\\
1&0&0&0&1&1&0&0\\
0&0&0&0&0&0&0&0\\
0&0&0&0&0&0&1&1\\
0&1&0&0&0&0&0&0\\
1&1&0&0&0&0&0&0\\
1&0&0&1&0&0&0&1\\
0&0&0&1&0&0&1&0\\
\end{smallmatrix}
\right)
$ &
\raisebox{-1.5cm}{\includegraphics[width=0.2\textwidth]{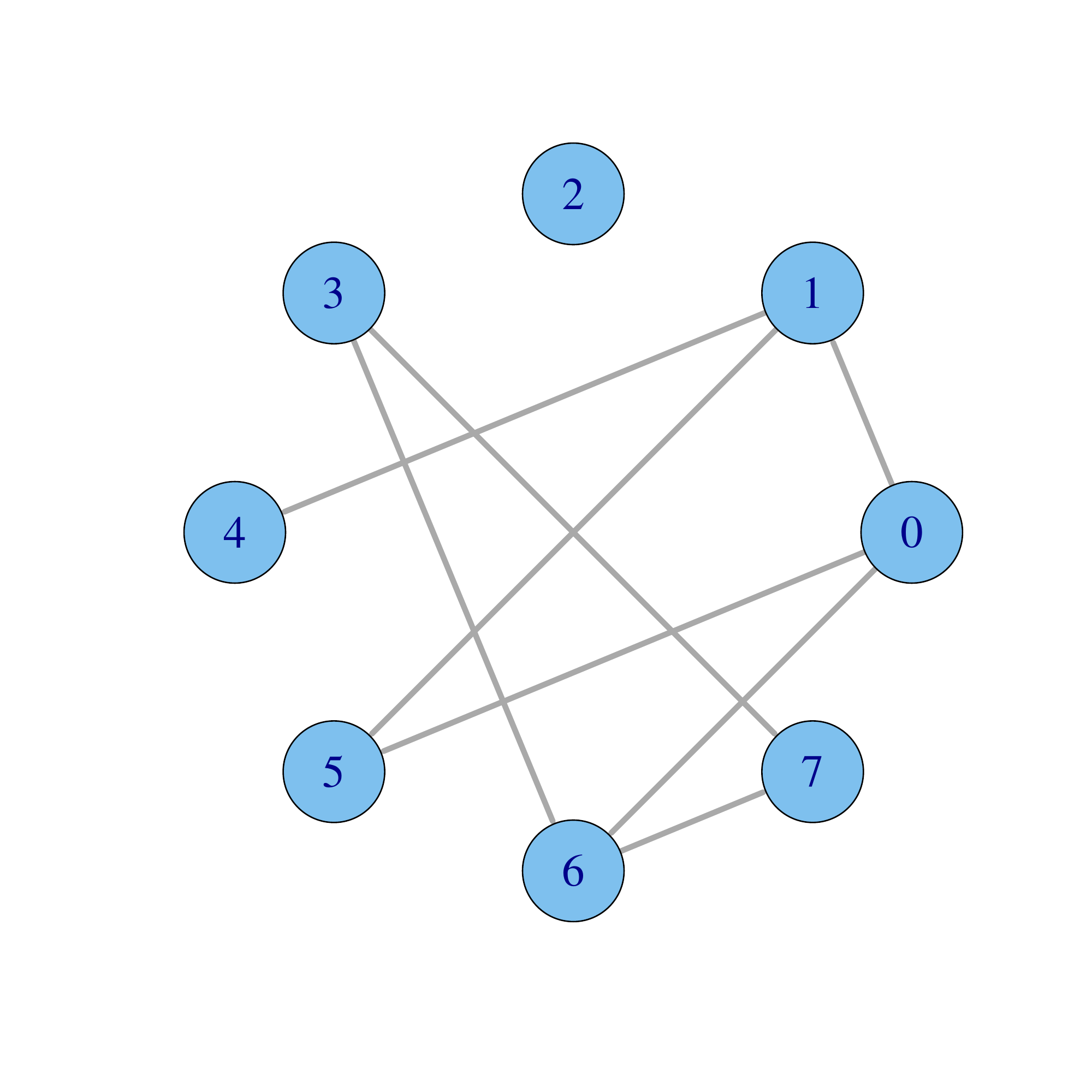}}
\end{tabular}
\end{center}
\caption{Adjacency matrix and graphical representation of $I_1$ and $I_2$.}
\label{fig:I1_I2}
\end{figure}
\section{The HIM kernel}
\label{sec:kernel}
Following \cite{cortes03positive}, a kernel can be naturally derived from a distance by means of a Gaussian (Radial Basis Function) map (see also \cite{bolla13spectral}).
Thus, given two graphs $x$ and $y$ on the same $n$ nodes and a positive real number $\gamma$, the HIM kernel can be defined as 
\begin{displaymath}
K(x,y) = e^{-\gamma\cdot\textrm{HIM}_\xi^2(x,y)}\ .
\end{displaymath}
Whenever a novel kernel is introduced, one has to check whether it is positively defined.

A function $\Psi\colon X\times X\to \mathbb{R}$ is a kernel of condionally negative type if
\begin{enumerate}
\item $\Psi(x,x)=0\quad\forall x\in X$;
\item $\Psi(x,y)=\Psi(y,x)\quad\forall x,y\in X$;
\item $\displaystyle{\sum_{i,j=1}^n} c_i c_j \Psi(x_i,x_j)\leq 0\quad \forall n\in\mathbb{N}, \forall x_1,\ldots,x_n\in X,\forall c_1,\ldots,c_n\in\mathbb{R}\;\textrm{such as}\; \displaystyle{\sum_{i=1}^n c_i}=0$.
\end{enumerate}
A variant of Schoenberg's theorem \cite{schoenberg38metric} (proved in \cite{ressel76short,bekka08kazhdan}) states that
\begin{theorem}
For a function $\Psi\colon X\times X\to \mathbb{R}$, the following are equivalent:
\begin{enumerate}
\item $\Psi$ is of conditionally negative type;
\item $K(x,y)=e^{-\gamma \Psi(x,y)}$ is a positive semidefinite kernel for all $\gamma\in\mathbb{R}_0^+$.
\end{enumerate}
\end{theorem}
The above theorem describes the correspondence between negative-type distances and positive definite kernel, which is also equivalent to $\ell_2^2$ embeddability \cite{berg84harmonic}.
Hence, $K(x,y)$ is a positive semidefinite kernel if and only if $\textrm{HIM}_\xi^2(x,y)$ is a symmetric function of conditionally negative type.
Although the square of many distances are condionally negative type functions, $\textrm{HIM}_\xi^2(x,y)$ cannot be proven to be of conditionally negative type (actually, it is probably not of negative type, as it is the case for many edit distances \cite{cortes03positive,martins06generative,neuhaus07bridging,li05class,cuturi09positive}, the HIM kernel $K$ is not positively defined in general for all $\gamma\in\mathbb{R}_0^+$.
Nevertheless, this problem can be overcomed by using Prop. 1.3.4 in \cite{schoelkopf97support} (see also \cite{bolla13spectral,li05class}):
\begin{theorem}
\label{th:k}
Suppose the data $x_1,\ldots,x_l$ and the kernel $k(\cdot,\cdot$) are such that the matrix
\begin{displaymath}
K_{ij} =k(x_i,x_j) 
\end{displaymath}
is positive. Then it is possible to construct a map $\Phi$ into a feature space $F$ such that
\begin{displaymath}
k(x_i, x_j) = \langle \Phi(x_i), \Phi(x_j) \rangle\ .
\end{displaymath}
Conversely, for a map $\Phi$Φ into some feature space $F$, the matrix $Ki_ij = \langle \Phi(x_i), \Phi(x_j) \rangle$ is positive.
\end{theorem}
Note that Th.~\ref{th:k} does not even require $x_1,\ldots,x_l$ to belong to a vector space.
This theorem implies that, even though the kernel is not positive definite, it is still possible to use it in Support Vector Machines or other algorithms requiring $k$ to correspond to a dot product in some space if the kernel matrix $K$ is positive for the given training data. 
This condition can be obtained by choosing a suitable value of $\gamma$: in the experiments shown hereafter, the HIM kernel is always positively defined on the given training data, leading to positive definite matrices, and thus posing no difficulties for the SVM optimizer, as in \cite{sonnenburg05large}.
\section{Applications}
\label{sec:apps}
\subsection{A minimal example}
\label{ssec:minimal}
Consider the two networks $I_1, I_2\in \pmb{N}(8)$ with corresponding adjacency matrices $A^{I_1}, A^{I_2}$ shown in Fig.~\ref{fig:I1_I2}.
The Hamming distance between $I_1$ and $I_2$ is 
\begin{displaymath}
\textrm{H}(I_1,I_2) = \frac{1}{N(N-1)} \sum_{1\leq i\not = j\leq N} \vert A^{I_1}_{ij} - A^{I_2}_{ij} \vert
= \frac{1}{56} \sum_{1\leq i\not= j \leq 8} 
\left(
\begin{smallmatrix}
0&0&0&0&1&1&1&1\\
0&0&0&0&1&1&1&1\\
0&0&0&0&0&1&1&0\\
0&0&0&0&1&1&1&1\\
1&1&0&1&0&0&0&0\\
1&1&1&1&0&0&0&0\\
1&1&1&1&0&0&0&0\\
1&1&0&1&0&0&0&0\\
\end{smallmatrix}
\right)
= \frac{28}{56}
= 0.5\ .
\end{displaymath}

\begin{figure}[!t]
\begin{center}
\begin{tabular}{ccc}
\includegraphics[width=0.33\textwidth]{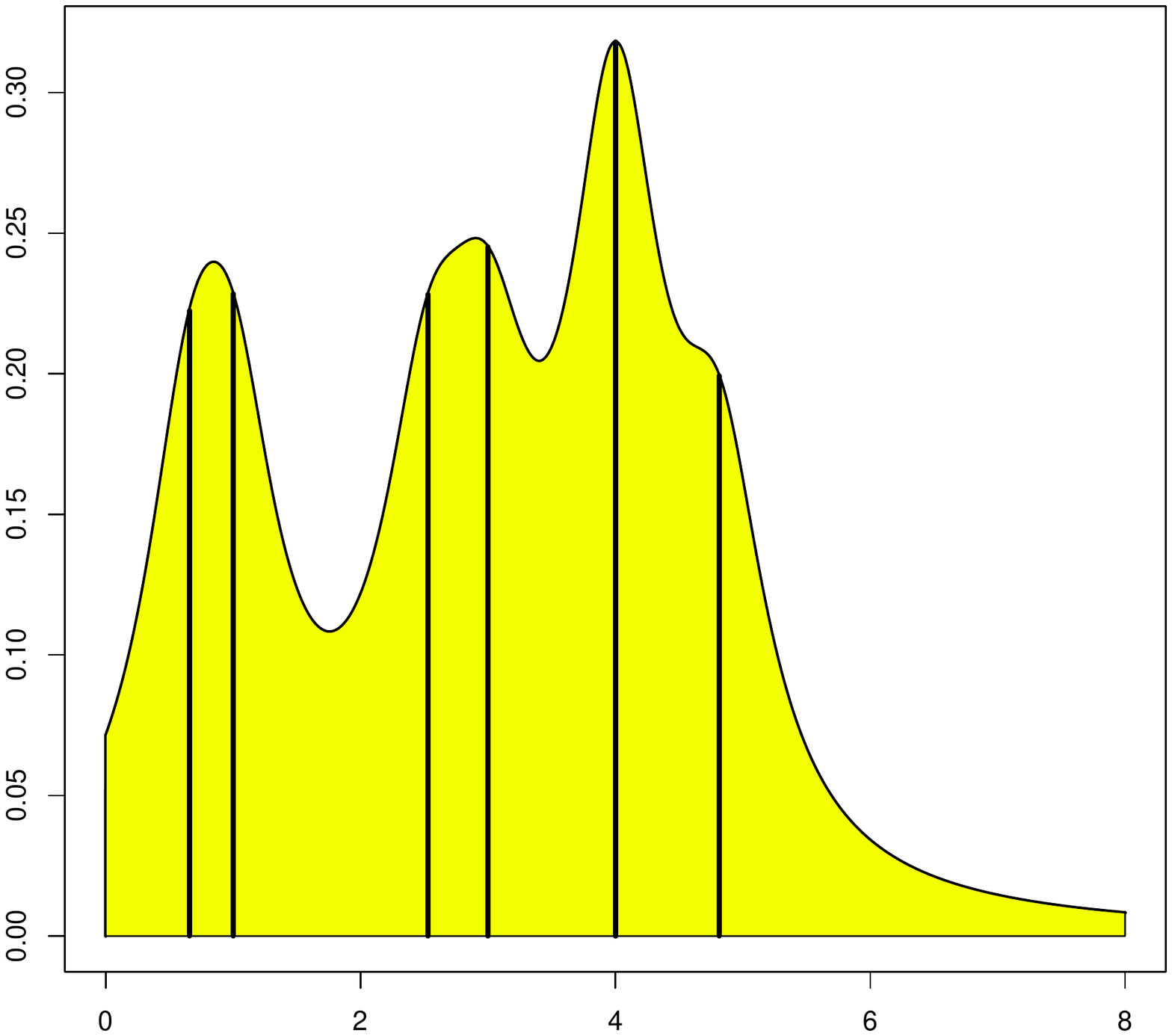}& \includegraphics[width=0.33\textwidth]{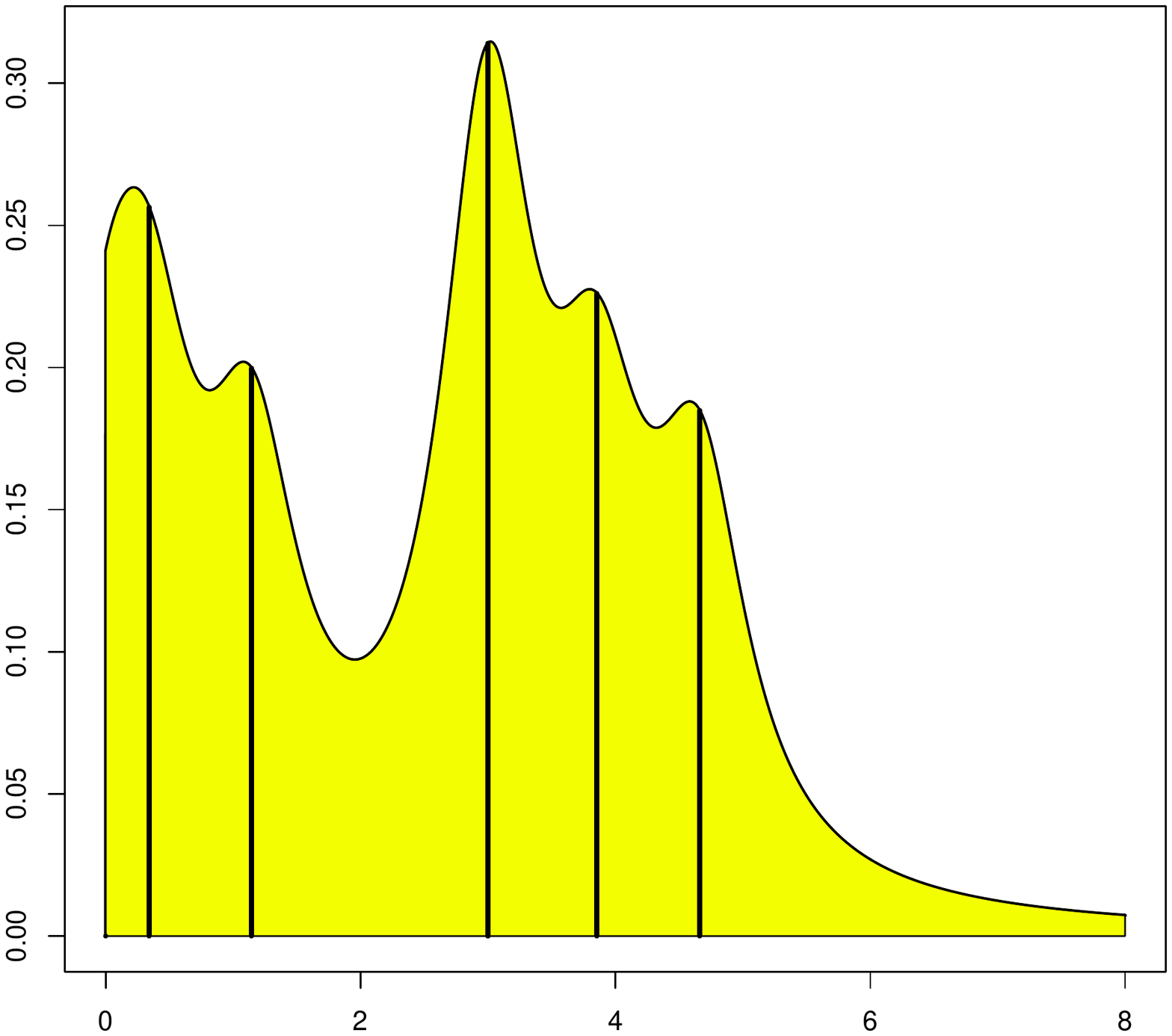}& \includegraphics[width=0.33\textwidth]{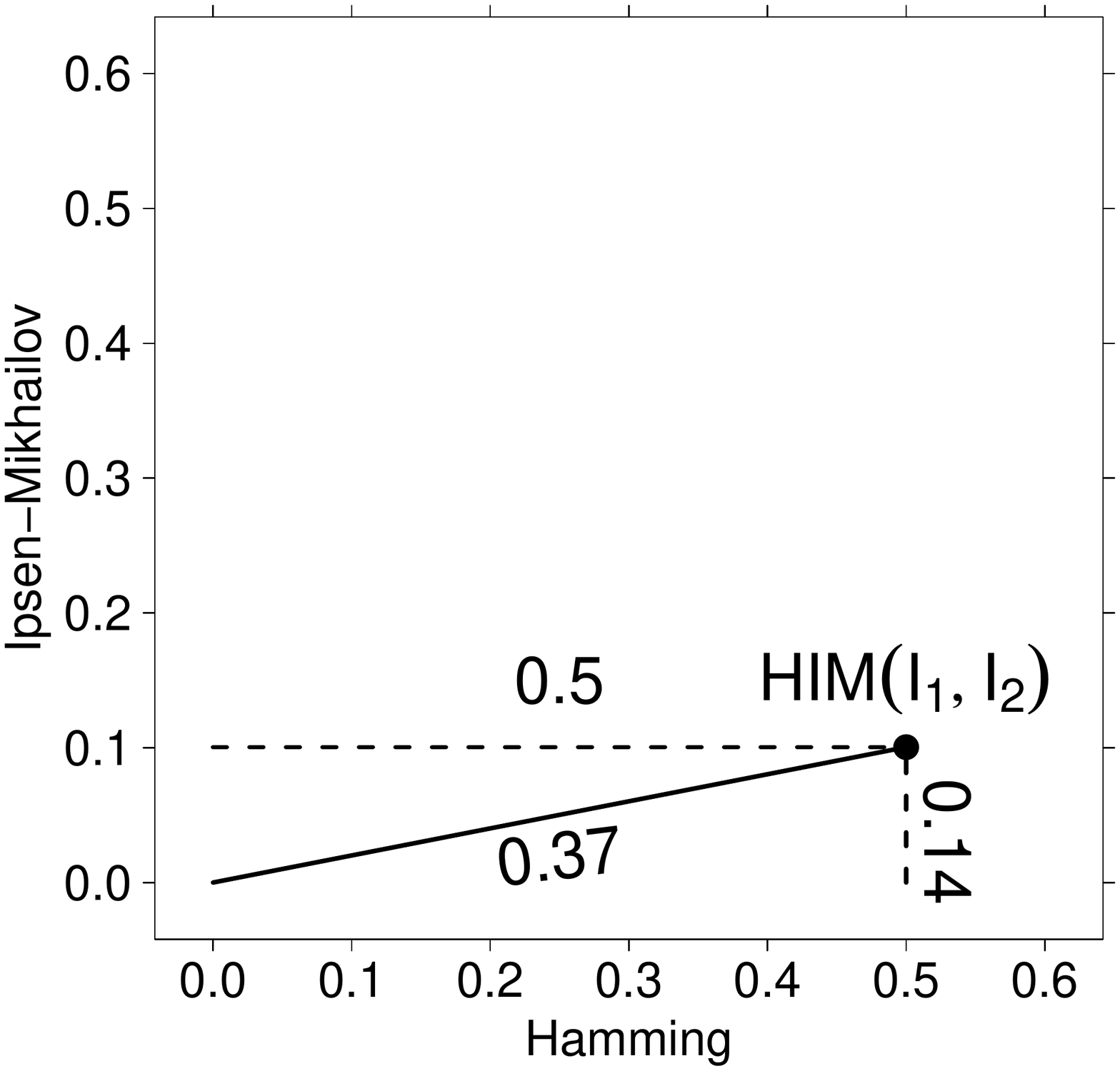}\\
\\
$\rho_{I_1}(\omega,\overline{\gamma})$ & $\rho_{I_2}(\omega,\overline{\gamma})$ & $\textrm{HIM}(I_1,I_2)$
\end{tabular}
\caption{Lorentzian distribution of the Laplacian spectra for $I_1$ (left) and $I_2$ (center) with vertical lines indicating eigenvalues, and $\textrm{HIM}(I_1,I_2)$ in the Hamming/Ipsen-Mikhailov space (right).}
\label{fig:lorentz}
\end{center}
\end{figure}

From the spectral point of view, the corresponding Laplacian matrices and eigenvalues are
\begin{align*}
L^{I_1} &= 
\left(
\begin{smallmatrix}
 3&-1& 0& 0&-1& 0& 0&-1\\
-1& 3& 0& 0& 0& 0&-1&-1\\
 0& 0& 2& 0& 0&-1&-1& 0\\
 0& 0& 0& 2&-1&-1& 0& 0\\
-1& 0& 0&-1& 2& 0& 0& 0\\
 0& 0&-1&-1& 0& 2& 0& 0\\
 0&-1&-1& 0& 0& 0& 3&-1\\
-1&-1& 0& 0& 0& 0&-1& 3\\
\end{smallmatrix}
\right)
\quad
&\textrm{spec}(L^{I_1}) &= [0,0.657077,1,2.529317,3,4,4,4.813607]
\\
L^{I_2} &= 
\left(
\begin{smallmatrix}
 3&-1& 0& 0& 0&-1&-1& 0\\
-1& 3& 0& 0&-1&-1& 0& 0\\
 0& 0& 0& 0& 0& 0& 0& 0\\
 0& 0& 0& 2& 0& 0&-1&-1\\
 0&-1& 0& 0& 1& 0& 0& 0\\
-1&-1& 0& 0& 0& 2& 0& 0\\
-1& 0& 0&-1& 0& 0& 3&-1\\
 0& 0& 0&-1& 0& 0&-1& 2\\
\end{smallmatrix}
\right)
\quad
&\textrm{spec}(L^{I_2}) &= [0,0,0.340321,1.145088,3,3,3.854912,4.659679]\ . 
\end{align*}
From the above spectra, we can compute the corresponding Lorentz distributions $\rho_{I_{\{1,2\}}}(\omega,\overline{\gamma})$, where $\overline{\gamma}=0.4450034$: their plots are shown in Fig.~\ref{fig:lorentz}.

The resulting Ipsen-Mikhailov distance is
\begin{displaymath}
\textrm{IM}(I_1,I_2)= \sqrt{\int_0^\infty \left[\rho_{I_1}(\omega,\overline{\gamma})-\rho_{I_2}(\omega,\overline{\gamma})\right]^2 \textrm{d}\omega}= 0.1004144\ ,
\end{displaymath}
so that the HIM distance results
\begin{displaymath}
\textrm{HIM}(I_1,I_2)=\frac{\sqrt{2}}{2}|| (\textrm{H}(I_1,I_2) , \textrm{IM}(I_1,I_2)) ||_{2} \approx 0.707168\sqrt{0.5^2+0.1004144^2} \approx 0.3606127\ .
\end{displaymath}
The situation can be graphically represented as in Fig.~\ref{fig:lorentz}: the two networks are quite different in terms of matching links, but their structures are not so diverse.

\subsection{Small networks}
\label{ssec:small}
Fixed the number of nodes $N$, there are exactly $2^\frac{N(N-1)}{2}$ different simple undirected unweighted networks, which can be grouped into isomorphism classes. As anticipated before, isomorphic graphs cannot be distiguished by spectral metrics, while their mutual Hamming distances are non zero, since their links are in different positions.
As an example, for $N=3$ there are 8 networks grouped in 4 isomorphism classes, for $N=4$ there are 11 isomorphism classes including a total of 64 graphs and for $N=5$ 34 classes with 1024 networks (for $N=6,7$, the number of classes is respectively 156 e 1044).

\begin{figure}[!t]
\begin{center}
\begin{tabular}{cccc}
\includegraphics[width=0.22\textwidth]{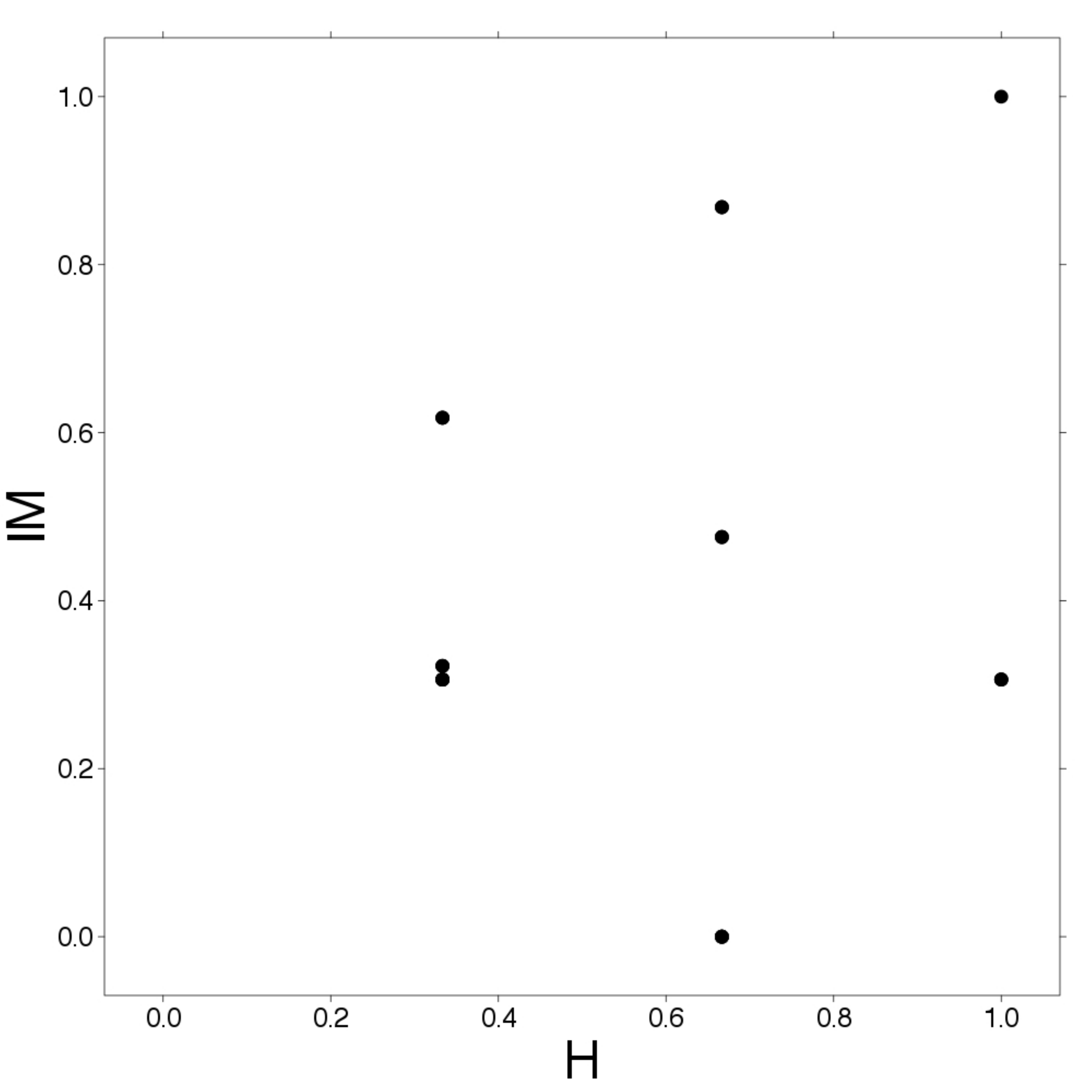} & 
\includegraphics[width=0.22\textwidth]{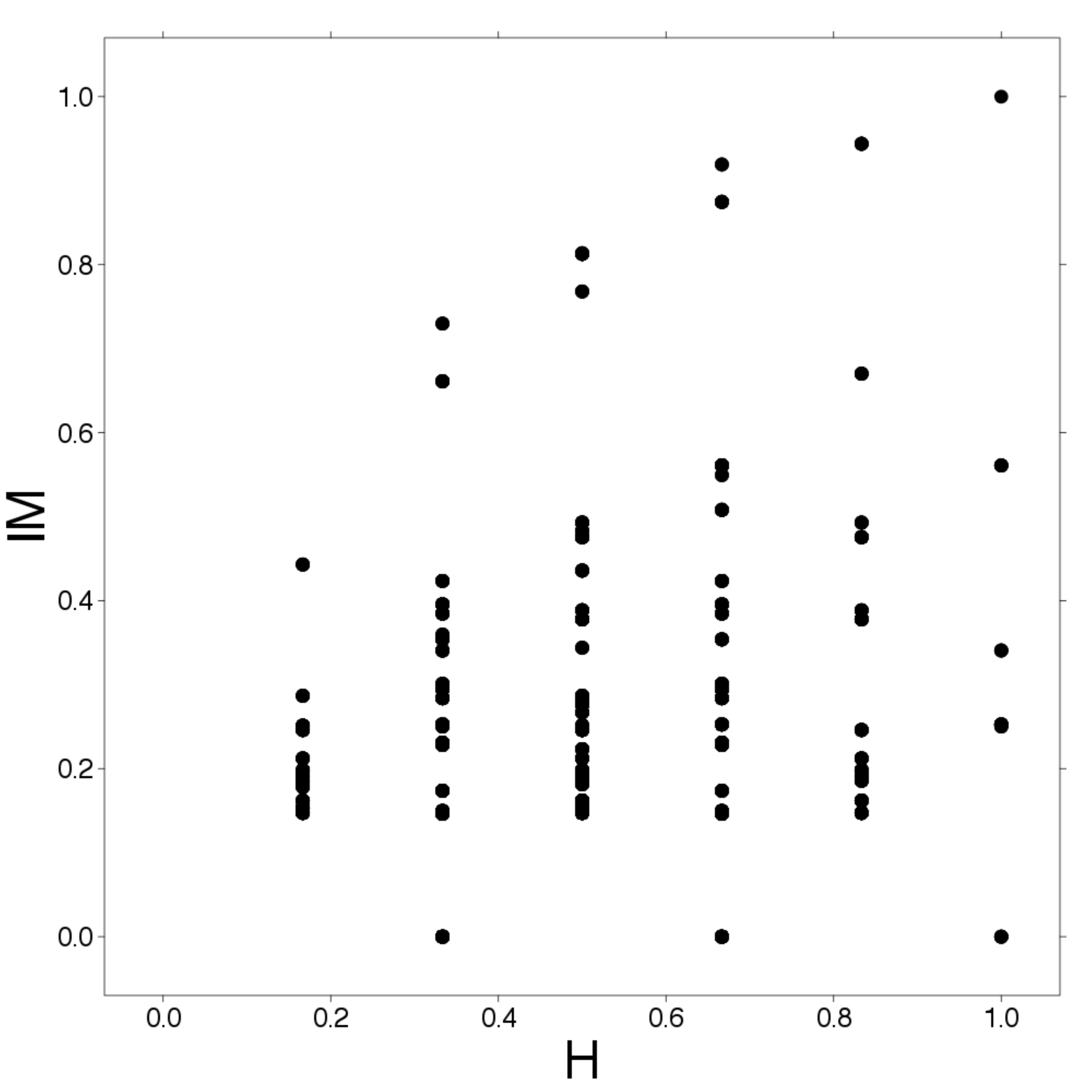} &
\includegraphics[width=0.22\textwidth]{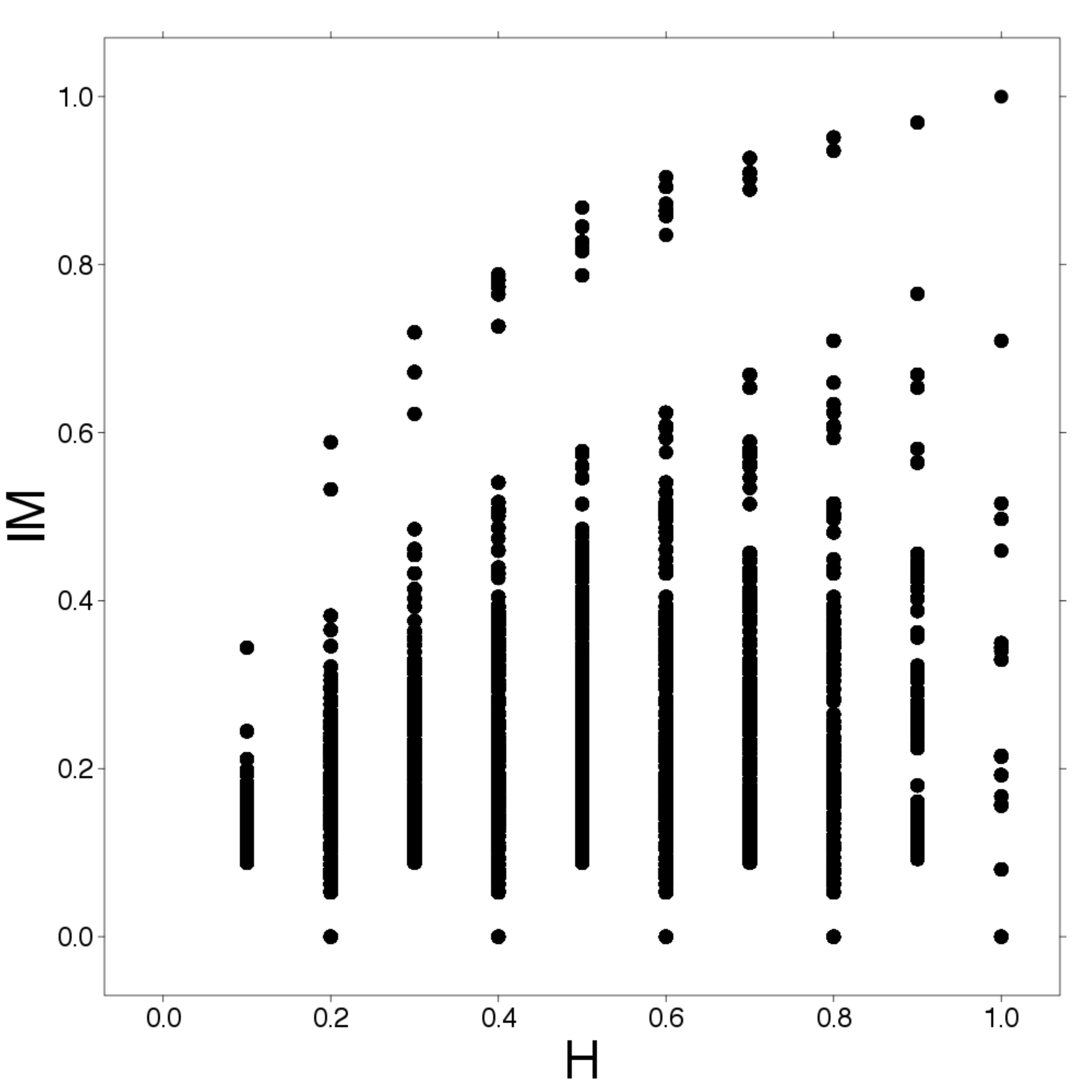} & 
\includegraphics[width=0.22\textwidth]{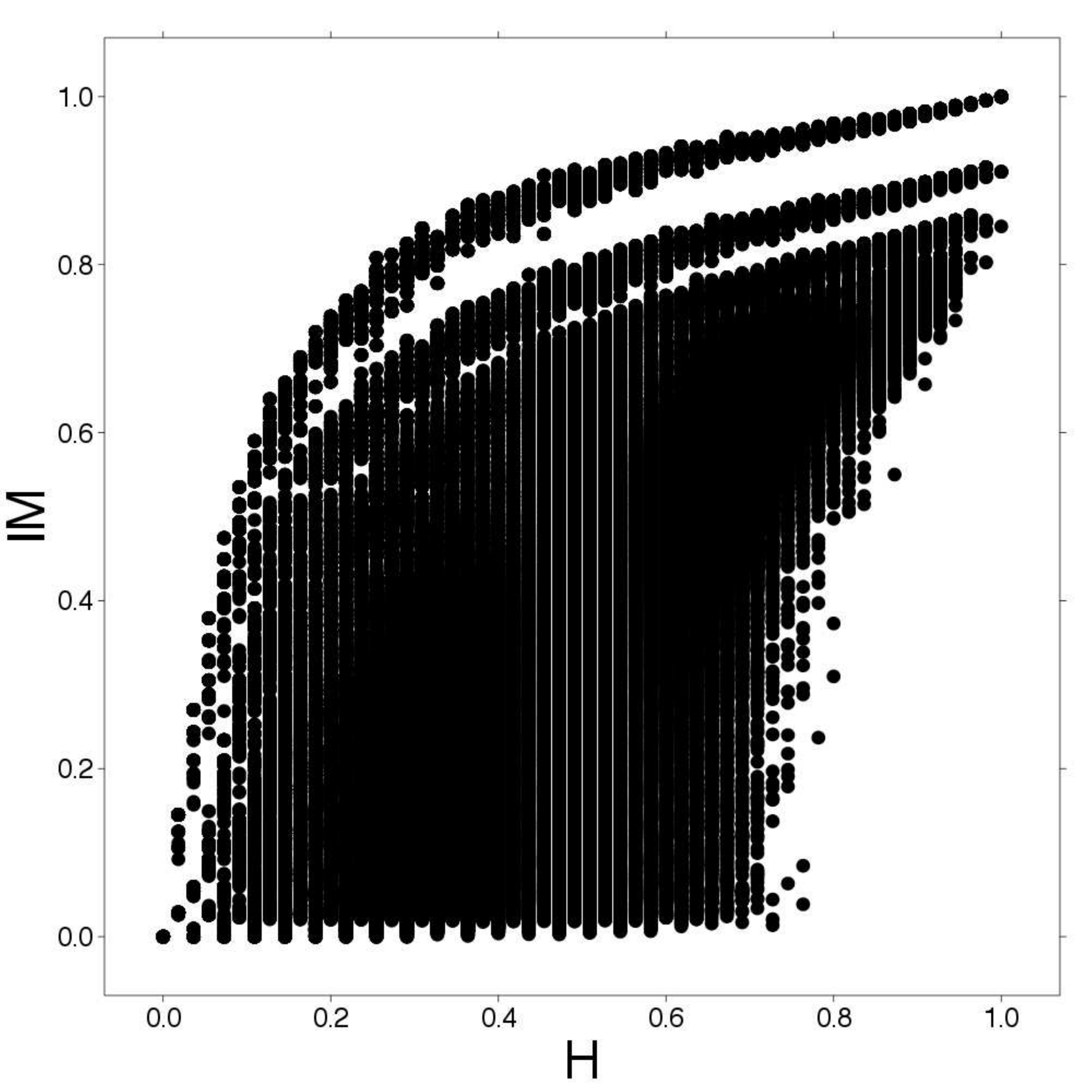} \\ 
(a) & (b) &
(c) & (d) \\
\end{tabular}
\caption{Mutual distances between (a) all 28 couples of networks with 3 nodes, (b) all 2016 couples of networks with 4 nodes, (c) all 523776 couples of networks with 5 nodes and (d) the 542361 mutual distances between a set of 1042 networks with 15 nodes.}
\label{fig:mutual}
\end{center}
\end{figure}

To give an overview of a broader situation, we compute a number of mutual distances between networks with a given number of nodes (all possible couples for $N=3,4,5$ and a subset of them for $N=15$) and we display the results in Fig.~\ref{fig:mutual}.
To select a good range of variability for the networks with 15 nodes, we select the empty graph, the full graph (with 105 nodes) and 10 different graphs with $i$ edges each, for $1\leq i\leq 104$.

\begin{figure}[!t]
\begin{center}
\begin{tabular}{cc|cc|cc}
\includegraphics[width=0.13\textwidth]{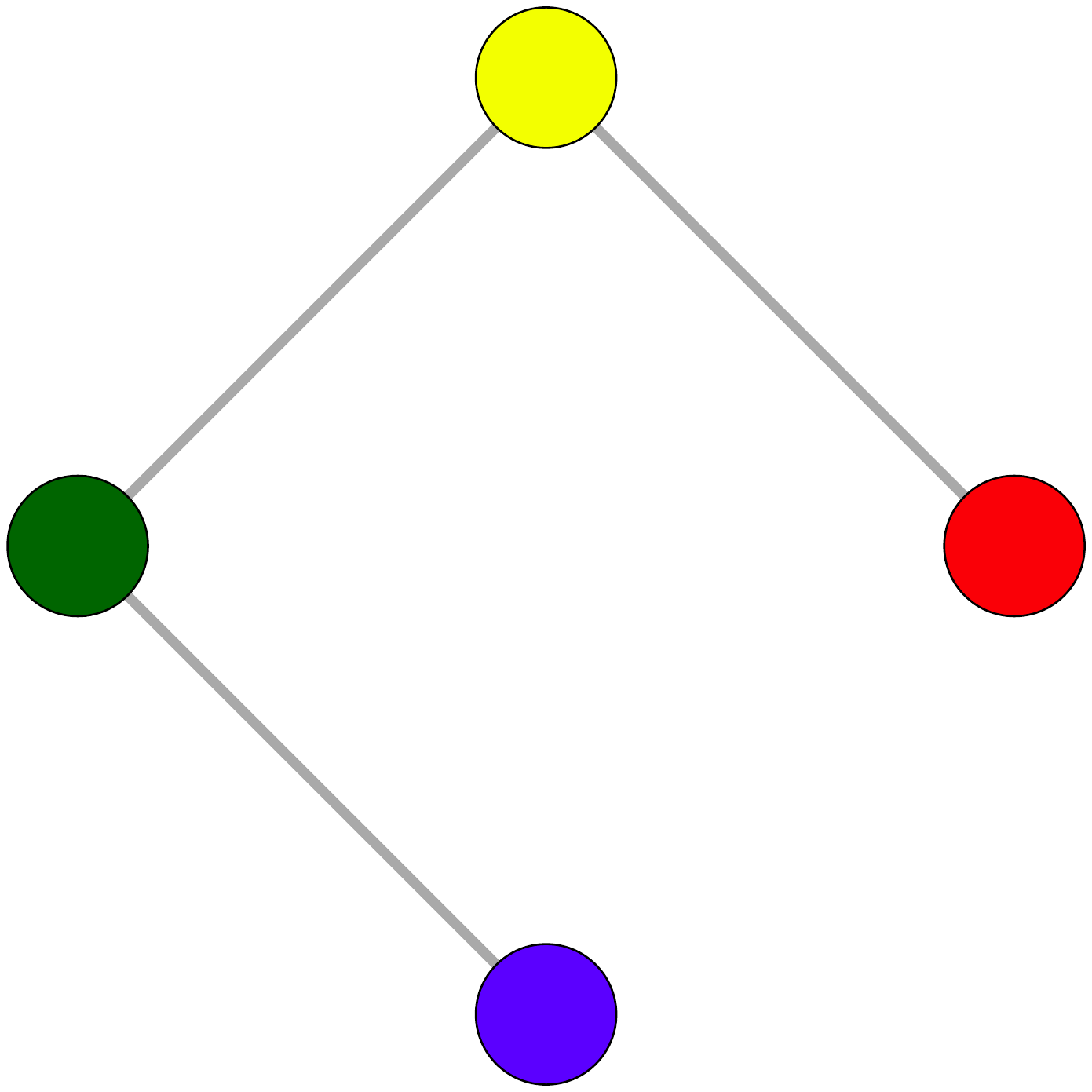} & 
\includegraphics[width=0.13\textwidth]{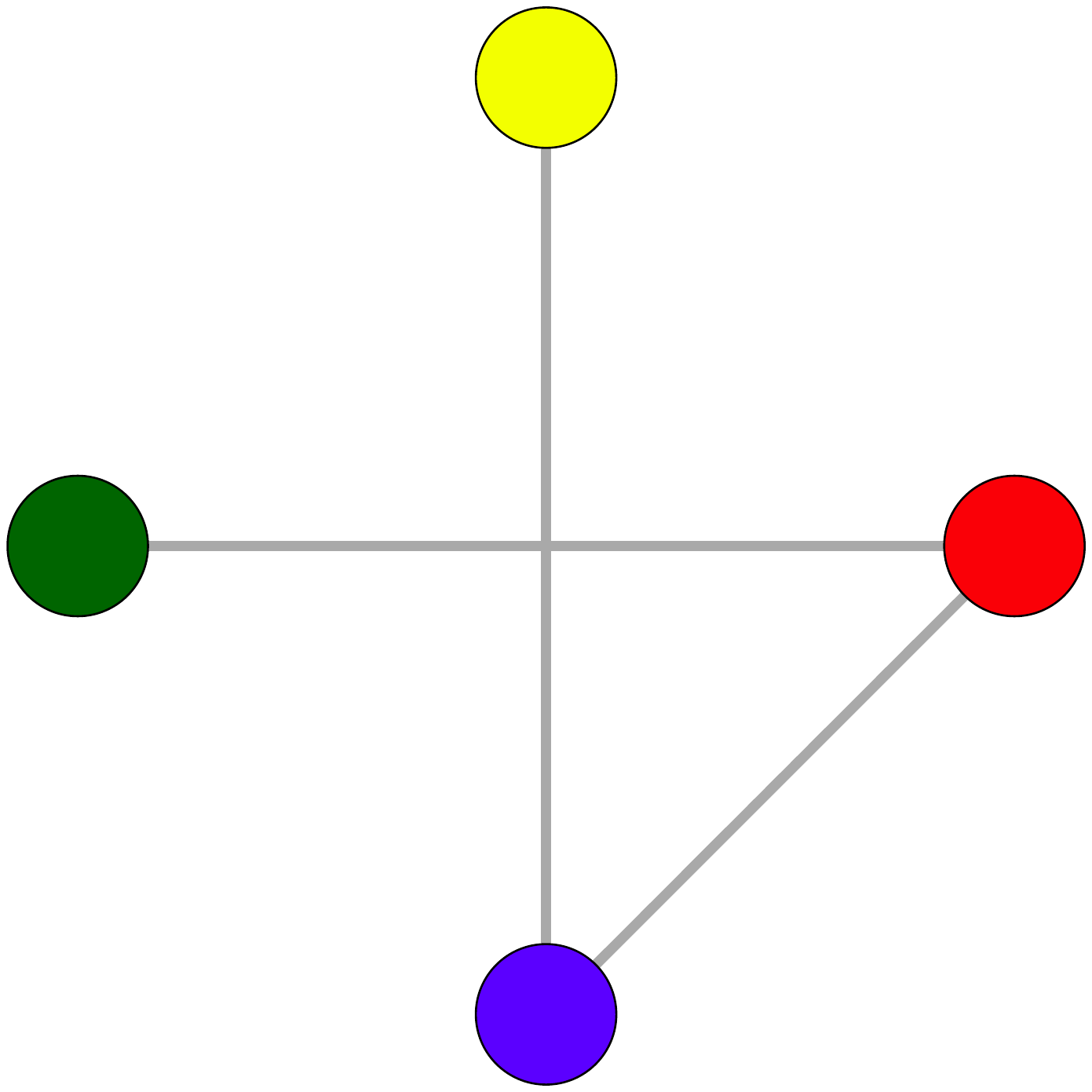} & 
\includegraphics[width=0.13\textwidth]{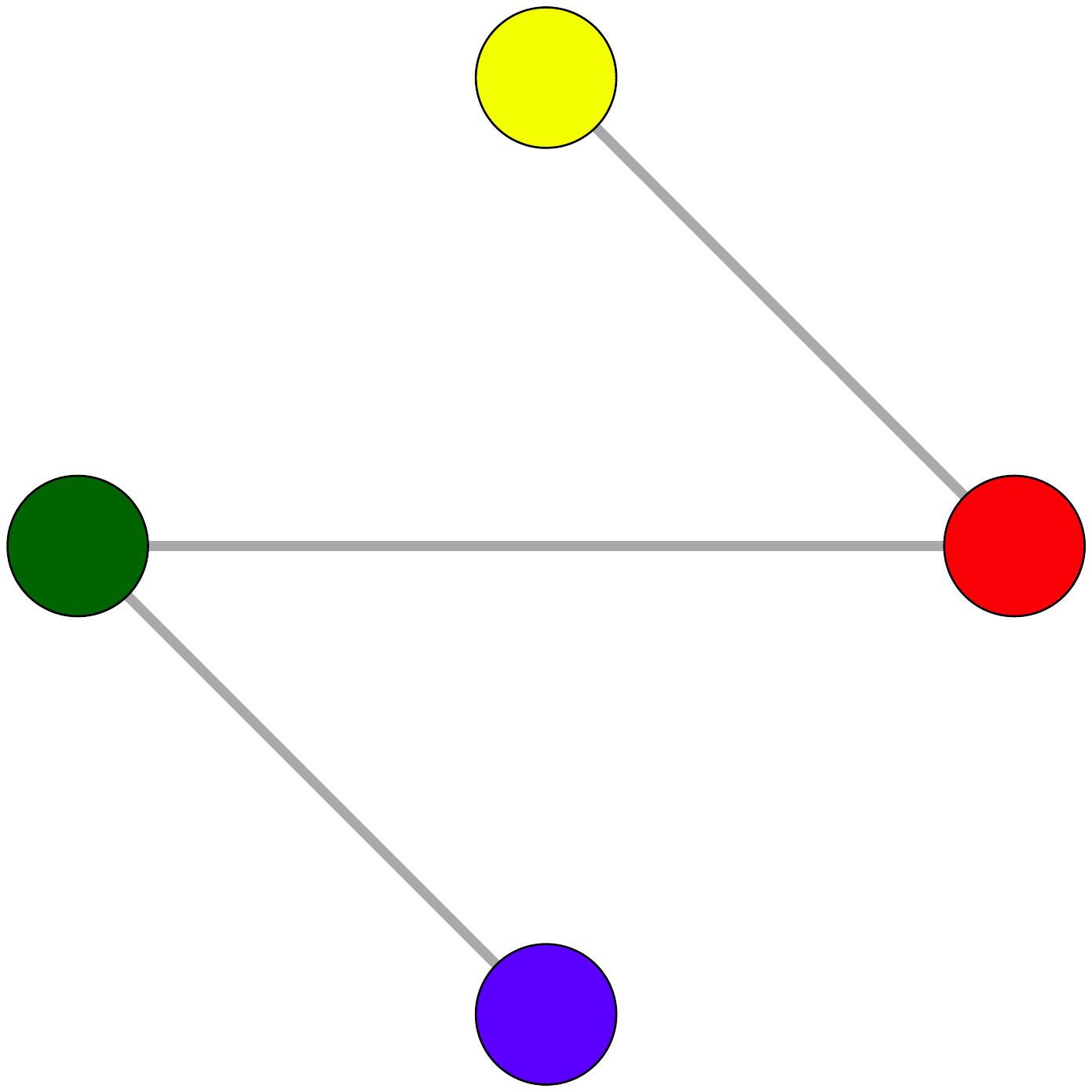} & 
\includegraphics[width=0.13\textwidth]{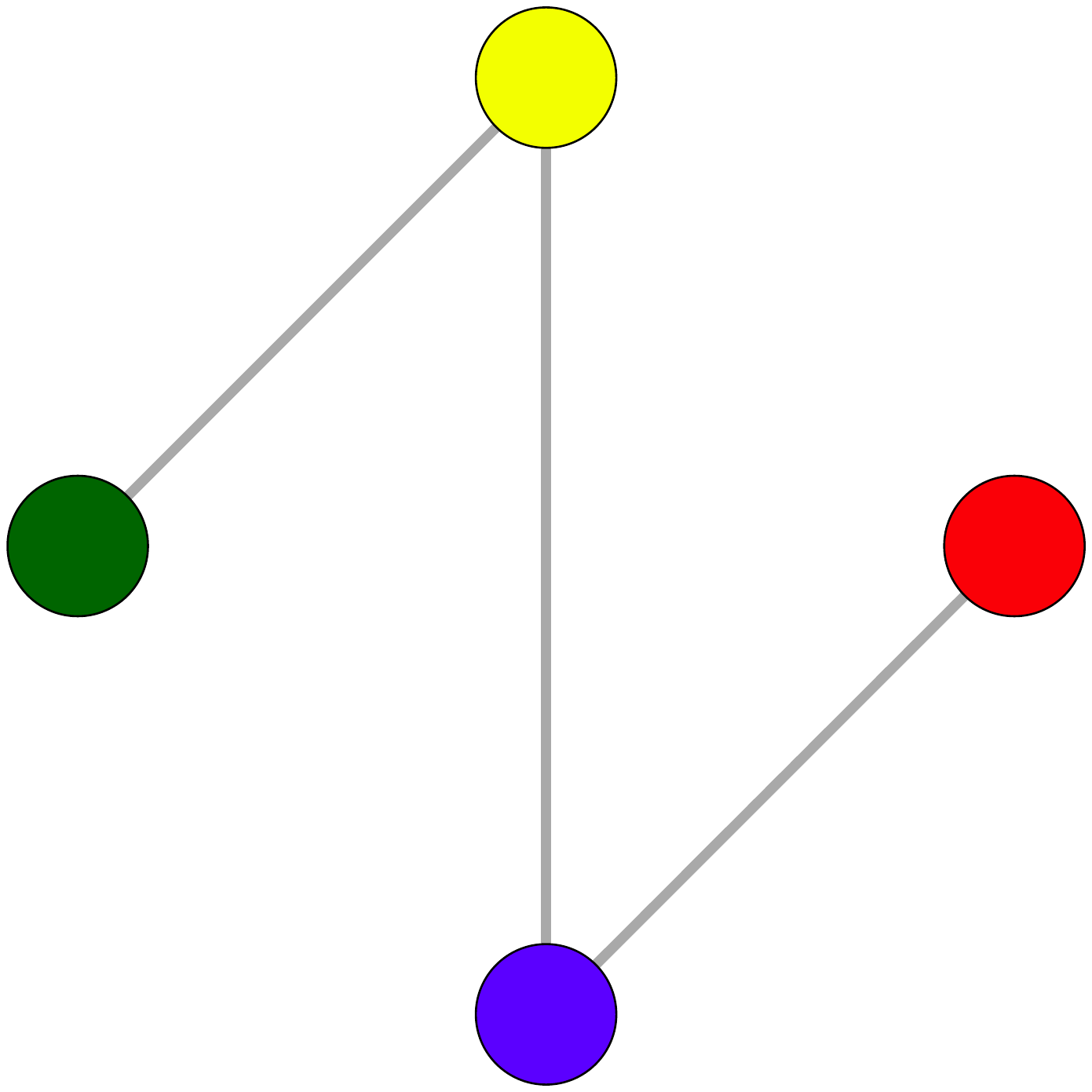} & 
\includegraphics[width=0.13\textwidth]{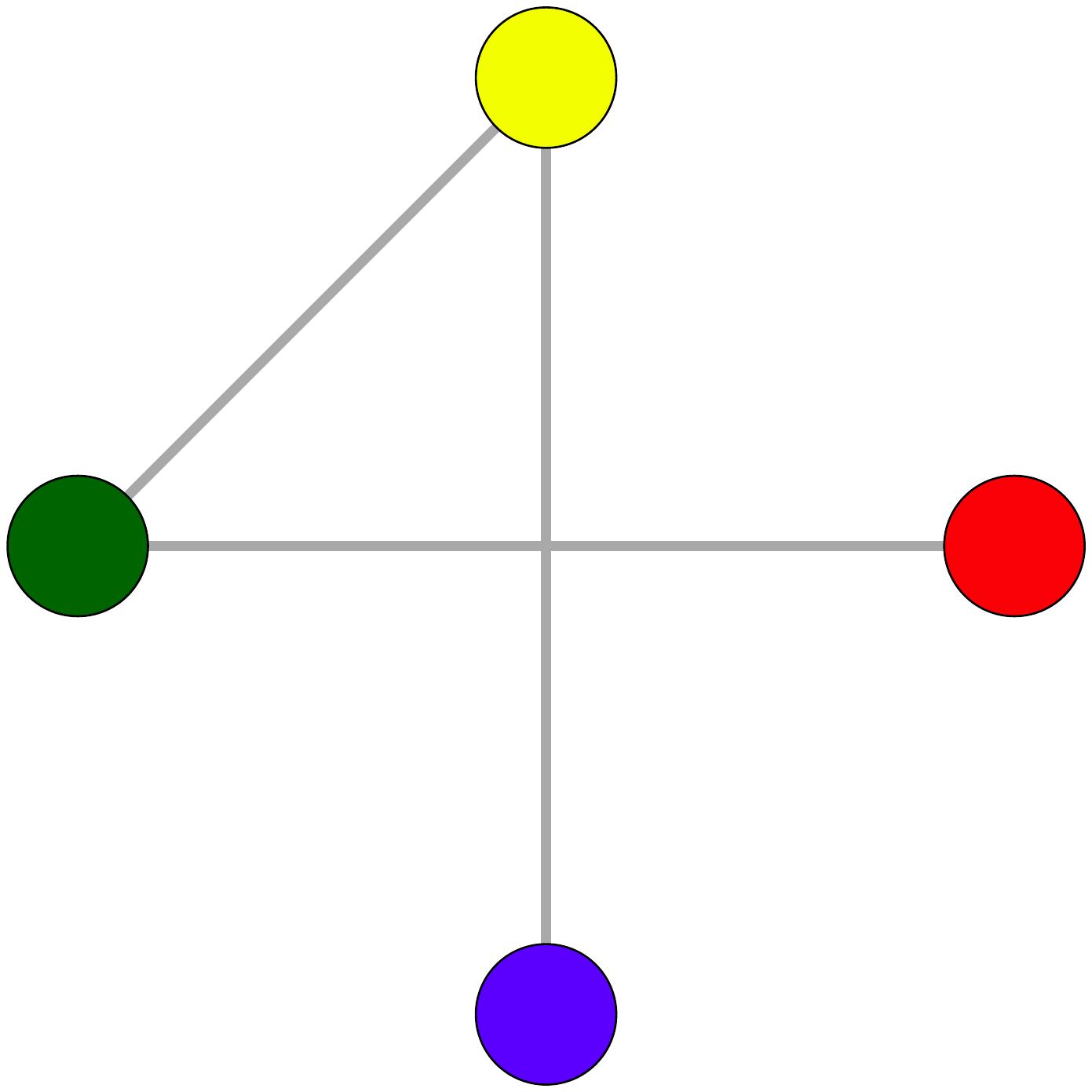} & 
\includegraphics[width=0.13\textwidth]{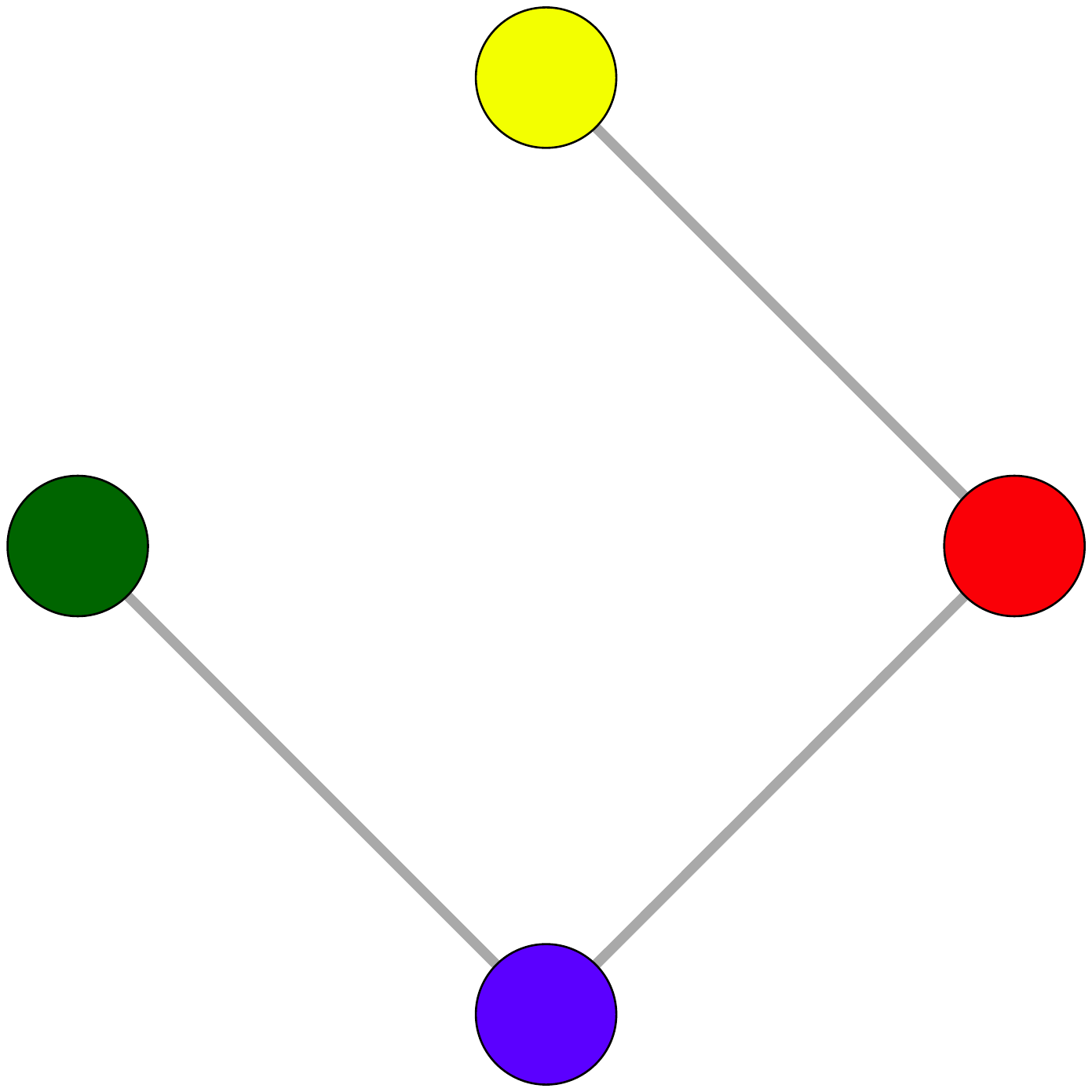} \\
\hline
\includegraphics[width=0.13\textwidth]{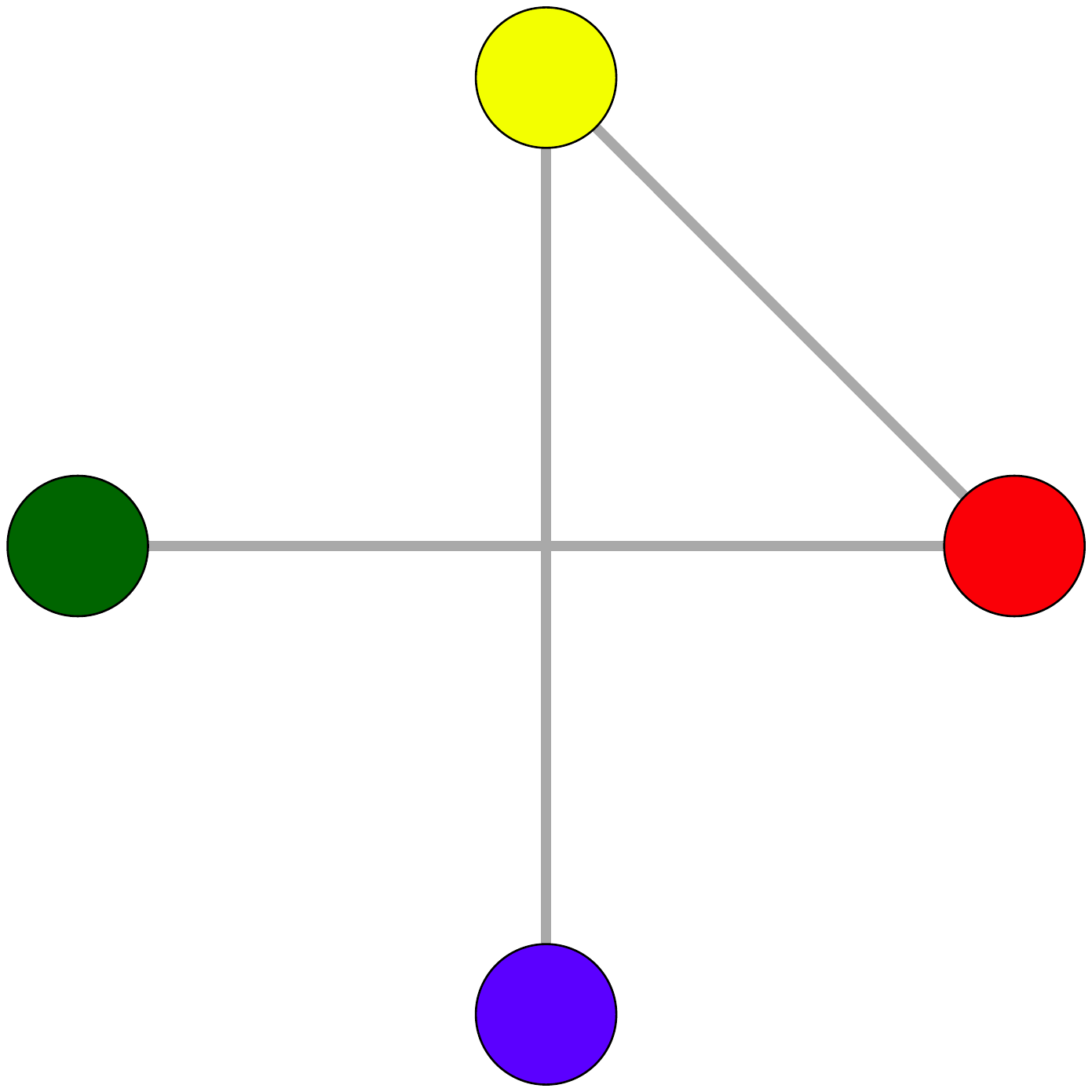} & 
\includegraphics[width=0.13\textwidth]{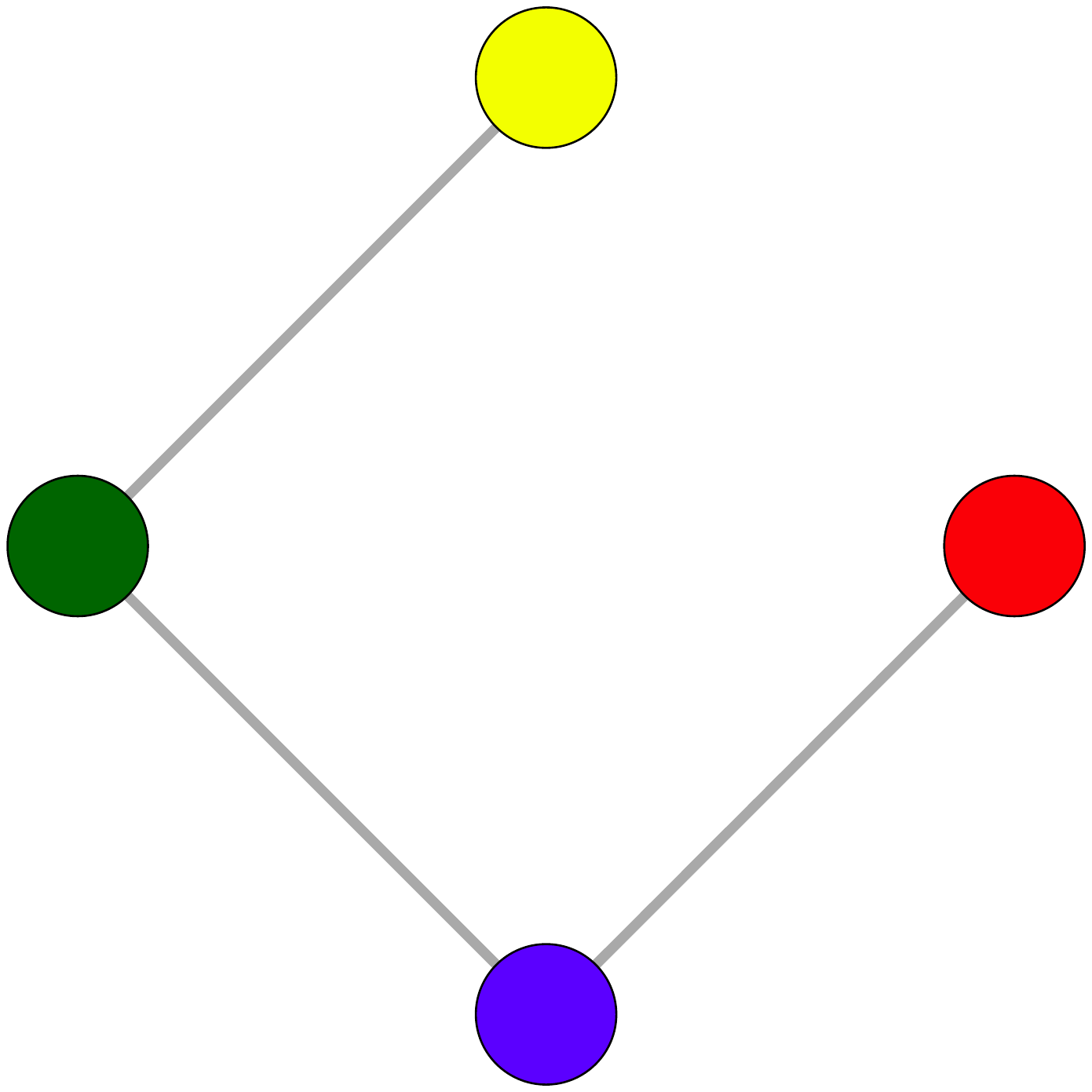} & 
\includegraphics[width=0.13\textwidth]{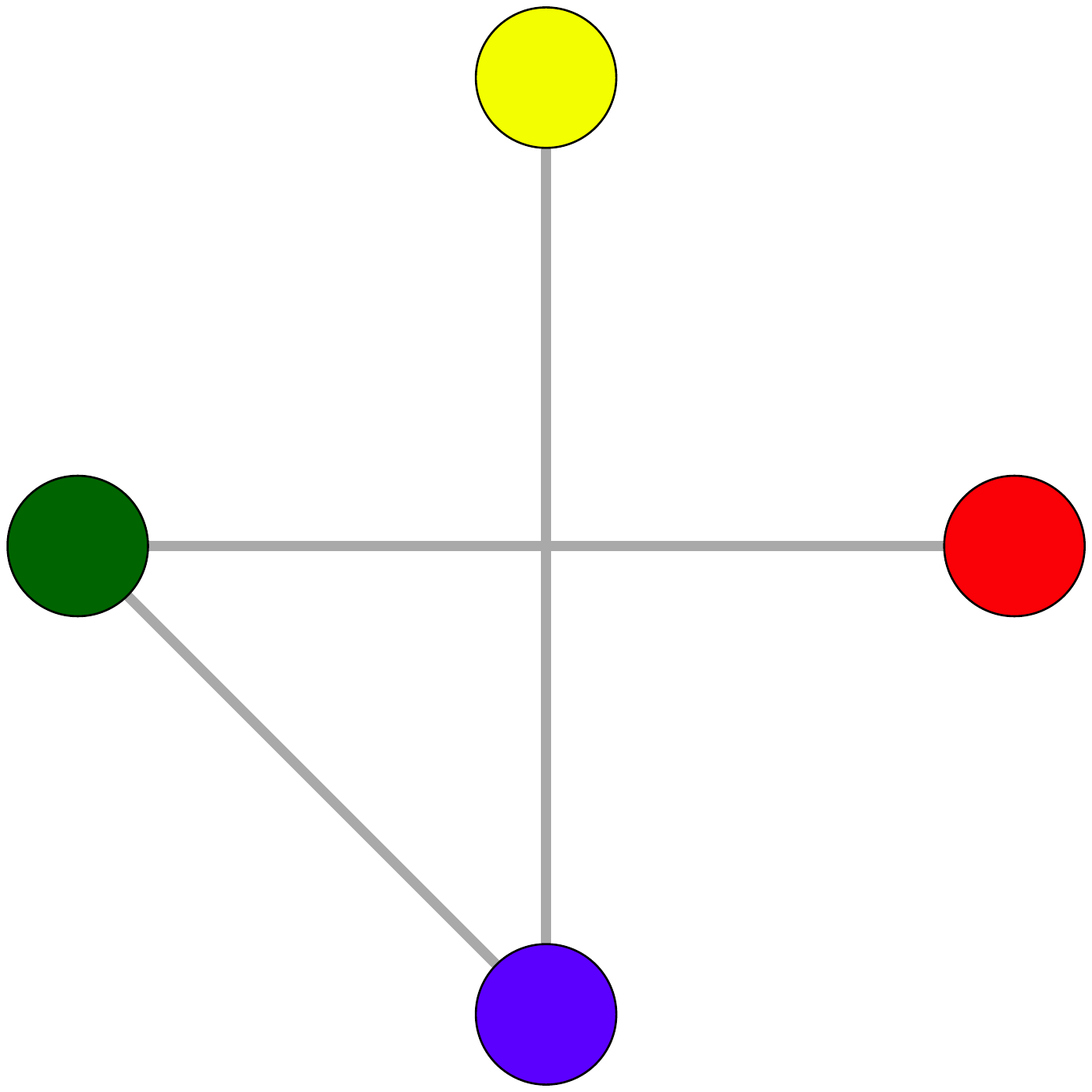} & 
\includegraphics[width=0.13\textwidth]{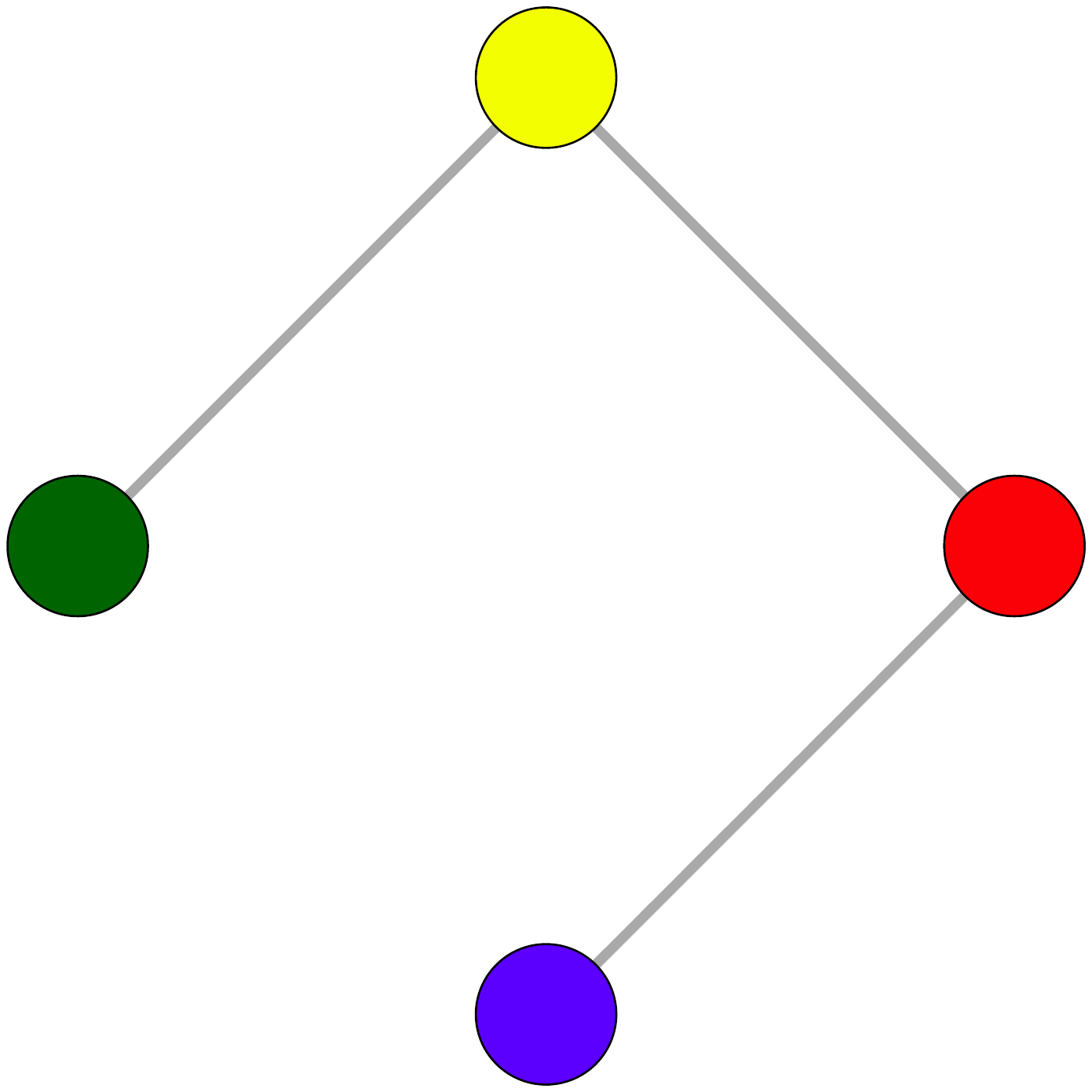} & 
\includegraphics[width=0.13\textwidth]{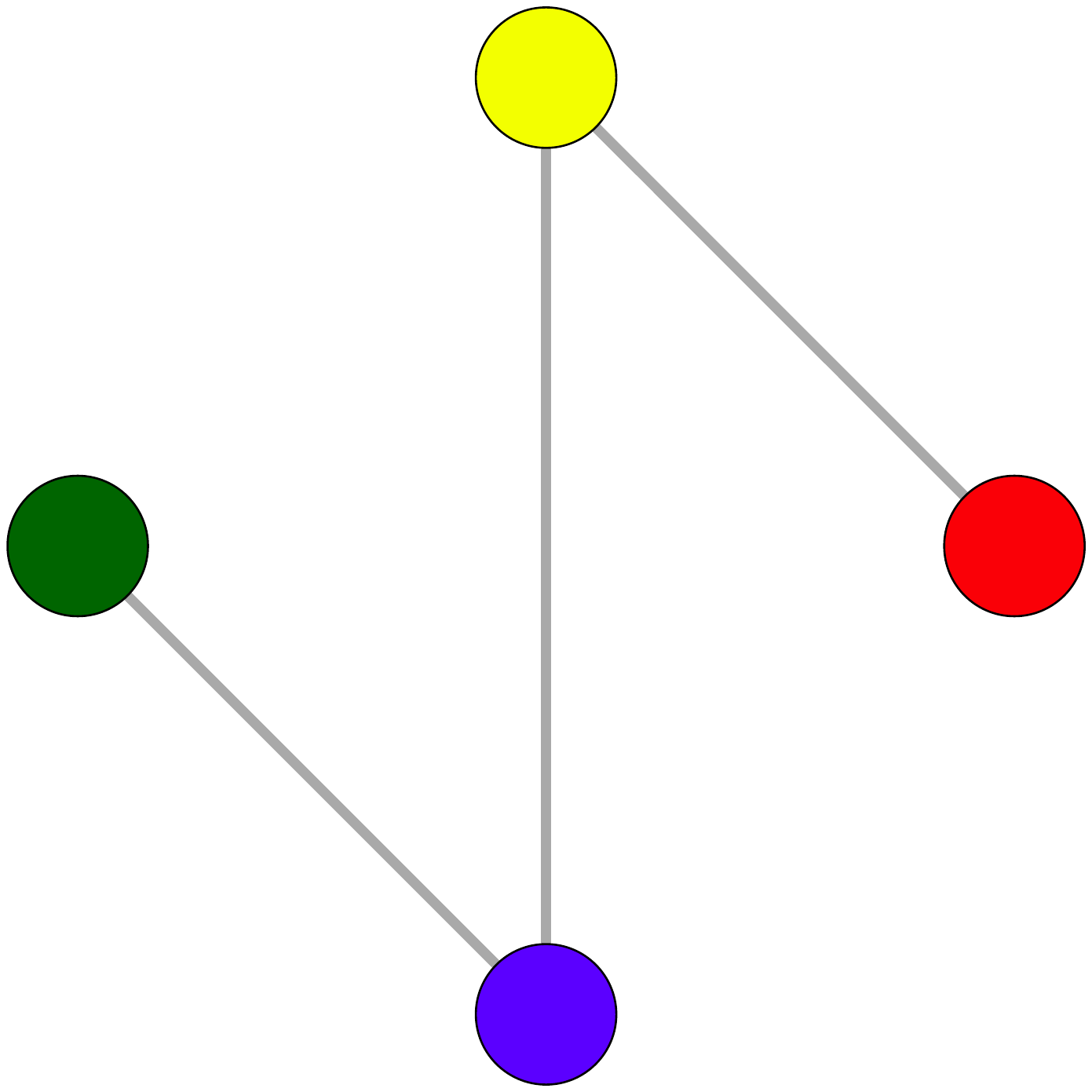} & 
\includegraphics[width=0.13\textwidth]{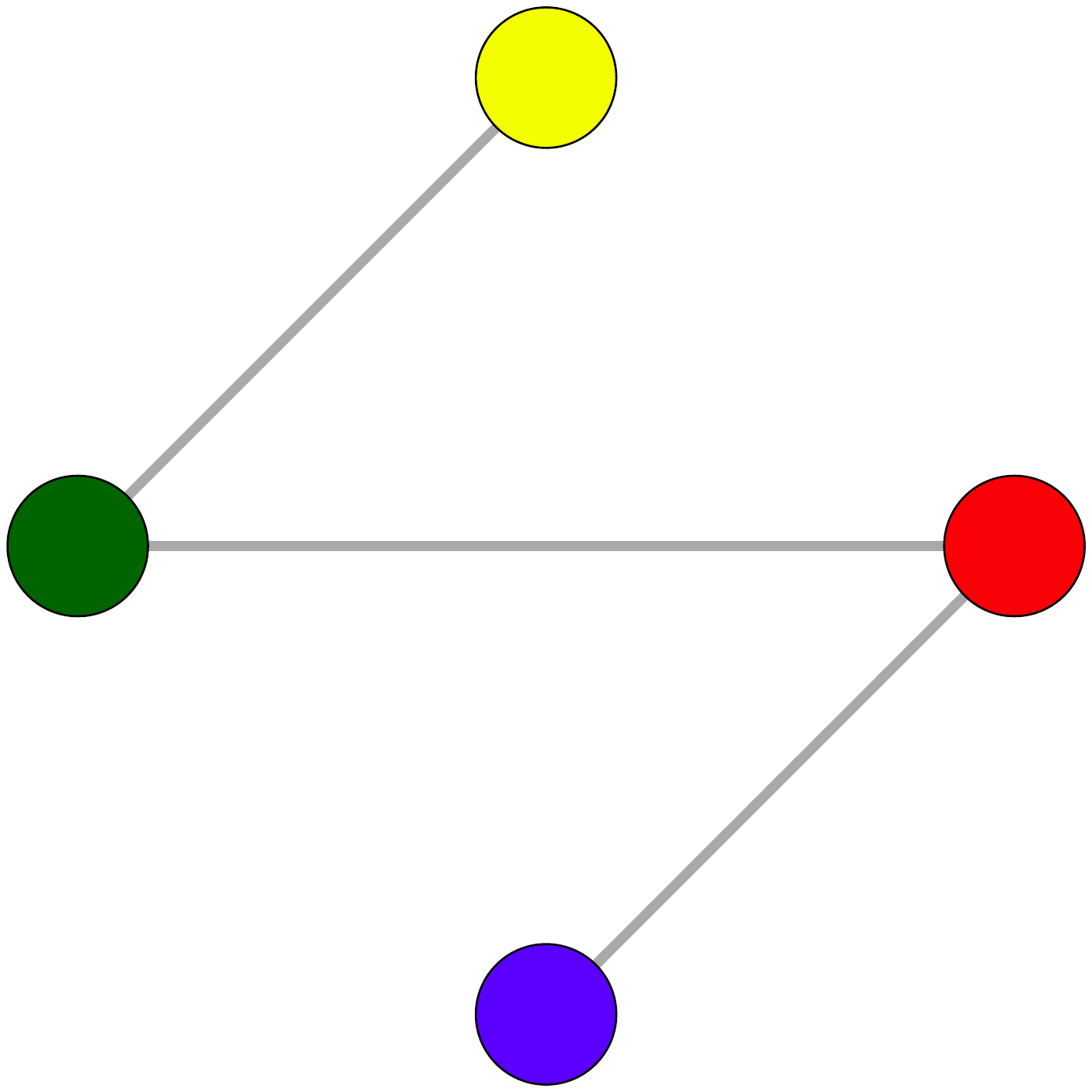} \\
\end{tabular}
\caption{The six pairs of networks on four nodes with Hamming distance one and Ipsen-Mikhailov distance zero.}
\label{fig:P}
\end{center}
\end{figure}

As shown by the plots, all possible situations can occur, apart from points in the northwest corner of zone II which are the rarest. 
For instance, the point $P(1,0)$ in Fig.~\ref{fig:mutual}(b) corresponds to 6 different pairs $(O_1,O_2)$ of networks with $4$ nodes with maximal Hamming distance and minimal spectral distance: as an example, we show one of these pairs in Fig.~\ref{fig:P}.

\subsection{Comparison with Matthews Correlation Coefficient}
\label{ssec:mcc}
When assessing performances in a link prediction task (for instance, in the series of DREAM challenges \cite{stolovitzky07dialogue,marbach10revealing,prill10towards}), the standard strategy following the machine learning approach, is to rely on functions of the confusion matrix, \textit{i.e.}, the table collecting the number of correct and wrong predictions with respect to the ground truth.
Classical measures of this kind are the pairs Sensitivity/Specificity and Precision/Recall, and the derived Area Under the Curve.

A reliable alternative is the Matthews Correlation Coefficient ($\textrm{MCC}$ for short) \cite{matthews75comparison}, summarizing into a single value the confusion matrix of a binary classification task.
This is a measure of common use in the machine learning community \cite{baldi00assessing}, recently accepted as an effective metric also for network comparison \cite{supper07reconstructing,stokic09fast}.
Also known as the $\phi$-coefficient, for a $2\times 2$ contingency table $\textrm{MCC}$ corresponds to the square root of the average $\chi^2$ statistic 
\begin{displaymath}
\textrm{MCC}=\sqrt{\chi^2 / N}\ ,
\end{displaymath}
where $N$ is the total number of observations.
In the binary case of two classes positive P and negative N, for the confusion matrix $\left(\begin{smallmatrix} \textrm{TP} & \textrm{FN} \\ \textrm{FP} & \textrm{TN}\end{smallmatrix}\right)$, where T and F stand for true and false respectively, the Matthews Correlation Coefficient has the following shape:
\begin{displaymath}
\textrm{MCC} = \frac{\textrm{TP}\cdot\textrm{TN}-\textrm{FP}\cdot\textrm{FN}}{\sqrt{\left(\textrm{TP}+\textrm{FP}\right)\left(\textrm{TP}+\textrm{FN}\right)\left(\textrm{TN}+\textrm{FP}\right)\left(\textrm{TN}+\textrm{FN}\right)}}\ .
\end{displaymath}
$\textrm{MCC}$ lives in the range $[-1,1]$, where $1$ is perfect classification, $-1$ is reached in the complete misclassification case while $0$ corresponds to coin tossing classification, and it is invariant for scalar multiplication of the whole confusion matrix.

\begin{figure}[!t]
\begin{tabular}{ccc}
\includegraphics[width=0.3\textwidth]{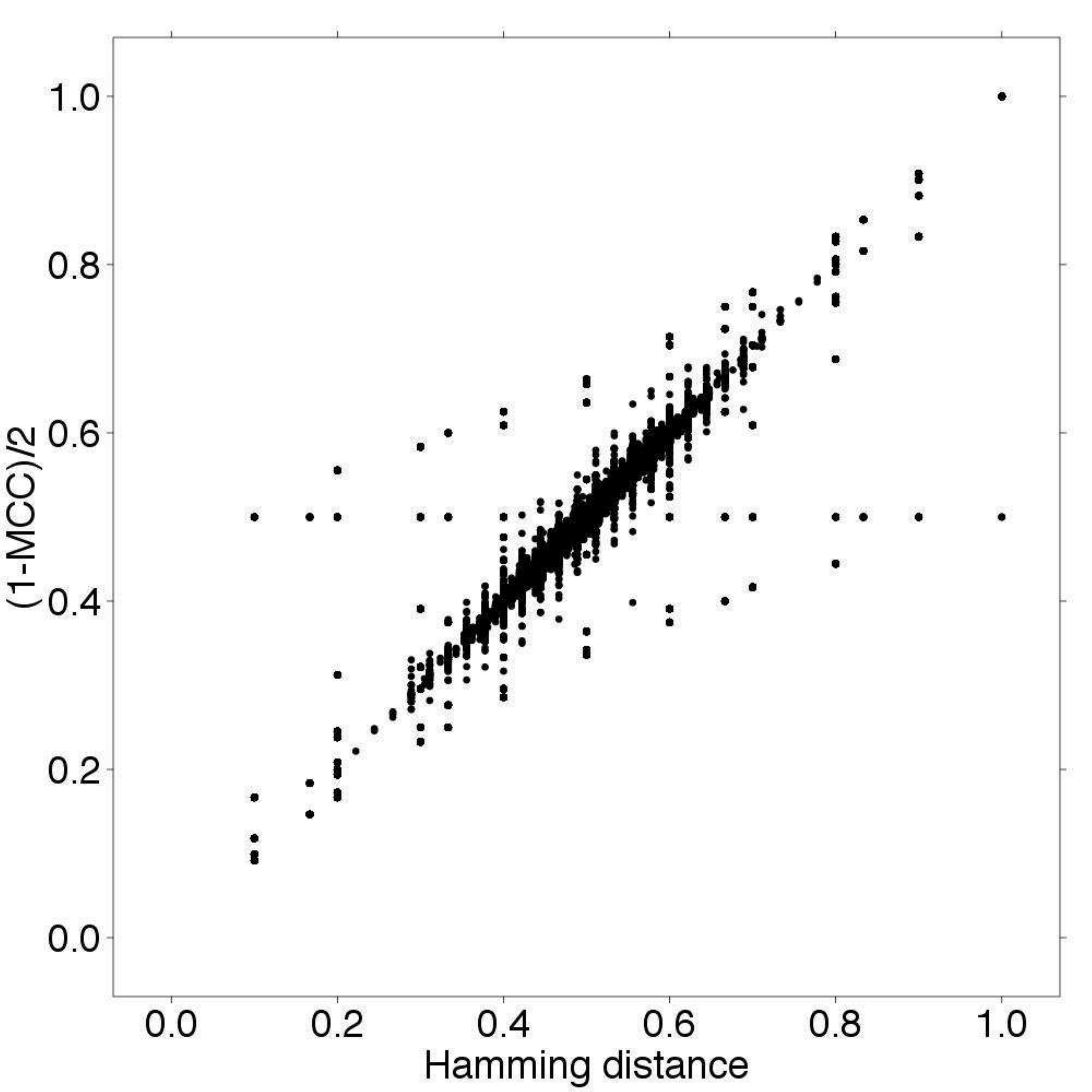} &
\includegraphics[width=0.3\textwidth]{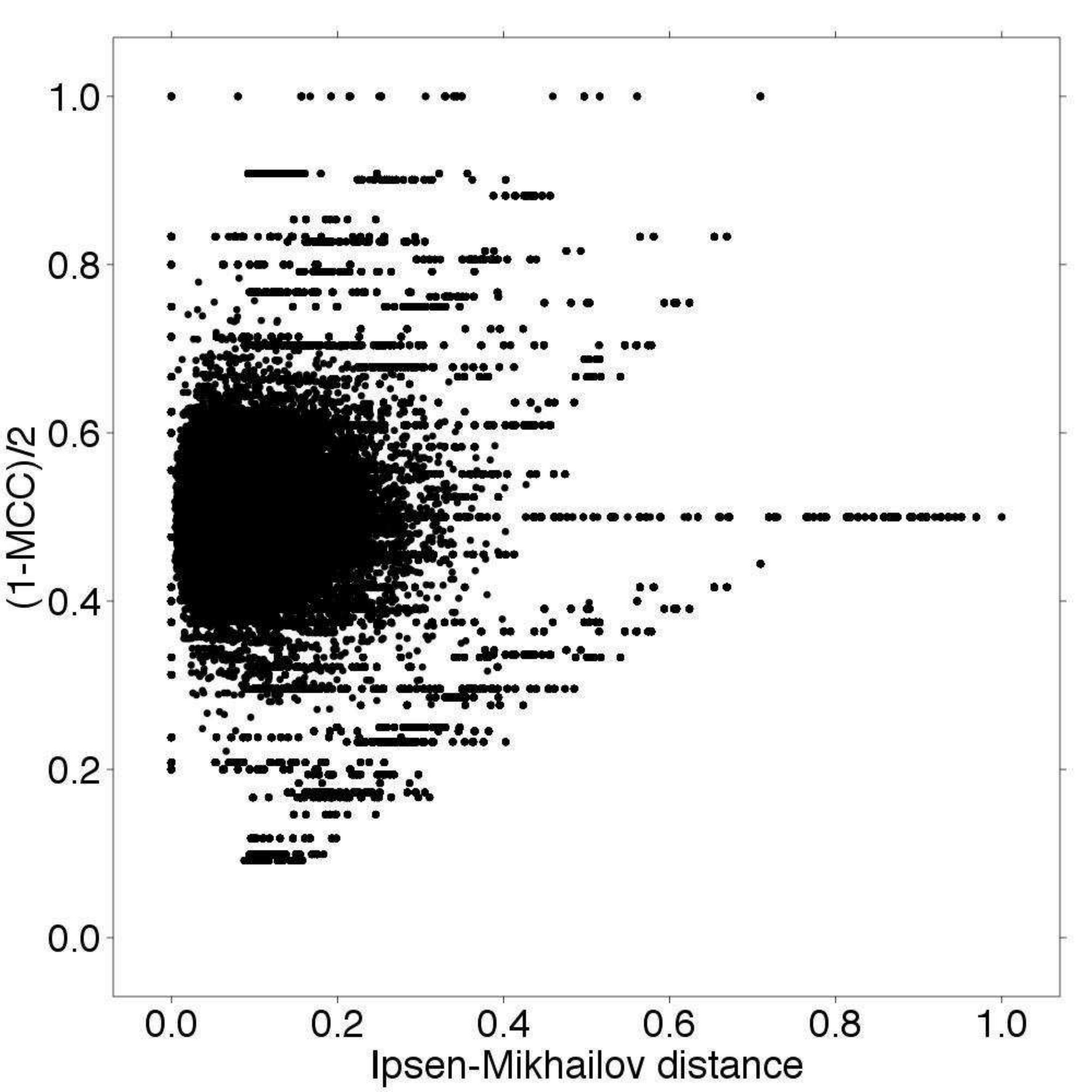} &
\includegraphics[width=0.3\textwidth]{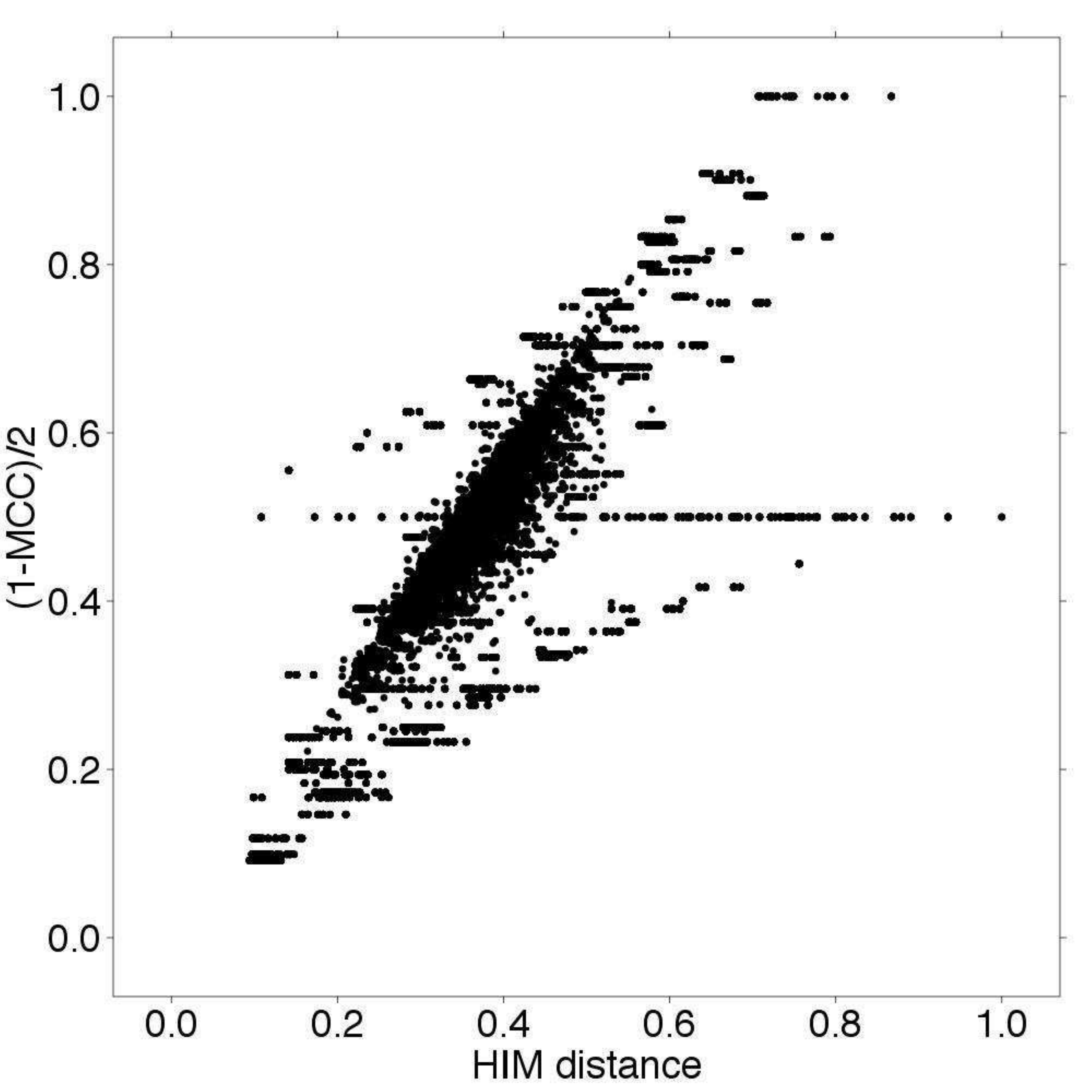} \\
(a) & (b) & (c)
\end{tabular}
\caption{Scatterplot of $\frac{1-\textrm{MCC}}{2}$ versus Hamming (a), Ipsen-Mikhailov (b) and HIM (c) distances when comparing 250,000 random pairs of networks of different size 3-100.}
\label{fig:mcc}
\end{figure}

Here we want to provide a quick comparison of $\textrm{MCC}$ and $\textrm{HIM}$ distances in a few cases.
First of all, some considerations on the extreme cases: 
\begin{itemize}
\item $\textrm{HIM}(G,H)=1$ only for $\{G,H\}=\{\mathcal{E}_N,\mathcal{F}_N\}$, which has $\textrm{MCC}=0$.
\item $\textrm{HIM}(G,H)=0$ only for $G=H$; in this case, $\textrm{MCC}(G,H)=0$ for $G=H\in \{\mathcal{E}_N,\mathcal{F}_N\}$, and $\textrm{MCC}(G,H)=1$ in all other cases.
\item $\textrm{MCC}(G,H)=1$ only for $G=H\not\in \{\mathcal{E}_N,\mathcal{F}_N\}$, and thus HIM=0.
\item The two values $\textrm{MCC}=0$ or $\textrm{MCC}=-1$ can correspond to a landscape of quite different pairs of networks, for which the $\textrm{HIM}$ distance can assume very diverse values.
\end{itemize}
To investigate the last case in the above list, we randomly generated 250,000 pairs of networks of different size, and we compared the $\textrm{MCC}$ with the H, IM and HIM distances: the corresponding scatterplots are shown in Fig.~\ref{fig:mcc}.
Since $\textrm{MCC}$ is a similarity measure, for a direct comparison we displayed it as the $[0,1]$-normalized dissimilarity measure $\frac{1-\textrm{MCC}}{2}$.

As predictable, since the confusion matrix is unaware of the network structure but it takes into account only matching and mismatching links, the $\textrm{MCC}$ is well correlated with the Hamming distance (Pearson Coefficient PC=0.92) and poorly correlated with the Ipsen-Mikhailov distance (PC=0.01), resulting in an good global correlation with the HIM distance (PC=0.79).
Nonetheless, the plots in Fig.~\ref{fig:mcc} show that the relevant variability of one measure for a given value of the other one supports the claim of a strong independency between $\textrm{MCC}$ and $\textrm{HIM}$.
Finally, as an example giving a quantitative basis to the last claim of the above list, for all the pairs with $\textrm{MCC}=0$ we obtain values of HIM ranging in $[0.11, 1]$, with median 0.37 and mean 0.39, while when $\textrm{MCC}=-1$ the range of the HIM values is $[0.71,0.86]$, with mean and median equal to 0.74.

\subsection{Dynamical networks}
\label{ssec:dynamic}
In what follows we show the evolution of the Hamming, the Ipsen-Mikhailov and the HIM distance during the evolution of the following dynamical processes $\mathbb{P}(i)$ moving through consecutive steps:
\begin{itemize}
\item{Random Addition} $\mathbb{P}_\textsc{ra}(i+1)$ is obtained from $\mathbb{P}_\textsc{ra}(i)$ by randomly adding a not already present link. 
\item{Random Removal} $\mathbb{P}_\textsc{rr}(i+1)$ is obtained from $\mathbb{P}_\textsc{rr}(i)$ by randomly removing an existing link.
\item{Sequential Addition} $\mathbb{P}_\textsc{sa}(i+1)$ is obtained from $\mathbb{P}_\textsc{sa}(i)$ by adding a new link in the same row and in the next available column of the last added link, if possible, or in the following row, starting from the first available column. The whole process starts from the first available row with the smallest index. As an example, if $\mathbb{P}_\textsc{sa}(0)=\mathcal{E}_5$, then the process evolves inserting ones in the adjacency matrix following the sequence $1\to 2\to 3\to\cdots 10$ in $\left( \begin{smallmatrix} & 1 & 2 & 3 & 4 \\ & & 5 & 6 & 7 \\ & & & 8 & 9 \\ & & & & 10 \end{smallmatrix} \right)$.
\item{Sequential Removal} $\mathbb{P}_\textsc{sr}$: as in $\mathbb{P}_\textsc{sa}$, but removing one link at each step.
\item{Highest Degree Addition} $\mathbb{P}_\textsc{hda}(i+1)$ is obtained from $\mathbb{P}_\textsc{hda}(i)$ by adding a previously not existing link connecting the node with the highest degree.
\item{Highest Degree Removal} $\mathbb{P}_\textsc{hdr}(i+1)$ is obtained from $\mathbb{P}_\textsc{hdr}(i)$ by removing an existing link connecting the node with the highest degree.
\end{itemize}
 
\begin{figure}[!b]
\begin{tabular}{cccc}
\includegraphics[width=0.25\textwidth]{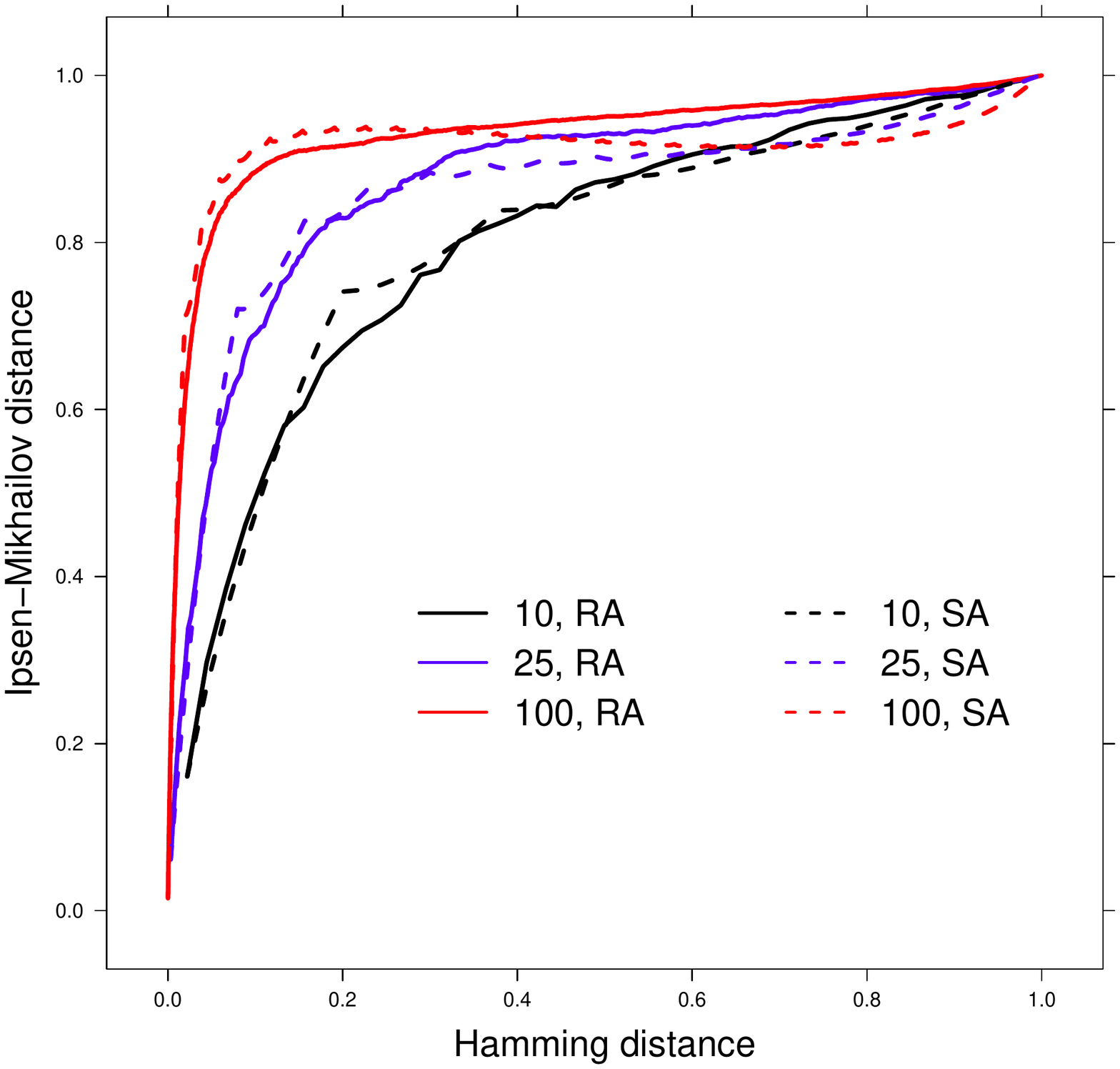} &
\includegraphics[width=0.25\textwidth]{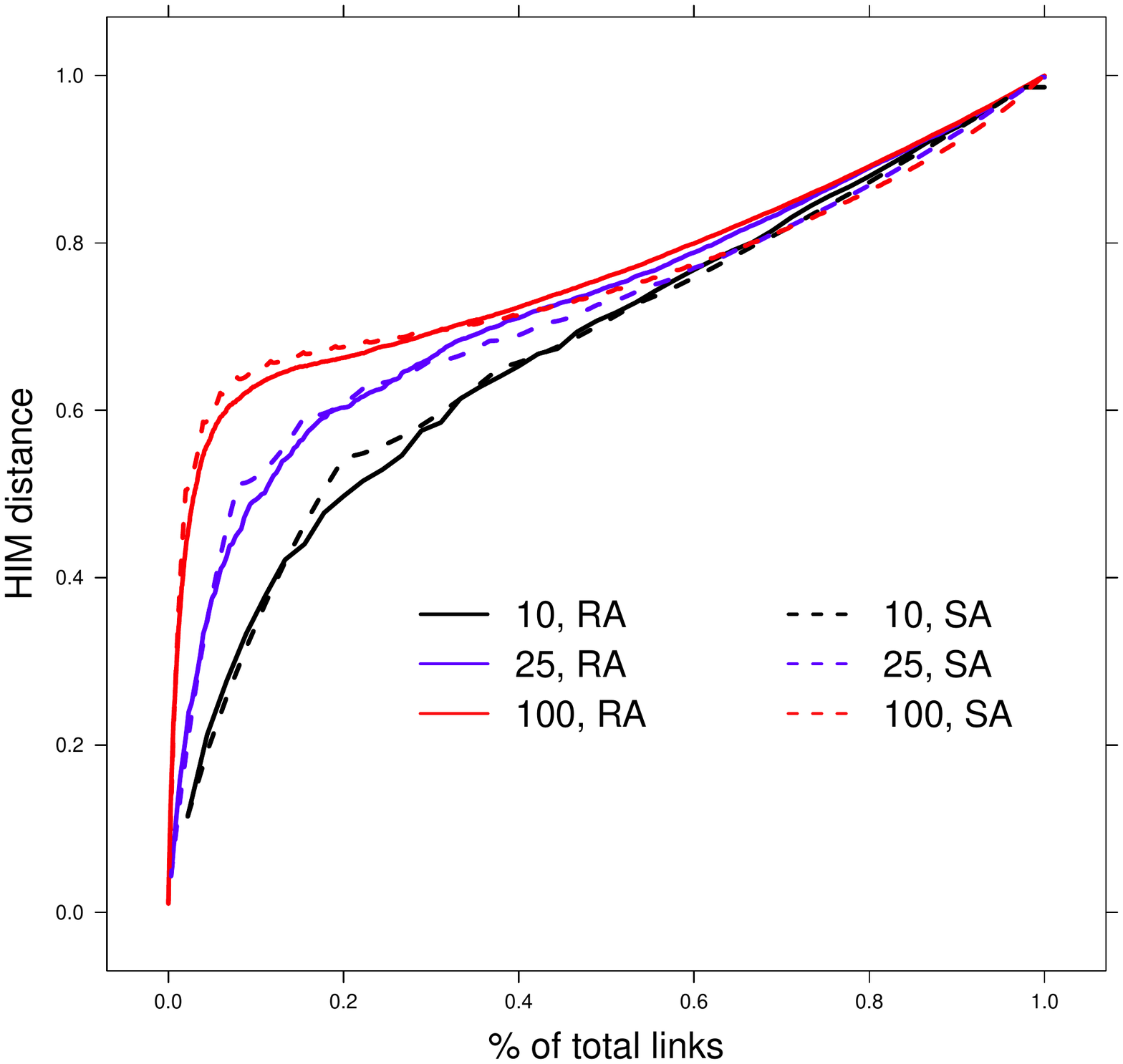} &\
\includegraphics[width=0.25\textwidth]{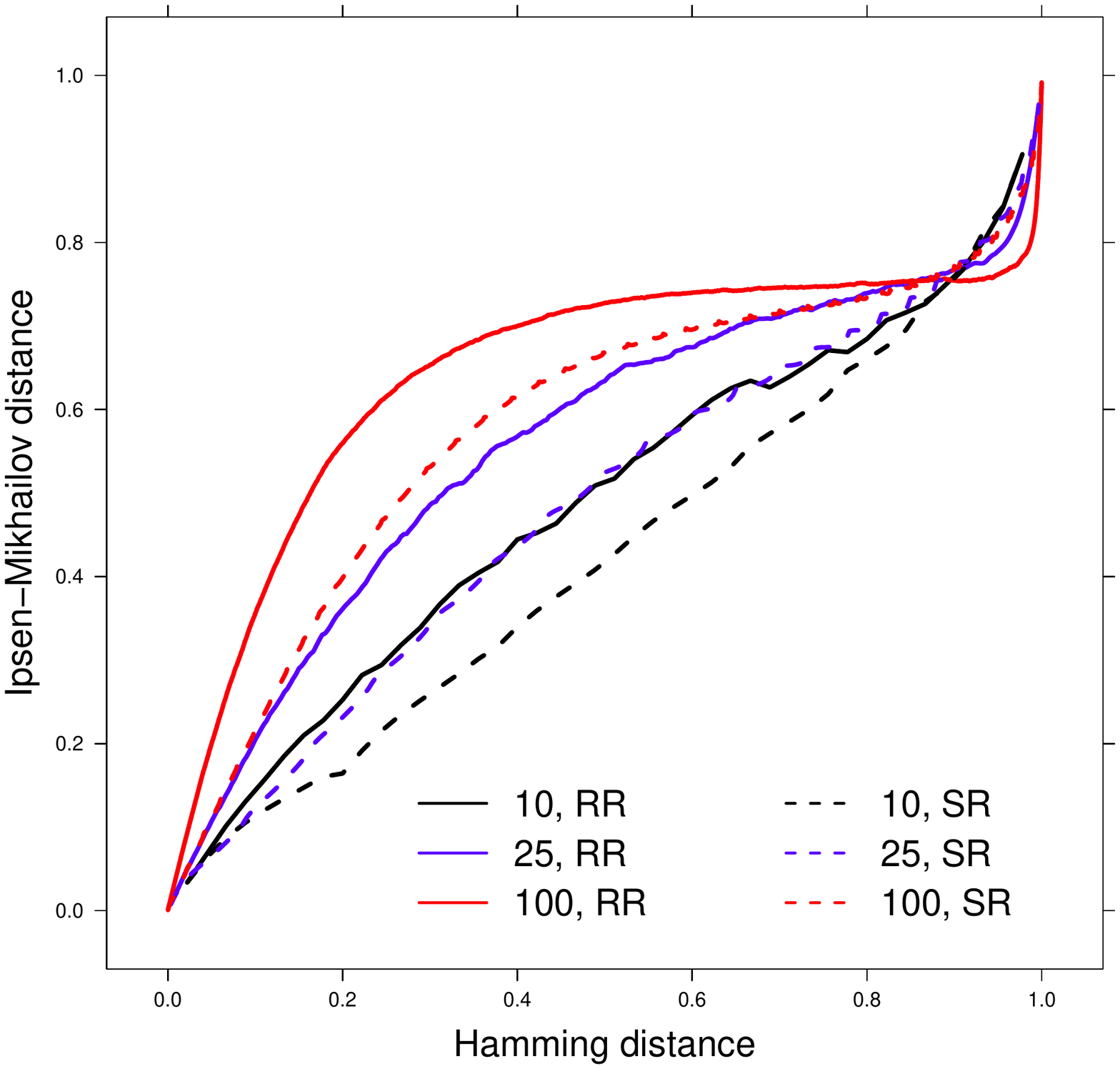} &
\includegraphics[width=0.25\textwidth]{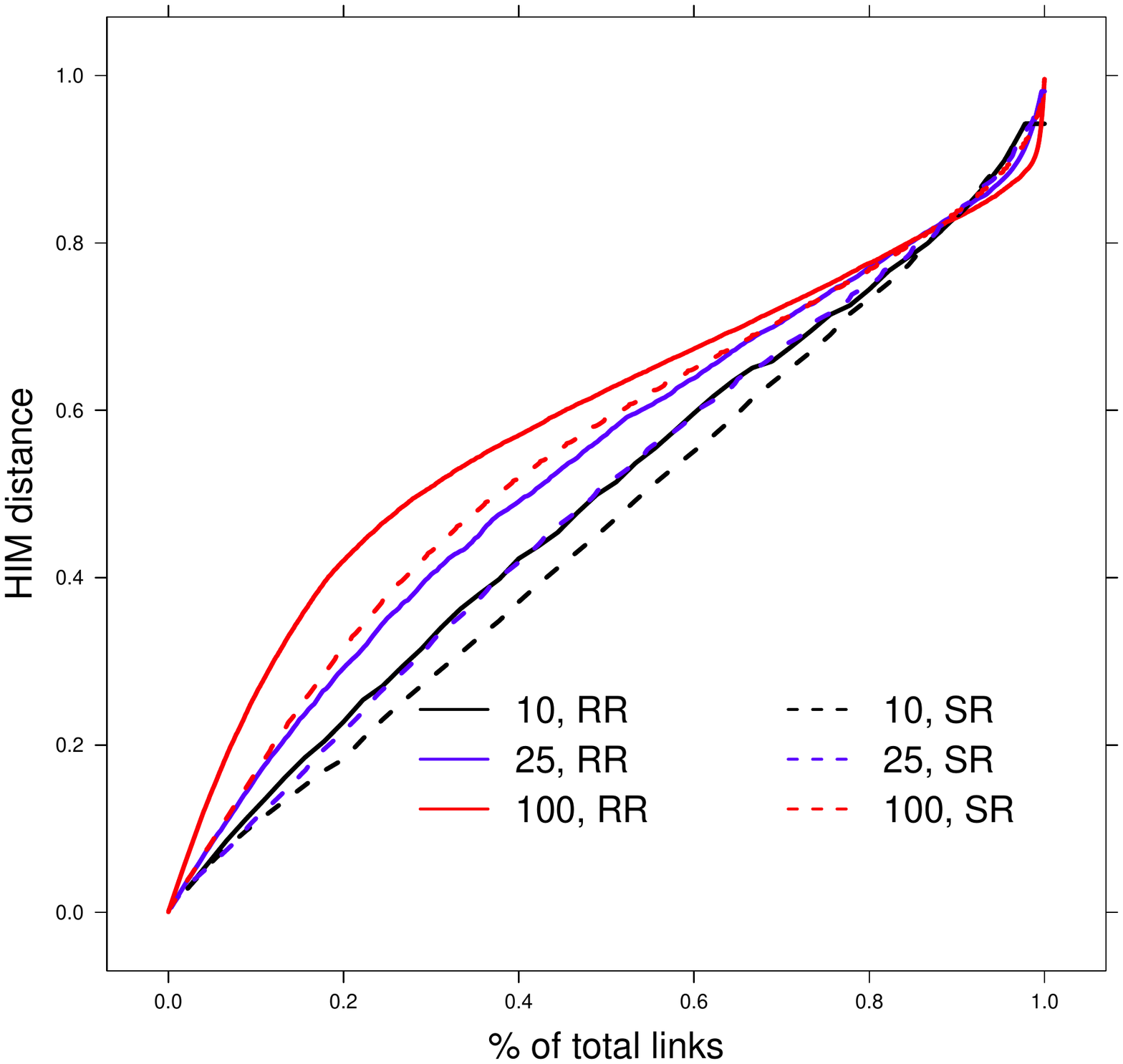} \\
(a) & (b) & (c) & (d) 
\end{tabular}
\caption{Distances between $\mathbb{P}_{\circ}(i)$ and $\mathbb{P}_{\circ}(0)$ for the processes evolving from the empty network to the clique (a and b) or vice versa (c and d), in the Hamming/Ipsen-Mikhailov space (a and c) or plot of the HIM distance as a function of the ratio of added/removed links (b and d), for $N=$ 10 (black), 25 (blue), 100 (red) nodes. Solid lines denote the average of distances for 100 runs of random evolution, while dashed lines denote the sequential processes $\mathbb{P}_\textsc{sa}$ and $\mathbb{P}_\textsc{sr}$. In all cases, the process evolves from the left-bottom corner to the right-top corner.}
\label{fig:process}
\end{figure}

As a first example, consider the processes $\mathbb{P}_\textsc{ra}$ and $\mathbb{P}_\textsc{sa}$ with the empty graph as starting network $\mathbb{P}_\textsc{ra}(0)=\mathbb{P}_\textsc{sa}(0)=\mathcal{E}_N$.
They both end at the N-nodes clique after $N_{\max}=\frac{N(N-1)}{2}$ steps: $\mathbb{P}_\textsc{ra}(N_{\max})=\mathbb{P}_\textsc{sa}(N_{\max})=\mathcal{F}_N$.
The corresponding inverse processes $\mathbb{P}_\textsc{rr}$ and $\mathbb{P}_\textsc{sr}$ evolve in the opposite direction: $\mathbb{P}_\textsc{rr}(0)=\mathbb{P}_\textsc{sr}(0)=\mathcal{F}_N$ and $\mathbb{P}_\textsc{rr}(N_{\max})=\mathbb{P}_\textsc{sr}(N_{\max})=\mathcal{E}_N$.
In Fig.~\ref{fig:process} we show the curves of $d(\mathbb{P}_{\circ}(i),\mathbb{P}_{\circ}(0))$ for $d$=H, IM and HIM in the cases $N=$10, 25 and 100 nodes. 

For the representation of the curves, we use two different spaces: the already introduced Hamming/Ipsen-Mikhailov space, with the metric H on the $x$ axis and the metric IM on the $y$ axis, and the Fraction-of-nodes/HIM space, with the ratio between the number of newly added or removed links over the total number $N_{\max}$ of links on the $x$ axis and the HIM distance on the $y$ axis; in this representation, the $i$-th step $\mathbb{P}_{\circ}(i)$ of a process has coordinates $\left(\frac{2i}{N(N-1)},\textrm{HIM}\left(\mathbb{P}_{\circ}(i),\mathbb{P}_{\circ}(0)\right)\right)$.
In all cases, since one edge is removed or added at each step, in both spaces the evolution of the processes proceeds from left to right in both graphs, and the trend of the curves representing the distances of the same process in the Hamming/Ipsen-Mikhailov space and in the Fraction-of-nodes/HIM space are similar, varying only for a scaling factor.

\begin{figure}[!t]
\begin{tabular}{cccc}
\includegraphics[width=0.25\textwidth]{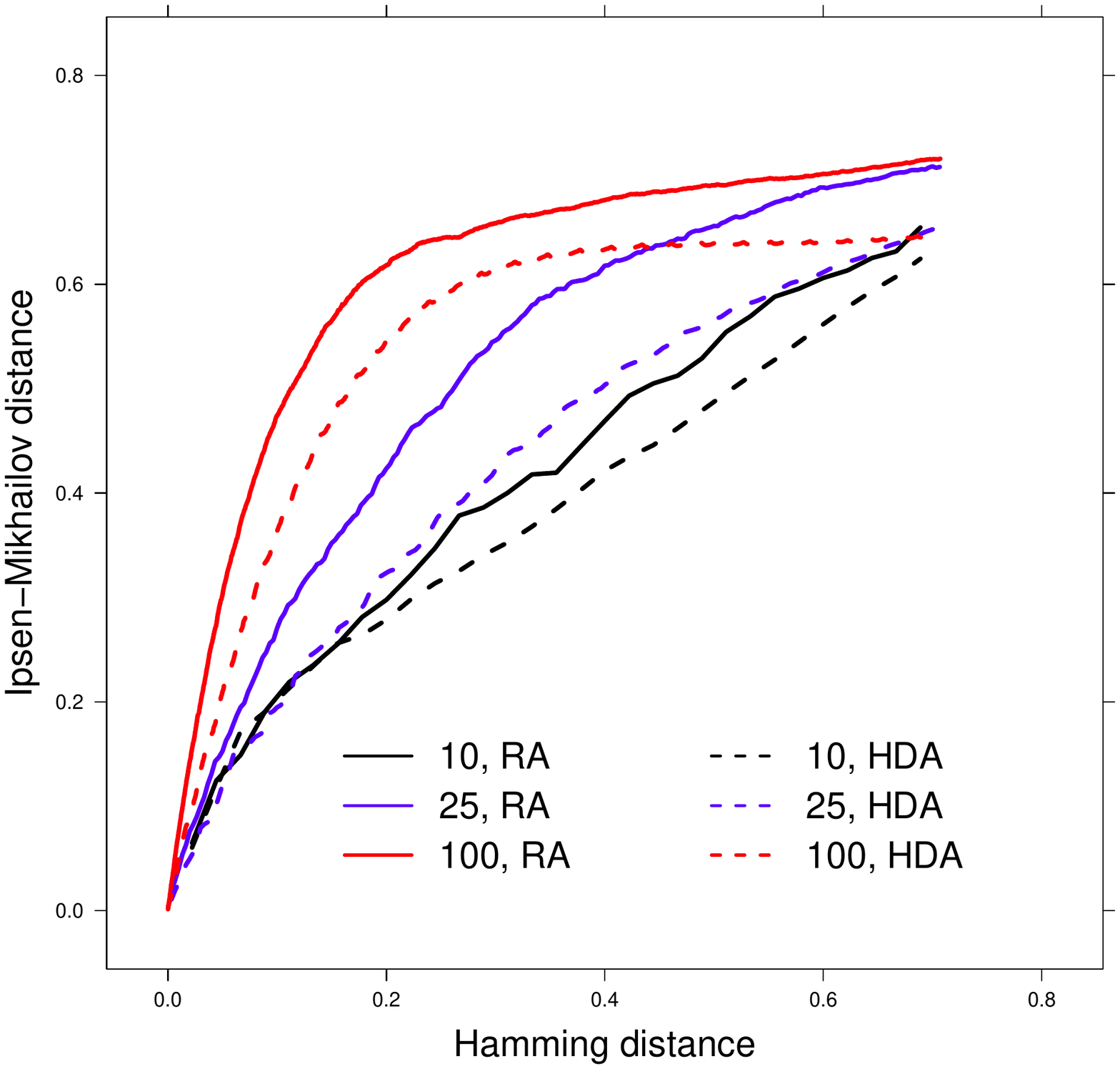} &
\includegraphics[width=0.25\textwidth]{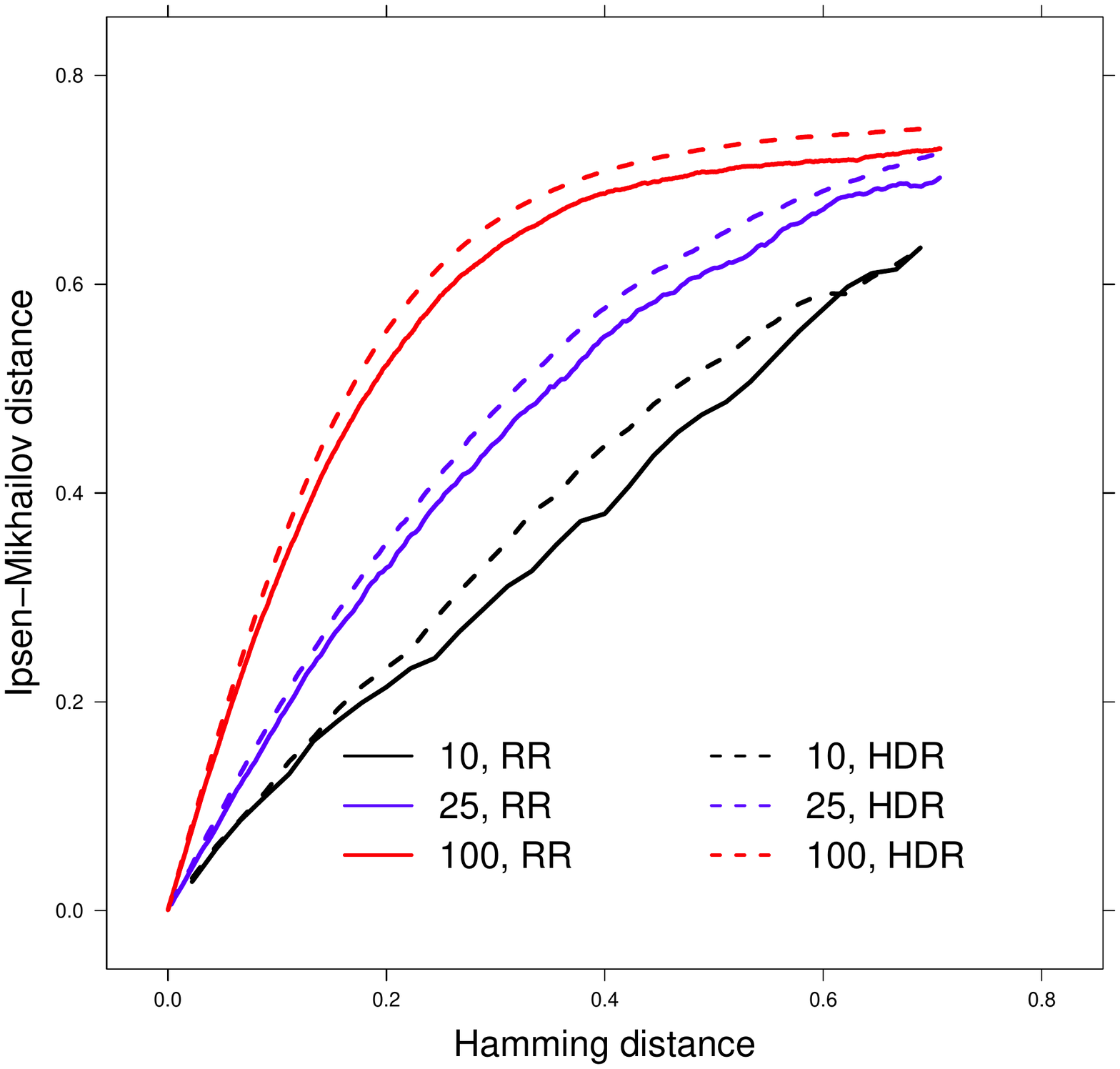} &
\includegraphics[width=0.25\textwidth]{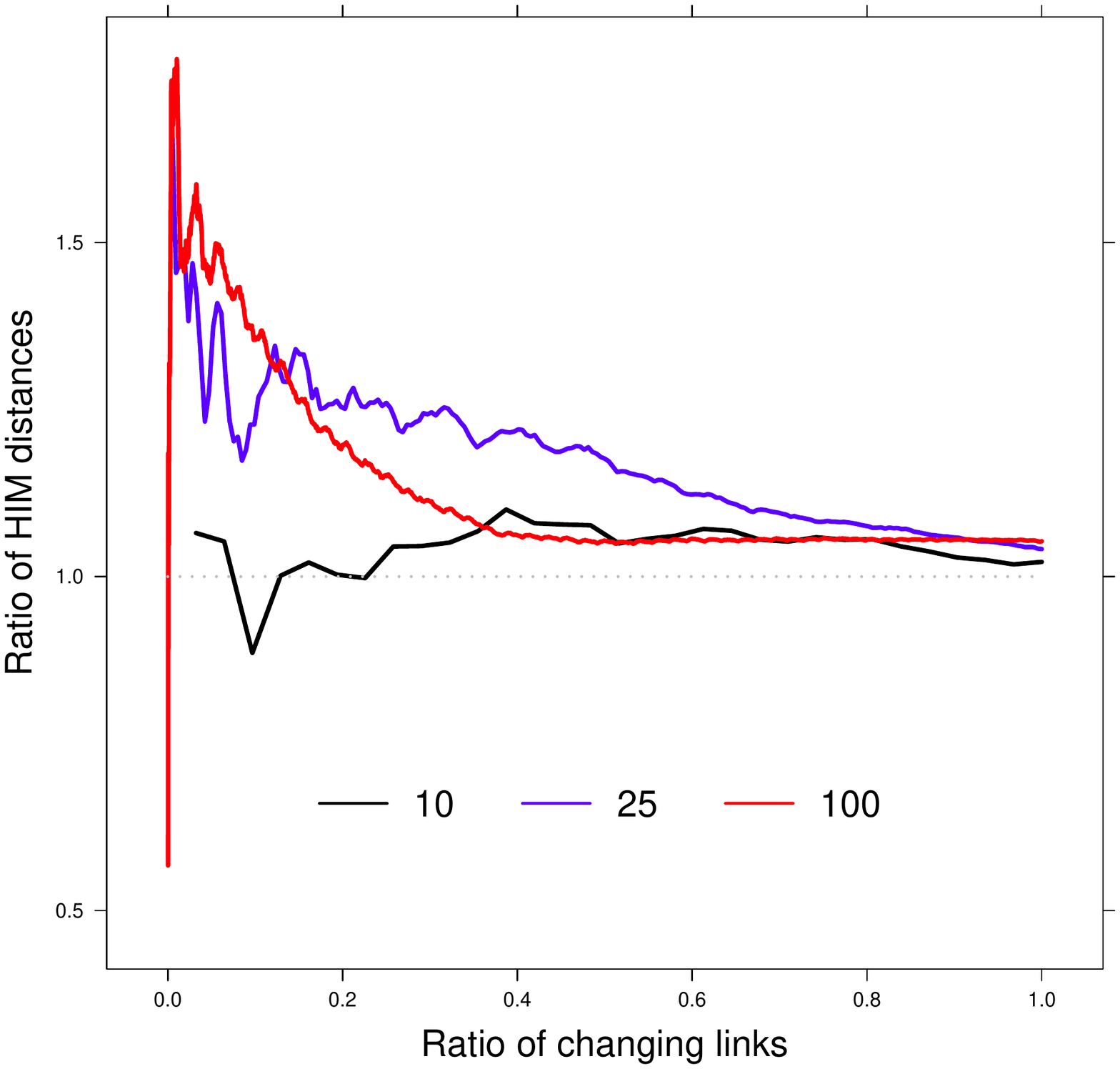} &
\includegraphics[width=0.25\textwidth]{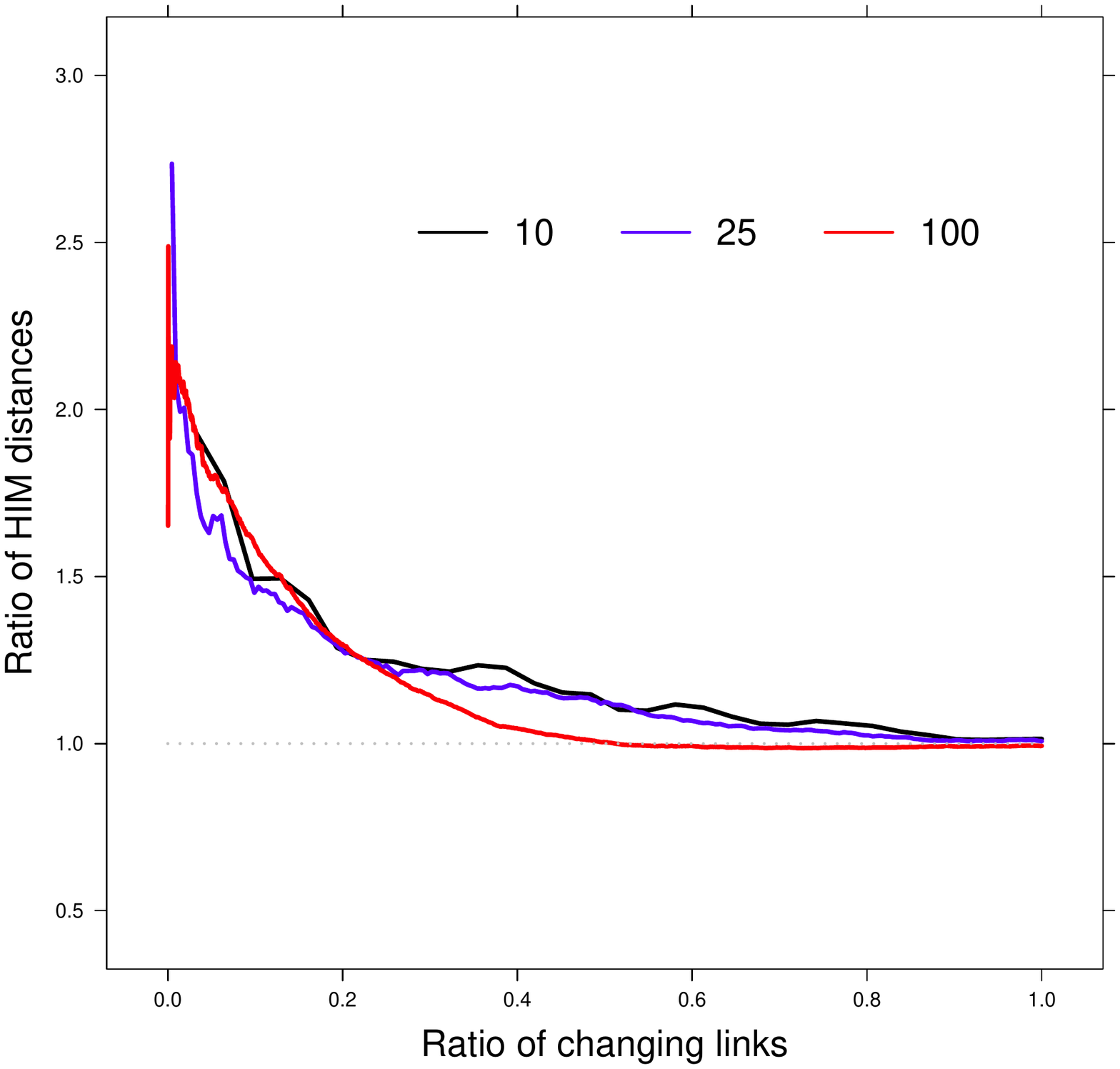} \\
(a) & (b) & (c) & (d)  
\end{tabular}
\caption{Plot of H, IM and HIM distances between $\mathbb{P}_{\circ}(i)$ and $\mathbb{P}_{\circ}(0)$ for $\circ=\textsc{ra},\textsc{hda}$ (a) and $\circ=\textsc{rr},\textsc{hdr}$ (b), for $N=$ 10 (black), 25 (blue), 100 (red) nodes. Solid lines denote the average of distances for 100 runs of $\mathbb{P}_\textsc{ra}$ and $\mathbb{P}_\textsc{rr}$, while dashed lines identify $\mathbb{P}_\textsc{hda}$ and $\mathbb{P}_\textsc{hdr}$. In all cases, the process evolves from the left-bottom corner to the right-top corner.
In panel (c), plot of $\frac{\textrm{HIM}(\mathbb{P}_\textsc{ra}(i),\mathbb{P}_\textsc{ra}(0))}{\textrm{HIM}(\mathbb{P}_\textsc{hda}(i),\mathbb{P}_\textsc{hda}(0))}$ as a function of $\frac{i}{N_{\max}}$ and, in panel (d), plot of $\frac{\textrm{HIM}(\mathbb{P}_\textsc{ra}(i),\mathbb{P}_\textsc{ra}(0))}{\textrm{HIM}(\mathbb{P}_\textsc{rr}(i),\mathbb{P}_\textsc{rr}(0))}$ as a function of $\frac{i}{N_{\max}}$.}
\label{fig:poisson}
\end{figure}

For the random processes $\mathbb{P}_\textsc{ra}$ and $\mathbb{P}_\textsc{rr}$ we show the means of the distances computed on 100 runs; no standard deviation or confidence intervals are plotted, because they are negligible at the scale of the plot. 
For instance, in the case $N=$25, the order of magnitude of the standard deviation for HIM at each step is $10^{-3}$, and the span of the 95\% boostrap confidence intervals is in the range of $10^{-4}$.
As a first observation, all curves are monotonically increasing and the bigger the graph, the larger the distances, but in the empty-to-clique case, where, in the second half of the process, $\mathbb{P}_\textsc{sa}$ induces distances which are smaller than $\mathbb{P}_\textsc{ra}$ and which are smaller for larger graphs.
The most interesting observation is the different shape of the curves between the empty-to-clique process and the clique-to-empty: for the same Hamming distance (or fraction of links), the corresponding Ipsen-Mikhailov (or HIM, respectively) distance is larger when the nodes are added rather than removed, because adding links quickly generates degree correlation.
Furthermore, in the empty-to-clique case, not much difference occurs between the random and the sequential process, while this difference is much wider (with the random one inducing larger distances) for the clique-to-empty case.
An analogous experiment was carried out within the family of Poissonian graphs, with Erd\"os-R\'enyi model \cite{erdos59random,erdos60evolution} $G(N,p)$.
In particular, for $N=$10, 25 and 100, let $S_N$ be a sparse network $G(N,p=0.05)$, with 2, 11 and 230 edges respectively and let $D_N$ be a dense network $G(N,p=0.9)$, with 39, 275 and 4462 edges respectively.
Consider the following four processes, of which we represent the initial $\lfloor\frac{N_{\max}}{\sqrt{2}}\rfloor=\lfloor\frac{\sqrt{2}}{2}\cdot\frac{N(N-1)}{2}\rfloor+1$ steps in Fig.~\ref{fig:poisson}:
\begin{itemize}
\item $\mathbb{P}_\textsc{ra}(i)$, for $i=0,\ldots, \lfloor\frac{\sqrt{2}}{2}\frac{N(N-1)}{2}\rfloor$, with $\mathbb{P}_\textsc{ra}(0)=S_N$, for $N$=10, 25, 100.
\item $\mathbb{P}_\textsc{rr}(i)$, for $i=0,\ldots, \lfloor\frac{\sqrt{2}}{2}\frac{N(N-1)}{2}\rfloor$, with $\mathbb{P}_\textsc{rr}(0)=D_N$, for $N$=10, 25, 100.
\item $\mathbb{P}_\textsc{hda}(i)$, for $i=0,\ldots, \lfloor\frac{\sqrt{2}}{2}\frac{N(N-1)}{2}\rfloor$, with $\mathbb{P}_\textsc{hda}(0)=S_N$, for $N$=10, 25, 100.
\item $\mathbb{P}_\textsc{hdr}(i)$, for $i=0,\ldots, \lfloor\frac{\sqrt{2}}{2}\frac{N(N-1)}{2}\rfloor$, with $\mathbb{P}_\textsc{hdr}(0)=D_N$, for $N$=10, 25, 100.
\end{itemize}
In this case, too, results on the random processes are averaged over 100 runs, with negligible confidence intervals.
To better highlight the differences of the resulting distances in the various processes, in Fig.~\ref{fig:poisson} (c) and (d) we show the ratio of some pairs of HIM distances as a function of the removed/added links.
In particular, in subfigure (c), for each step $i$, we show the quotient of HIM distances for $\mathbb{P}_\textsc{ra}$ over $\mathbb{P}_\textsc{hda}$, in the three cases $N=$10, 25 and 100.
The three curves show that HIM distances for $\mathbb{P}_\textsc{ra}$ are larger than the HIM distances for $\mathbb{P}_\textsc{hda}$ for $N$=25 and 100, and their difference is higher in the first steps of the process $i< 0.3 N_{\max}$, while they tend to get closer as far as the processes evolve. 
In the other cases $\mathbb{P}_\textsc{rr}$ and $\mathbb{P}_\textsc{hdr}$ (not shown here), the differences are smaller and they converge faster to one, but in this case the process $\mathbb{P}_\textsc{hdr}$ accounts for the smaller values of HIM distances.
In the plot (c) of Fig.~\ref{fig:poisson}, we show the curves for $\frac{\textrm{HIM}(\mathbb{P}_\textsc{ra}(i),\mathbb{P}_\textsc{ra}(0))}{\textrm{HIM}(\mathbb{P}_\textsc{rr}(i),\mathbb{P}_\textsc{rr}(0))}$ as a function of $\frac{i}{N_{\max}}$.
All the three curves are monotonically decreasing and converging to one after the first stages of the processes, yielding that, for all values of $N$, adding links produces higher values of HIM distance.
In the case of the evolution targeting higher degree nodes first (not shown here), the trend is the same, only scaled down to smaller ratios.

\begin{figure}[!t]
\begin{center}
\begin{tabular}{cccc}
\includegraphics[width=0.25\textwidth]{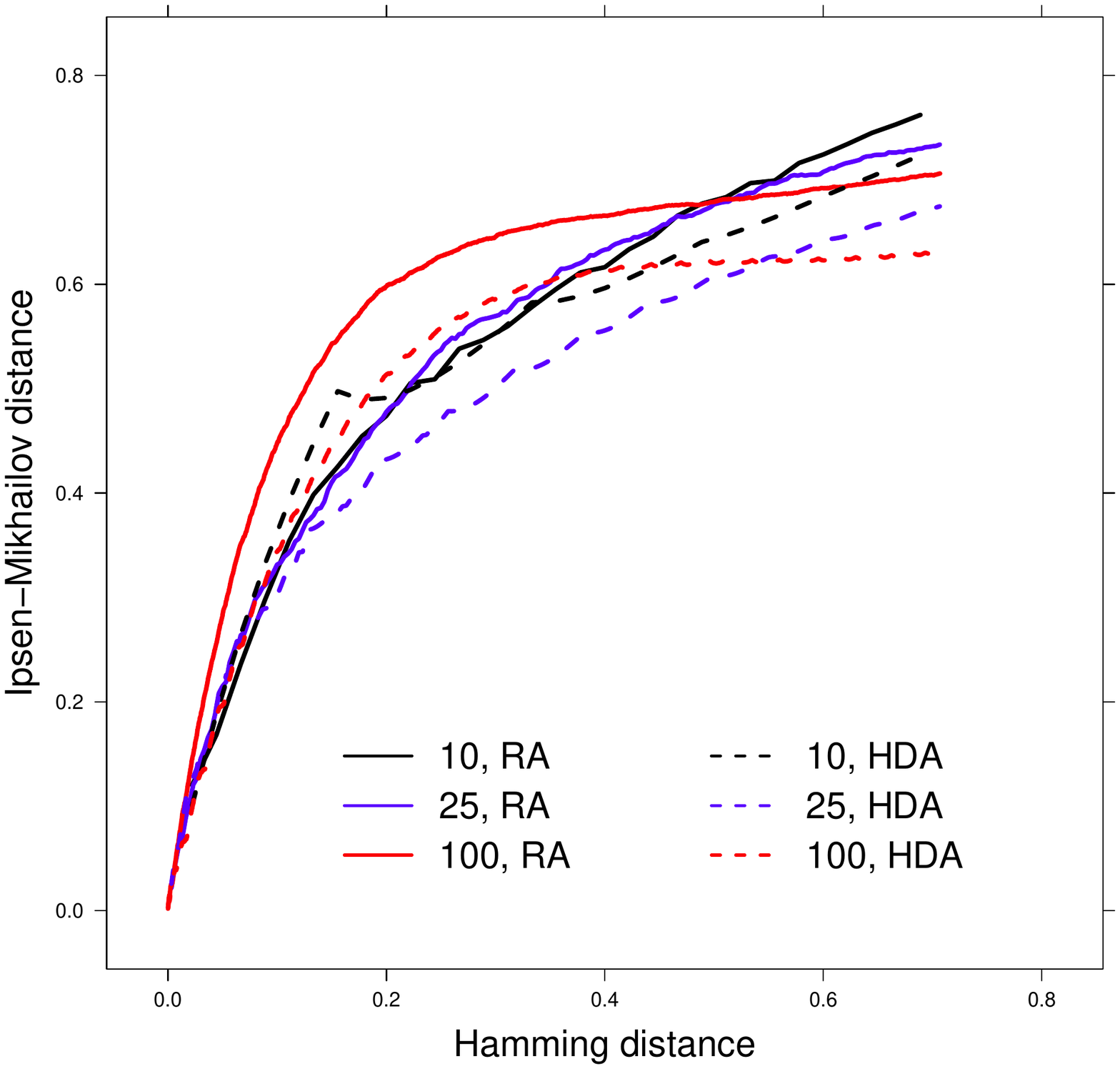} &
\includegraphics[width=0.25\textwidth]{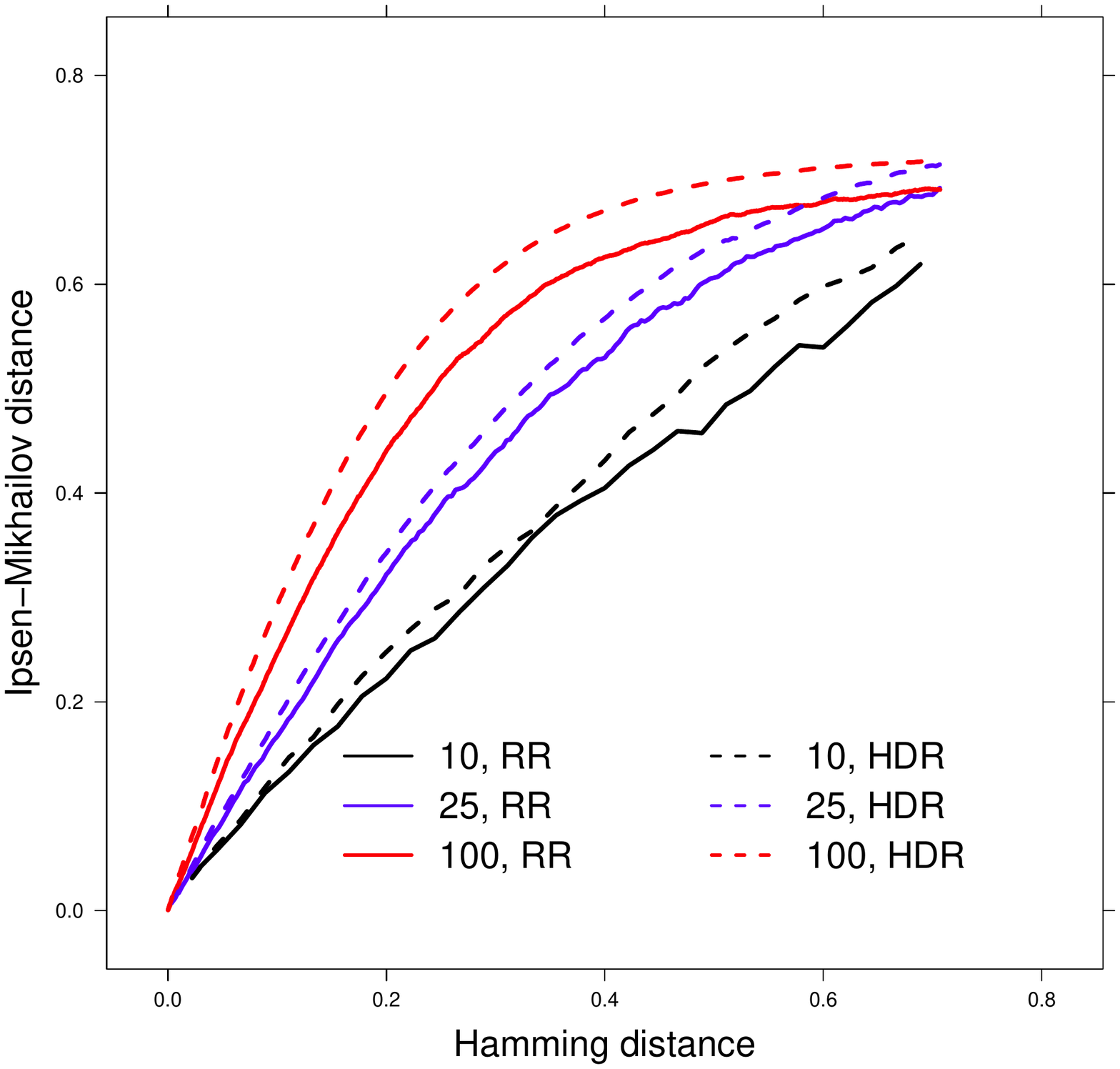} &
\includegraphics[width=0.25\textwidth]{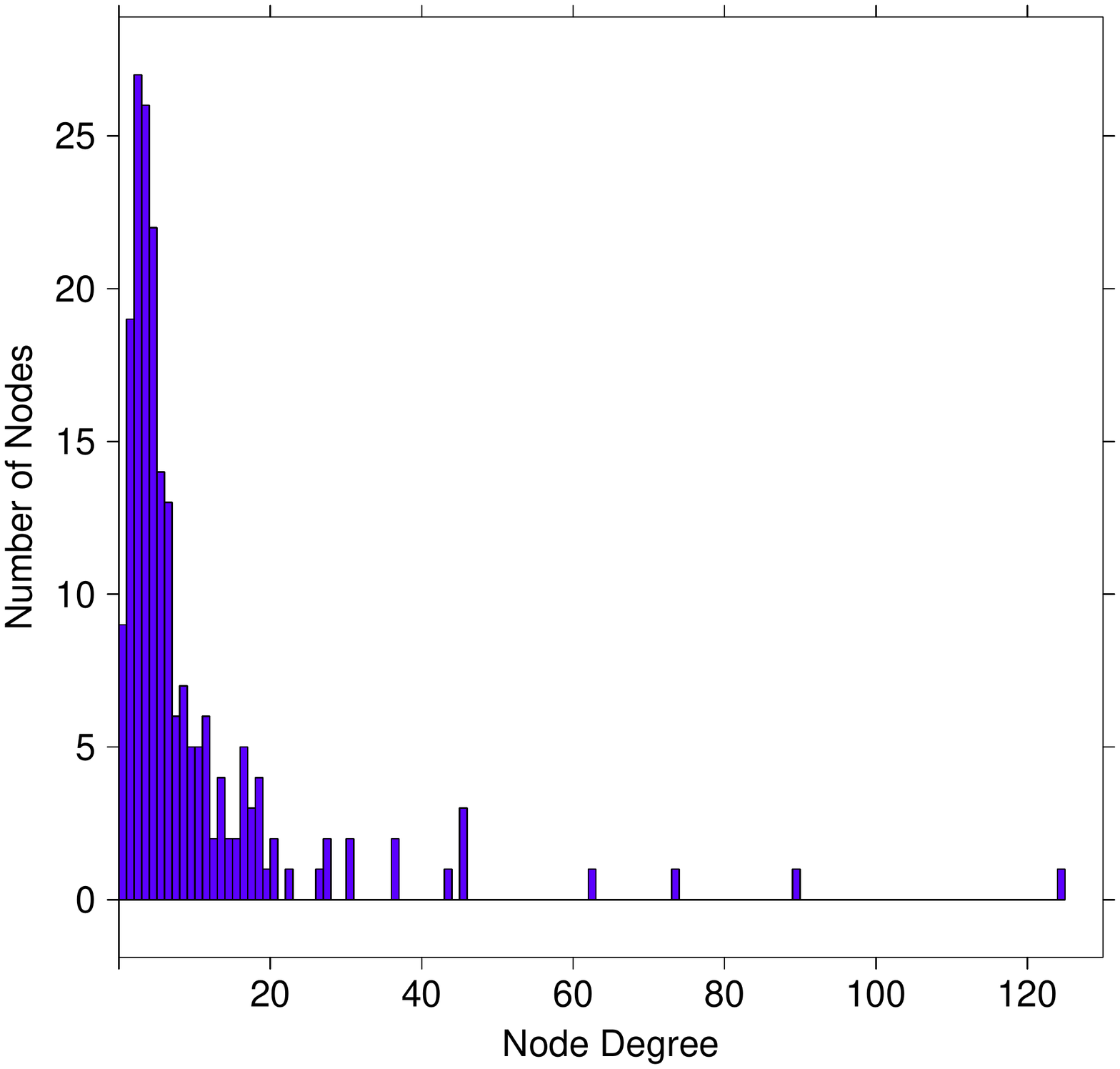} &
\includegraphics[width=0.25\textwidth]{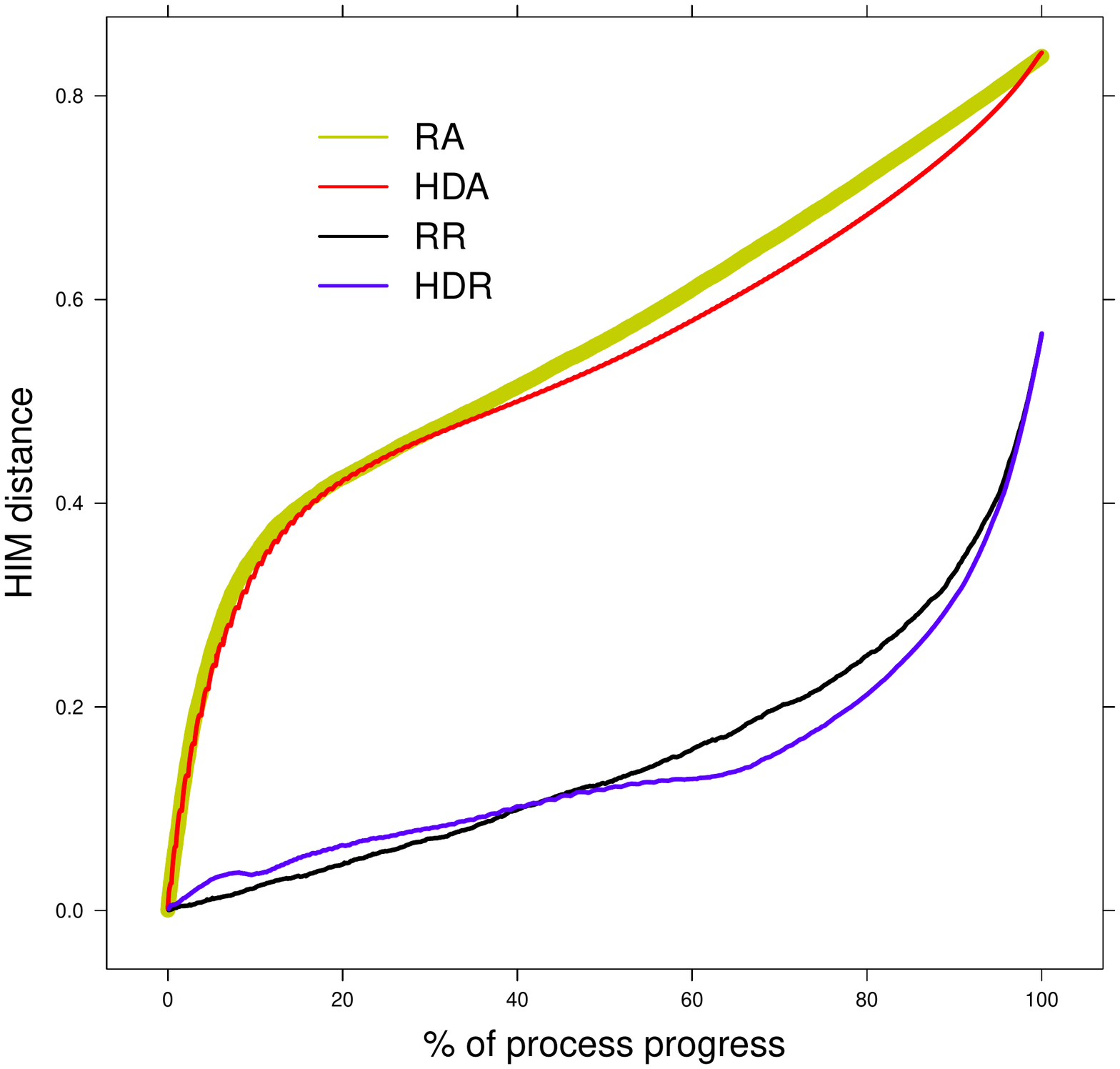} \\
(a) & (b) & (c) & (d)
\end{tabular}
\end{center}
\caption{Plot of Hamming versus Ipsen-Mikhailov distance for $\mathbb{P}_{\circ}(i)$ and $\mathbb{P}_{\circ}(0)=SB_{200}$ for $\circ=\textsc{ra},\textsc{hda}$ (a) and $\circ=\textsc{rr},\textsc{hdr}$ (b), for $N=$ 10 (black), 25 (blue), 100 (red) nodes. Solid lines denote the average of distances for 100 runs of $\mathbb{P}_\textsc{ra}$ and $\mathbb{P}_\textsc{rr}$, while dashed lines identify $\mathbb{P}_\textsc{hda}$ and $\mathbb{P}_\textsc{hdr}$. In all cases, the process evolves from the left-bottom corner to the right-top corner.
(c) Histogram of the node degrees of the network ${SB}_{200}$. (d) Hamming distances versus percentage of processes steps for $\mathbb{P}_\textsc{ra}$, $\mathbb{P}_\textsc{hda}$, $\mathbb{P}_\textsc{rr}$ and $\mathbb{P}_\textsc{hdr}$.}
\label{fig:scalefree}
\end{figure}

The final examples consider processes having scale free networks as starting graphs.
For $N$=10, 25 and 100, let ${SS}_N$ be a scale free sparse network generated following the Albert-Barabasi model \cite{barabasi99emergence}, with power law exponent 2.3 and with 9, 24 and 99 edges rispectively, and $SD_N$ a dense network with the same exponent 2.3 but with 35, 300 and 4150 edges respectively.
The same four processes of the previous case were tested for the initial $\lfloor\frac{N_{\max}}{\sqrt{2}}\rfloor=\lfloor\frac{\sqrt{2}}{2}\cdot\frac{N(N-1)}{2}\rfloor+1$ steps:
\begin{itemize}
\item $\mathbb{P}_\textsc{ra}(i)$, for $i=0,\ldots, \lfloor\frac{\sqrt{2}}{2}\frac{N(N-1)}{2}\rfloor$, with $\mathbb{P}_\textsc{ra}(0)=SS_N$, for $N$=10, 25, 100.
\item $\mathbb{P}_\textsc{rr}(i)$, for $i=0,\ldots, \lfloor\frac{\sqrt{2}}{2}\frac{N(N-1)}{2}\rfloor$, with $\mathbb{P}_\textsc{rr}(0)=SD_N$, for $N$=10, 25, 100.
\item $\mathbb{P}_\textsc{hda}(i)$, for $i=0,\ldots, \lfloor\frac{\sqrt{2}}{2}\frac{N(N-1)}{2}\rfloor$, with $\mathbb{P}_\textsc{hda}(0)=SS_N$, for $N$=10, 25, 100.
\item $\mathbb{P}_\textsc{hdr}(i)$, for $i=0,\ldots, \lfloor\frac{\sqrt{2}}{2}\frac{N(N-1)}{2}\rfloor$, with $\mathbb{P}_\textsc{hdr}(0)=SD_N$, for $N$=10, 25, 100.
\end{itemize}
The corresponding curves are plotted in Fig.~\ref{fig:scalefree} (a) and (b).
We recall here that scalefree networks are not invariant for percolation, \textit{i.e.}, they do not remain scalefree when links are randomly removed or added.
However, the evolution of the processes is not very different from the Erd{\"o}s-R{\'e}nyi case, especially for the processes removing links as shown in Fig.~\ref{fig:scalefree}, panel (b).
Some differences emerge for the processes adding links, and a few peculiarities that are also present in the Poissonian case here become more evident.
In particular, for all $N$ and for both $\mathbb{P}_\textsc{ra}$ and $\mathbb{P}_\textsc{hda}$ the derivative of the curves are larger than those in panel (b), and it is not true anymore that the larger the number of nodes, the larger the distances.
For instance, in the case $N=100$, both the processes quickly modify the network structure, resulting in a fast increment of the Ipsen-Mikhailov distance for $i<0.2N_{\max}$, while later the curves grow at a much smaller rate.
To better study this behaviour, a larger starting network $SB_{200}$ was generated following the scale free model in \cite{goh01universal}, with 200 nodes and 1000 edges, power law exponent 2.001 and degree distribution as in the histogram of Fig.~\ref{fig:scalefree}, panel (c).
The following processes were started from $SB_{200}$, and they were carried on until they reach either the empty network or the clique:
\begin{itemize}
\item $\mathbb{P}_\textsc{ra}(i)$, with $\mathbb{P}_\textsc{ra}(0)={SB}_{200}$ and $\mathbb{P}_\textsc{ra}(18900)=\mathcal{F}_{200}$.
\item $\mathbb{P}_\textsc{hda}(i)$, with $\mathbb{P}_\textsc{hda}(0)={SB}_{200}$ and $\mathbb{P}_\textsc{hda}(18900)=\mathcal{F}_{200}$.
\item $\mathbb{P}_\textsc{rr}(i)$, with $\mathbb{P}_\textsc{rr}(0)={SB}_{200}$ and $\mathbb{P}_\textsc{rr}(1000)=\mathcal{E}_{200}$.
\item $\mathbb{P}_\textsc{hdr}(i)$, with $\mathbb{P}_\textsc{hdr}(0)={SB}_{200}$ and $\mathbb{P}_\textsc{hdr}(1000)=\mathcal{E}_{200}$.
\end{itemize}
The curves corresponding to HIM distances from $SB_{200}$ in the four aforementioned processes are plotted in Fig.~\ref{fig:scalefree}, panel (d), versus the percentage of progress of the process, \textit{i.e.}, $100\cdot \frac{i}{N_{\circ}}$, with $0\leq i \leq N_{\circ}$ and $N_{\textsc{ra}}=N_{\textsc{hda}}=18900$, $N_{\textsc{rr}}=N_{\textsc{hdr}}=1000$. 
The HIM distance for the processes $\mathbb{P}_\textsc{rr}$ and $\mathbb{P}_\textsc{hdr}$ are monotonically and similarly increasing when evolving from $SB_{200}$ to $\mathcal{E}_{200}$, slower at the beginning and much faster in the last steps of the process.
The two other processes instead show the same effect previously noted: $\textrm{HIM}(\mathbb{P}_\textsc{ra}(i),SB_{200})$ and $\textrm{HIM}(\mathbb{P}_\textsc{hdaa}(i),SB_{200})$ change rapidly in the initial 10\% of the processes, yielding a fast increase in the Ipsen-Mikhailov distance, due to the quick modification in the network structure. After this initial period, the growth of both curves proceed with a smaller derivatives until they reach their maximum at the end of the process.

\subsection{Graph families}
\label{ssec:families}
\begin{figure}[!t]
\begin{center}
\begin{tabular}{ccccc}
\raisebox{2cm}{\textcolor{blue}{BA}} & 
\includegraphics[height=5cm]{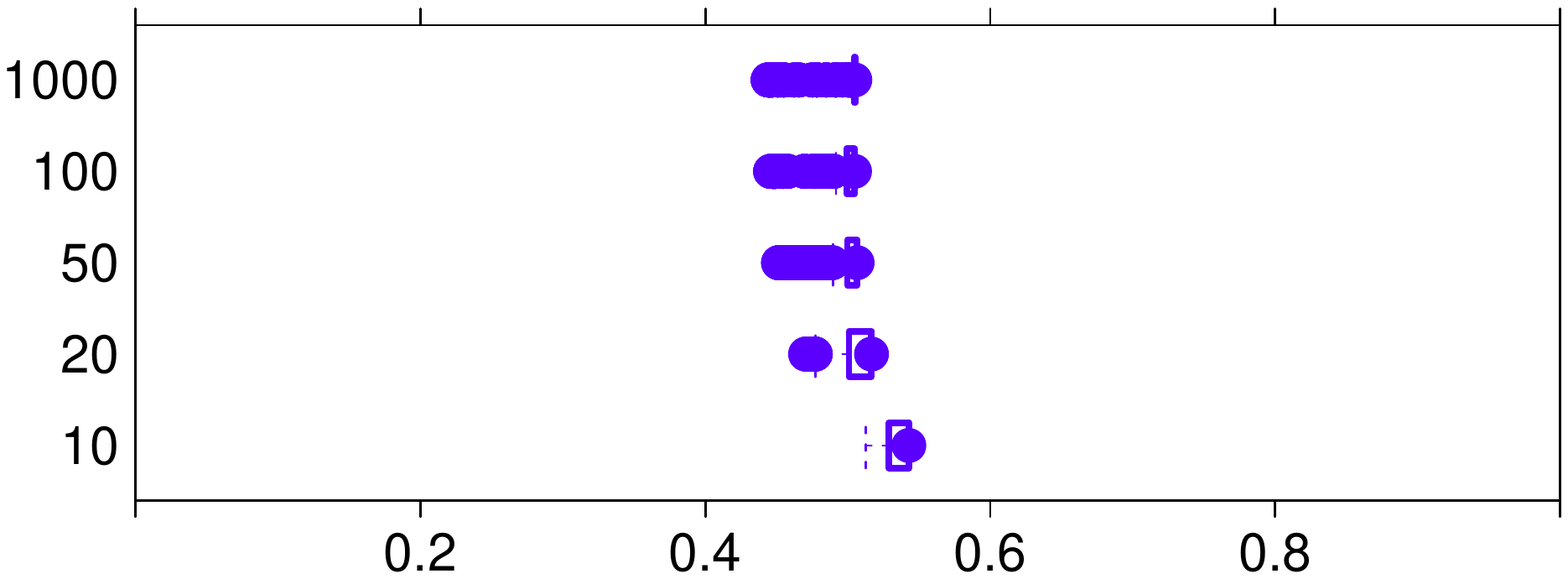} &
\includegraphics[height=5cm]{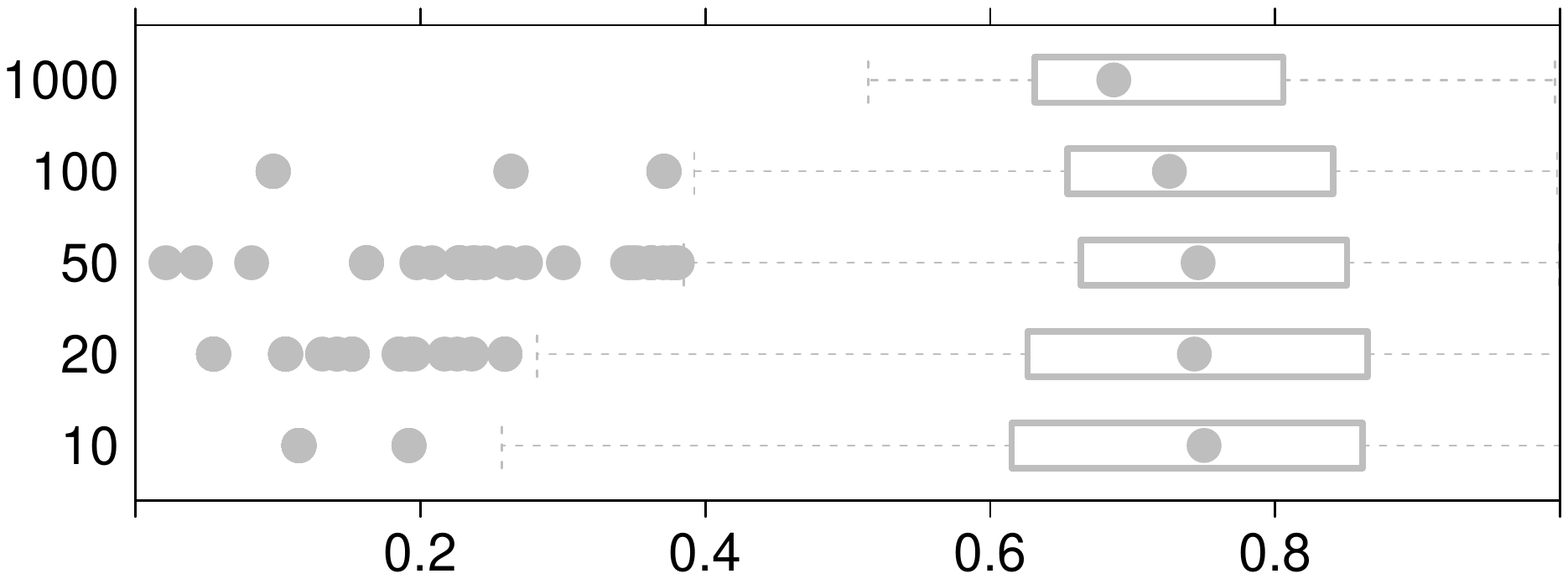} &
\raisebox{2cm}{\textcolor{gray}{PL}} \\
\raisebox{2cm}{\textcolor{green}{ER}} &
\includegraphics[height=5cm]{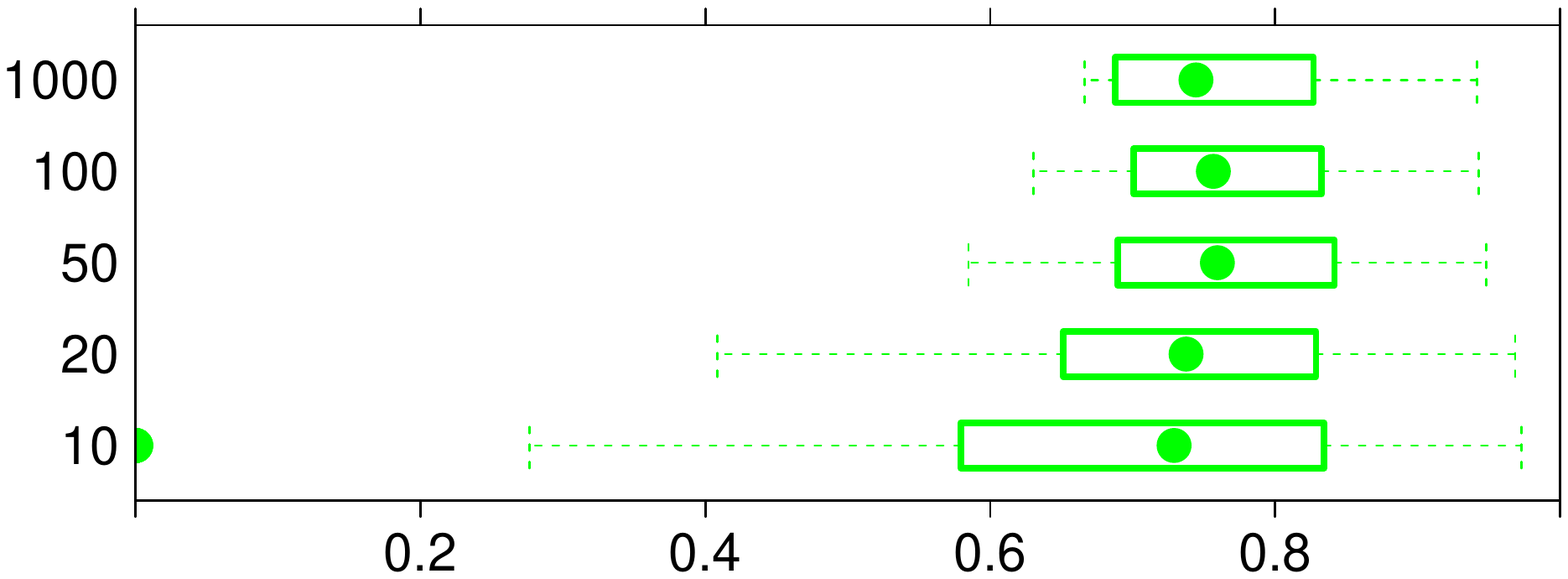} & 
\includegraphics[height=5cm]{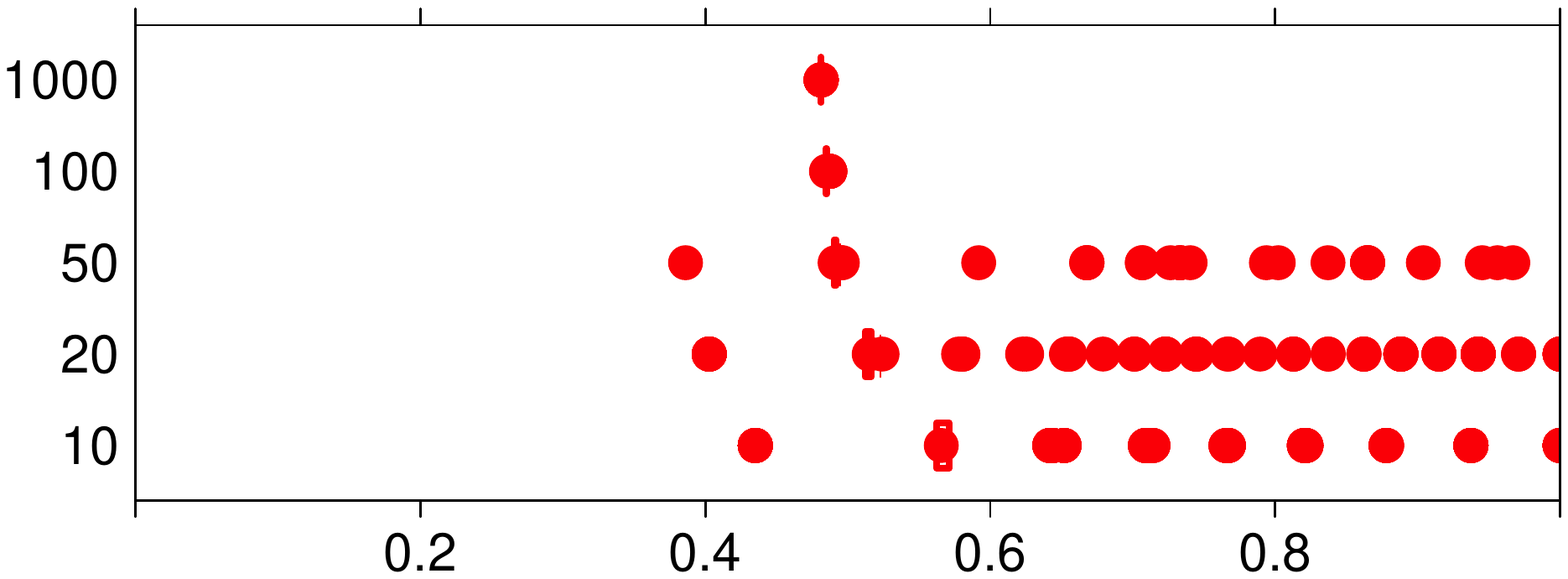} &
\raisebox{2cm}{\textcolor{red}{KR}} \\
\raisebox{2cm}{\textcolor{orange}{WS}} &
\includegraphics[height=5cm]{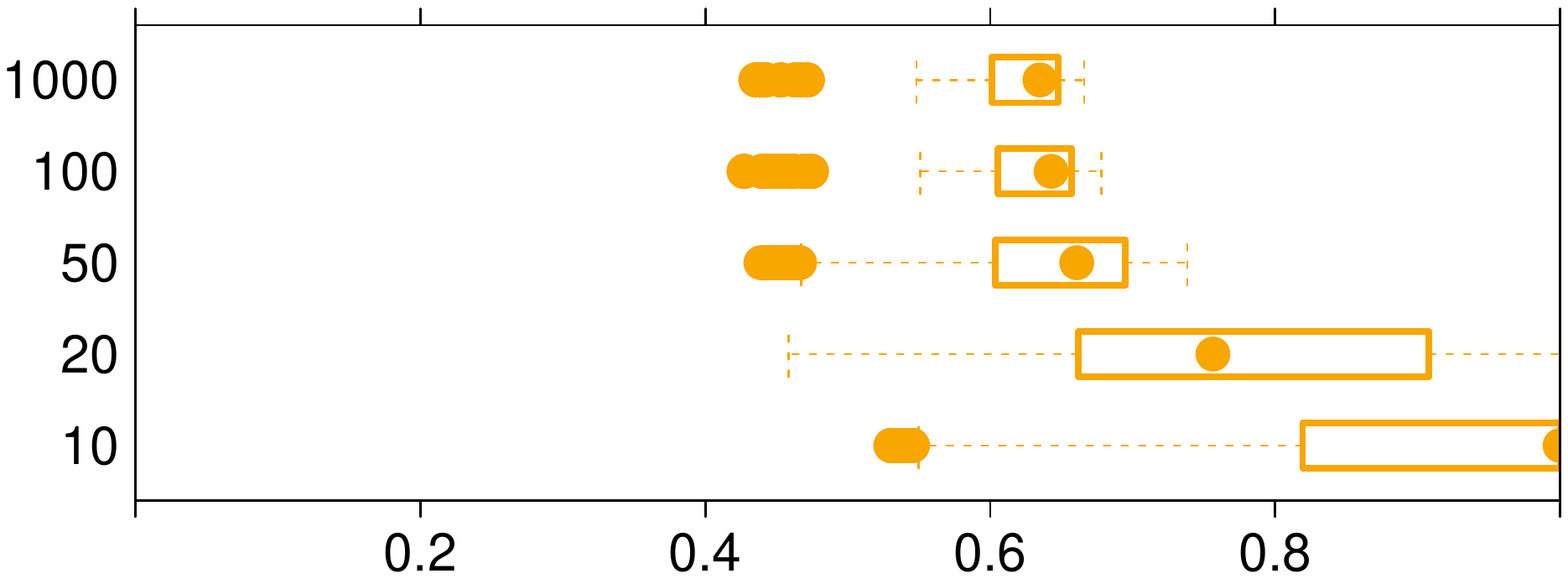} & 
\includegraphics[height=5cm]{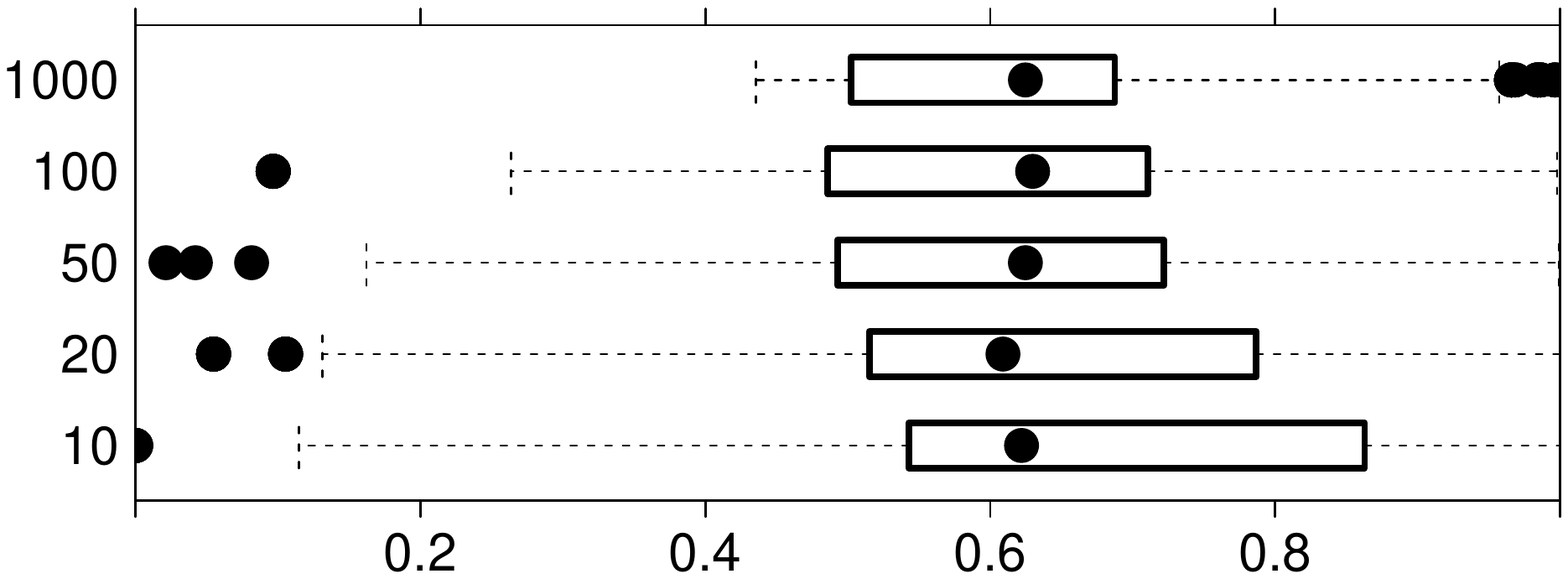} &
\raisebox{2cm}{\textcolor{black}{All}} \\
\end{tabular}
\end{center}
\caption{Boxplots of the $\textrm{HIM}(\circ,\mathcal{E}_N)$ for all combinations of node size and network type.} 
\label{fig:boxplots}
\end{figure}

\begin{figure}[!t]
\begin{center}
\begin{tabular}{cc}
\includegraphics[height=6cm]{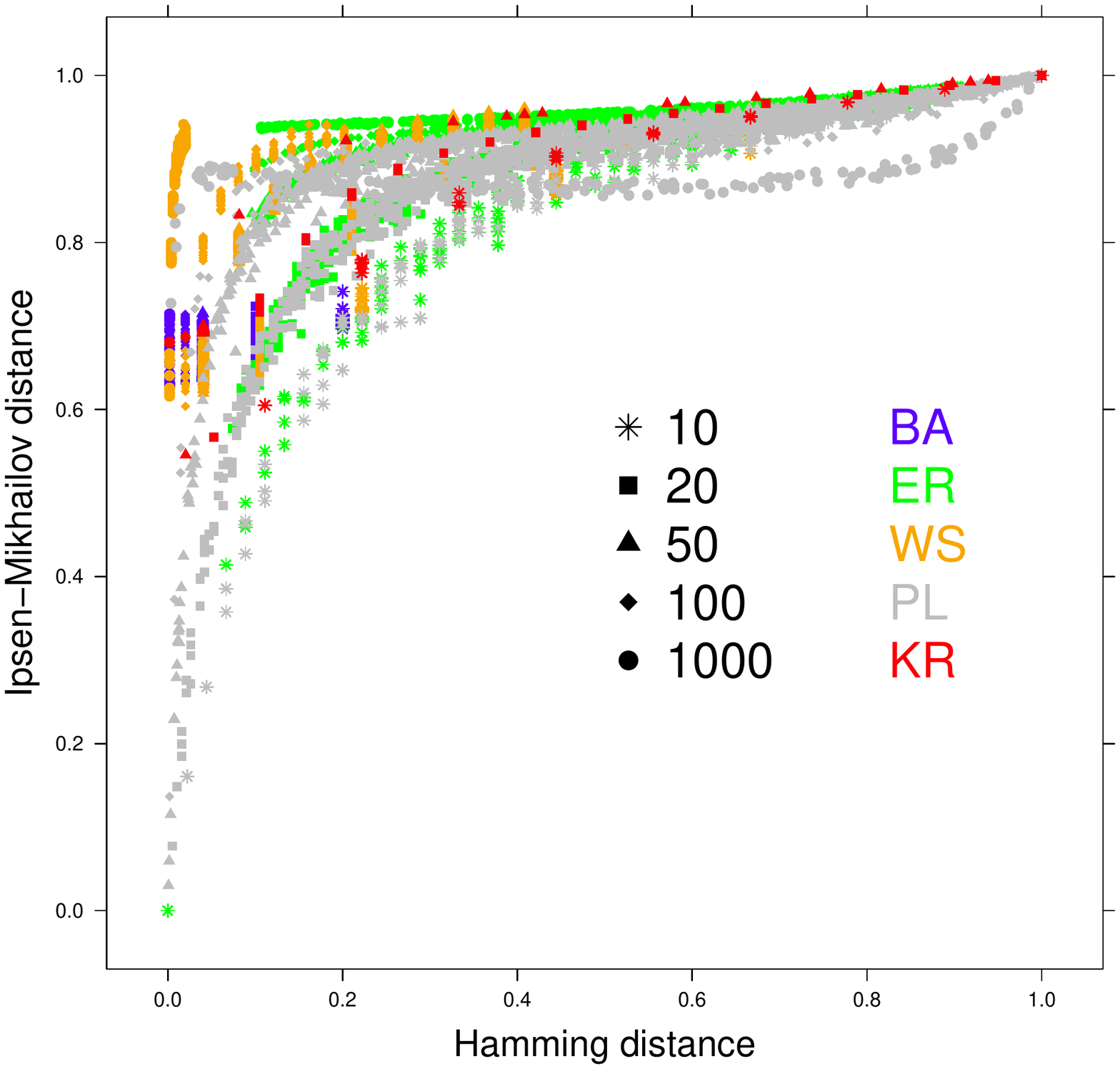} & \includegraphics[height=6cm]{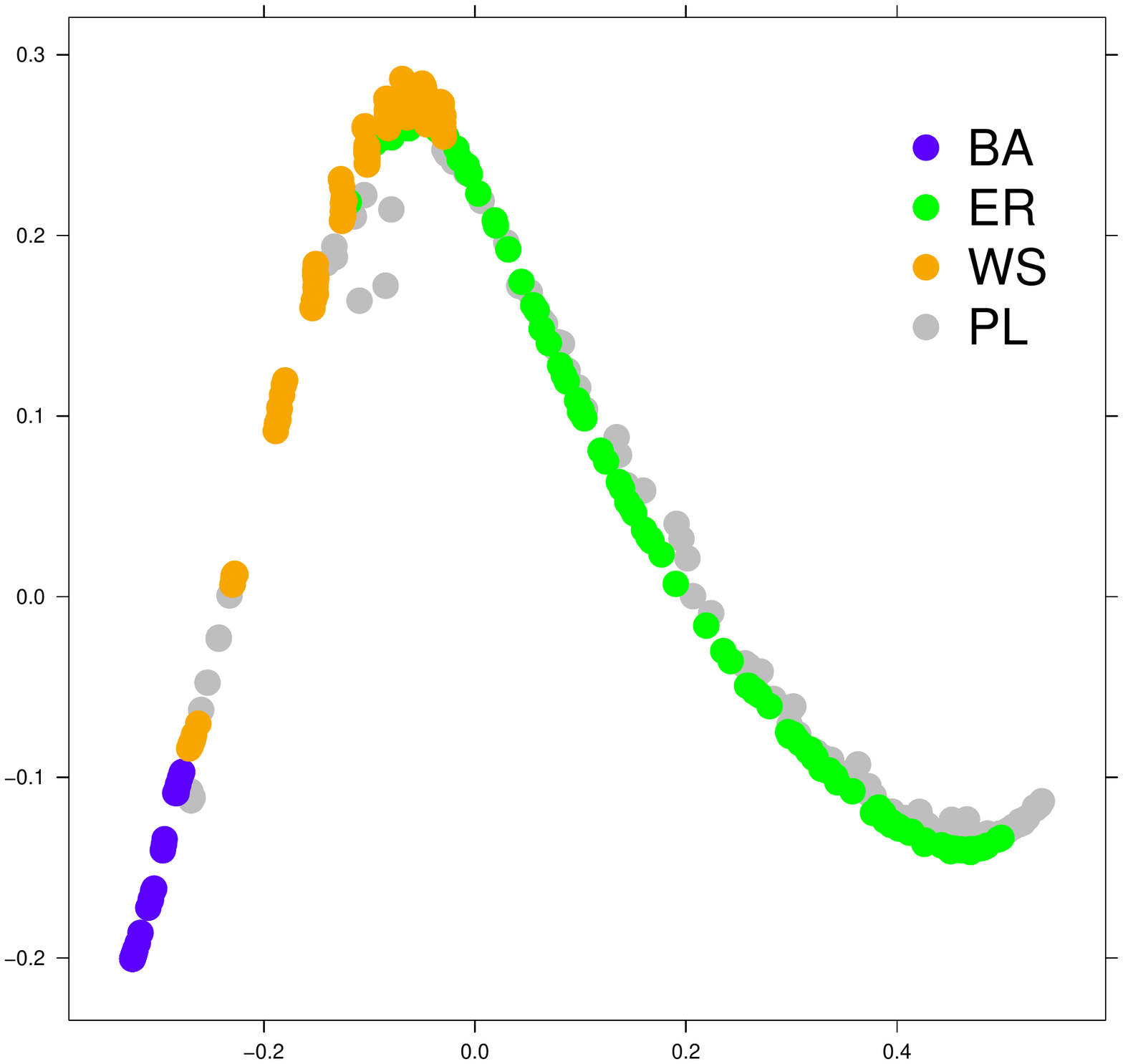} \\
(a) & (b) 
\end{tabular}
\end{center}
\caption{(a) Scatterplot in the Hamming/Ipsen-Mikhailov space of the $\textrm{HIM}(\circ,\mathcal{E}_N)$ for all combinations of node size and network type; (b) Multidimensional Scaling of the mutual HIM distances of 400 networks with 100 nodes in the BA, ER, WS and PL families.}
\label{fig:families}
\end{figure}

In this section we investigate the distribution of the distances from the empty network of a set of graphs randomly extracted from five families.
In particular, for each $N$= 10, 20, 50, 100 and 1000 we extracted 1000 networks on $N$ nodes from each of the following class of graphs:
\begin{itemize}
\item{BA} Barabasi-Albert model \cite{barabasi99emergence}, with power of preferential attachment extracted from the uniform distribution between 0.1 and 10.
\item{ER} Erd\"os-R\'enyi model \cite{erdos59random,erdos60evolution}, with link probability extracted from the uniform distribution between 0.1 and 0.9. 
\item{WS} Watts-Strogatz model \cite{watts98collective}, with neighborhood within which the vertices of the lattice will be connected uniformly sampled in $\{1,\ldots,10\}$ and rewiring probability extracted from the uniform distribution between 0.1 and 0.9.
\item{PL} Scale-free random graphs from vertex fitness scores \cite{goh01universal}, with number of edges uniformly sampled between 1 and $\frac{N(N-1)}{2}$ and power law exponent of the degree distribution extracted from the uniform distribution between 2.005 and 3.
\item{KR} Random regular graphs, with all possible values of node degree.
\end{itemize}
In Tab.~\ref{tab:families} we list mean $\mu$ and standard deviation $\sigma$ of $\textrm{HIM}(\circ,\mathcal{E}_N)$ for all combinations of node size and network type: note that we do not report the corresponding median, because its distance from the mean $\mu$ is always smaller than 0.02 nor the bootstrap confidence intervals, whose range is always smaller than 0.02 from either side of the mean.
In Fig.~\ref{fig:boxplots} we also show the corresponding boxplots, while in Fig.~\ref{fig:families}(a) we display the scatterplot in the Hamming/Ipsen-Mikhailov space of all the aforementioned distances.
In the Hamming/Ipsen-Mikhailov space all the BA nets are confined in the narrow rectangle $[0,0.2]\times [0.6,9.75]$, while all other classes of graphs span a much wider area.
In particular, the points corresponding to distances of the PL nets occupy densely all the upper left triangle of the H/IM plane, and the same happens, with $H>0.1$, also for the ER networks, while WS and KR points lie in the upper rectangle $[0,1]\times [0.6,1]$.
Thus, different PL networks show very different stucture, while the BA nets are very homogeneous.
Notably, no point occurs in the lower right corner of the H/IM space.
Moreover, in average, the standard deviation decreases inversely with the network size, showing larger homogeneity in bigger networks.
Finally, to better highlight the difference among the diverse families, we randomly extracted 100 networks with 100 nodes for the four families BA, ER, WS and PL and we computed the mutual distances between all possible pairs of these 400 graphs.
A few statistics of these HIM distances are reported in Tab.~\ref{tab:stats}, while the planar multidimensional scaling plot \cite{cox01multidimensional} is displayed in Fig.~\ref{fig:families}(b).
Apart from the PL networks, the three families BA, ER and WS can be mutually well separated as shown in the multidimensional plot; moreover, the graphs in the BA and in the WS families are mutually quite similar, as supported by the small interclass mean HIM distance.
On the other side, the PL networks have essentially the same distance from all other groups, so they cannot be easily distiguished.
\begin{table}[!b]
\caption{Mean $\mu$ and standard deviation $(\sigma)$ of the HIM distances between the 100-nodes graphs in the four families BA, ER, WS and PL; each family includes 100 networks.}
\label{tab:stats}
\begin{center}
\begin{tabular}{c|cccc}
& \textcolor{blue}{BA} &\textcolor{green}{ER} &\textcolor{orange}{WS} & \textcolor{gray}{PL} \\
\hline\\
\textcolor{blue}{BA}& 0.05 (0.06) & 0.69 (0.10) & 0.47 (0.13) & 0.65 (0.17) \\
\textcolor{green}{ER} & 		 & 0.50 (0.13) & 0.56 (0.15) & 0.51 (0.14) \\
\textcolor{orange}{WS}& 		 &	       & 0.29 (0.12) & 0.56 (0.18) \\
\textcolor{gray}{PL}& 		 &	       &	     & 0.50 (0.17) \\
\end{tabular}
\end{center}
\end{table}

\begin{table}[!b]
\caption{Mean $\mu$ and standard deviation $\sigma$ of HIM distances $\textrm{HIM}(\circ,\mathcal{E}_N)$ from the empty network for all combinations of network type T and network size N, and cumulatively across node sizes and graph classes.}
\label{tab:families}
\begin{center}
\begin{tabular}{c|cc|cc|cc|cc|cc|cc}
\multicolumn{1}{r|}{$\to$N}& \multicolumn{2}{c|}{10} & \multicolumn{2}{c|}{20} & \multicolumn{2}{c|}{50} & \multicolumn{2}{c|}{100} & \multicolumn{2}{c}{1000} & \multicolumn{2}{c}{All}\\
\multicolumn{1}{l|}{$\downarrow$ T}& $\mu$ & $\sigma$ & $\mu$ & $\sigma$ & $\mu$ & $\sigma$ & $\mu$ & $\sigma$ & $\mu$ & $\sigma$ & $\mu$ & $\sigma$\\
\hline
\textcolor{blue}{BA} & 0.53 &0.01 &0.51 &0.01 &0.50 &0.02 &0.49 &0.02 &0.50 &0.02 & 0.51 & 0.02	\\
\textcolor{green}{ER} & 0.69 &0.18 &0.73 &0.12 &0.77 &0.09 &0.77 &0.08 &0.76 &0.08 & 0.74 & 0.12\\
\textcolor{orange}{WS} &0.91 &0.14 &0.76 &0.15 &0.64 &0.08 &0.62 &0.06 &0.62 &0.05 & 0.71 & 0.16\\
\textcolor{gray}{PL}& 0.72 &0.19 &0.72 &0.19 &0.75 &0.15 &0.74 &0.14 &0.72 &0.11 & 0.73 & 0.16\\
\textcolor{red}{KR} & 0.60 &0.11 &0.54 &0.10 &0.50 &0.05 &0.49 &0.00 & 0.48 &0.00 & 0.52 & 0.08\\
All & 0.69 & 0.19 & 0.65 & 0.17 & 0.63 & 0.15 & 0.62 & 0.14 & 0.62 & 0.13 & 0.64 & 0.16 \\
\end{tabular}
\end{center}
\end{table}
\subsection{The \textit{D. melanogaster} development dataset}
\label{ssec:droso}
\begin{figure}[!t]
\begin{center}
\begin{tabular}{ccc}
\includegraphics[height=0.3\textwidth]{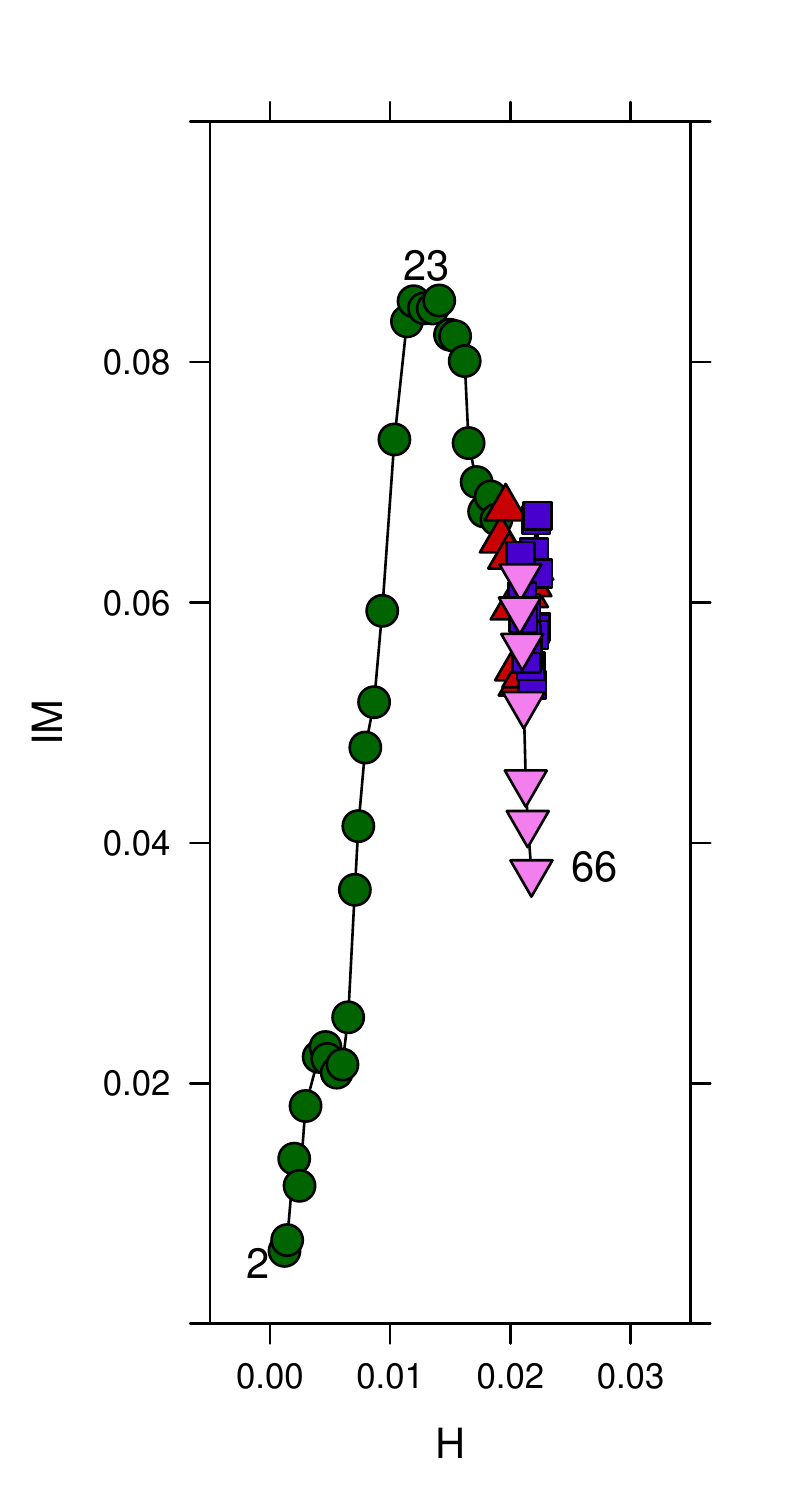} & 
\includegraphics[height=0.3\textwidth]{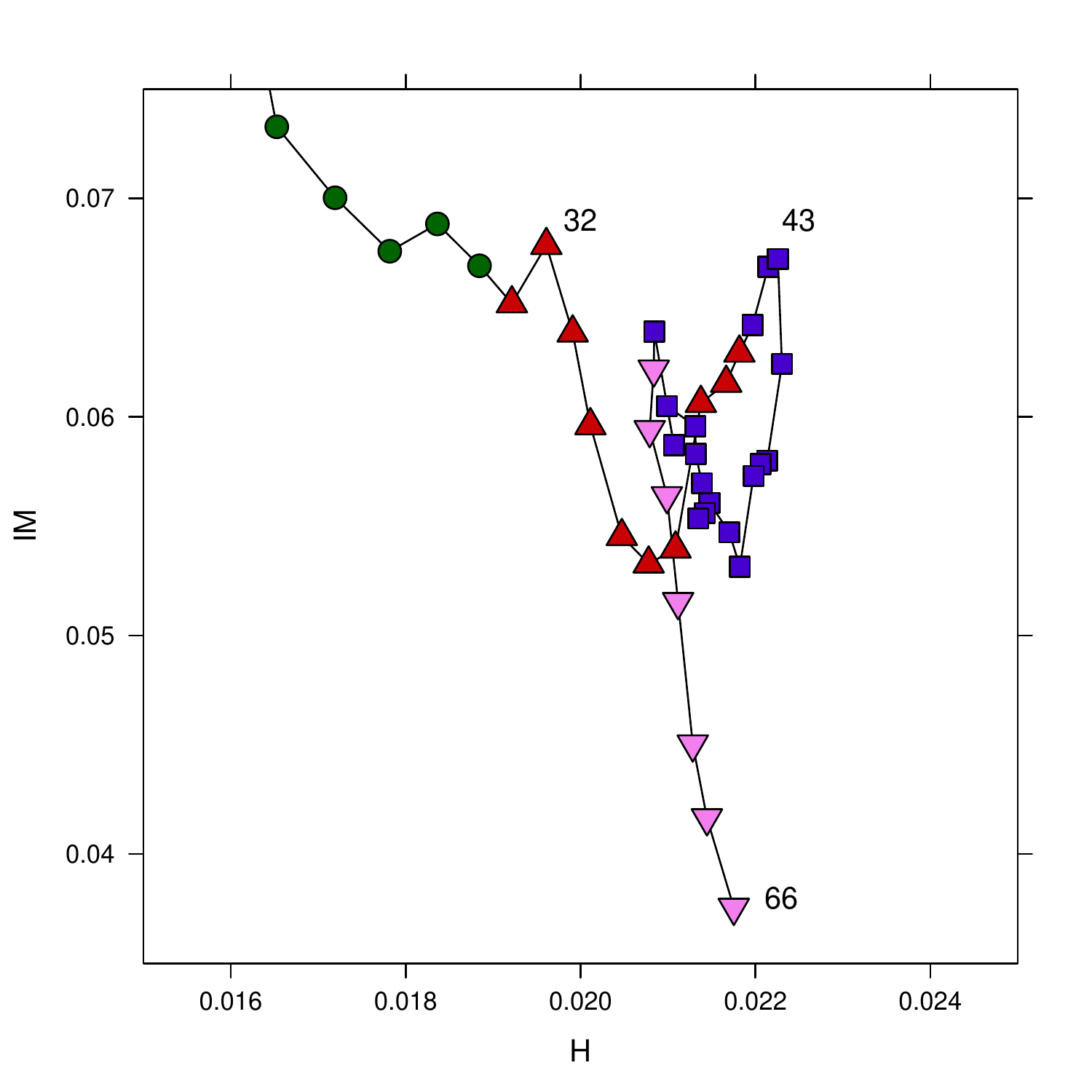} & 
\includegraphics[height=0.3\textwidth]{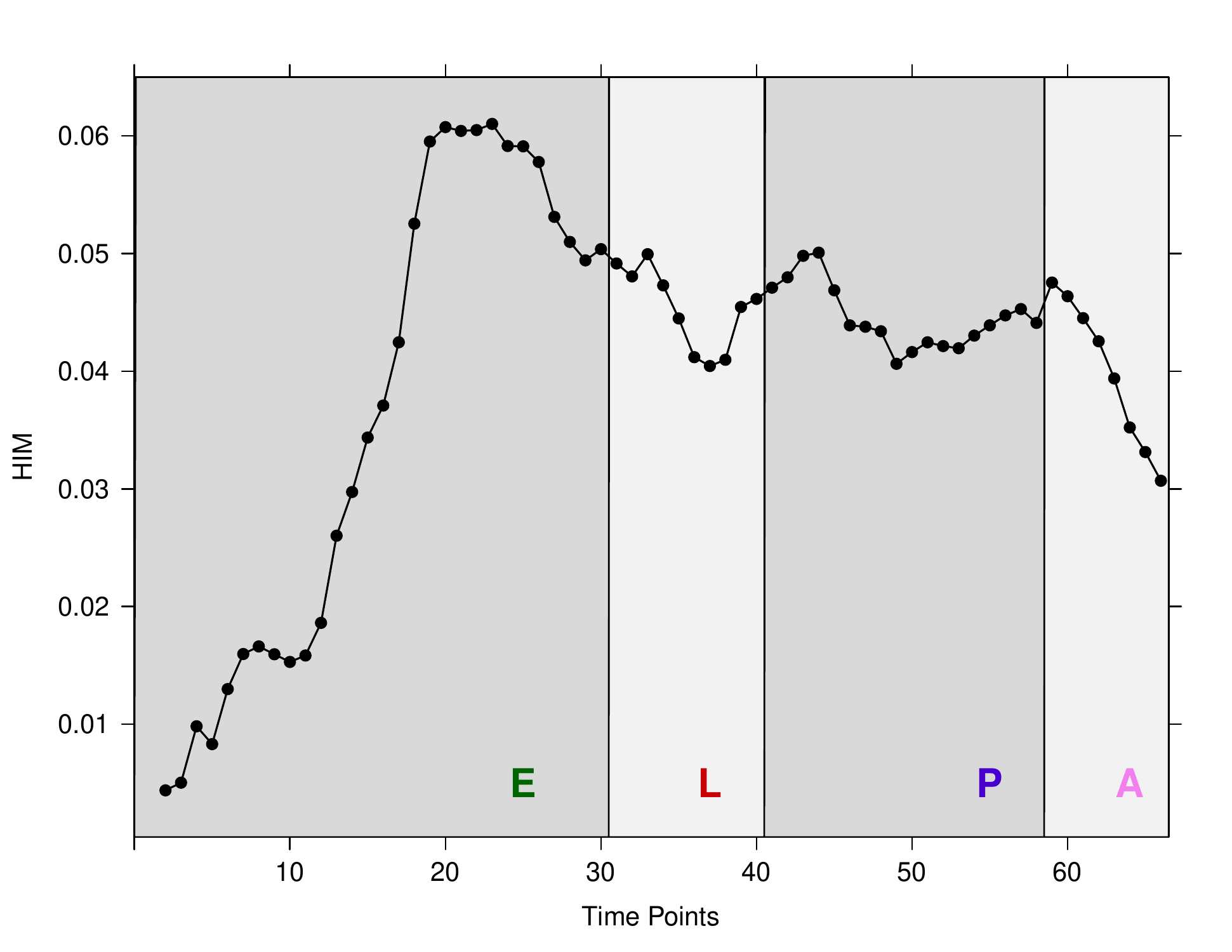} \\
(a) & (b) & (c)
\end{tabular}
\caption{(a) Evolution of distances of the \textit{D. melanogaster} development gene network time series in the Hamming/Ipsen-Mikhailov space with zoom (b) on the final timepoints and (c) evolution of same HIM distances along 66 time points in the 4 stages Embryonic (E), Larval (L), Pupal (P) and Adult (A).}
\label{fig:time}
\end{center}
\end{figure}

In \cite{kolar10estimating}, the authors used the Keller algorithm to infer the gene regulatory networks of \textit{Drosophila melanogaster} from a time series of gene expression data measured during its full life cycle, originally published in \cite{arbeitman02gene}.
They followed the dynamics of 588 development genes along 66 time points spanning through four different stages (Embryonic -- time points 1-30, Larval -- t.p. 31-40, Pupal -- t.p. 41-58, Adult -- t.p. 59-66), constructing a time series of inferred networks $N_i$, publicly available at \url{http://cogito-b.ml.cmu.edu/keller/downloads.html}.
Hereafter we evaluate the structural differences between $N_i$ and the initial network $N_1$, as measured by the HIM distance: the resulting plot is displayed in Fig.~\ref{fig:time}.
The largest variations, both between consecutive terms and with respect to the initial network $N_1$, occur in the embrional stage (E): in particular, the HIM distance grows until time points 23, then the following networks start getting closer again to $N_1$, showing that the interactions of the selected 588 genes in the adult stage are more similar to the corresponding net of interaction in the embrional stage, rather than in the other two stages.
Moreover, while Hamming distance ranges between $0$ and $0.0223$, the Ipsen-Mikhailov distance has $0.0851$ as its maximum, indicating an higher variability of the networks in terms of structure rather than matching links.
Finally, using a Support Vector Machine with HIM kernel built in the \textit{kernlab} package in R, a 5-fold Cross Validation with $\gamma=10^3$ and $C=1$ reached accuracy 0.97 in discriminating Embryonic and Adult networks from Larval and Pupal, while, in the same setup, we reach perfect separation between Embryonic and Adult stages for all values of $\gamma$ larger than 1000.
\subsection{The HCC dataset} 
\label{ssec:hcc}
\begin{figure}[!t]
\begin{center}
\begin{tabular}{ccc}
\includegraphics[clip,width=0.18\textwidth, angle=-90]{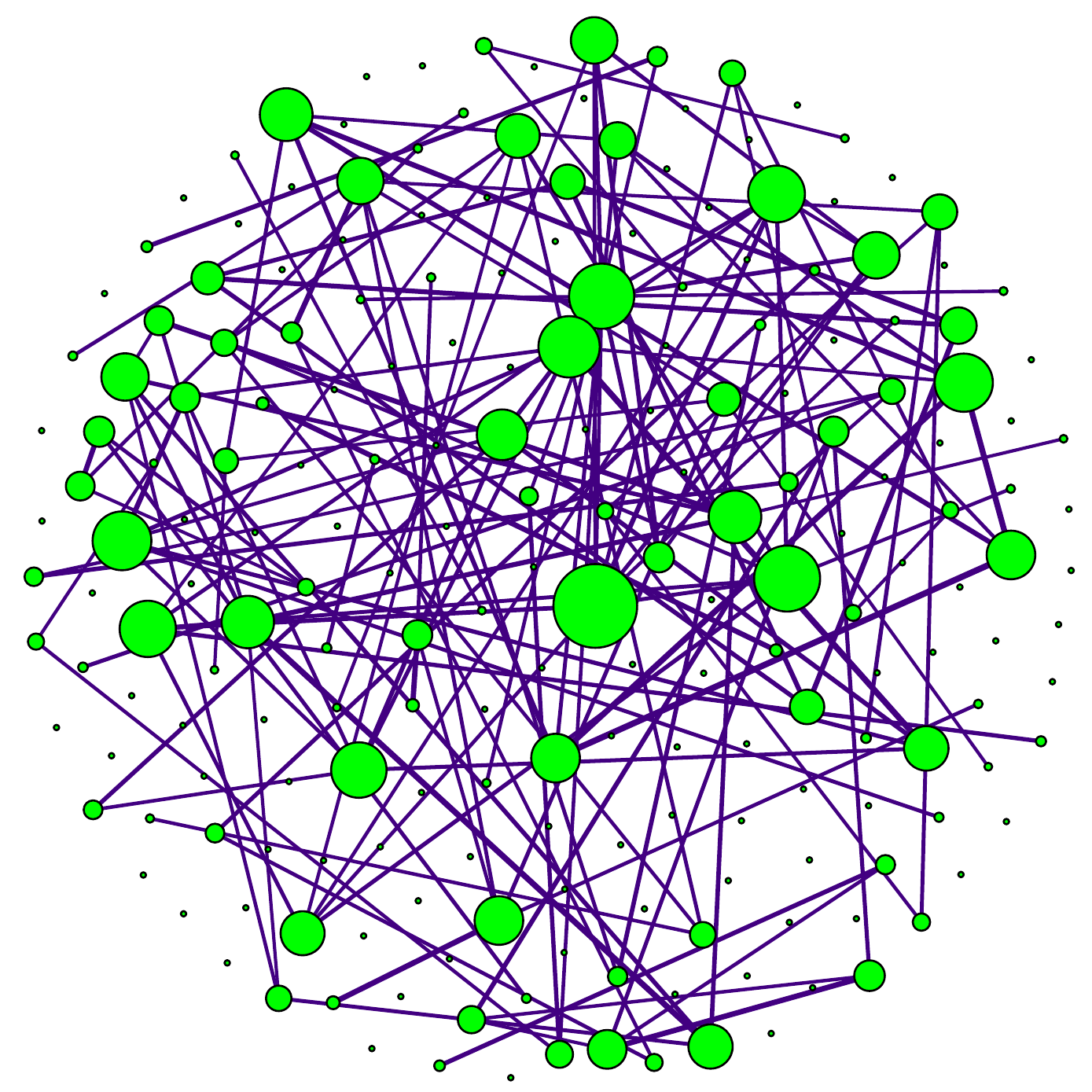} &
\includegraphics[clip,width=0.18\textwidth, angle=-90]{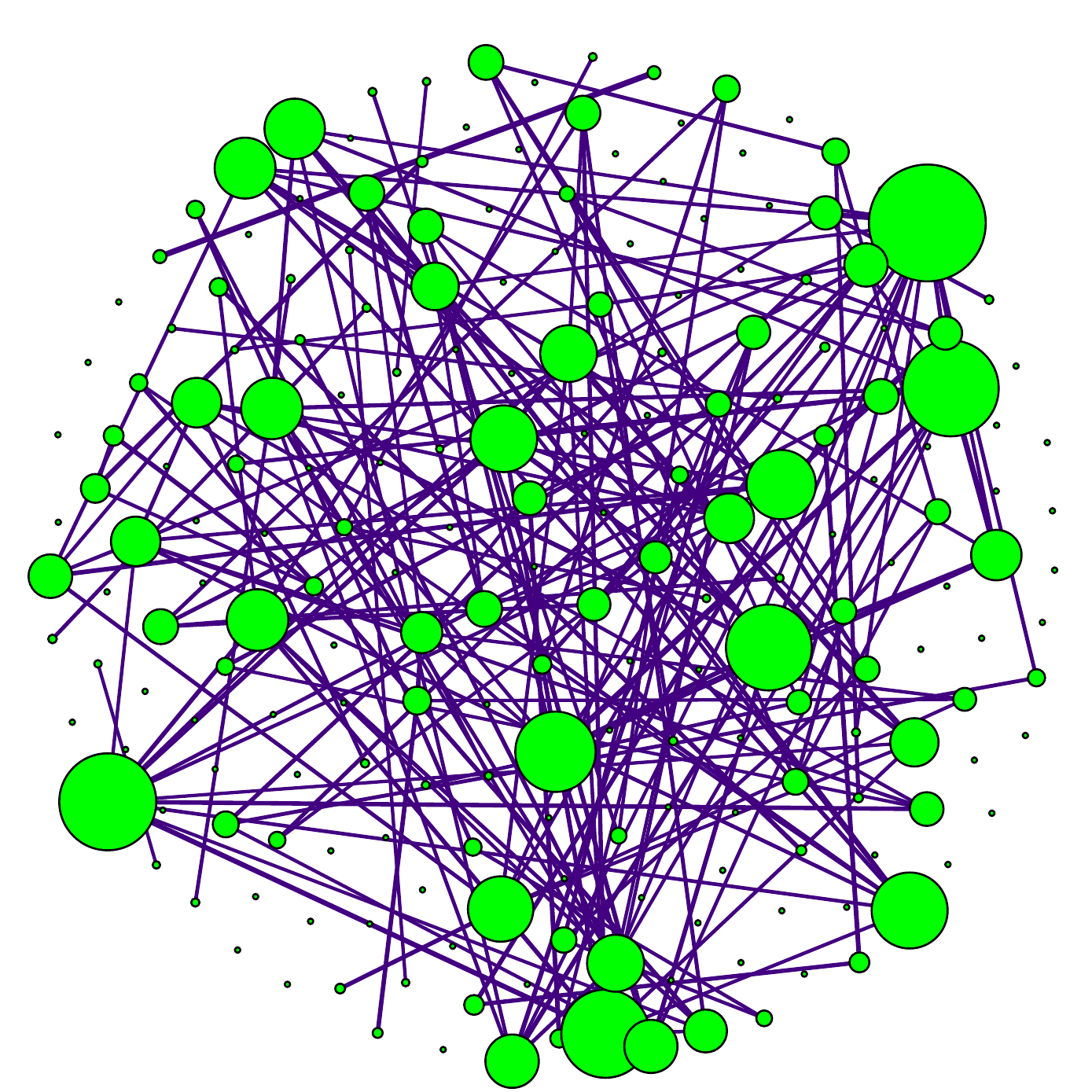} &
\multirow{4}{*}{\includegraphics[clip,width=0.6\textwidth]{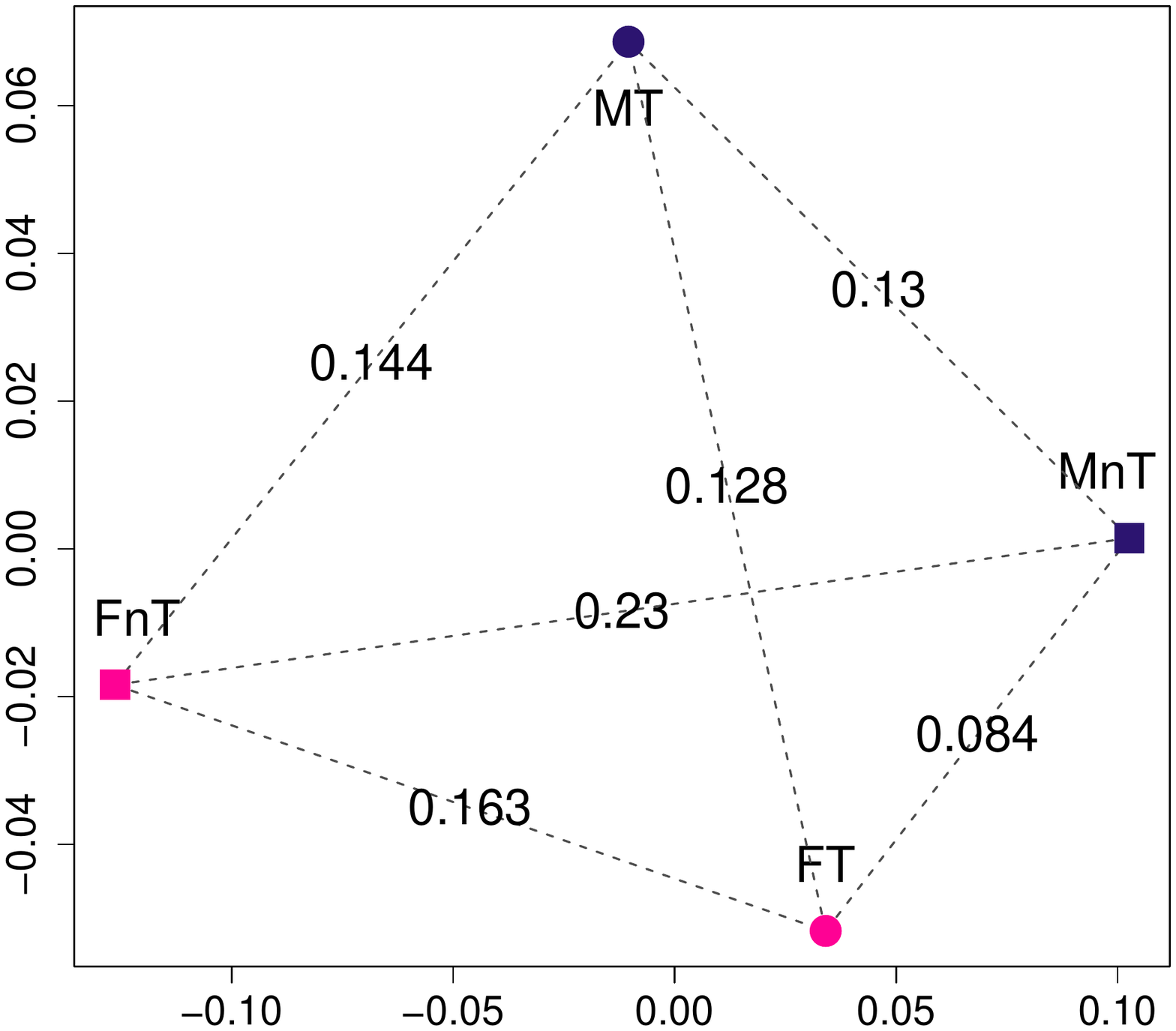}} \\
\\
MT \textcolor{myblue}{\CircleSolid} & MnT \textcolor{myblue}{\SquareSolid} \\
\includegraphics[clip,width=0.18\textwidth, angle=-90]{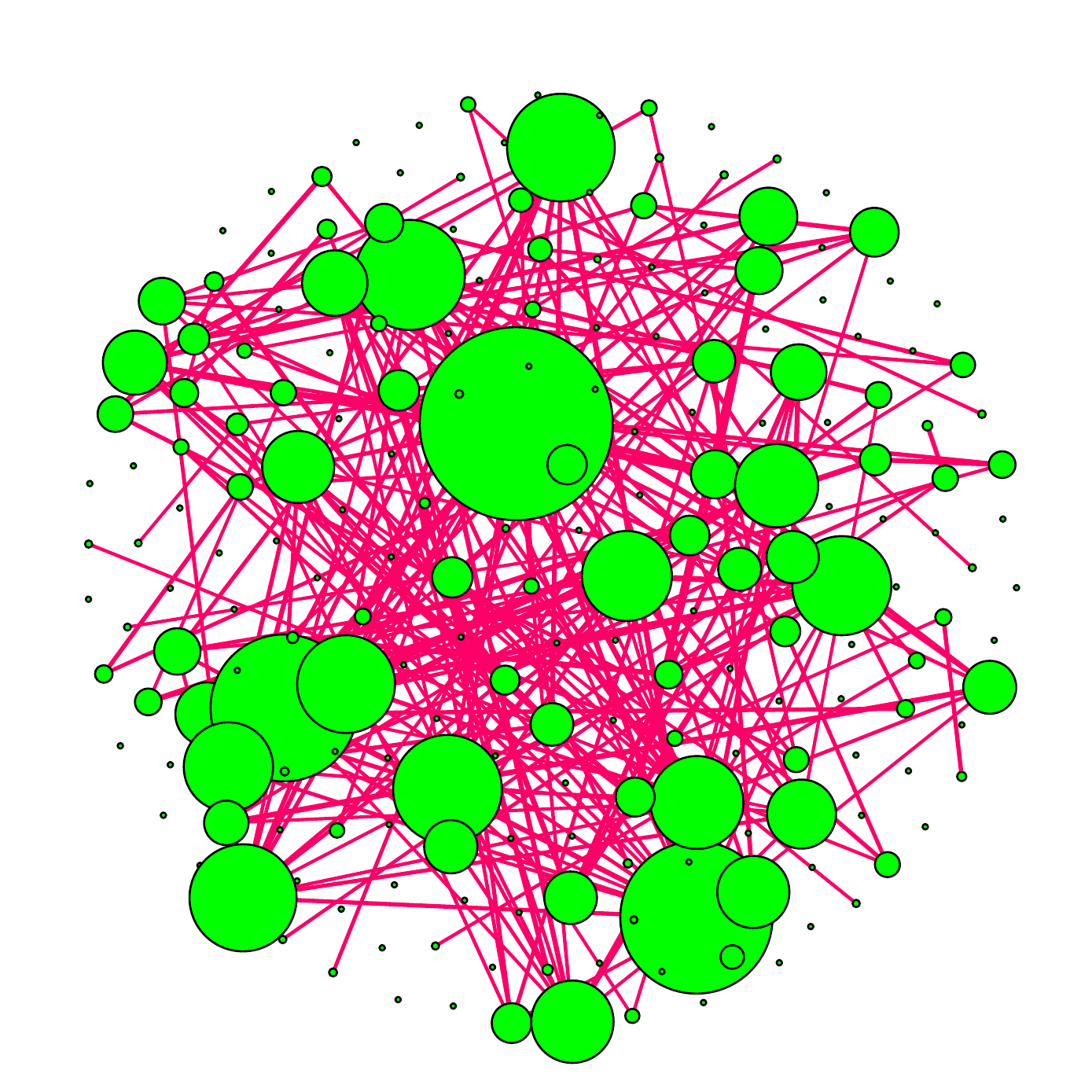} &
\includegraphics[clip,width=0.18\textwidth, angle=-90]{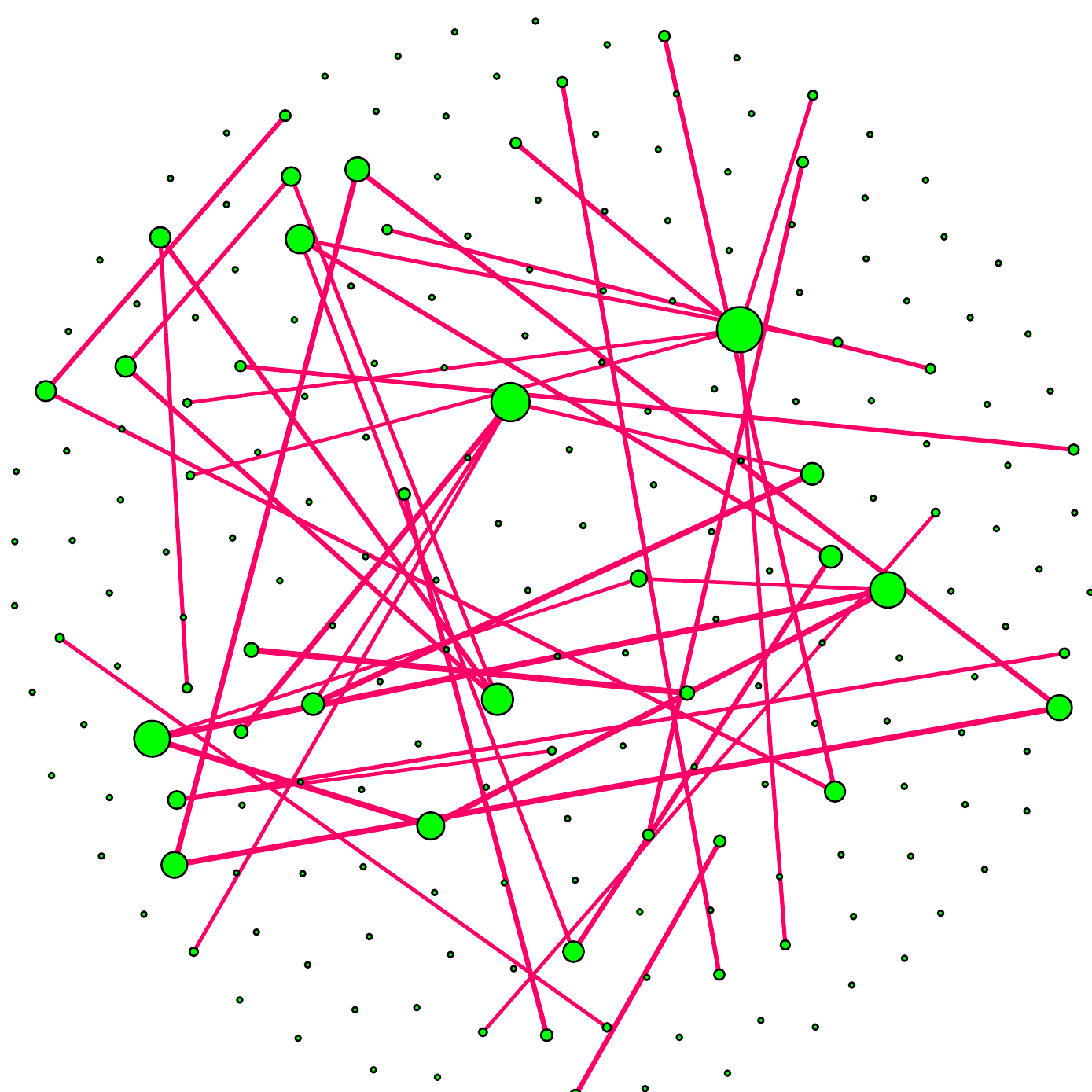} \\
\\
FT \textcolor{mypurple}{\CircleSolid} & FnT \textcolor{mypurple}{\SquareSolid} \\
\end{tabular}
\end{center}
\caption{(left) Fruchterman-Reingold layout of the networks MT, MnT, FT, FnT where the male networks are inferred by one random extraction of 30 samples out of the whole cohort of 210 patients. Node size is proportional to node degree. (right) Multidimensional scaling plot, with mutual HIM distances, of the networks FT, FnT, MT, MnT.}
\label{fig:hcc_him}
\end{figure}

Publicly available at the Gene Expression Omnibus (GEO) \url{http://www.ncbi.nlm.nih.gov/geo}, at the Accession Number GSE6857, the HepatoCellular Carcinoma (HCC) dataset \cite{budhu08identification,ji09microrna} collects 482 tissue samples from 241 patients affected by HCC, a well-studied pathology \cite{law11emerging,gu12gene} where the impact of microRNA (miRNA) is notably relevant \cite{volinia10reprogramming,bandyopadhyay10development}.
For each patient, a sample from cancerous hepatic tissue and a sample from surrounding non-cancerous hepatic tissue are available, hybridized on the Ohio State University CCC MicroRNA Microarray Version 2.0 platform collecting the signals of 11,520 probes of 250 non-redundant human and 200 mouse miRNA.
After a preprocessing phase including imputation of missing values \cite{troyanskaya01missing} and discarding probes corresponding to non-human (mouse and controls) miRNA, we consider the dataset HCC of 240+240 paired samples described by 210 human miRNA,
with the cohort consisting of 210 male and 30 female patients.
We thus parted the whole dataset HCC into four subsets combining the sex and disease status phenotypes, collecting respectively the cancer tissue for the male patients (MT), the cancer tissue for the female patients (FT) and the corresponding two datasets including the non cancer tissues (MnT, FnT).
Then we first generated the four co-expression networks on the 210 miRNA as vertices, inferred via absolute Pearson's correlation and corresponding to the combinations of the two binary phenotypes, and we computed all mutual HIM distances.
In particular, to show the possible effects due to the different sample size, we computed 30 instances of the MT and MnT networks, inferred using only 30 matching samples and then averaging all the mutual HIM distances.
One instance of MT and MnT is displayed as an hairball in Fruchterman-Reingold layout \cite{fruchterman91graph} together with the nets FT and FnT.
The corresponding two-dimensional scaling plot \cite{cox01multidimensional} in the right panel of the Figure~\ref{fig:hcc_him}. 
The four networks are widely separated, with orthogonal separations for the two phenotypes, but the values of the  HIM distances between the network support the known different development of HCC in male and female: for instance, the FT network is closer to the MnT net (HIM=0.08), rather than to the MT and FnT (HIM=0.13 and 0.16, respectively). 
Note that the largest distance (HIM=0.23) is detected between the two non-tumoral networks MnT, FnT.
An expanded version of the example is shown in \cite{jurman12stability,filosi13stability}, where more networks are generated from the same dataset using different inference algorithms and a stability analysis is performed.
\subsection{The Gulf Dataset} 
\label{ssec:gulf}
\begin{figure}[!b]
\centering
\includegraphics[width=0.7\textwidth]{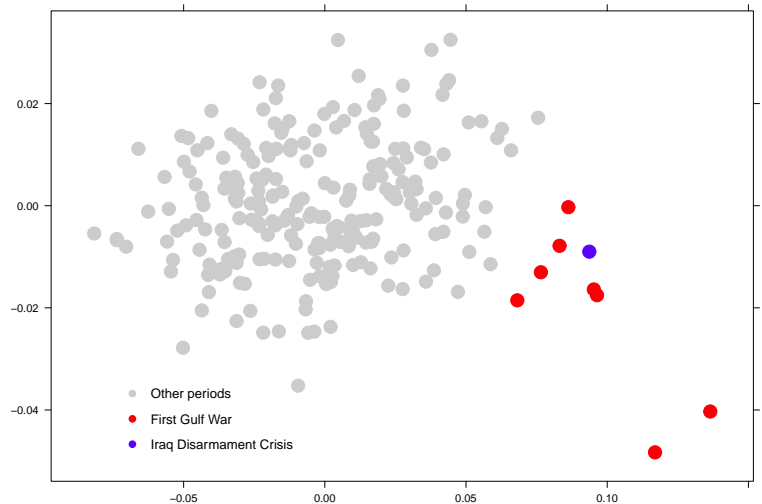}
\caption{Planar HIM distance based multidimensional scaling plot of the monthly Gulf Dataset. Red dots corresponds to the First Gulf War months (July 1990 - April 1991), while grey points correspond to months outside that temporal window and the blue point corresponds to February 1998, the month of the Iraq disarmament crisis.}
\label{fig:gulf}
\end{figure}

Part of the Kansas Event Data System, available at \url{http://vlado.fmf.uni-lj.si/pub/networks/data/KEDS/}, the Gulf Dataset collects, on a monthly bases, political events between pairs of countries focusing on the Gulf region and the Arabian peninsula for the period 15 April 1979 to 31 March 1999, for a total of 240 months. Political events belong to 66 classes (including for instance ''pessimist comment'', ''meet'', ''formal protest'', ''military engagement'', etc.) and involve 202 countries. 
This dataset formally translates into a time series of 240 unweighted and undirected graphs with 202 nodes, for which we computed all the mutual $\frac{240\cdot 239}{2}$ HIM distances.
These distances are then used to project the 240 networks on a plane through a multidimensional scaling \cite{cox01multidimensional}: the resulting plot is displayed in Fig.~\ref{fig:gulf}. 
The months corresponding to the First Gulf War months (July 1990 - April 1991) are close together and confined in the lower left corner of the plane, showing both a mutual high degree of homogeneity and, at the same time, a relevant difference to the graphs of all other months.
This shows that, at the onset of the conflict, the diplomatic relations worldwide changed consistently and their structure remained very similar throughout the whole event.
Note that the blue point (closer to the war-like period) corresponds to February 1998, the time of Iraq disarmament crisis: Iraqi President Saddam Hussein negotiates a deal with U.N. Secretary General Kofi Annan, allowing weapons inspectors to return to Baghdad, preventing military action by the United States and Britain.
\subsection{The International Trade Network data}
\label{ssec:wtn}
As an application of the HIM distance on directed and weighted networks, we show four examples based on the International Trade Network (ITN) data, version 4.1, by Gledisch \cite{gleditsch02expanded} available at \url{http://privatewww.essex.ac.uk/∼ksg/exptradegdp.html}, collecting estimates of trade flows between independent states (1948-2000) and GDP per capita of independent states (1950-2000).
As noted by \cite{fronczak12statistical}, due to differences in reporting procedures between countries, incongruences occur between exports from $i$ to $j$ and imports from $i$ to $j$: to avoid this issues, in our analysis we only use the figures reported as export in the dataset.

In what follows, we extract four sets of countries, and we study the evolution of their trade subnetworks during the aforementioned period.
In each example, chosen the set of $N$ countries $C_1,\ldots C_N$, we construct, for every year, the weighted directed network having $C_1,\ldots C_N$ as nodes.
A link between country $C_i$ and country $C_j$ represents the export from $C_i$ to $C_j$, and its weight $w_{ij}$ corresponds to the volume of the export flow. 
Then we compute all mutual HIM distances among these networks, first rescaling link weights in the unit interval.
Finally, using these $\frac{N(N-1)}{2}$ HIM distances we construct a planar classical Multidimensional Scaling plot, transforming the networks in a set of points such that the distances between the points are approximately equal to the mutual HIM dissimilarities, using the methods in \cite{gower66some,mardia78some,cailliez83analytical,cox01multidimensional} as implemented in R.
The aim here is to connect the structural changes in yearly trade networks with time periods and events having a role in explaing such changes.
Note that in \cite{fronczak12statistical}, the authors show that bilateral trade fulfills fluctuation-response theorem \cite{fronczak06fluctuation}, stating that the average relative change in import (export) between two countries is a sum of relative changes in their GDPs.
This result yields that directed connections, \textit{i.e.}, bilateral trade volumes, are only characterized by the product of the trading countries’ GDPs.
 
\begin{figure}[!t]
\begin{center}
\begin{tabular}{cc}
\raisebox{5cm}{
\begin{tabular}{ccccc}
\multicolumn{5}{c}{\includegraphics[height=6cm]{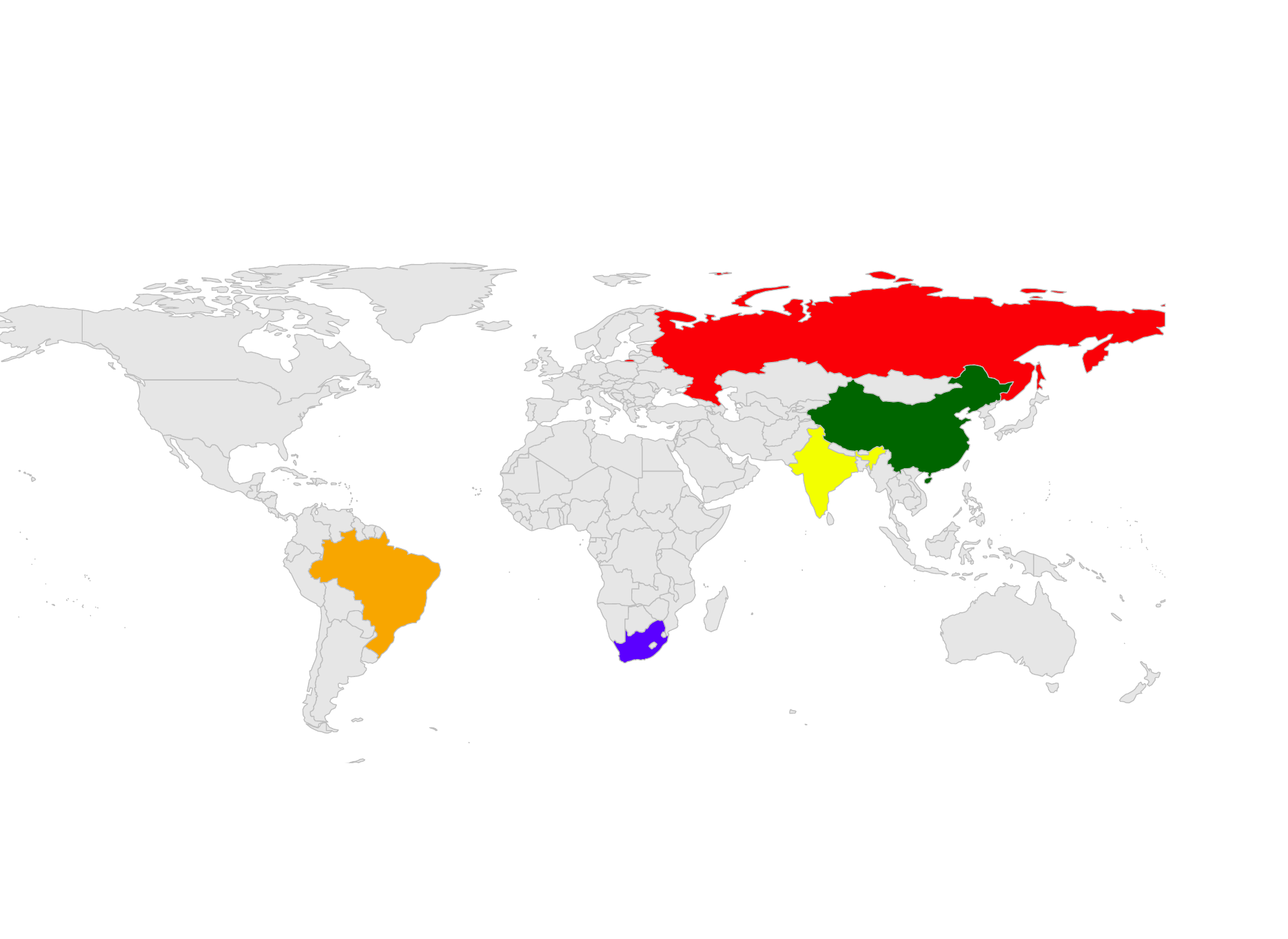}}\\
\includegraphics[height=0.75cm]{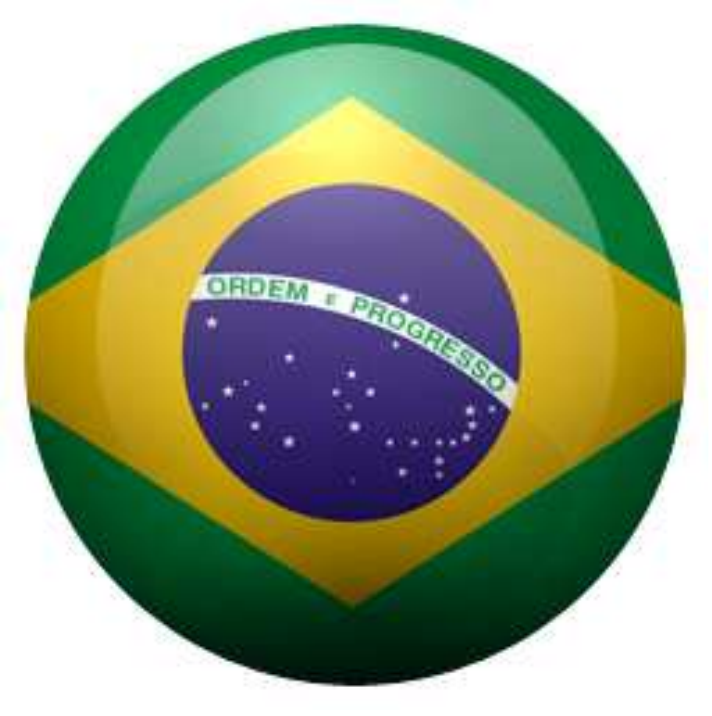}& 
\includegraphics[height=0.75cm]{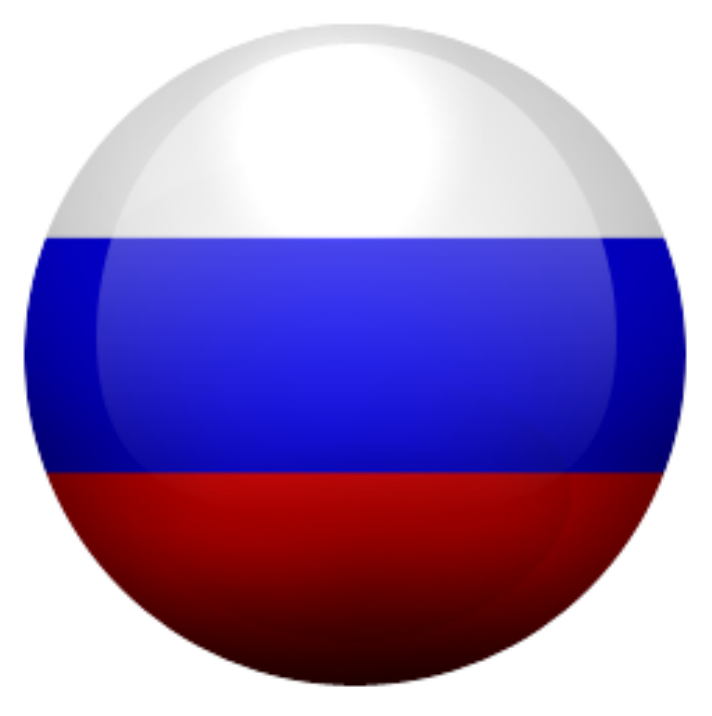}&
\includegraphics[height=0.75cm]{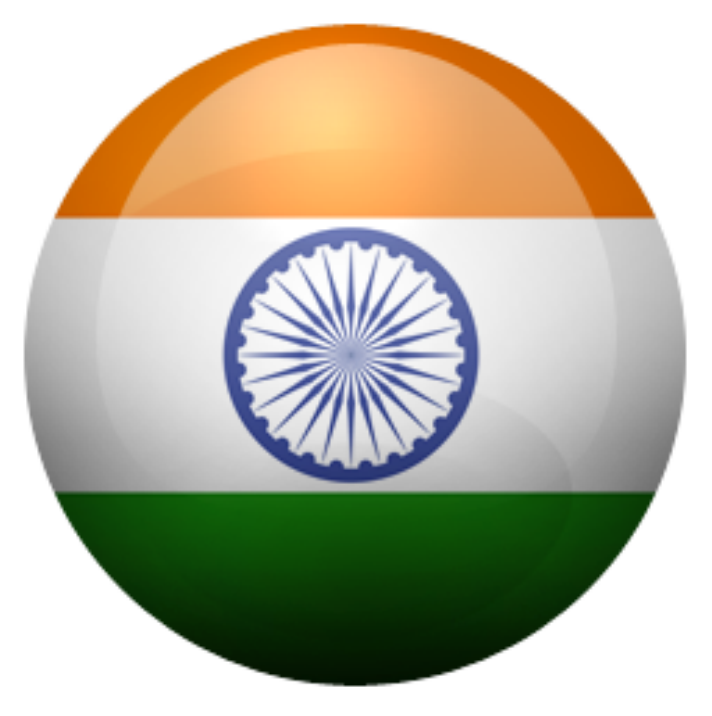}&
\includegraphics[height=0.75cm]{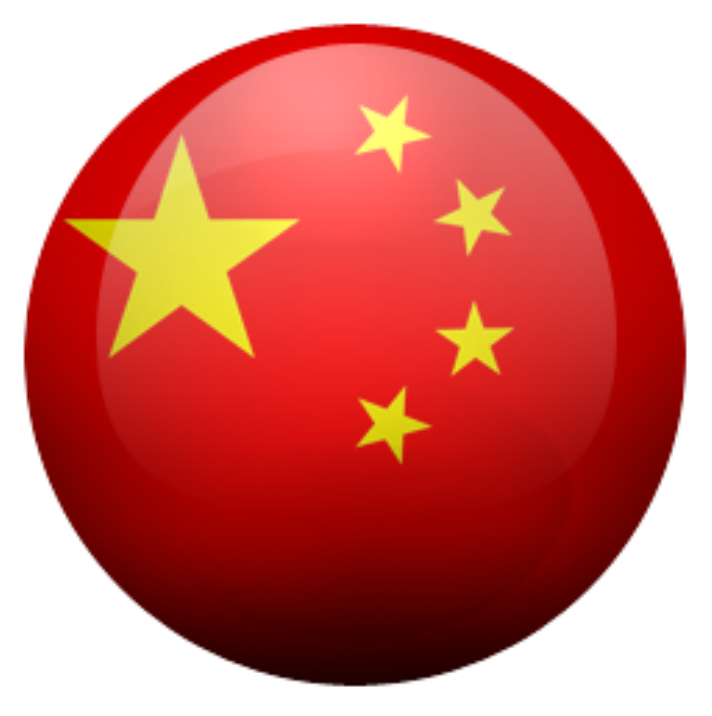}&
\includegraphics[height=0.75cm]{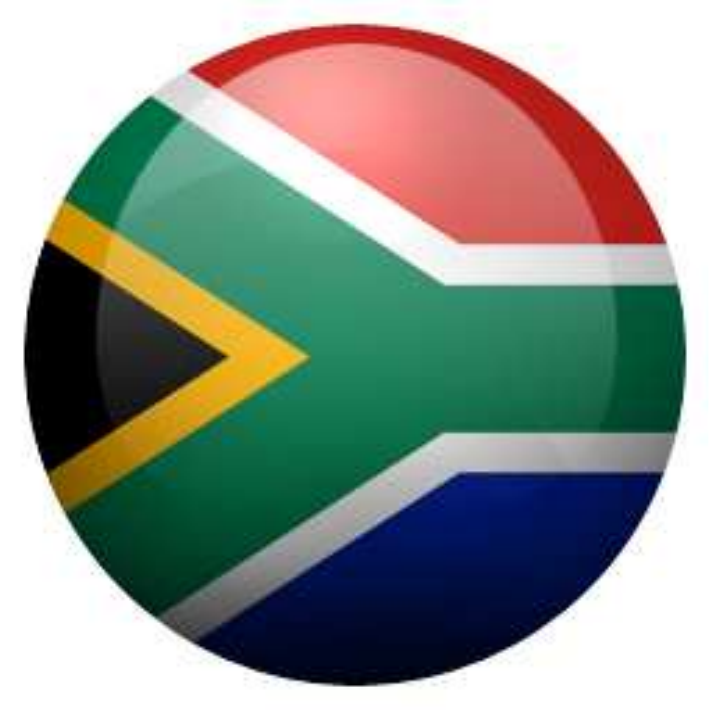}\\
\tiny{BRA} & 
\tiny{RUS} &
\tiny{IND} & 
\tiny{CHN} &
\tiny{SAF} \\
\end{tabular}
}
& 
\raisebox{0cm}{\includegraphics[height=8cm]{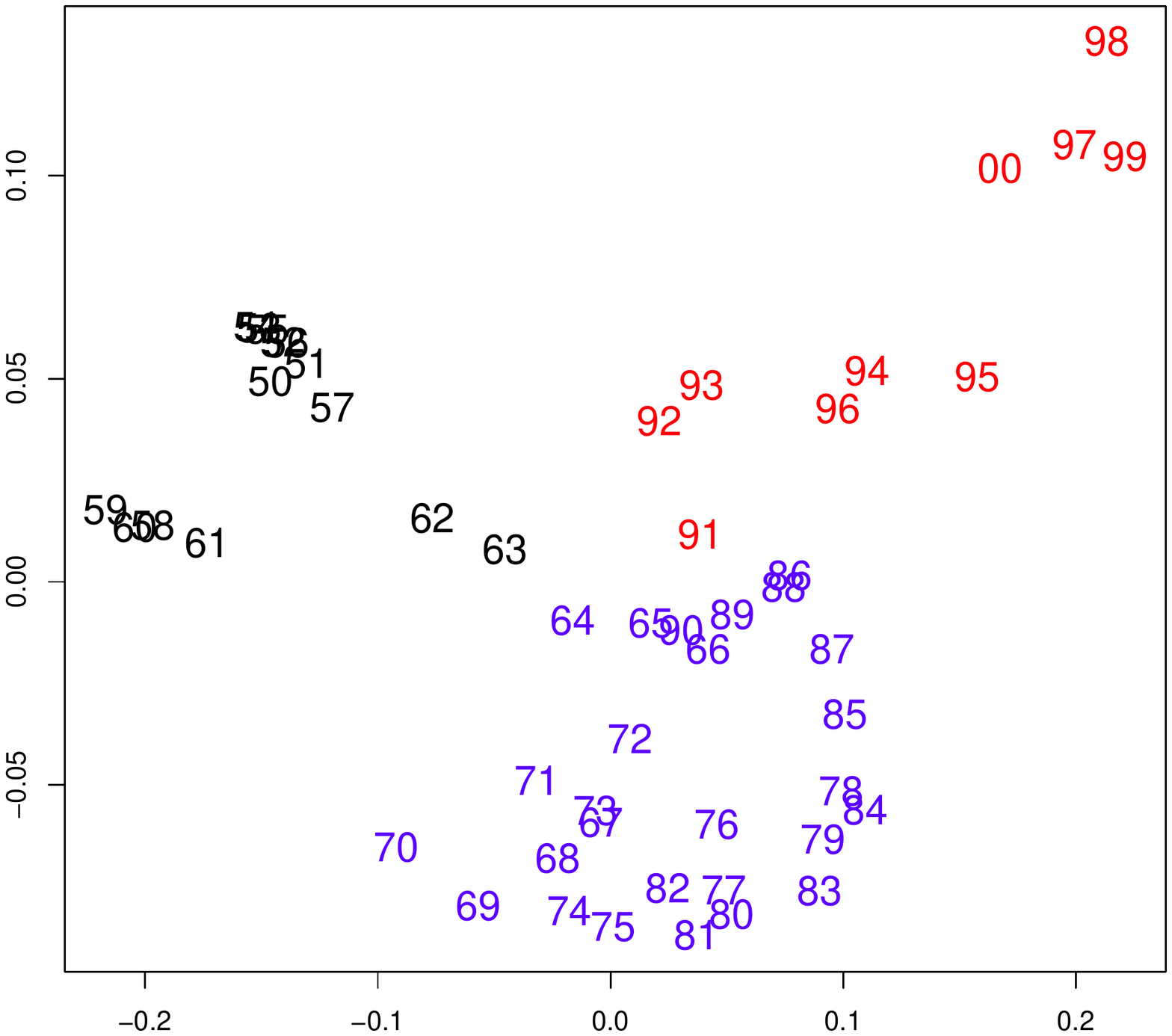}} \\
(a) & (b)
\end{tabular}
\caption{(a) Maps and flags of the BRICS countries. (b) Multidimensional Scaling of the HIM distances among the intertrade networks of the BRICS countries in the periods 1950--1963 (black), 1964--1990 (blue), 1991--2000 (red).}
\label{fig:brics}
\end{center}
\end{figure}

As a first example we present the BRICS countries case.
Introduced in 2001, the acronym BRICS collects the five nations Brazil, Russia, India, China and South Africa (Fig.~\ref{fig:brics}(a)) which, although developing or newly industrialized countries, are distinguished by their large and fast-growing economies and by their significant influence on regional and global affairs.
To such aim, in Fig.~\ref{fig:brics} we show the bidimensional scaling of their trade networks for the years 1950--2000, with the HIM matrix as the distance constraint. 
As shown by the plot, three groups of years can be clearly divided, thus yielding that the corresponding networks are similar within each group, but diverse across different groups: the early years recovering after WWII (until about 1963), the seventies and eighties, where the economies of the involved countries started to develop, and the nineties, where their growth begun to accelerate.
 
\begin{figure}[!t]
\begin{center}
\begin{tabular}{cc}
\raisebox{-0cm}{\includegraphics[height=7cm]{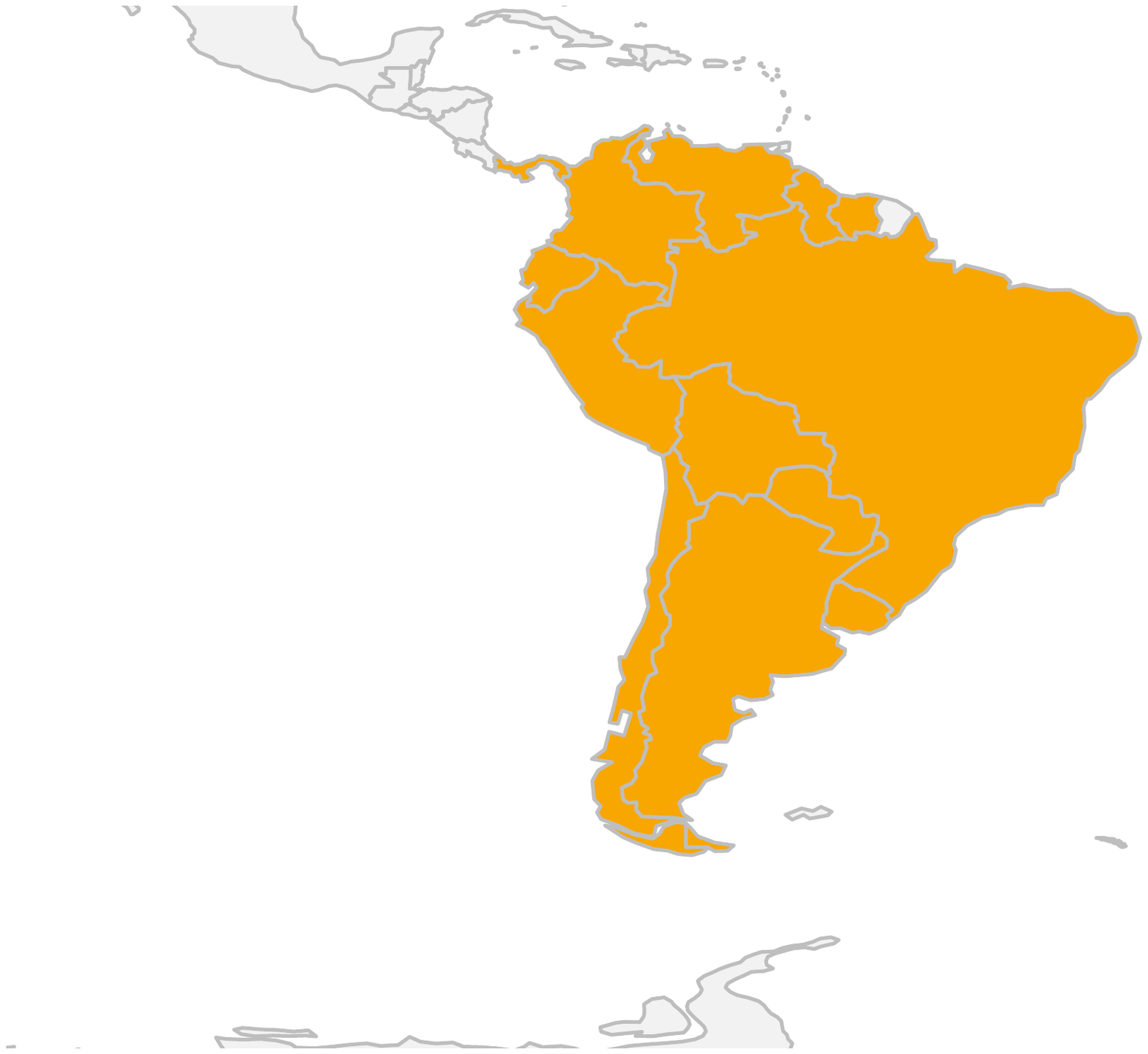}} & 
\raisebox{0cm}{\includegraphics[height=7cm]{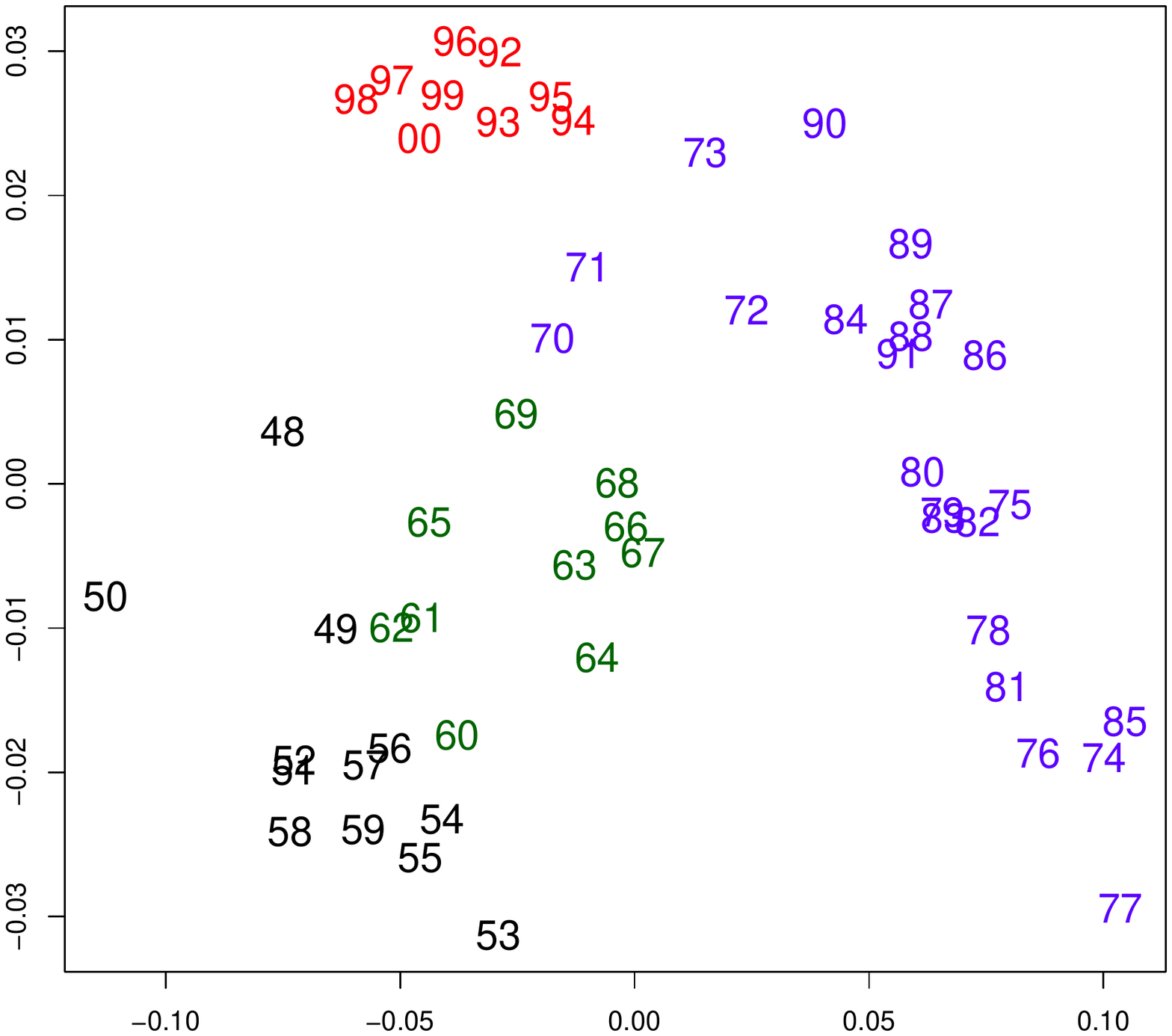}} \\
\multicolumn{2}{c}{
\begin{tabular}{ccccccc}
\tiny{PAN} & 
\tiny{COL} & 
\tiny{VEN} & 
\tiny{GUY} & 
\tiny{SUR} & 
\tiny{ECU} & 
\tiny{PER} \\ 
\includegraphics[height=0.7cm]{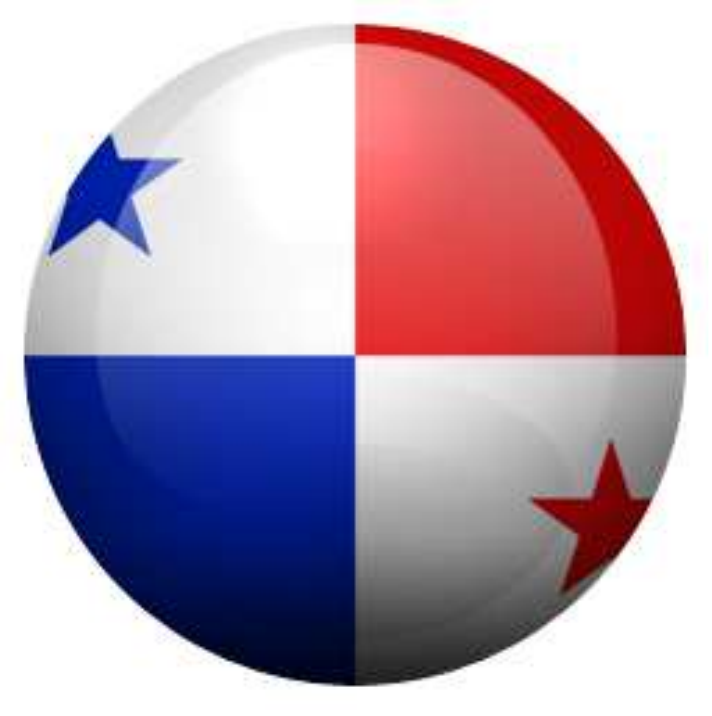} &
\includegraphics[height=0.7cm]{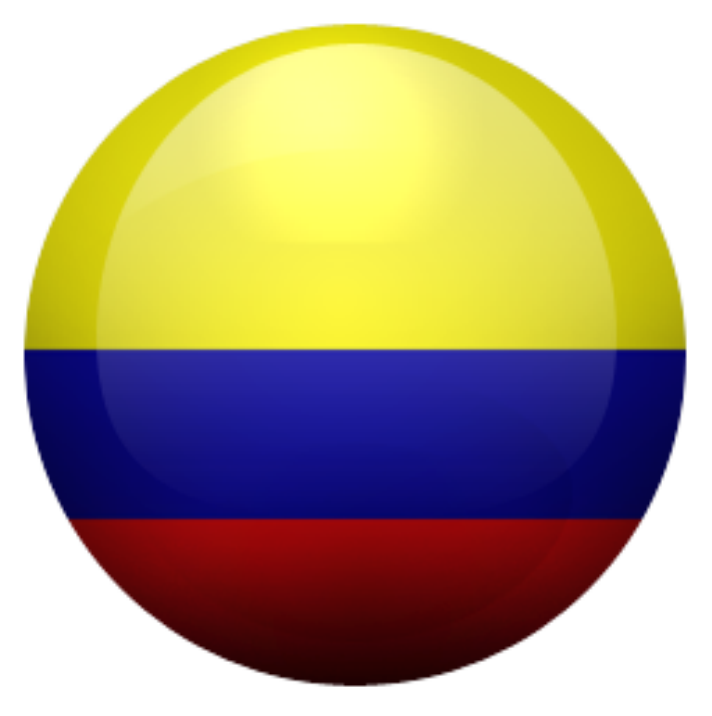} &
\includegraphics[height=0.7cm]{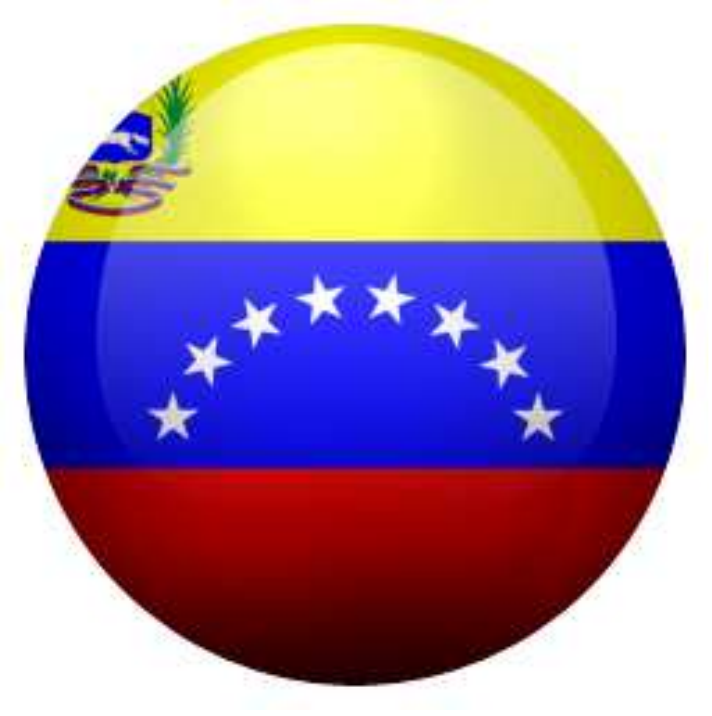} &
\includegraphics[height=0.7cm]{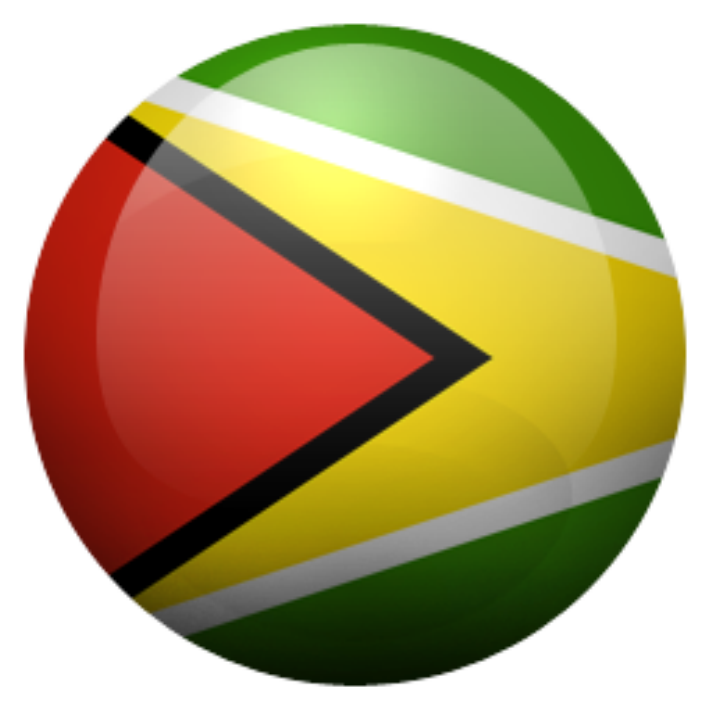} &
\includegraphics[height=0.7cm]{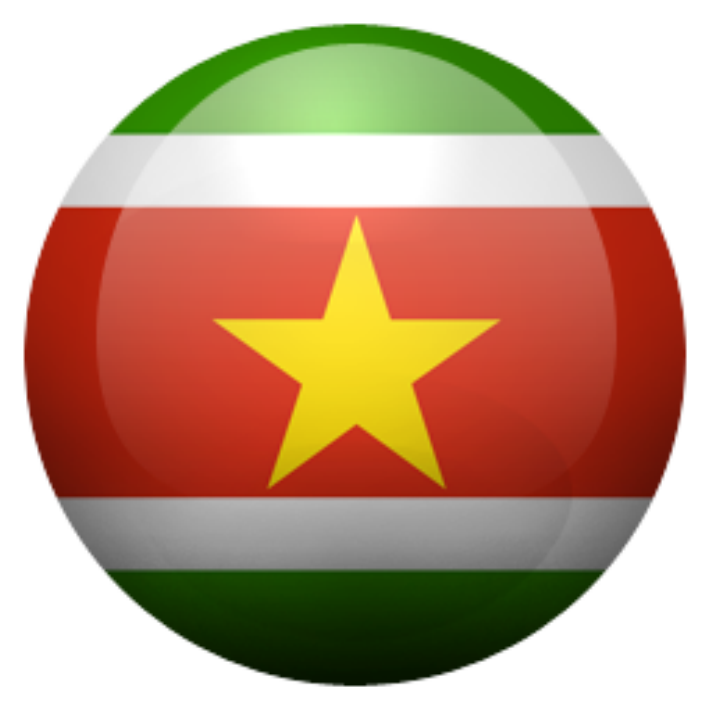} &
\includegraphics[height=0.7cm]{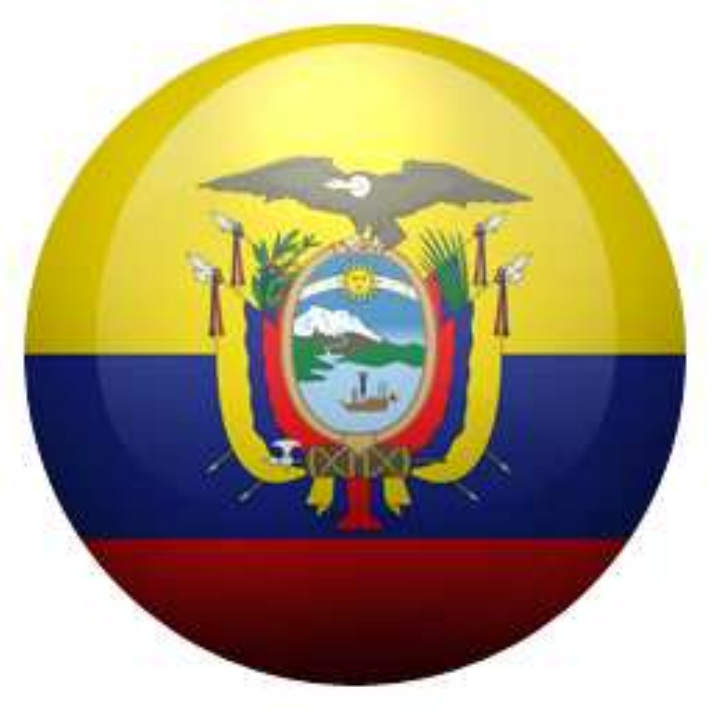} &
\includegraphics[height=0.7cm]{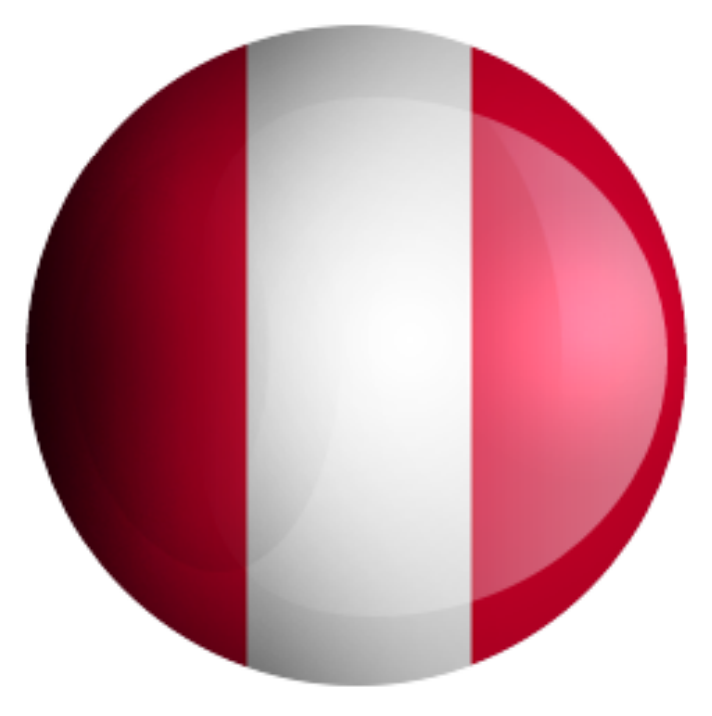} \\
\includegraphics[height=0.7cm]{./br.pdf} &
\includegraphics[height=0.7cm]{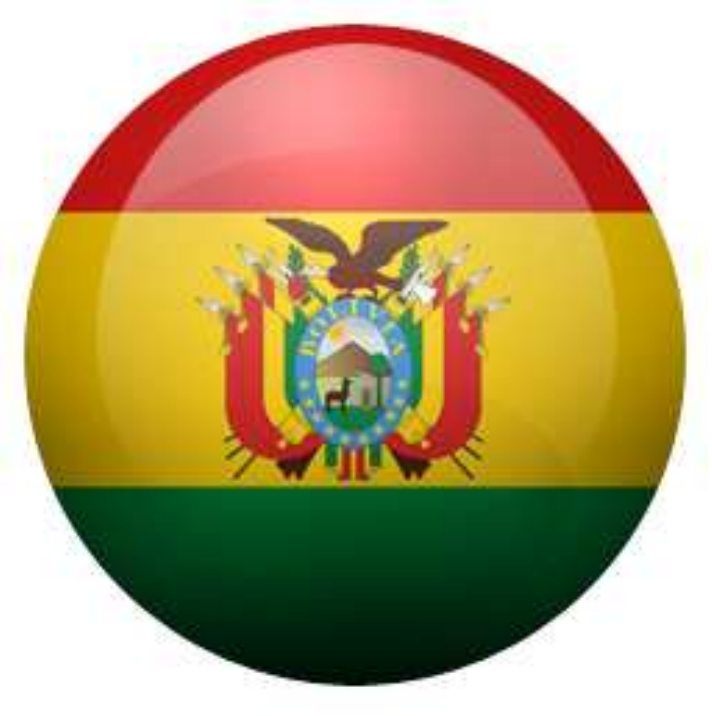} &
\includegraphics[height=0.7cm]{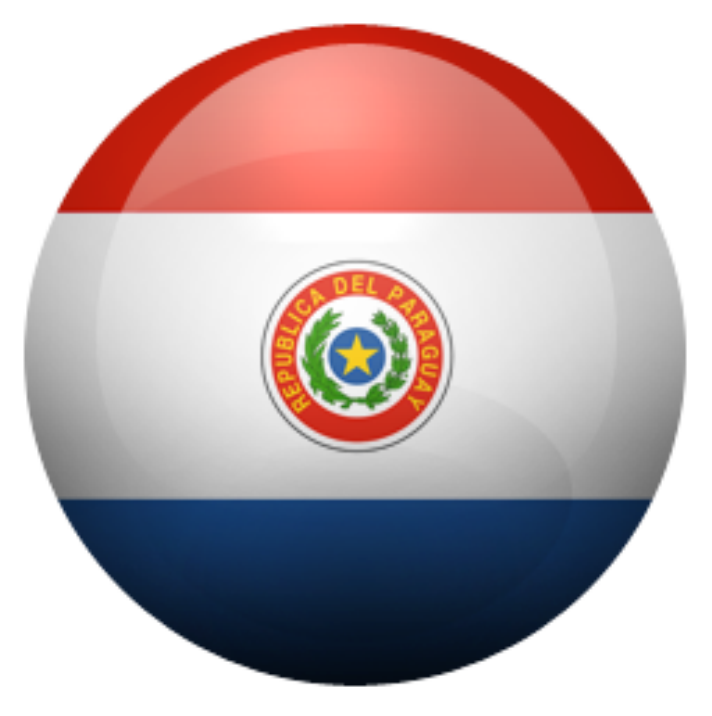} &
\includegraphics[height=0.7cm]{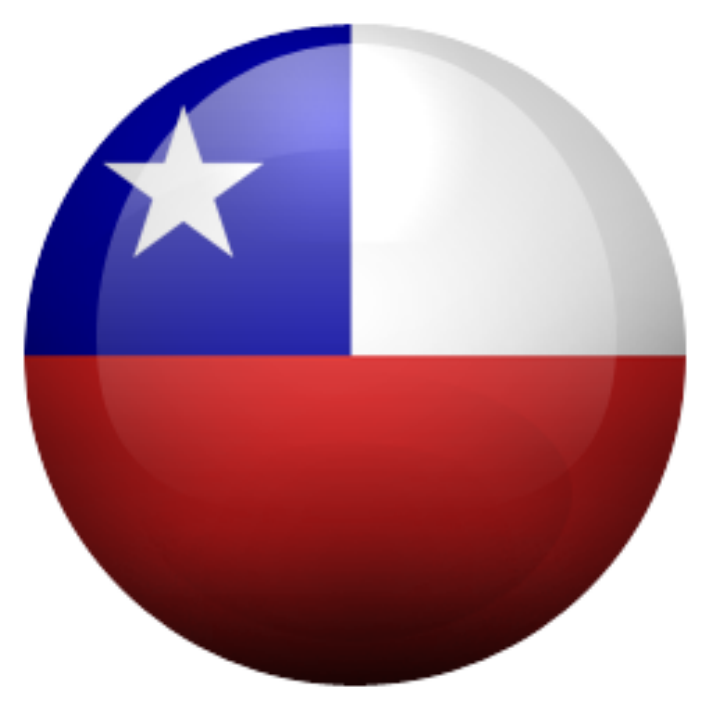} &
\includegraphics[height=0.7cm]{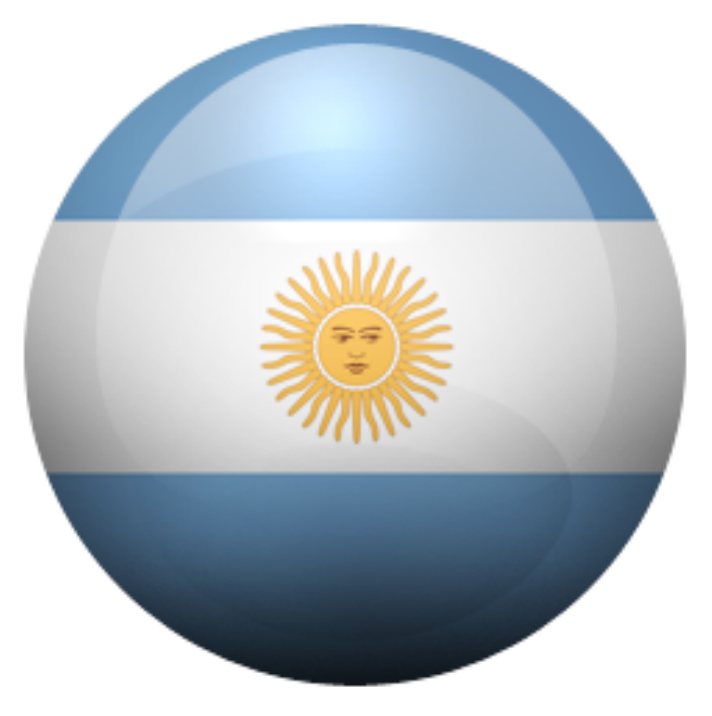} &
\includegraphics[height=0.7cm]{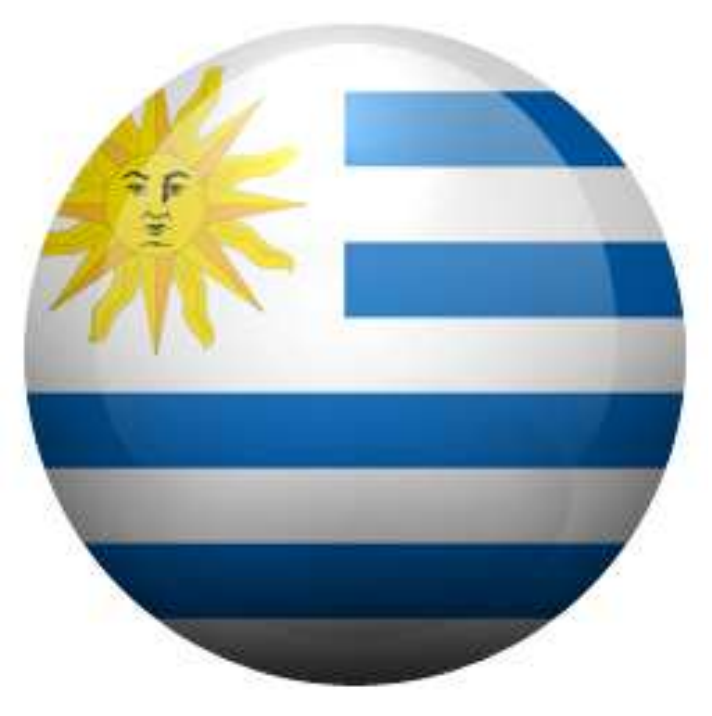} \\
\tiny{BRA} & 
\tiny{BOL} & 
\tiny{PAR} & 
\tiny{CHI} & 
\tiny{ARG} & 
\tiny{URU}\\
\end{tabular} 
}
\end{tabular}
\end{center}
\caption{Maps (left) and flags (bottom) of the Latin America countries. (right) Multidimensional Scaling of the HIM distances among the intertrade networks of the countries in the Latin America in the periods 1948--1959 (black), 1960--1969 (green), 1970--1990 (blue), 1991--2000 (red).}
\label{fig:southamerica}
\end{figure}

A very similar situation occurs in the regional trade network among the South American countries (Fig.~\ref{fig:southamerica}), where the global behaviour is essentially controlled by the two local giants Brazil and Argentina, and for which the larger differences between the nets can be appreciated between the economic growth of the 90s and the suffering economies in the late 70s / early 80s due to the struggling political situations.

\begin{figure}[!t]
\begin{center}
\begin{tabular}{cc}
\raisebox{0cm}{\includegraphics[height=6cm]{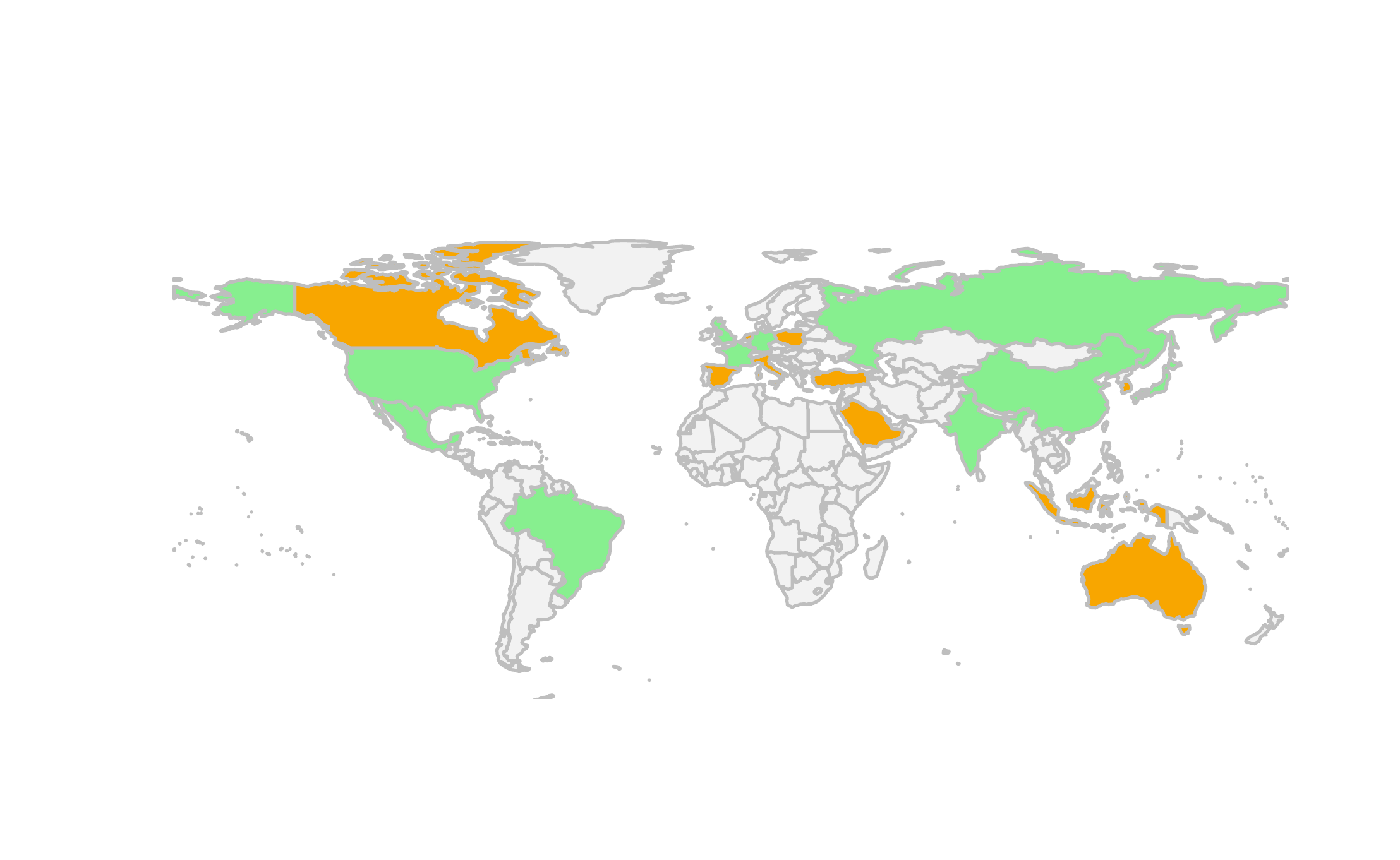}} & 
\raisebox{0cm}{\includegraphics[height=7cm]{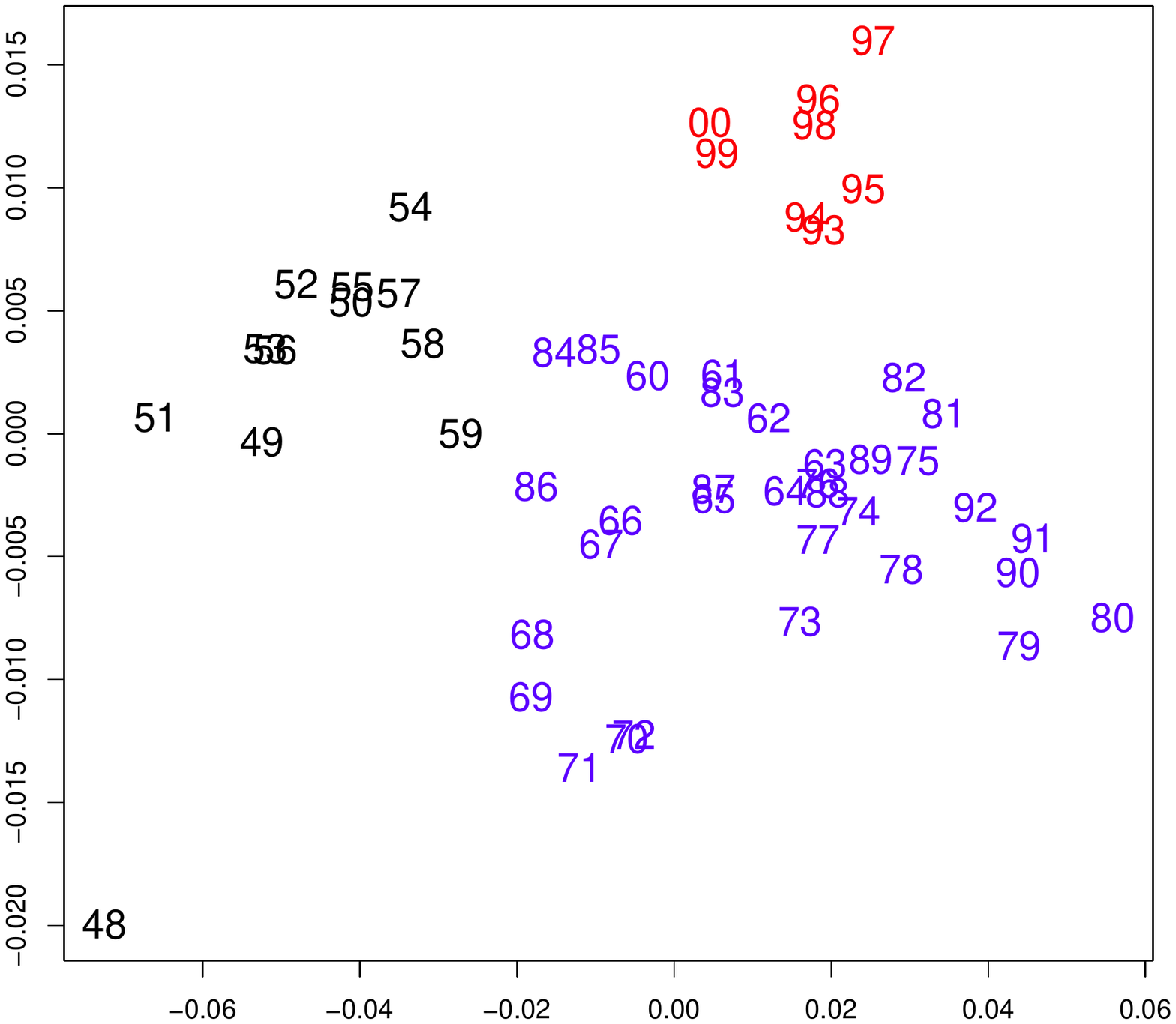}} \\
\multicolumn{2}{c}{
\begin{tabular}{cccccccccc}
\tiny{1 USA} & 
\tiny{2 CHN} & 
\tiny{3 IND} & 
\tiny{4 JPN} & 
\tiny{5 RUS} & 
\tiny{6 GER} & 
\tiny{7 FRA} & 
\tiny{8 BRA} & 
\tiny{9 UK} & 
\tiny{10 MEX}\\
\includegraphics[height=0.7cm]{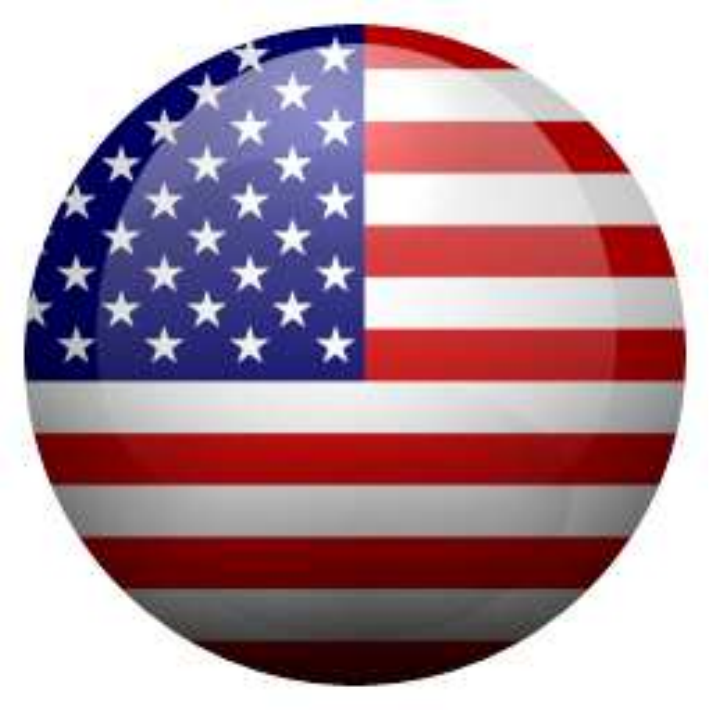} &
\includegraphics[height=0.7cm]{./cn.pdf} &
\includegraphics[height=0.7cm]{./in.pdf} &
\includegraphics[height=0.7cm]{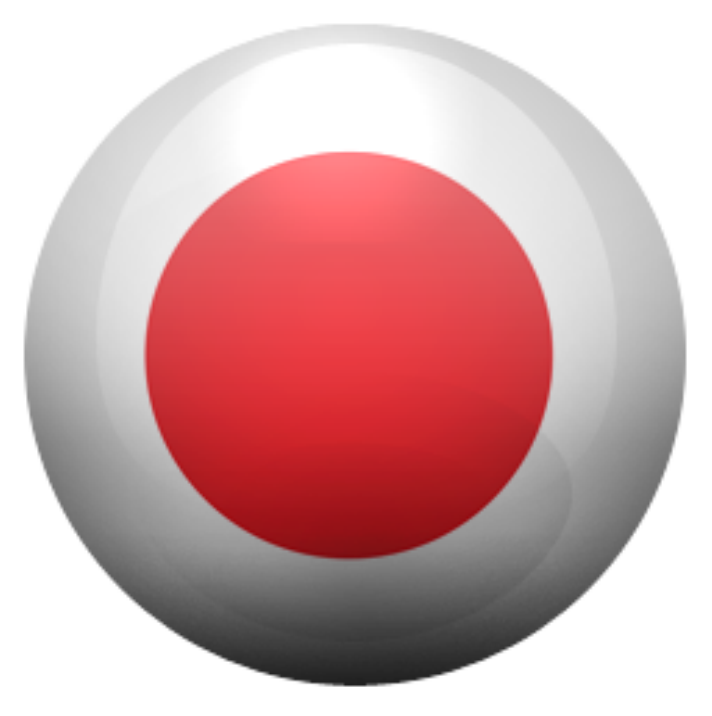} &
\includegraphics[height=0.7cm]{./ru.pdf} &
\includegraphics[height=0.7cm]{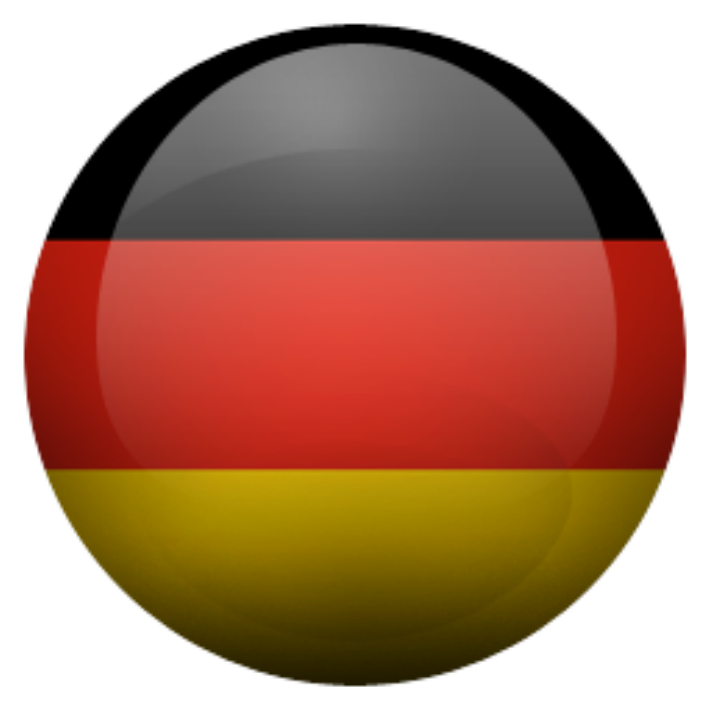} &
\includegraphics[height=0.7cm]{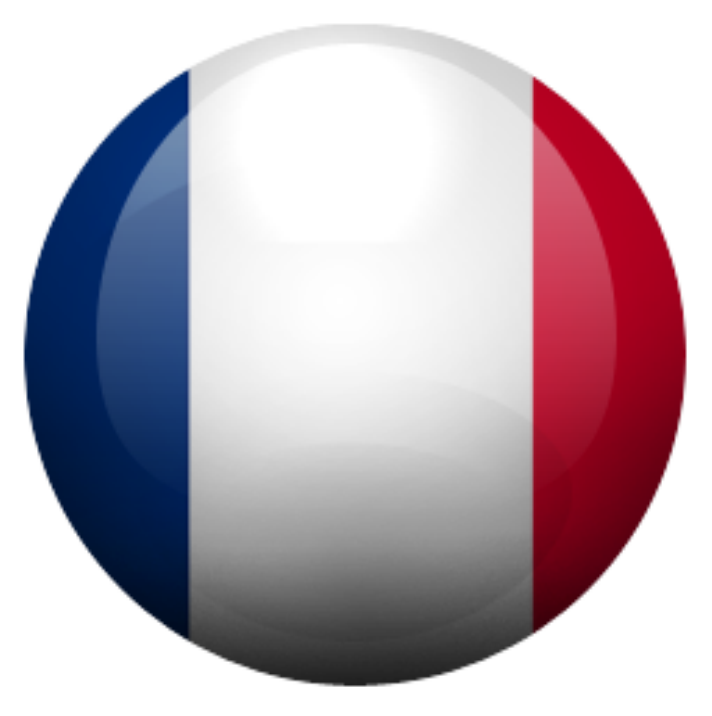} &
\includegraphics[height=0.7cm]{./br.pdf} &
\includegraphics[height=0.7cm]{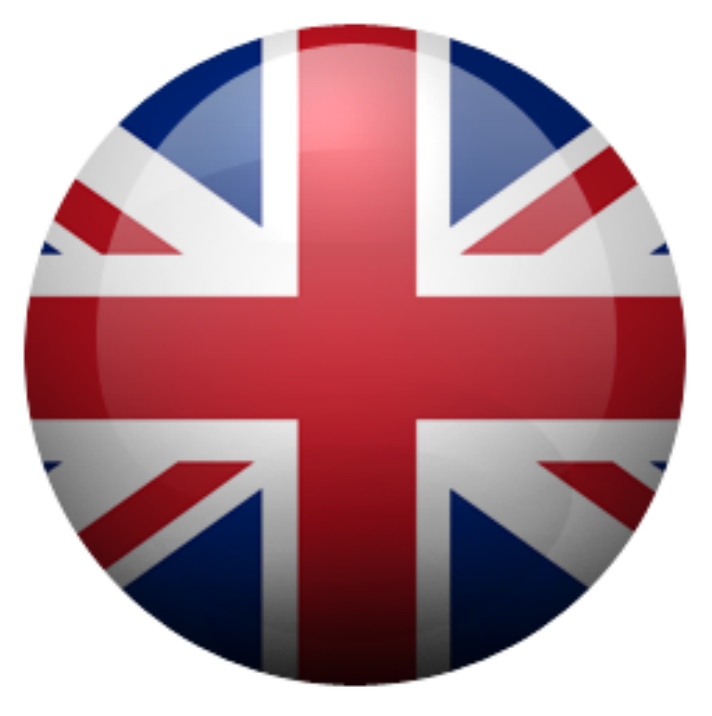} &
\includegraphics[height=0.7cm]{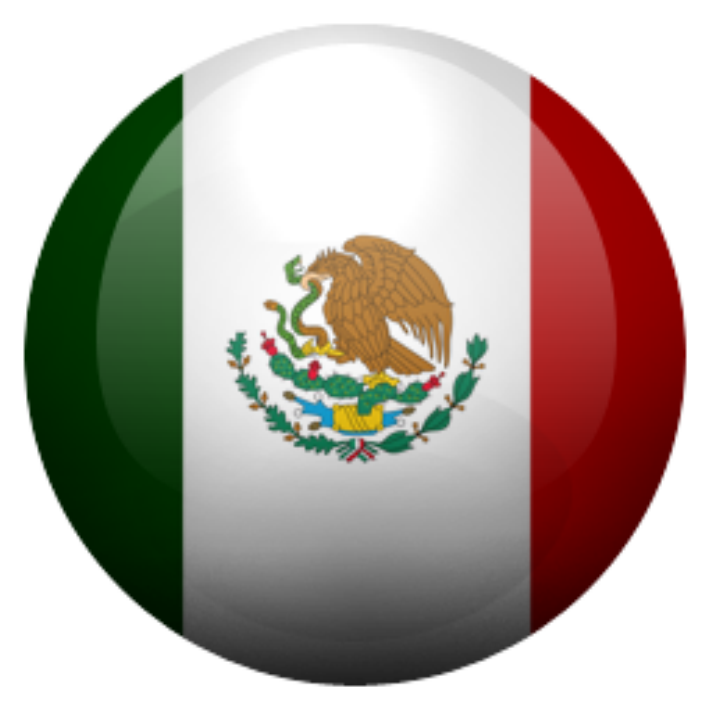} \\
\includegraphics[height=0.7cm]{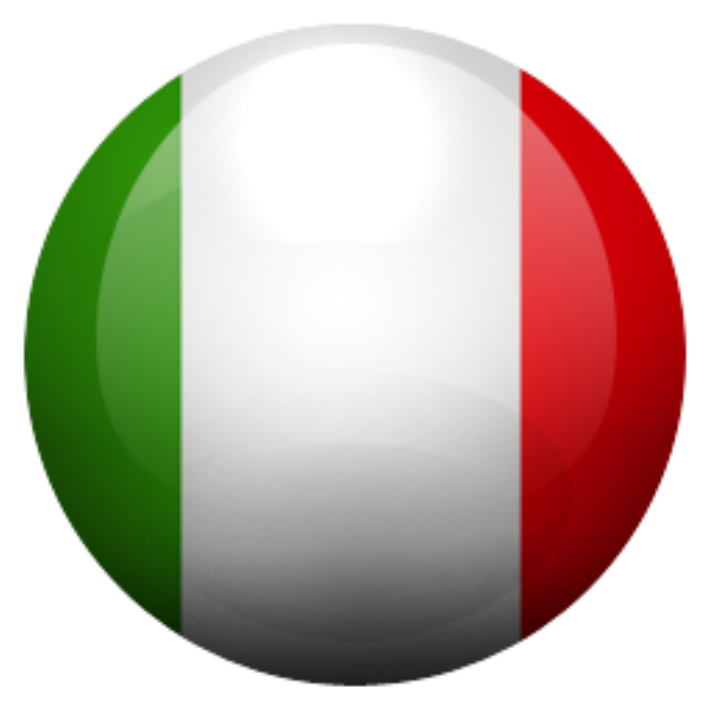} &
\includegraphics[height=0.7cm]{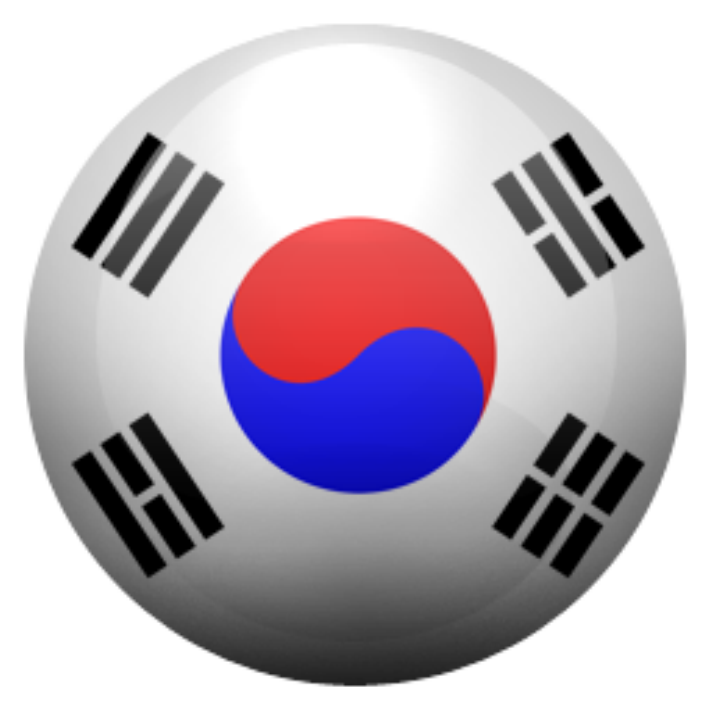} &
\includegraphics[height=0.7cm]{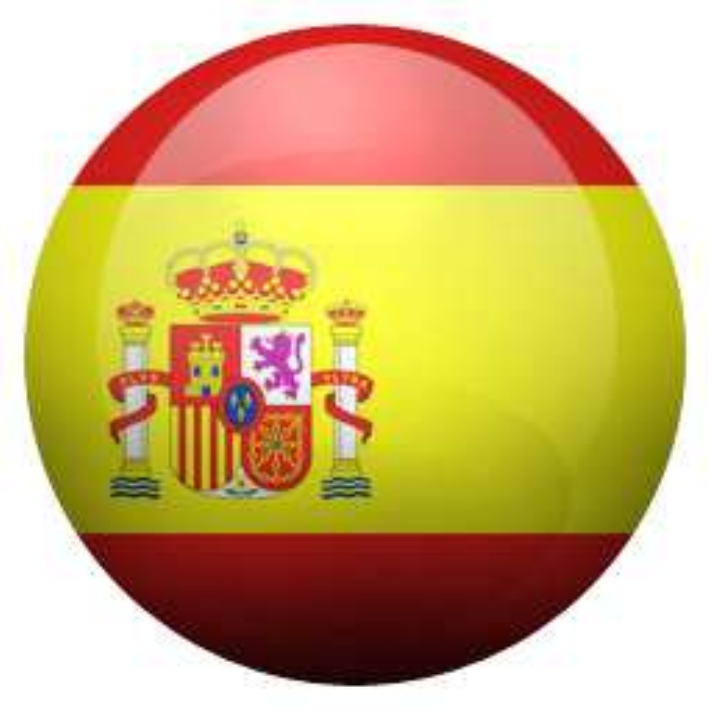} &
\includegraphics[height=0.7cm]{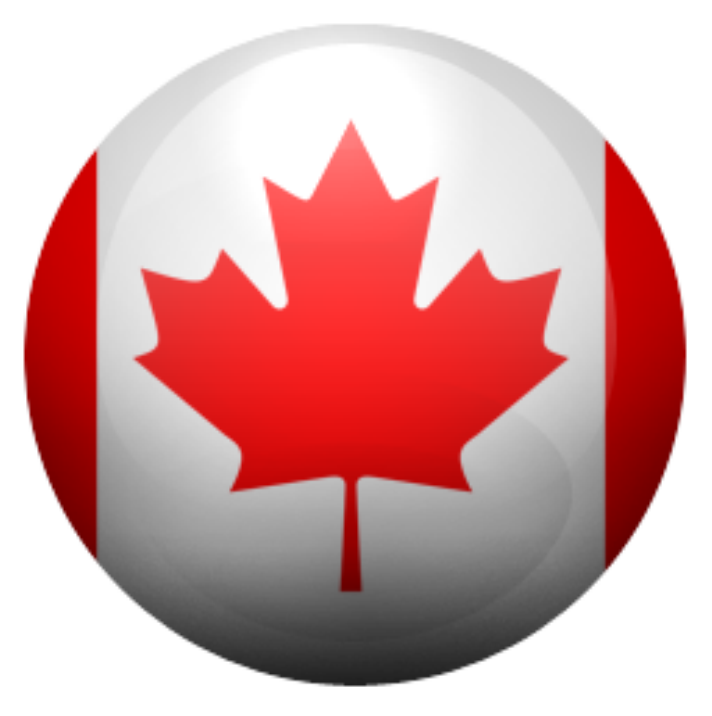} &
\includegraphics[height=0.7cm]{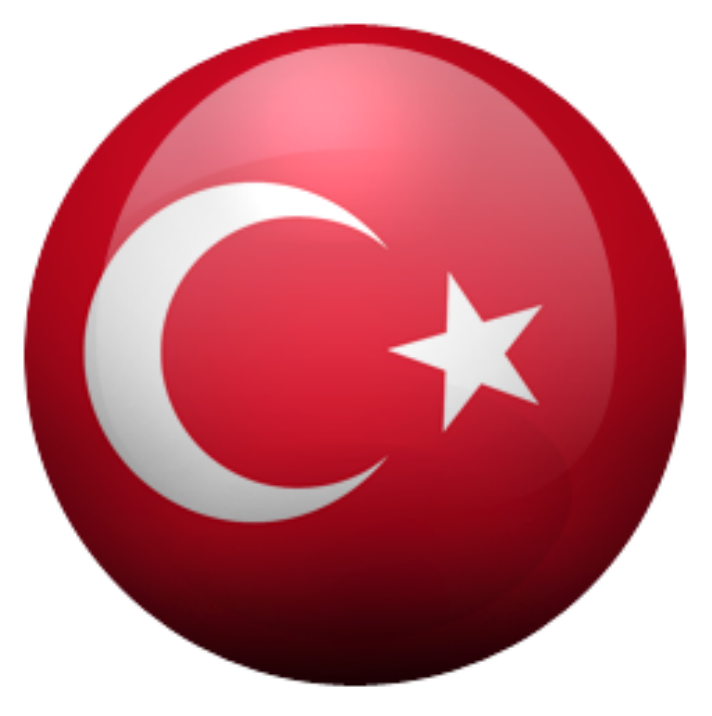} &
\includegraphics[height=0.7cm]{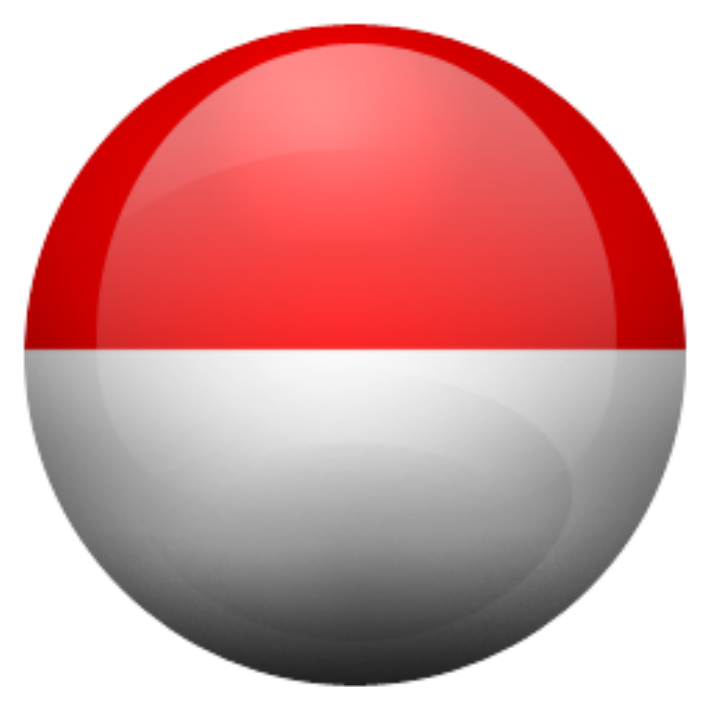} &
\includegraphics[height=0.7cm]{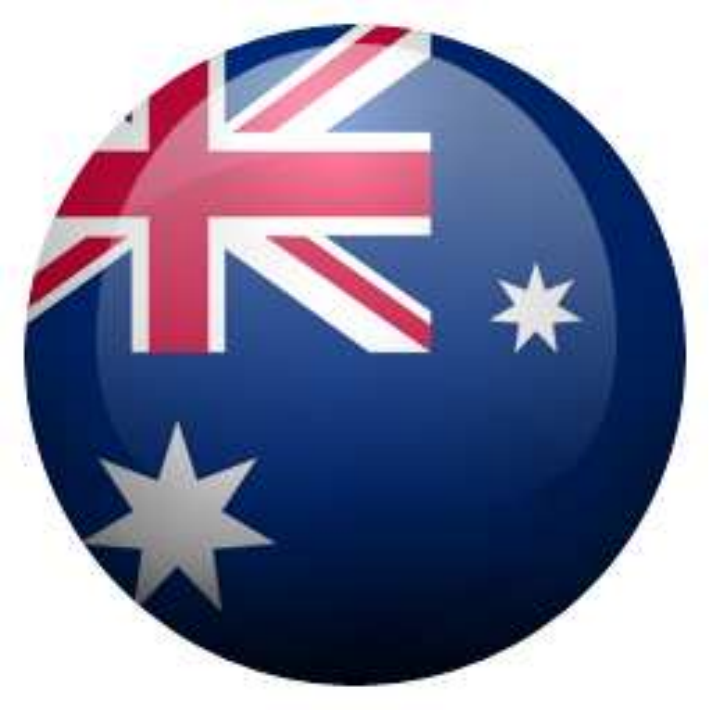} &
\includegraphics[height=0.7cm]{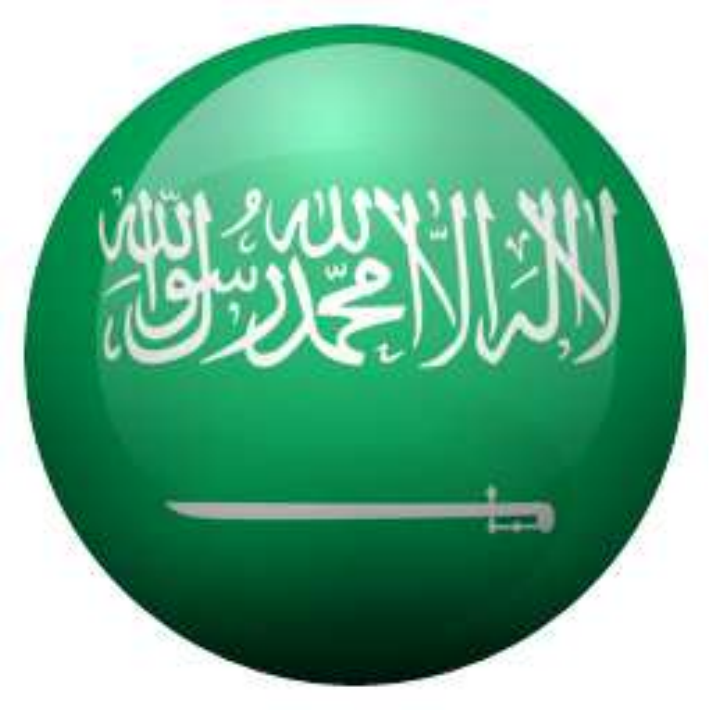} &
\includegraphics[height=0.7cm]{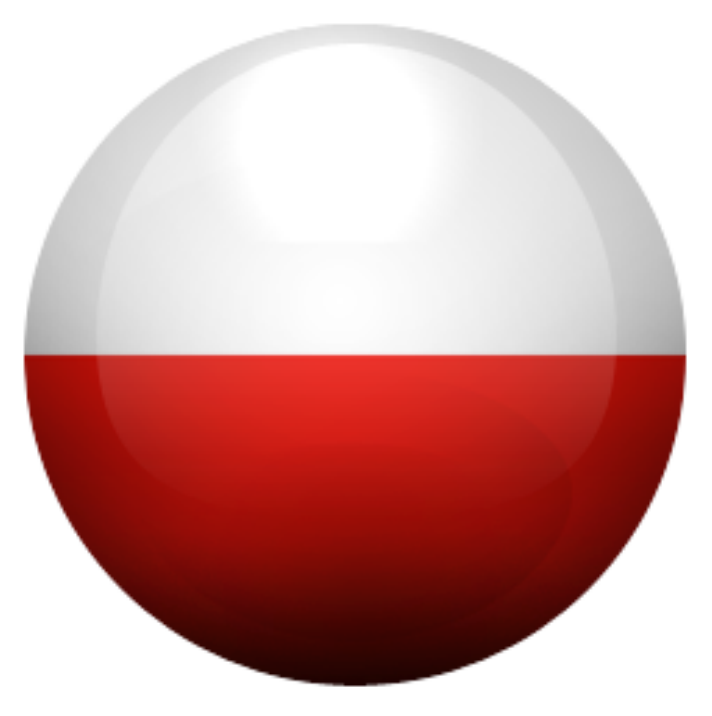} &
\includegraphics[height=0.7cm]{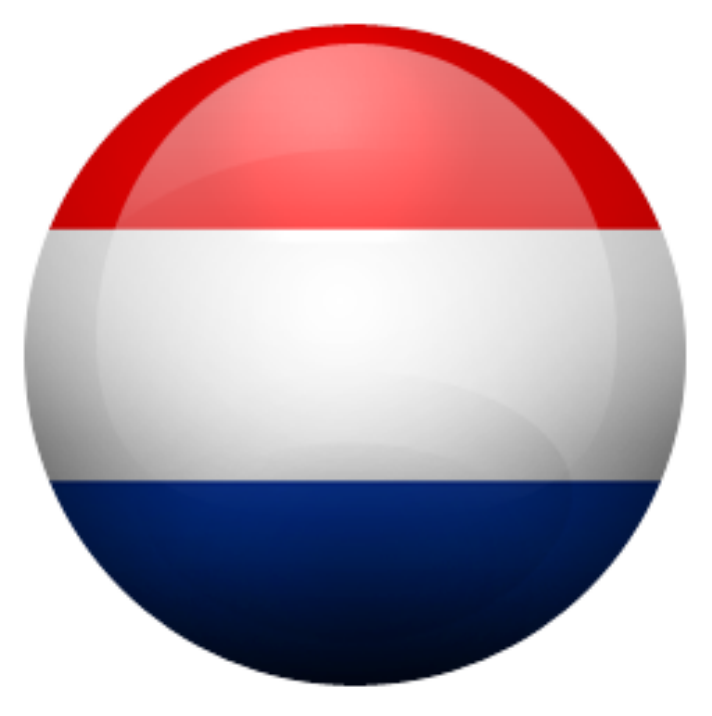} \\
\tiny{11 ITA} & 
\tiny{12 KOR} & 
\tiny{13 SPA} & 
\tiny{14 CAN} & 
\tiny{15 TUR} & 
\tiny{16 IDN} & 
\tiny{17 AUS} & 
\tiny{18 SAU} & 
\tiny{19 POL} & 
\tiny{20 NED}\\
\end{tabular} 
}
\end{tabular}
\end{center}
\caption{(left) Maps and flags (bottom) of the Top20 countries: in green, the top-10 economies, in orange the countries ranking 11--20. (right) Multidimensional Scaling of the HIM distances among the intertrade networks of the Top20 countries in the periods 1948--1959 (black), 1960--1992 (blue), 1993--2000 (red).}
\label{fig:top20}
\end{figure}

Not much different is the case of the larger trade subnetworks of the top 20 world economies ranked by Gross Domestic Product 2012 (PPP) (Top20 for short) as listed by the World Bank \url{http://data.worldbank.org} and shown in Fig.~\ref{fig:top20}, with the notable difference that the networks for the 60s are more homogeneous to those of the 70s and 80s, supporting a faster recovery of these economies after WWII than the BRICS or the South American countries.
Again, the 90s are remarkably separated by the previous periods, as a consequence of the fact that economic growth for high-income countries such as the United States, Japan, Singapore, Hong Kong, Taiwan, South Korea and Western Europe was steady and coupled with "an unprecedented extension and intensification of globalization in terms of the international integration of capital and product markets" \cite{crafts06world}, thus causing a structural evolution of the trade networks for these countries, whose economies account for approximately 85\% of the gross world product (GWP), 80 percent of world trade (including EU intra-trade), and two-thirds of the world population.

\begin{figure}[!b]
\begin{center}
\begin{tabular}{cccc}
\includegraphics[height=3cm]{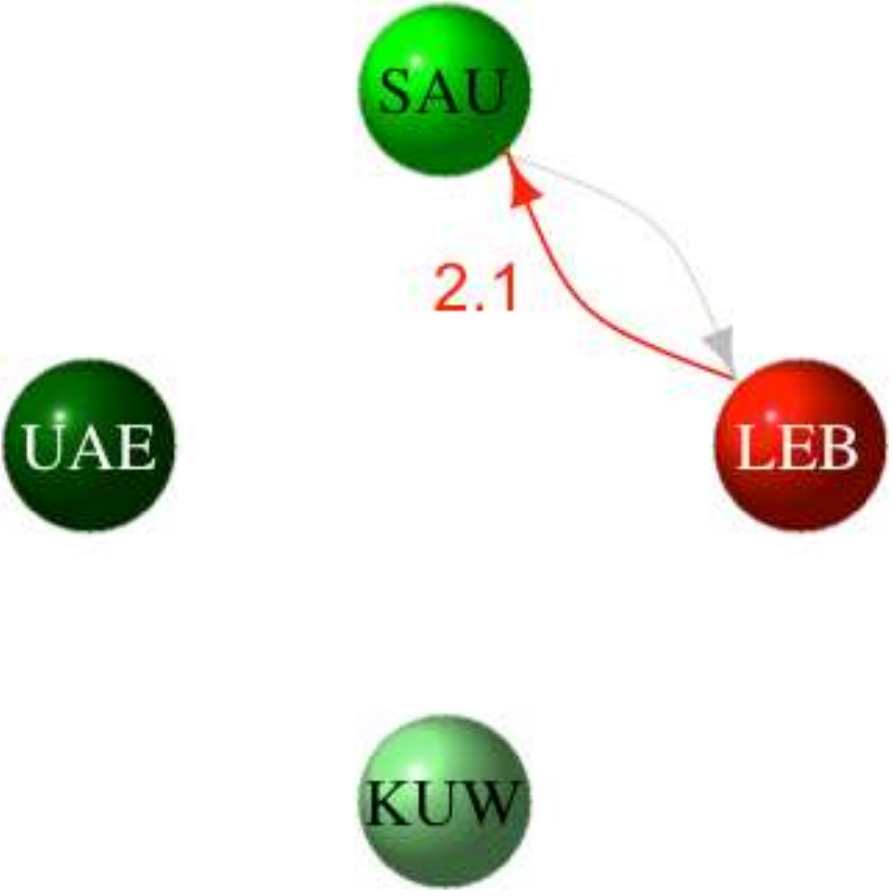}&
\includegraphics[height=3cm]{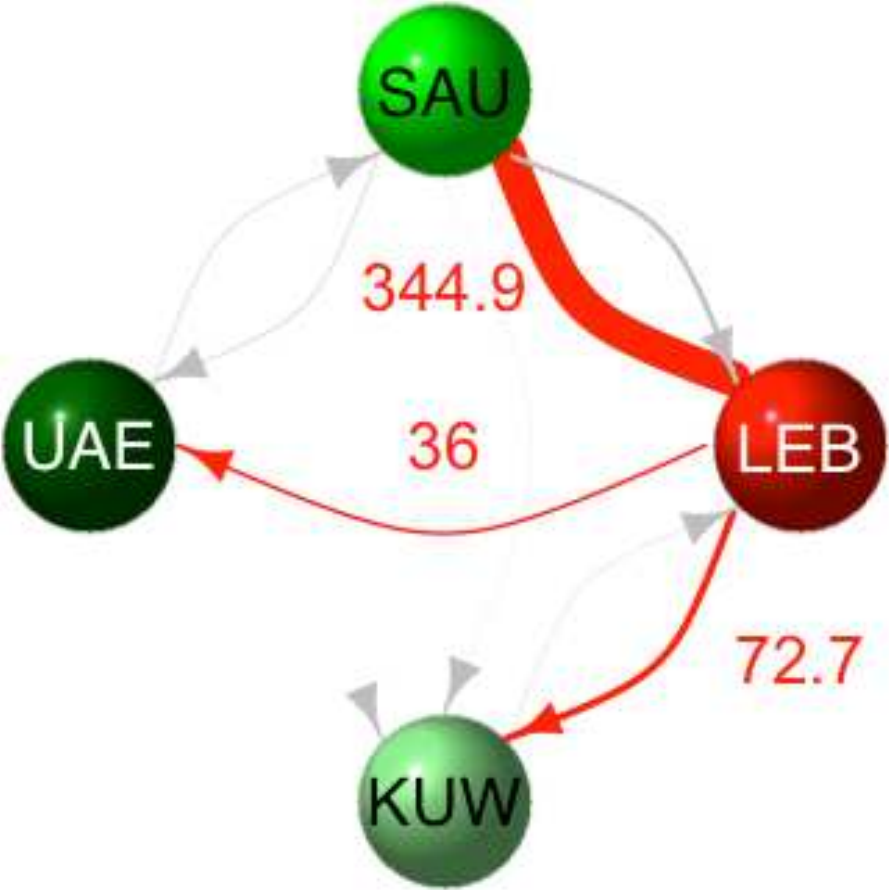}&
\includegraphics[height=3cm]{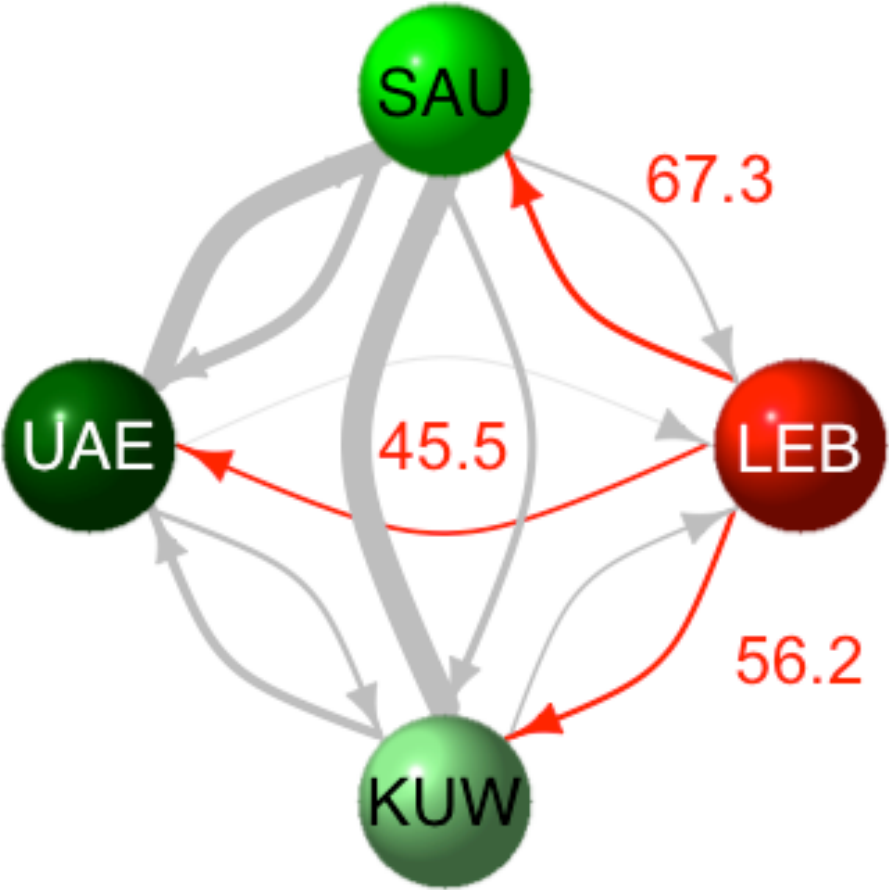}&
\includegraphics[height=3cm]{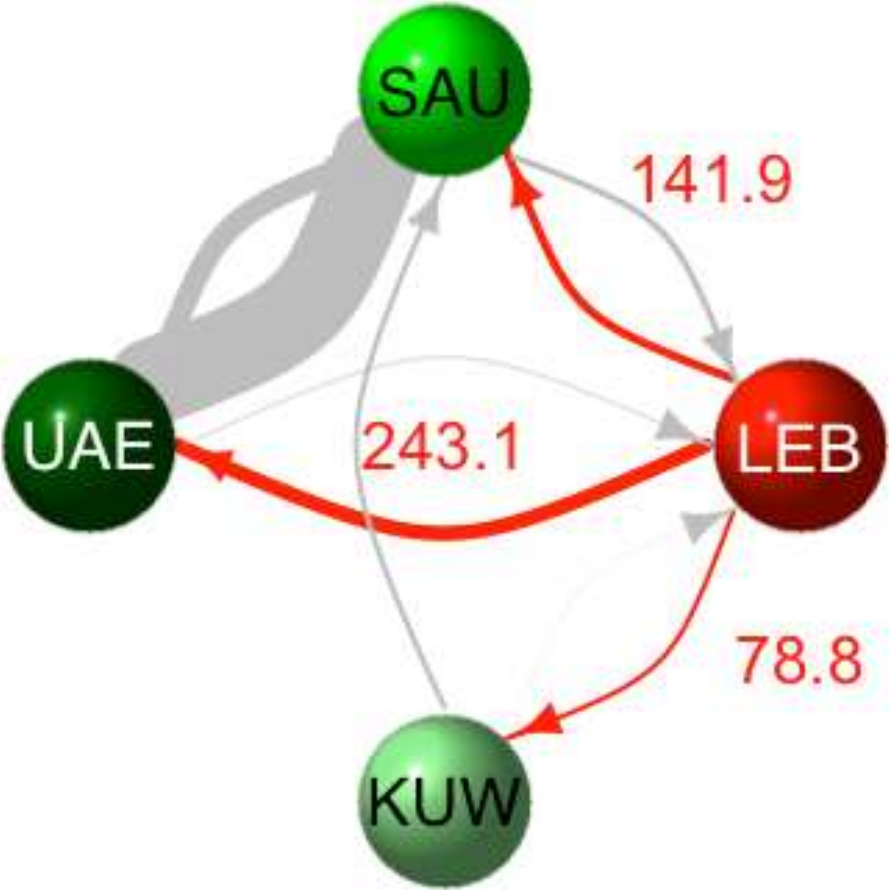}\\
1951 & 1974 & 1985 & 1996
\end{tabular}
\end{center}
\caption{The trade network between Lebanon, Saudi Arabia, Kuwait and United Arab Emirates in 1951, 1974, 1985 and 1996; red links indicate Lebanon export flow, with the corresponding volume figure. Edge width is proportional to export flow volume.}
\label{fig:leb_nets}
\end{figure}

\begin{figure}[!t]
\begin{center}
\begin{tabular}{cc}
\begin{tabular}{c}
\begin{tabular}{cccc}
\includegraphics[height=0.75cm]{./sa.pdf}&
\includegraphics[height=0.75cm]{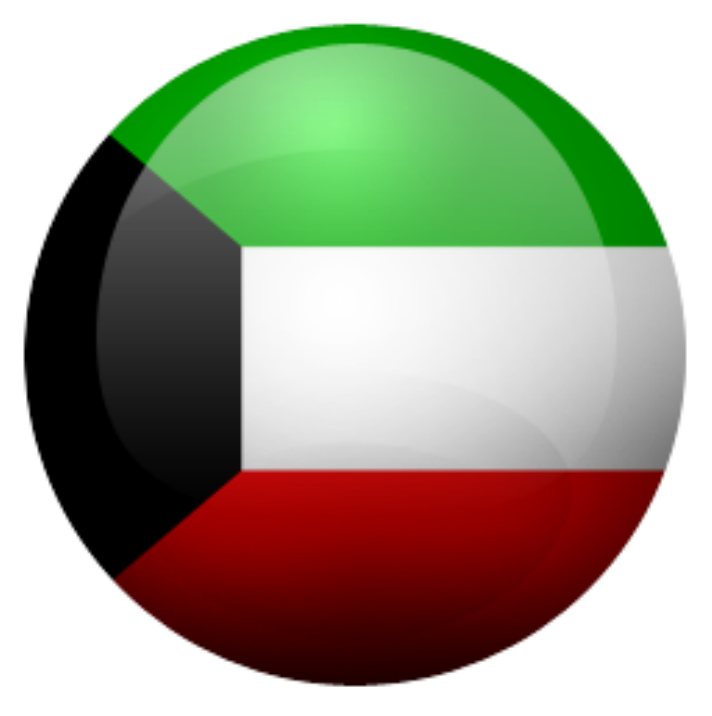}&
\includegraphics[height=0.75cm]{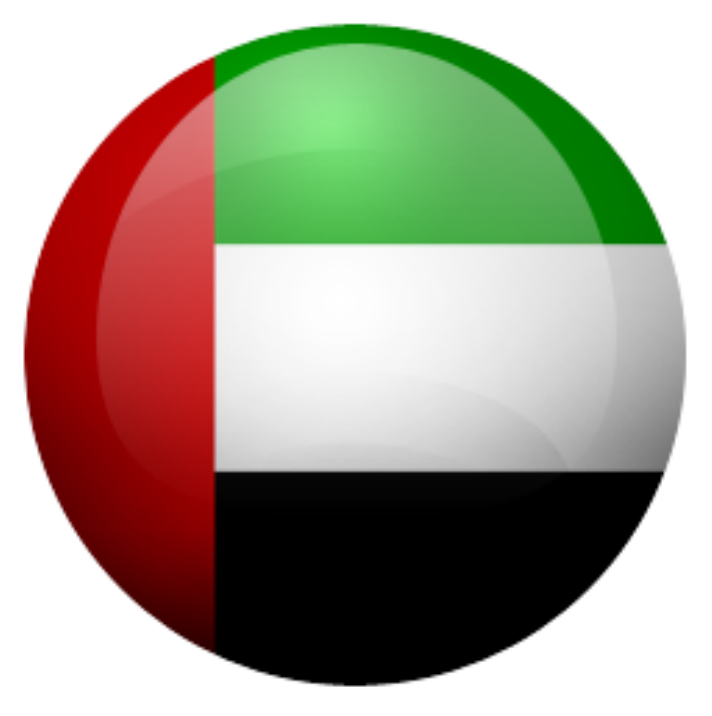}&
\includegraphics[height=0.75cm]{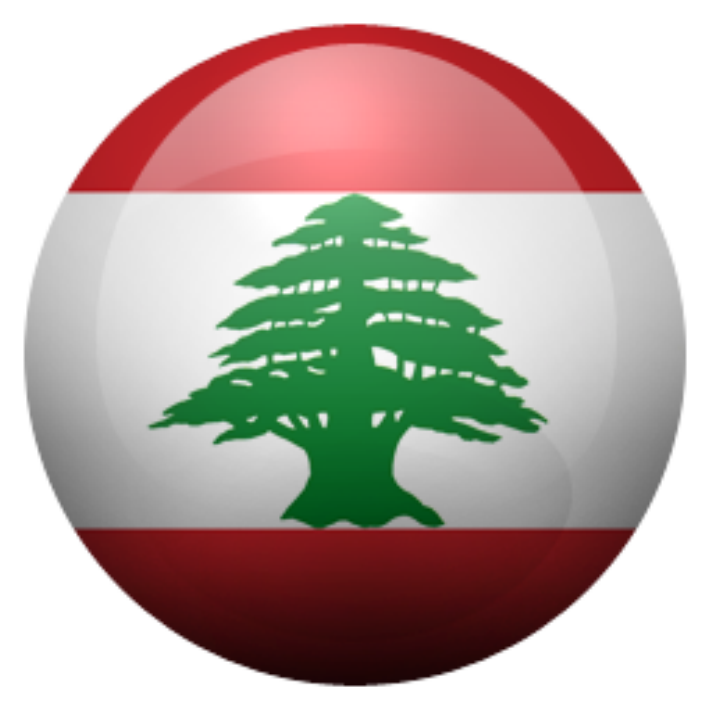}\\
\tiny{SAU} &
\tiny{KUW} &
\tiny{UAE} & 
\tiny{LEB} 
\end{tabular}
\\
\raisebox{6cm}{\includegraphics[height=4cm]{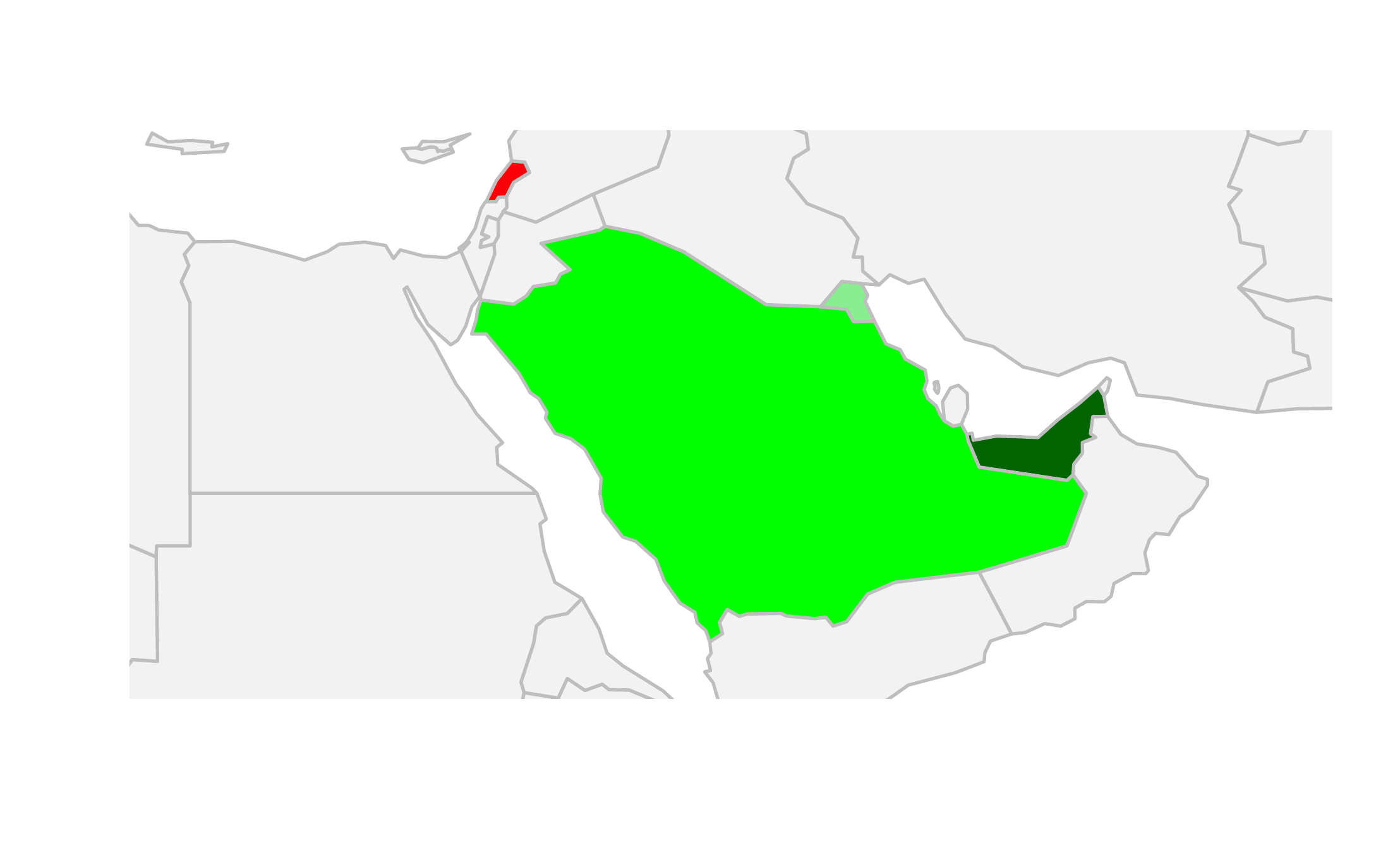}}\\
\end{tabular}
&
\raisebox{0cm}{\includegraphics[height=8cm]{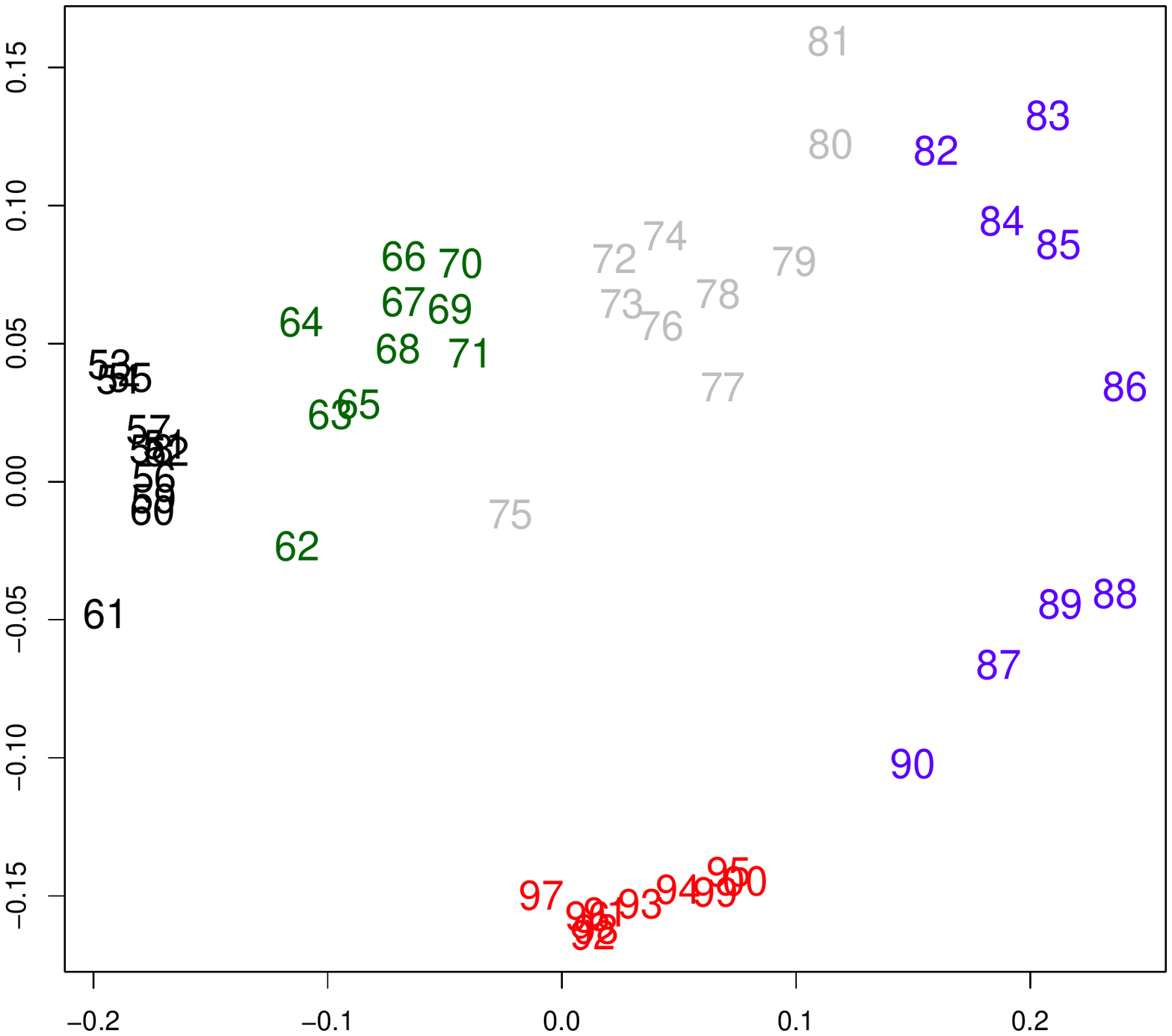}} \\
\end{tabular}
\caption{(left) Maps and flags of the countries in the Lebanon trade net . (top right) Multidimensional Scaling of the HIM distances among the intertrade networks of the Lebanon trade net countries in the periods 1950--1961 (black), 1962--1971 (green), 1971--1981 (gray), 1982--1990 (blue), 1991--2000 (red).}
\label{fig:lebanon}
\end{center}
\end{figure}

\begin{figure}[!t]
\begin{center}
\includegraphics[height=10cm]{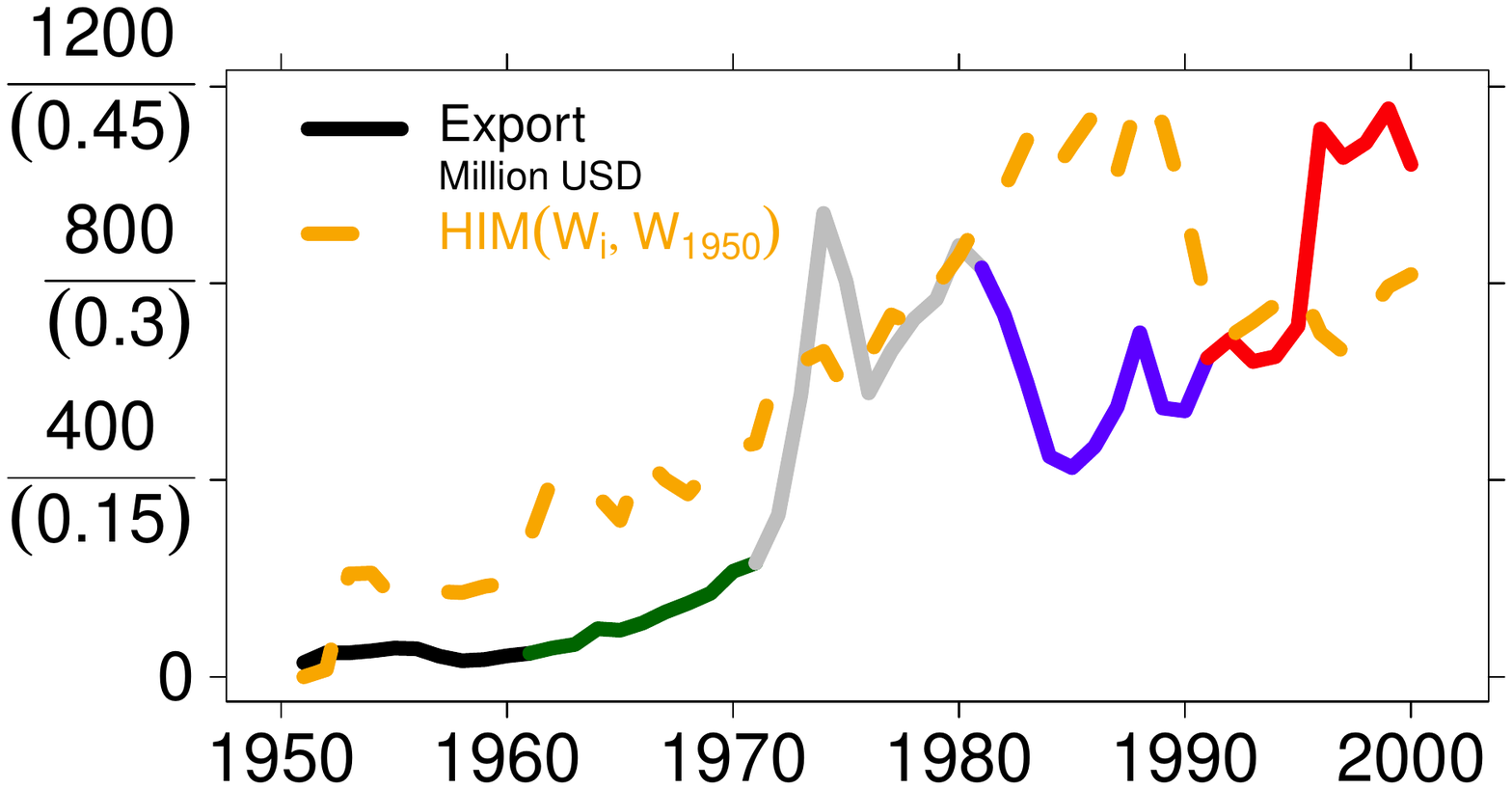}
\end{center}
\caption{Lebanon export flow in the years 1950--2000 (Million of USD, solid multicolor line left y axis) and curve of $\textrm{HIM}(W_i,W_{1950})$ (dashed orange line).}
\label{fig:lebanonex}
\end{figure}

We conclude with a more local example: between 1975 and 1990, the civil war heavily damaged Lebanon's economic infrastructure, reducing the role of the country as the major West Asian banking hub.
The following period of relative peace stimulated economic recovery also through an increasing flow of  manufactured and farm exports.
In this last example we consider the trading network $W$ between Lebanon and its three major economic partners, Saudi Arabia, Kuwait and United Arab Emirates.
In Fig.~\ref{fig:leb_nets} we show 4 examples of the trade networks with the Lebanon export figures.
In the bidimensional scaling plot of Fig.~\ref{fig:lebanon} the trajectory emerges of the evolution of the $W$ graphs across the different decades 50s, 60s, 70s, 80 and 90s, even more clearly than in the previous cases. 
Here the rightmost points in the plot, corresponding to the years 1977--1990, in the middle of the civil war in Lebanon, where a contraction of the trading flow was recorded.
Finally, in the plot of Fig.~\ref{fig:lebanonex}, we show the relation between the volume of export flow of Lebanon and the curve representing the HIM distance of $W_i$ from $W_{1950}$ for $i\in\{1951,\ldots,2000\}$. 
Pearson correlation between the two curves is 0.71, and their shape shows that the trade network is following the trend of the other curve with a temporal shift of about a decade.

\subsection{The MEG Biomag 2010 competition 1 dataset} 
\label{ssec:meg}
The challenge dataset for the 2010 Biomag competition was derived from \cite{vangerven09attention} and consisted in monitoring 4 subjects by a MEG in a set of trials where a fixation cross was presented to the subject and after that at regular intervals a cue indicated which direction, either left or right, they had to covertly attend to during the next 2500ms. After this period, a target in the indicated direction appeared. Brain activity was recorded from 500ms before cue offset to 2500ms after cue offset through 274 sensors at 300Hz; a total of 128 trials per condition were collected, 256 total trials per subject. 
MEG data were first preprocessed as explained in \cite{vangerven09attention}: the raw signals of each trial are independently decomposed with a multitaper frequency transformation in the 5-40  Hz interval with 2  Hz bin width. 
The results of the frequency transforms are used to construct a coherence network for each trial, which is successively rescaled such that its eigenvalues are between +1 and -1. 
After rescaling, on a separated instance of the dataset, the eigenvalues of the network are subjected to a network deconvolution procedure as explained in \cite{feizi13network}.
Finally, an Elastic Net \cite{zou05regularization} linear regression using the Lasso \cite{tibshirani96regression} in two phases with the mixed $\ell 1\ell 2$ algorithm \cite{demol09elastic,mosci10solving}, resulting in a final dataset of 252 covariance networks on 274 nodes, equally distributed between label ''right'' and ''left''.
A suite of Support Vector Machines from \textit{mlpy} \url{http://mlpy.fbk.eu} \cite{albanese12mlpy} with different $\textrm{HIM}_\xi$ kernels were tested, together with the linear kernel L-SVM, Random Forest RF \cite{breiman01random} and Elastic-Net EN as a baseline, on a set of 100 MonteCarlo resampling of stratified training (84+84 networks) and test (42+42) sets of both the deconvolved and the original dataset, yielding the performance shonw in Tab.~\ref{tab:acc}.

\begin{table}[!t]
\centering
\caption {Average Classification Accuracies for Deconvolved and non Deconvolved 11 Hz Networks, standard errors in brackets. As a baseline, authors in \cite{kia13discrete} reach 0.67 accuracy using the Elastic Net with summary statistics of spatio-temporal activations and 0.73 using 2-D DCT basis.\\}
\begin{tabular}{l|c|c|}
~       & Non-Deconvolved & Deconvolved\\ \hline
L-SVM & ~ 0.65 (0.02)           & 0.72 (0.02) ~          \\
$\textrm{HIM}_0$-SVM & ~ 0.67 (0.03)           & \textbf{0.74 (0.03)} ~          \\
$\textrm{HIM}_{+\infty}$-SVM & ~  0.56 (0.04)     & 0.48 (0.05) ~         \\ 
$\textrm{HIM}$-SVM & ~  0.61 (0.04)     & 0.63 (0.04) ~         \\ 
RF & ~  0.71 (0.01)     & 0.70 (0.02) ~          \\
EN  & ~     0.71 (0.03)    & \textbf{0.74 (0.03)} ~          \\
\end{tabular}
\label{tab:acc}
\end{table}
As a general consideration, the deconvolution procedure helps improving classification.
The better accuracy is reached by the kernel with only the Hamming component ($\textrm{HIM}_0$), which performs better of baseline methods, while the Ipsen-Mikhailov component ($\textrm{HIM}_{+\infty}$) performs very poorly, on both versions of the datasets. 
All intermediate values of $\xi$ (including $\xi=1$ reported in the table) gives decreasing performance for increasing values of $\xi$, implying that the topological features of the graph are not useful for classification in this task, maybe because of the symmetric nature of the task.

The obtained results of this off-the-shelf method are comparable with the range of performances obtained by far more complex and properly targeted approaches \cite{bahramisharif10covert,signoretto12classification,kia13discrete}, representing a promising starting point for an effective use of the HIM kernel, for instance coupled with other graph kernels or some feature selection techniques.
An extended version of this example can be found at \cite{furlanello13sparse}.
\section{Conclusions}
\label{sec:conclusion}
We introduced $\textrm{HIM}_{\xi}$, a novel family of distances between graphs with same nodes, even directed and weighted, aimed at combining the local and global aspects of the comparison between networks, \textit{i.e.}, the difference between matching vertices and the difference between the spectral structure.
After unveiling definitions and properties, we provided a range of applications in several fields, from functional genomics to economics, to show the usefulness of the proposed solution.
In particular, we underlined the effectiveness of the HIM metrics when used as a kernel functions for classification purposes, \textit{e.g.}, in Support Vector Machines, applied to heterogenous data in diverse areas.
A final comment on the computational feasibility: the costly part when computing the HIM distance is the extraction of the spectrum from the Laplacian matrices of the two compared graph. 
This task is both CPU intensive and requiring a fair amount of RAM, but allows for a wide parallelization: nonetheless, huge graphs should be dealt with HPC facilities.
As an example, the size of the largest graphs we compared (using a Python implementation making use of the NumPy library) is about 40,000 nodes: on a workstation with 48 Intel Xeon CPU E5649 at 2.53GHz and powered by 48Gb RAM we were able to run 4 parallel processes which took about 36 hours to compute the mutual distances between a set of 45 networks, for a total of 990 comparisons.

\bibliographystyle{unsrt}
\bibliography{jurman12glocal}

\begin{thebibliography}{100}

\bibitem{barabasi12network}
A.-L. Barab{\'a}si.
\newblock The network takeover.
\newblock {\em Nature Physics}, 8:14--16, 2012.

\bibitem{sharan06modeling}
R.~Sharan and T.~Ideker.
\newblock Modeling cellular machinery through biological network comparison.
\newblock {\em Nature Biotechnology}, 24(4):427--433, 2006.

\bibitem{ideker12differential}
T.~Ideker and N.J. Krogan.
\newblock Differential network biology.
\newblock {\em Molecular Systems Biology}, 8:565, 2012.

\bibitem{yoon12comparative}
B.-J. Yoon, X.~Qian, and S.M.E. Sahraeian.
\newblock {Comparative Analysis of Biological Networks}.
\newblock {\em {IEEE} Signal Processing Magazine}, 29(1):22--34, 2012.

\bibitem{csermely13structure}
P.~Csermely, T.~Korcsm{\'a}ros, H.J.M. Kiss, G.~London, and R.~Nussinov.
\newblock {Structure and dynamics of biological networks: a novel paradigm of
  drug discovery. A comprehensive review}.
\newblock {\em Pharmacology and Therapeutics}, 138:333--408, 2013.

\bibitem{chuang07network}
H.-Y. Chuang, E.~Lee, Y.-T. Liu, D.~Lee, and T.~Ideker.
\newblock Network-based classification of breast cancer metastasis.
\newblock {\em Molecular Systems Biology}, 3:140, 2007.

\bibitem{yang13network}
B.~Yang, J.~Zhang, Y.~Yin, and Y.~Zhang.
\newblock {Network-Based Inference Framework for Identifying Cancer Genes from
  Gene Expression Data}.
\newblock {\em BioMed Research International}, 2013:Article ID 401649, 2013.

\bibitem{pavlopoulos11using}
G.~Pavlopoulos, M.~Secrier, C.~Moschopoulos, T.~Soldatos, S.~Kossida, J.~Aerts,
  R.~Schneider, and P.~Bagos.
\newblock Using graph theory to analyze biological networks.
\newblock {\em BioData Mining}, 4(1):10, 2011.

\bibitem{barla12machine}
A.~Barla, G.~Jurman, R.~Visintainer, M.~Squillario, M.~Filosi, S.~Riccadonna,
  and C.~Furlanello.
\newblock {A Machine Learning Pipeline for Discriminant Pathways
  Identification}.
\newblock In E.~Biganzoli, A.~Vellido, F.~Ambrogi, and R.~Tagliaferri, editors,
  {\em {Computational Intelligence Methods for Bioinformatics and
  Biostatistics}}, volume 7548 of {\em Lecture Notes in Computer Science},
  pages 36--48. Springer, 2012.

\bibitem{barla13machine}
A.~Barla, G.~Jurman, R.~Visintainer, M.~Squillario, M.~Filosi, S.~Riccadonna,
  and C.~Furlanello.
\newblock {A Machine Learning Pipeline for Discriminant Pathways
  Identification}.
\newblock In N.K. Kasabov, editor, {\em {Springer Handbook of
  Bio-/Neuroinformatics}}, chapter~53, page 1200. Springer, Berlin, 2013.

\bibitem{xiao08structure}
Y.~Xiao, H.~Dong, W.~Wu, M.~Xiong, W.~Wang, and B.~Shi.
\newblock {Structure-based Graph Distance Measures of High Degree of
  Precision}.
\newblock {\em Pattern Recognition}, 41(12):3547--3561, 2008.

\bibitem{dehmer13discrimination}
M.~Dehmer and A.~Mowshowitz.
\newblock {The Discrimination Power of Structural SuperIndices}.
\newblock {\em Plos ONE}, 8(7):e70551, 2013.

\bibitem{entringer76distance}
R.C. Entringer, D.E. Jackson, and D.A. Snyder.
\newblock Distance in graphs.
\newblock {\em Czechoslovak Mathematical Journal}, 26(2):283--296, 1976.

\bibitem{zhu11classifying}
L.~Zhu, W.K. Ng, and S.~Han.
\newblock Classifying graphs using theoretical metrics: A study of feasibility.
\newblock In J.~Xu, G.~Yu, S.~Zhou, and R.~Unland, editors, {\em {Database
  Systems for Adanced Applications}}, volume 6637 of {\em Lecture Notes in
  Computer Science}, pages 53--64. Springer, 2011.

\bibitem{aliakbary13learning}
S.~Aliakbary, S.~Motallebi, J.~Habibi, and A.~Movaghar.
\newblock {Learning an Integrated Distance Metric for Comparing Structure of
  Complex Networks}.
\newblock arXiv:1307.3626v1 [cs.SI], 2013.

\bibitem{chen12discovery}
Z.~Chen.
\newblock {\em {Discovery of Informative and Predictive Patterns in Dynamic
  Networks of Complex Systems}}.
\newblock PhD thesis, North Carolina State University, 2012.

\bibitem{chen11identifying}
L.~Chen, J.~Xuan, R.~Riggins, R.~Clarke, and Y.~Wang.
\newblock {Identifying cancer biomarkers by network-constrained support vector
  machines}.
\newblock {\em BMC Systems Biology}, 5(1):161, 2011.

\bibitem{thorat13survey}
D.R. Thorat and S.S. Sonawane.
\newblock {A survey of graph classification approaches}.
\newblock In {\em Proceedings of International Conference on Computer Science
  and Information Technology ICSIT 2013}, pages 169--178. IRNet Explore, 2013.

\bibitem{mahe04extension}
P.~Mah{\'e}, N.~Ueda, T.~Akutsu, J.-L. Perret, and J.-P. Vert.
\newblock Extensions of marginalized graph kernels.
\newblock In {\em Proceedings of the 21 International Conference on Machine
  learning ICML '04}, page~70. ACM, 2004.

\bibitem{gaertner06short}
T.~G{\"a}rtner, Q.V. Le, and A.J. Smola.
\newblock {A Short Tour of Kernel Methods for Graphs}.
\newblock Technical report, Fraunhofer IAIS and National ICT Australia, 2006.

\bibitem{gaertner07kernel}
T.~G{\"a}rtner, T.~Horvath, Q.V. Le, A.J. Smola, and S.~Wrobel.
\newblock {Kernel methods for graphs}.
\newblock In D.J. Cook and L.B. Holder, editors, {\em {Mining Graph Data}},
  page 500. Wiley, 2007.

\bibitem{borgwardt07graph}
K.M. Borgwardt.
\newblock {\em {Graph kernels}}.
\newblock PhD thesis, Ludwig Maximilians Universitaet Muenchen, 2007.

\bibitem{ketkar09empirical}
N.S. Ketkar, L.B. Holder, and D.J. Cook.
\newblock {Empirical comparison of graph classification algorithms}.
\newblock In {\em {IEEE Symposium on Computational Intelligence and Data Mining
  CIDM '09}}, pages 259--266. IEEE, 2009.

\bibitem{vishwanathan10graph}
S.V.N. Vishwanathan, N.N. Schraudolph, R.~Kondor, and K.M. Borgwardt.
\newblock {Graph Kernels}.
\newblock {\em Journal of Machine Learning Research}, 11:1201--1242, 2010.

\bibitem{tsuda10graph}
K.~Tsuda and H.~Saigo.
\newblock {Graph Classification}.
\newblock In C.C. Aggarwal and H.~Wang, editors, {\em {Managing and Mining
  Graph Data}}, volume~40 of {\em Advances in Database Systems}, pages
  337--363. Springer, 2010.

\bibitem{vert05supervised}
J.-P. Vert and Y.~Yamanishi.
\newblock {Supervised graph inference}.
\newblock In L.K. Saul, Y.~Weiss, and L.~Bottou, editors, {\em Advances in
  Neural Information Processing Systems 17 NIPS 2004}, pages 1433--1440. MIT
  Press, 2005.

\bibitem{vert03graph}
J.-P. Vert and M.~Kanehisa.
\newblock {Graph-driven features extraction from microarray data using
  diffusion kernels and kernel CCA}.
\newblock In S.~Becker, S.~Thrun, and K.~Obermayer, editors, {\em Advances in
  Neural Information Processing Systems 15 NIPS 2002}, pages 1425--1432. MIT
  Press, 2003.

\bibitem{shervashidze11weisfeiler}
N.~Shervashidze, P.~Schweitzer, E.J. van Leeuwen, M.~Kurt, and K.M. Borgwardt.
\newblock Weisfeiler-lehman graph kernels.
\newblock {\em Journal of Machine Learning Research}, 12:2539--2561, 2011.

\bibitem{jie13integration}
B.~Jie, D.~Zhang, W.~Gao, Q.~Wang, C.-Y. Wee, and D.~Shen.
\newblock {Integration of Network Topological and Connectivity Properties for
  Neuroimaging Classification}.
\newblock IEEE Transactions on Biomedical Engineering, in press.

\bibitem{richiardi13machine}
J.~Richiardi, S.~Achard, H.~Bunke, and D.~Van De~Ville.
\newblock {Machine Learning with Brain Graphs: Predictive Modeling Approaches
  for Functional Imaging in Systems Neuroscience}.
\newblock {\em IEEE Signal Processing Magazine}, 30(3):58--70, 2013.

\bibitem{su13discriminative}
L.~Su, L.~Wang, H.~Shen, G.~Feng, and D.~Hu.
\newblock {Discriminative analysis of non-linear brain connectivity in
  schizophrenia: an fMRI Study}.
\newblock {\em Frontiers in Human Neurosciences}, 7:702, 2013.

\bibitem{dougherty10validation}
E.R. Dougherty.
\newblock Validation of gene regulatory networks: scientific and inferential.
\newblock {\em Briefings in Bioinformatics}, 12(3):245--252, 2010.

\bibitem{tun06metabolic}
K.~Tun, P.~Dhar, M.~Palumbo, and A.~Giuliani.
\newblock Metabolic pathways variability and sequence/networks comparisons.
\newblock {\em BMC Bioinformatics}, 7(1):24, 2006.

\bibitem{iwayama12characterizing}
K.~Iwayama, Y.~Hirata, K.~Takahashi, K.~Watanabe, K.~Aihara, and H.~Suzuki.
\newblock {Characterizing global evolutions of complex systems via intermediate
  network representations}.
\newblock {\em Nature Scientific Report}, 2:srep00423, 2012.

\bibitem{morris08specification}
M.~Morris, M.S. Handcock, and D.R. Hunter.
\newblock {Specification of Exponential-Family Random Graph Models: Terms and
  Computational Aspects}.
\newblock {\em Journal of Statistical Software}, 24(4):1--24, 2008.

\bibitem{ipsen02evolutionary}
M.~Ipsen and A.S. Mikhailov.
\newblock Evolutionary reconstruction of networks.
\newblock {\em Physics Review E}, 66(4):046109, 2002.

\bibitem{rajendran13analysis}
K.~Rajendran and I.G. Kevrekidis.
\newblock {Analysis of data in the form of graphs}.
\newblock arXiv:1306.3524v1 [physics.data-an], 2013.

\bibitem{jurman11introduction}
G.~Jurman, R.~Visintainer, and C.~Furlanello.
\newblock {An introduction to spectral distances in networks}.
\newblock {\em Frontiers in Artificial Intelligence and Applications},
  226:227--234, 2011.

\bibitem{kivela13multilayer}
M.~Kivel{\"a}, A.~Arenas, M.~Barthelemy, J.P. Gleeson, Y.~Moreno, and M.A.
  Porter.
\newblock {Multilayer Networks}.
\newblock arXiv:1309.7233v1 [physics.soc-ph], 2013.

\bibitem{dedomenico13mathematical}
M.~{De Domenico}, A.~Sol{\'e}-Ribalta, E.~Cozzo, M.~Kivel{\"a}, Y.~Moreno, M.A.
  Porter, S.~G{\'o}mez, and A.~Arenas.
\newblock {Mathematical Formulation of Multi-Layer Networks}.
\newblock arXiv:1307.4977v1 [physics.soc-ph], 2013.

\bibitem{sole-ribalta13spectral}
A.~Sol{\'e}-Ribalta, M.~{De Domenico}, N.E. Kouvaris, A.~D{\'\i}az-Guilera,
  S.~G{\'o}mez, and A.~Arenas.
\newblock {Spectral properties of the Laplacian of multiplex networks}.
\newblock arXiv:1307.2090v1 [physics.soc-ph], 2013.

\bibitem{sanchezgarcia13dimensionality}
R.J. S{\'a}nchez-Garc{\'\i}a, E.~Cozzo, and Y.~Moreno.
\newblock {Dimensionality reduction and spectral properties of multiplex
  networks}.
\newblock arXiv:1311.1759v1 [physics.soc-ph], 2013.

\bibitem{cortes03positive}
C.~Cortes, P.~Haffner, and M.~Mohri.
\newblock {Positive Definite Rational Kernels}.
\newblock In {\em Proceedings of The 16th Annual Conference on Computational
  Learning Theory COLT 2003}, pages 41--56. Springer, 2003.

\bibitem{shawe-taylor04kernel}
J.~Shawe-Taylor and N.~Cristianini.
\newblock {\em {Kernel Methods for Pattern Analysis}}.
\newblock Cambridge University Press, 2004.

\bibitem{kloft11lp}
M.~Kloft, U.~Brefeld, S.~Sonnenburg, and A.~Zien.
\newblock {{$\ell_p$}-Norm Multiple Kernel Learning}.
\newblock {\em Journal of Machine Learning Research}, 12:953--997, 2011.

\bibitem{R2013}
{R Core Team}.
\newblock {\em R: A Language and Environment for Statistical Computing}.
\newblock R Foundation for Statistical Computing, Vienna, Austria, 2013.

\bibitem{liu11controllability}
Y.-Y. Liu, J.-J. Slotine, and A.-L. Barab{\'a}si.
\newblock Controllability of complex networks.
\newblock {\em Nature}, 473(7346):167--173, 2011.

\bibitem{chung97spectral}
F.~Chung.
\newblock {\em {Spectral Graph Theory}}.
\newblock American Mathematical Society, 1997.

\bibitem{spielman09spectral}
D.A. Spielman.
\newblock {Spectral Graph Theory: The Laplacian (Lecture 2)}.
\newblock Lecture notes, 2009.

\bibitem{tonjes09perturbation}
R.~T\"{o}njes and B.~Blasius.
\newblock {Perturbation Analysis of Complete Synchronization in Networks of
  Phase Oscillators}.
\newblock arXiv:0908.3365, 2009.

\bibitem{atay06network}
F.M. Atay, T.~B{\i}y{\i}ko{\u g}lu, and J.~Jost.
\newblock {Network synchronization: Spectral versus statistical properties}.
\newblock {\em Physica D Nonlinear Phenomena}, 224:35--41, 2006.

\bibitem{rodgers05eigenvalue}
G.J. Rodgers, K.~Austin, B.~Kahng, and D.~Kim.
\newblock Eigenvalue spectra of complex networks.
\newblock {\em Journal of Physics A: Mathematical and General}, 38(43):9431,
  2005.

\bibitem{jost02evolving}
J.~Jost and M.P. Joy.
\newblock {Evolving Networks with distance preferences}.
\newblock {\em Phys. Rev. E}, 66:036126, 2002.

\bibitem{jost07dynamical}
J.~Jost.
\newblock {Dynamical Networks}.
\newblock In J.~Feng, J.~Jost, and M.~Qian, editors, {\em {Networks: From
  Biology to Theory}}, pages 35--64. Springer-Verlag, 2007.

\bibitem{almendral07dynamical}
J.A. Almendral and A.~D\'\i{}az-Guilera.
\newblock Dynamical and spectral properties of complex networks.
\newblock {\em New J. Phys.}, 9:187, 2007.

\bibitem{haemers04enumeration}
W.H. Haemers and E.~Spence.
\newblock Enumeration of cospectral graphs.
\newblock {\em Eur. J. Comb.}, 25(2):199--211, 2004.

\bibitem{deza09encyclopedia}
E.~Deza.
\newblock {\em Encyclopedia of Distances}.
\newblock Springer Verlag, 2009.

\bibitem{bolla13spectral}
M.~Bolla.
\newblock {\em {Spectral Clustering and Biclustering: Learning Large Graphs and
  Contingency Tables}}.
\newblock Wiley, 2013.

\bibitem{schoenberg38metric}
I.J. Schoenberg.
\newblock {Metric spaces and positive defined functions}.
\newblock {\em Transactions of the American Mathematical Society},
  44(3):522--536, 1938.

\bibitem{ressel76short}
P.~Ressel.
\newblock {A short proof of Schoenberg's theorem}.
\newblock {\em Proceedings of the American Mathematical Society}, 57(1):66--68,
  1976.

\bibitem{bekka08kazhdan}
B.~Bekka, P.~de~la Harpe, and A.~Valette.
\newblock {\em {Kazhdan's property (T)}}, volume~11 of {\em New Mathematical
  Monographs}.
\newblock Cambridge University Press, 2008.

\bibitem{berg84harmonic}
C.~Berg, J.P.R. Christensen, and P.~Ressel.
\newblock {\em {Harmonic Analysis on Semigroups: Theory of Positive Definite
  and Related Functions}}, volume 100 of {\em Graduate Texts in Mathematics}.
\newblock Springer, 1984.

\bibitem{martins06generative}
A.~Martins.
\newblock {Generative kernels}.
\newblock Technical report, Instituto Superior Técnico, UTL Lisbon, 2006.

\bibitem{neuhaus07bridging}
M.~Neuhaus and H.~Bunke.
\newblock {\em {Bridging the Gap between Graph Edit Distance and Kernel
  Machines}}, volume~68 of {\em Series in Machine Perception and Artificial
  Intelligence}.
\newblock World Scientific, 2007.

\bibitem{li05class}
H.~Li and T.~Jiang.
\newblock {A Class of Edit Kernels for SVMs to Predict Translation Initiation
  Sites in Eukaryotic mRNAs}.
\newblock {\em Journal of Computational Biology}, 12(6):702--718, 2005.

\bibitem{cuturi09positive}
M.~Cuturi.
\newblock {Positive Definite Kernels in Machine Learning}.
\newblock arXiv:0911.5367v2 [stat.ML], 2009.

\bibitem{schoelkopf97support}
B.~Sch{\"o}lkopf.
\newblock {\em {Support Vector Learning}}.
\newblock Oldenbourg, Muenchen, 1997.

\bibitem{sonnenburg05large}
S.~Sonnenburg, G.~R\"{a}tsch, and B.~Sch\"{o}lkopf.
\newblock {Large Scale Genomic Sequence SVM Classifiers}.
\newblock In {\em Proceedings of the 22nd International Conference on Machine
  Learning ICML '05}, pages 848--855. ACM, 2005.

\bibitem{stolovitzky07dialogue}
G.~Stolovitzky, D.~Monroe, and A.~Califano.
\newblock {Dialogue on Reverse-Engineering Assessment and Methods: The DREAM of
  High-Throughput Pathway Inference}.
\newblock {\em Annals of the New York Academy of Sciences}, 1115:11--22, 2007.

\bibitem{marbach10revealing}
D.~Marbach, R.J. Prill, T.~Schaffter, C.~Mattiussi, D.~Floreano, and
  G.~Stolovitzky.
\newblock Revealing strengths and weaknesses of methods for gene network
  inference.
\newblock {\em PNAS}, 107(14):6286--6291, 2010.

\bibitem{prill10towards}
R.J. Prill, D.~Marbach, J.~Saez-Rodriguez, P.K. Sorger, L.G. Alexopoulos,
  X.~Xue, N.D. Clarke, G.~Altan-Bonnet, and G.~Stolovitzky.
\newblock {Towards a Rigorous Assessment of Systems Biology Models: The DREAM3
  Challenges}.
\newblock {\em PLoS ONE}, 5(2):e9202, 02 2010.

\bibitem{matthews75comparison}
B.W. Matthews.
\newblock {Comparison of the predicted and observed secondary structure of T4
  phage lysozyme}.
\newblock {\em Biochimica et Biophysica Acta - Protein Structure},
  405(2):442--451, 1975.

\bibitem{baldi00assessing}
P.~Baldi, S.~Brunak, Y.~Chauvin, C.A.F. Andersen, and H.~Nielsen.
\newblock {Assessing the accuracy of prediction algorithms for classification:
  an overview}.
\newblock {\em Bioinformatics}, 16(5):412--424, 2000.

\bibitem{supper07reconstructing}
J.~Supper, C.~Spieth, and A.~Zell.
\newblock {Reconstructing Linear Gene Regulatory Networks}.
\newblock In E.~Marchiori, J.H. Moore, and J.C. Rajapakse, editors, {\em
  Proceedings of the 5th European Conference on Evolutionary Computation,
  Machine Learning and Data Mining in Bioinformatics, EvoBIO2007, LNCS 4447},
  pages 270--279. Springer-Verlag, 2007.

\bibitem{stokic09fast}
D.~Stokic, R.~Hanel, and S.~Thurner.
\newblock A fast and efficient gene-network reconstruction method from multiple
  over-expression experiments.
\newblock {\em BMC Bioinformatics}, 10(1):253, 2009.

\bibitem{erdos59random}
P.~Erd{\"o}s and A.~R{\'e}nyi.
\newblock {On Random Graphs. I}.
\newblock {\em Publicationes Mathematicae}, 6:290--–297, 1959.

\bibitem{erdos60evolution}
P.~Erd{\"o}s and A.~R{\'e}nyi.
\newblock {On the evolution of random graphs}.
\newblock {\em Publications of the Mathematical Institute of the Hungarian
  Academy of Sciences}, 5:17--–61, 1960.

\bibitem{barabasi99emergence}
A.-L. Barabasi and R.~Albert.
\newblock {Emergence of scaling in random networks}.
\newblock {\em Science}, 286:509--–512, 1999.

\bibitem{goh01universal}
K.-I. Goh, B.~Kahng, and D.~Kim.
\newblock {Universal behaviour of load distribution in scale-free networks}.
\newblock {\em Physical Review Letters}, 87(27):278701, 2001.

\bibitem{watts98collective}
D.J. Watts and S.H. Strogatz.
\newblock {Collective dynamics of ‘small world’ networks}.
\newblock {\em Nature}, 393:440--442, 1998.

\bibitem{cox01multidimensional}
T.F. Cox and M.A.A. Cox.
\newblock {\em {Multidimensional Scaling}}.
\newblock Chapman and Hall, 2001.

\bibitem{kolar10estimating}
M.~Kolar, L.~Song, A.~Ahmed, and E.P. Xing.
\newblock {Estimating time-varying networks}.
\newblock {\em Ann. Appl. Stat.}, 4(1):94--123, 2010.

\bibitem{arbeitman02gene}
M.~Arbeitman, E.~Furlong, F.~Imam, E.~Johnson, B.~Null, B.~Baker, M.~Krasnow,
  M.~Scott, R.~Davis, and K.~White.
\newblock {Gene expression during the life cycle of \emph{Drosophila
  melanogaster}}.
\newblock {\em Science}, 297:2270–--2275, 2002.

\bibitem{budhu08identification}
A.~Budhu, H.-L. Jia, M.~Forgues, C.-G. Liu, D.~Goldstein, A.~Lam, K.~A.
  Zanetti, Q.-H. Ye, L.-X. Qin, C.~M. Croce, Z.-Y. Tang, and X.~W. Wang.
\newblock {Identification of Metastasis-Related MicroRNAs in Hepatocellular
  Carcinoma}.
\newblock {\em Hepatology}, 47(3):897--907, 2008.

\bibitem{ji09microrna}
J.~Ji, J.~Shi, A.~Budhu, Z.~Yu, M.~Forgues, S.~Roessler, S.~Ambs, Y.~Chen, P.S.
  Meltzer, C.M. Croce, L.-X. Qin, K.~Man, C.-M. Lo, J.~Lee, I.O.L. Ng, J.~Fan,
  Z.-Y. Tang, H.-C. Sun, and X.W. Wang.
\newblock {MicroRNA Expression, Survival, and Response to Interferon in Liver
  Cancer}.
\newblock {\em New England Journal of Medicine}, 361:1437--1447, 2009.

\bibitem{law11emerging}
P.~T.-Y. Law and N.~Wong.
\newblock {Emerging roles of microRNA in the intracellular signaling networks
  of hepatocellular carcinoma}.
\newblock {\em Journal of Gastroenterology and Hepatology}, 26(3):437--449,
  2011.

\bibitem{gu12gene}
Z.~Gu, C.~Zhang, and J.~Wang.
\newblock Gene regulation is governed by a core network in hepatocellular
  carcinoma.
\newblock {\em BMC Systems Biology}, 6(1):32, 2012.

\bibitem{volinia10reprogramming}
S.~Volinia, M.~Galasso, S.~Costinean, L.~Tagliavini, G.~Gamberoni, A.~Drusco,
  J.~Marchesini, N.~Mascellani, M.E. Sana, R.~Abu~Jarour, C.~Desponts,
  M.~Teitell, R.~Baffa, R.~Aqeilan, M.V. Iorio, C.~Taccioli, R.~Garzon,
  G.~Di~Leva, M.~Fabbri, M.~Catozzi, M.~Previati, S.~Ambs, T.~Palumbo,
  M.~Garofalo, A.~Veronese, A.~Bottoni, P.~Gasparini, C.C. Harris, R.~Visone,
  Y.~Pekarsky, A.~de~la Chapelle, M.~Bloomston, M.~Dillhoff, L.Z. Rassenti,
  T.J. Kipps, K.~Huebner, F.~Pichiorri, D.~Lenze, S.~Cairo, M.-A. Buendia,
  P.~Pineau, A.~Dejean, N.~Zanesi, S.~Rossi, G.A. Calin, C.-G. Liu,
  J.~Palatini, M.~Negrini, A.~Vecchione, A.~Rosenberg, and C.M. Croce.
\newblock {Reprogramming of miRNA networks in cancer and leukemia}.
\newblock {\em Genome Research}, 20(5):589--599, 2010.

\bibitem{bandyopadhyay10development}
S.~Bandyopadhyay, R.~Mitra, U.~Maulik, and M.Q. Zhang.
\newblock {Development of the human cancer microRNA network}.
\newblock {\em Silence}, 1:6, 2010.

\bibitem{troyanskaya01missing}
O.G. Troyanskaya, M.~Cantor, G.~Sherlock, P.O. Brown, T.~Hastie, R.~Tibshirani,
  D.~Botstein, and R.B. Altman.
\newblock {Missing value estimation methods for DNA microarrays}.
\newblock {\em Bioinformatics}, 17(6):520--525, 2001.

\bibitem{fruchterman91graph}
T.M.J. Fruchterman and E.M. Reingold.
\newblock {Graph Drawing by Force-directed Placement}.
\newblock {\em Software - Practice and Experience}, 21(11):1129--1164, 1991.

\bibitem{jurman12stability}
G.~Jurman, M.~Filosi, R.~Visintainer, S.~Riccadonna, and C.~Furlanello.
\newblock {Stability Indicators in Network Reconstruction}.
\newblock arXiv:1209.1654v1 [q-bio.MN], 2012.

\bibitem{filosi13stability}
M.~Filosi, R.~Visintainer, S.~Riccadonna, G.~Jurman, and C.~Furlanello.
\newblock {Stability Indicators in Network Reconstruction}.
\newblock submitted, 2013.

\bibitem{gleditsch02expanded}
K.S. Gleditsch.
\newblock {Expanded Trade and GDP Data}.
\newblock {\em Journal of Conflict Resolution}, 46:712--724, 2002.

\bibitem{fronczak12statistical}
A.~Fronczak and P.~Fronczak.
\newblock {Statistical mechanics of the international trade network}.
\newblock {\em Physical Review E}, 85:056113, 2012.

\bibitem{gower66some}
J.C. Gower.
\newblock {Some distance properties of latent root and vector methods used in
  multivariate analysis}.
\newblock {\em Biometrika}, 53:325--–328, 1966.

\bibitem{mardia78some}
K.V. Mardia.
\newblock {Some properties of classical Multidimensional Scaling}.
\newblock {\em Communications on Statistics – Theory and Methods},
  A7:1233--–1241, 1978.

\bibitem{cailliez83analytical}
F.~Cailliez.
\newblock {The analytical solution of the additive constant problem}.
\newblock {\em Psychometrika}, 48:343--–349, 1983.

\bibitem{fronczak06fluctuation}
A.~Fronczak, P.~Fronczak, and J.A. Ho{\l}yst.
\newblock {Fluctuation-dissipation relations in complex networks}.
\newblock {\em Physical Review E}, 73:016108, 2006.

\bibitem{crafts06world}
N.~Crafts.
\newblock {The World Economy In The 1990s: A Long Run Perspective}.
\newblock In P.W. Rhode and G.~Toniolo, editors, {\em {The Global Economy in
  the 1990s: A Long-Run Perspective}}. Cambridgre Academic Press, 2006.

\bibitem{vangerven09attention}
{van Gerven M. and Jensen O.}
\newblock {Attention modulations of posterior alpha as a control signal for
  two-dimensional brain-computer interfaces}.
\newblock {\em Journal of Neuroscience Methods}, 179(1):78--84, 2009.

\bibitem{feizi13network}
S.~Feizi, D.~Marbach, M.~M{\'e}dard, and M.~Kellis.
\newblock Network deconvolution as a general method to distinguish direct
  dependencies in networks.
\newblock {\em Nature Biotechnology}, 31:726--733, 2013.

\bibitem{zou05regularization}
H.~Zou and T.~Hastie.
\newblock {Regularization and Variable Selection Via the Elastic Net}.
\newblock {\em Journal of the Royal Statistical Society: Series B},
  67(2):301--320, 2005.

\bibitem{tibshirani96regression}
R.~Tibshirani.
\newblock {Regression Shrinkage and Selection Via the Lasso}.
\newblock {\em Journal of the Royal Statistical Society. Series B},
  58(1):267--288, 1996.

\bibitem{demol09elastic}
C.~De~Mol, E.~De~Vito, and L.~Rosasco.
\newblock {Elastic-Net regularization in learning theory}.
\newblock {\em Journal of Complexity}, 25(2):201--230, 2009.

\bibitem{mosci10solving}
S.~Mosci, L.~Rosasco, M.~Santoro, A.~Verri, and S.~Villa.
\newblock Solving structured sparsity regularization with proximal methods.
\newblock In {\em Machine Learning and Knowledge Discovery in Databases}, pages
  418--433. Springer, 2010.

\bibitem{albanese12mlpy}
D.~Albanese, R.~Visintainer, S.~Merler, S.~Riccadonna, G.~Jurman, and
  C.~Furlanello.
\newblock {\textit{mlpy}: Machine Learning Python}.
\newblock arXiv:1202.6548 [cs.MS], 2012.

\bibitem{breiman01random}
L.~Breiman.
\newblock {Random Forests}.
\newblock {\em Machine Learning}, 45(1):5--32, 2001.

\bibitem{kia13discrete}
S.M. Kia, E.~Olivetti, and P.~Avesani.
\newblock {Discrete Cosine Transform for MEG Signal Decoding}.
\newblock In {\em Proceedings of the International Workshop on Pattern
  Recognition in Neuroimaging PRNI}, pages 132--135. IEEE, 2013.

\bibitem{bahramisharif10covert}
A.~Bahramisharif, M.~Van~Gerven, T.~Heskes, and O.~Jensen.
\newblock {Covert attention allows for continuous control of brain--computer
  interfaces}.
\newblock {\em European Journal of Neuroscience}, 31(8):1501--1508, 2010.

\bibitem{signoretto12classification}
M.~Signoretto, E.~Olivetti, L.~De~Lathauwer, and J.A~K Suykens.
\newblock {Classification of Multichannel Signals With Cumulant-Based Kernels}.
\newblock {\em IEEE Transactions on Signal Processing}, 60(5):2304--2314, 2012.

\bibitem{furlanello13sparse}
T.~Furlanello, M.~Cristoforetti, C.~Furlanello, and G.~Jurman.
\newblock {Sparse Predictive Structure of Deconvolved Functional Brain
  Networks}.
\newblock arXiv:1310.6547 [q-bio.NC], 2013.

\end{thebibliography}
\appendix
\section{Uniqueness of $\pmb{\overline{\gamma}}$}
\label{sec:appendix}
Fix the number $N$ of nodes, and consider the two extremal networks $\mathcal{E}_N$ and $\mathcal{F}_N$, whose Laplacian spectrum is respectively
\begin{displaymath}
\textrm{spec}(L(\mathcal{E}_N)) = ( \underbrace{0,\cdots,0}_N)
\quad\textrm{and}\quad
\textrm{spec}(L(\mathcal{F}_N)) = ( 0,\underbrace{N,\cdots,N}_{N-1})\ ,
\end{displaymath}
so that $\omega_i=0$ for the empty network and $\omega_i=\sqrt{N}$ for the fully connected network, for $i=1,\ldots,N-1$.

The Lorentz distribution for the empty network is thus
\begin{displaymath}
\begin{split}
\rho_{\mathcal{E}_N}(\omega,\gamma) &= K\sum_{i=1}^{N-1} \frac{\gamma}{\gamma^2+(\omega-\omega_i)^2} \\
&= \frac{K\gamma (N-1)}{\gamma^2+\omega^2}\ ,
\end{split}
\end{displaymath}
where $K$ can be computed as
\begin{displaymath}
\begin{split}
K &= \frac{1}{\displaystyle{\int_0^{+\infty} \frac{\gamma (N-1)}{\gamma^2+\omega^2} \textrm{d}\omega}} \\
&= \frac{1}{(N-1)\left[ \arctan\left(\frac{\omega}{\gamma}\right)\right]_0^{+\infty}} \\
&= \frac{1}{\displaystyle{\frac{\pi}{2}}(N-1)} \\
&= \frac{2}{(N-1)\pi}\ ,
\end{split}
\end{displaymath}
so that
\begin{displaymath}
\label{eq:rhoE}
\begin{split}
\rho_{\mathcal{E}_N}(\omega,\gamma) &= \frac{K\gamma (N-1)}{\gamma^2+\omega^2} \\
&= \frac{2\gamma}{\pi(\gamma^2+\omega^2)}\ .
\end{split}
\end{displaymath}

For the fully connected network we have
\begin{displaymath}
\begin{split}
\rho_{\mathcal{F}_N}(\omega,\gamma) &= K\sum_{i=1}^{N-1} \frac{\gamma}{\gamma^2+(\omega-\omega_i)^2} \\
&= K\sum_{i=1}^{N-1} \frac{\gamma}{\gamma^2+(\omega-\sqrt{N})^2} \\
&= \frac{\gamma K (N-1) }{\gamma^2+\omega^2+N-2\omega\sqrt{N}}\ ,
\end{split}
\end{displaymath}
where $K$ is
\begin{displaymath}
\begin{split}
K &= \frac{1}{\gamma (N-1) \displaystyle{\int_0^{+\infty} \frac{\textrm{d}\omega}{\gamma^2+\omega^2+N-2\omega\sqrt{N}}}} \\
&= \frac{1}{\frac{\gamma (N-1)}{\gamma} \left[ \arctan\left(\frac{\omega-\sqrt{N}}{\gamma}\right)\right]_0^{+\infty}} \\
&= \frac{1}{(N-1)\left( \frac{\pi}{2} + \arctan\left(\frac{\sqrt{N}}{\gamma}\right)\right)}\ ,
\end{split}
\end{displaymath}
so that
\begin{displaymath}
\begin{split}
\rho_{\mathcal{F}_N}(\omega,\gamma) &= \frac{\gamma K (N-1) }{\gamma^2+\omega^2+N-2\omega\sqrt{N}}\\
&= \frac{1}{(N-1)\left( \frac{\pi}{2} + \arctan\left(\frac{\sqrt{N}}{\gamma}\right)\right)} \cdot \frac{\gamma (N-1) }{\gamma^2+\omega^2+N-2\omega\sqrt{N}}\\
&= \frac{\gamma}{ \left( \frac{\pi}{2} + \arctan\left(\frac{\sqrt{N}}{\gamma}\right)\right) \left(\gamma^2+\omega^2+N-2\omega\sqrt{N}\right)} \ .
\end{split}
\end{displaymath}

Thus, we expand Eq.~\ref{eq:gamma_implicit} as follows:
\begin{equation}
\label{eq:gamma_explicit1}
\begin{split}
1 &= \epsilon_\gamma(\mathcal{E}_N, \mathcal{F}_N) \\
&= \sqrt{\int_0^\infty \left(\rho_{\mathcal{E}_N}(\omega,\gamma)-\rho_{ \mathcal{F}_N }(\omega,\gamma)\right)^2 \textrm{d}\omega}\\
&= \sqrt{\displaystyle{\int_0^\infty \left(  
\frac{2\gamma}{\pi(\gamma^2+\omega^2)} -
\frac{\gamma}{ \left( \frac{\pi}{2} + \arctan\left(\frac{\sqrt{N}}{\gamma}\right)\right) \left(\gamma^2+\omega^2+N-2\omega\sqrt{N}\right)}
  \right)^2 \textrm{d}\omega}
} 
\\
&= \sqrt{
\displaystyle{ \int_0^\infty  A^2  \textrm{d}\omega} +
\displaystyle{ \int_0^\infty B^2  \textrm{d}\omega} -
2\displaystyle{ \int_0^\infty AB  \textrm{d}\omega} 
}\ ,
\end{split}
\end{equation}
where
\begin{displaymath}
\begin{split}
A & = \frac{2\gamma}{\pi(\gamma^2+\omega^2)} \\
B &= \frac{\gamma}{ \left( \frac{\pi}{2} + \arctan\left(\frac{\sqrt{N}}{\gamma}\right)\right) \left(\gamma^2+\omega^2+N-2\omega\sqrt{N}\right)}\ .
\end{split}
\end{displaymath}
The three terms in Eq.~\ref{eq:gamma_explicit1} can be expanded as follows:
\begin{equation}
\label{eq:A}
\begin{split}
\displaystyle{\int_0^{+\infty}} A^2  \textrm{d}\omega &= \displaystyle{\int_0^{+\infty}} \left(\frac{2\gamma}{\pi(\gamma^2+\omega^2)}\right)^2  \textrm{d}\omega  \\
&= \frac{4\gamma^2}{\pi^2} \displaystyle{\int_0^{+\infty}}  \frac{\textrm{d}\omega}{(\gamma^2+\omega^2)^2} \\
&= \frac{4\gamma^2}{\pi^2} \frac{1}{2\gamma^3} \left[ \frac{\gamma\omega}{\gamma^2+\omega^2} + \arctan\left(\frac{\omega}{\gamma} \right) \right]_0^{+\infty} \\
&= \frac{2}{\gamma\pi^2}\left[\frac{\pi}{2}\right] \\
&= \frac{1}{\pi\gamma}\ ;
\end{split}
\end{equation}

\begin{equation}
\label{eq:B}
\begin{split}
\displaystyle{\int_0^{+\infty}} B^2  \textrm{d}\omega &= \displaystyle{\int_0^{+\infty}} \left( \frac{\gamma}{ \left( \frac{\pi}{2} + \arctan\left(\frac{\sqrt{N}}{\gamma}\right)\right) \left(\gamma^2+\omega^2+N-2\omega\sqrt{N}\right)} \right)^2  \textrm{d}\omega  \\
&= \displaystyle{\int_0^{+\infty}} \frac{\gamma^2}{ \left( \frac{\pi}{2} + \arctan\left(\frac{\sqrt{N}}{\gamma}\right)\right)^2 \left(\gamma^2+\omega^2+N-2\omega\sqrt{N}\right)^2}   \textrm{d}\omega  \\
&= \frac{\gamma^2}{\left( \frac{\pi}{2} + \arctan\left(\frac{\sqrt{N}}{\gamma}\right)\right)^2} \displaystyle{\int_0^{+\infty} \frac{\textrm{d}\omega}{\left(\gamma^2+\omega^2+N-2\omega\sqrt{N}\right)^2 } } \\
&= \frac{\gamma^2}{2\gamma^3\left( \frac{\pi}{2} + \arctan\left(\frac{\sqrt{N}}{\gamma}\right)\right)^2} \left[ \frac{\gamma(\omega-\sqrt{N})}{\gamma^2+(\omega-\sqrt{N})^2} + \arctan\left( \frac{\omega-\sqrt{N}}{\gamma} \right)  \right]_0^{+\infty} \\
&= \frac{1}{2\gamma\left( \frac{\pi}{2} + \arctan\left(\frac{\sqrt{N}}{\gamma}\right)\right)^2} \left( \frac{\pi}{2} + \frac{\gamma\sqrt{N}}{\gamma^2+N} + \arctan\left( \frac{\sqrt{N}}{\gamma} \right)  \right) \ ;
\end{split}
\end{equation}

\begin{equation}
\label{eq:AB}
\begin{split}
-2\displaystyle{\int_0^{+\infty}} AB  \textrm{d}\omega &= -2 \displaystyle{\int_0^{+\infty}} \frac{2\gamma}{\pi(\gamma^2+\omega^2)} \frac{\gamma}{ \left( \frac{\pi}{2} + \arctan\left(\frac{\sqrt{N}}{\gamma}\right)\right) \left(\gamma^2+\omega^2+N-2\omega\sqrt{N}\right)}  \textrm{d}\omega  \\
&= \frac{-2\cdot\gamma\cdot 2\gamma}{\pi\left( \frac{\pi}{2} + \arctan\left(\frac{\sqrt{N}}{\gamma}\right)\right)  } \displaystyle{\int_0^{+\infty}} \frac{\textrm{d}\omega}{ (\gamma^2+\omega^2) \left(\gamma^2+\omega^2+N-2\omega\sqrt{N}\right)  } \\
&= \frac{-4\gamma}{ \left( \frac{\pi}{2} + \arctan\left(\frac{\sqrt{N}}{\gamma}\right)\right) \pi ( 4\gamma^2+N) } \left[ 
\frac{\gamma}{\sqrt{N}} \log\frac{\gamma^2+\omega^2}{\gamma^2+\omega^2+N-2\omega\sqrt{N}} + \right. \\
&\phantom{=} \left. \arctan\left( \frac{\omega-\sqrt{N}}{\gamma}\right) + \arctan\left( \frac{\omega}{\gamma} \right)
\right]_0^{+\infty}\\
&= \frac{-4\gamma}{ \left( \frac{\pi}{2} + \arctan\left(\frac{\sqrt{N}}{\gamma}\right)\right) \pi ( 4\gamma^2+N) } \left[ 
\frac{\pi}{2} +\frac{\pi}{2} - \frac{\gamma}{\sqrt{N}} 
\log\frac{\gamma^2}{\gamma^2+N} + \right. \\
&\phantom{=} \left. \arctan\left( \frac{\sqrt{N}}{\gamma}\right) 
\right]\ .
\end{split}
\end{equation}

Plugging Eqs.~\ref{eq:A},\ref{eq:B},\ref{eq:AB} into Eq.~\ref{eq:gamma_explicit1}, we obtain:
\begin{displaymath}
\begin{split}
1 &= \epsilon_\gamma(\mathcal{E}_N, \mathcal{F}_N) \\
&= \frac{1}{\pi\gamma} + \frac{1}{2\gamma\left( \frac{\pi}{2} + \arctan\left(\frac{\sqrt{N}}{\gamma}\right)\right)^2} \left( \frac{\pi}{2} + \frac{\gamma\sqrt{N}}{\gamma^2+N} + \arctan\left( \frac{\sqrt{N}}{\gamma} \right)  \right) - \\ 
&\phantom{=} \frac{-4\gamma}{ \left( \frac{\pi}{2} + \arctan\left(\frac{\sqrt{N}}{\gamma}\right)\right) \pi ( 4\gamma^2+N) } \left[
\pi  - \frac{\gamma}{\sqrt{N}} 
\log\frac{\gamma^2}{\gamma^2+N} +  \arctan\left( \frac{\sqrt{N}}{\gamma}\right) \right] \ .
\end{split}
\end{displaymath}

\begin{figure}[!t]
\begin{center}
\begin{tabular}{cc}
\includegraphics[width=0.4\textwidth]{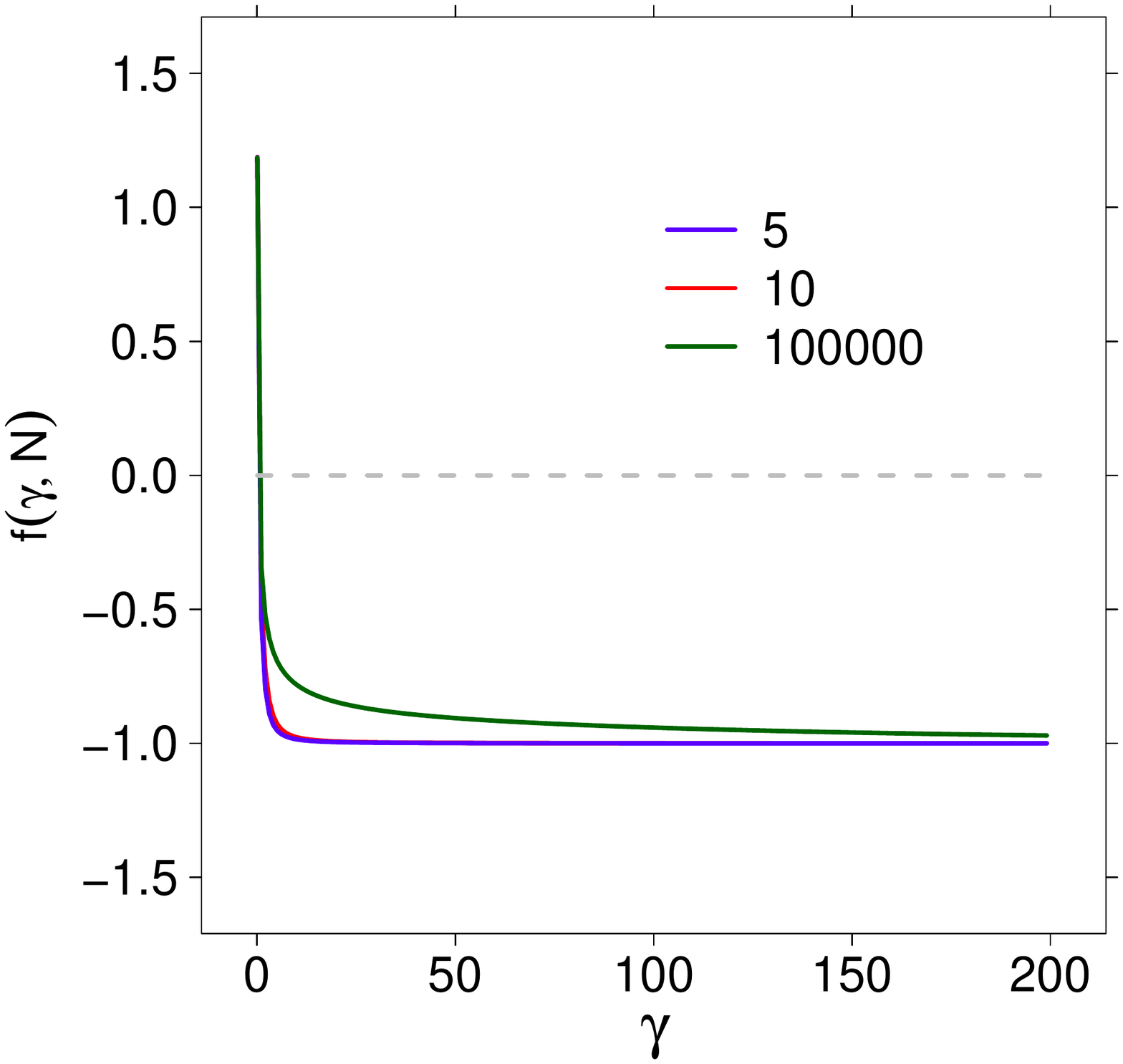}&
\includegraphics[width=0.4\textwidth]{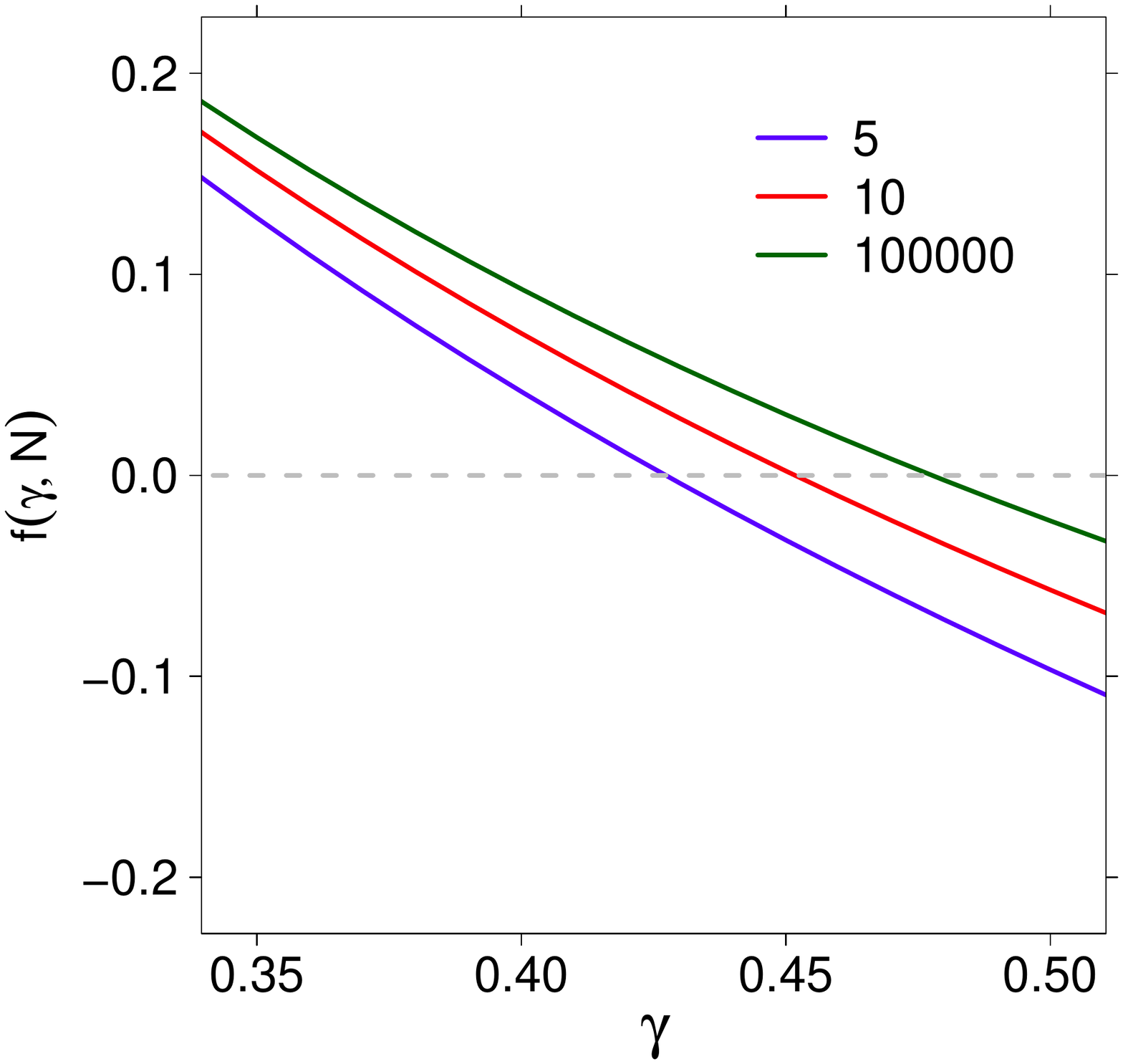}\\
(a) & (b) \\
\end{tabular}
\end{center}
\caption{(a) Behaviour of $f(\gamma,N)$ for $N$=5, 10 and 10000, in the interval $\gamma\in (0,200]$ and (b) zoomed in the interval $\gamma\in [0.35,0.5]$ with the solutions of Eq.~\ref{eq:gamma_implicit}.}
\label{fig:gammaN}
\end{figure}

Consider now the function $f(N,\gamma)=\epsilon_\gamma(\mathcal{E}_N, \mathcal{F}_N) -1$: for a fixed value of $N$, it is a monotonically decreasing function of $\gamma$, so the equation Eq.~\ref{eq:gamma_implicit} has an unique solution $\overline{\gamma}$. 
In Fig.~\ref{fig:gammaN} (a) and (b) we display the situation for $N$= 5,10 and 100000. 

\section{Uniqueness of $\pmb{\overline{\gamma}^\uparrow}$}
\label{sec:appendixb}
The spectra of the laplacian matrices of the two extremal graphs $\hat{\mathcal{E}}_N^\uparrow$ and $\hat{\mathcal{F}}_N^\uparrow$ are now
\begin{displaymath}
\textrm{spec}(L(\hat{\mathcal{E}}_N^\uparrow)) = ( \underbrace{0,\cdots,0}_{2N})
\quad\textrm{and}\quad
\textrm{spec}(L(\hat{\mathcal{F}}_N^\uparrow)) = ( 0,\underbrace{N-2,\cdots,N-2}_{N-1},\underbrace{N,\cdots,N}_{N-1},2N-2)\ .
\end{displaymath}

It follows that
\begin{displaymath}
K_{\hat{\mathcal{E}}_N^\uparrow} = \frac{2}{(2N-1)\pi}
\end{displaymath}
and
\begin{displaymath}
K_{\hat{\mathcal{F}}_N^\uparrow} = \frac{1}{ (2N-1)\frac{\pi}{2} + (N-1)\left(\arctan\frac{\sqrt{N-2}}{\gamma} + \arctan\frac{\sqrt{N}}{\gamma}\right) +\arctan\frac{\sqrt{2N-2}}{\gamma}} \ .
\end{displaymath}

Thus the equation 
\begin{displaymath}
\epsilon_{\gamma}(\hat{\mathcal{E}}^\uparrow,\hat{\mathcal{F}}^\uparrow) = 1
\end{displaymath}
(whose solution is the normalizing factor $\overline{\gamma}^\uparrow$) reads as follows:
\begin{equation}
1 = \sqrt{
\int_0^{+\infty}
\left[
\frac{2\gamma}{\gamma^2+\omega^2} -
\frac{
\gamma 
\left(
\frac{N-1}{\gamma^2+(\omega-\sqrt{N-2})^2} + 
\frac{N-1}{\gamma^2+(\omega-\sqrt{N})^2} + 
\frac{1}{\gamma^2+(\omega-\sqrt{2N-2})^2} 
\right)
}
{
(2N-1)\frac{\pi}{2} + (N-1)\left(\arctan\frac{\sqrt{N-2}}{\gamma} + \arctan\frac{\sqrt{N}}{\gamma}\right) +\arctan\frac{\sqrt{2N-2}}{\gamma}
}
\right]^2
\textrm{d}\omega
}\ .
\label{eq:gammahat}
\end{equation}

Introduce now a few shorthands: define, for $T,U\in\mathbb{R}$, the following integral
\begin{displaymath}
\int_0^{+\infty}{ \frac{\textrm{d}\omega}{(\gamma^2+(\omega-\sqrt{T})^2)(\gamma^2+(\omega-\sqrt{U})^2)}} = 
\begin{cases}
M(T) & \text{if $T=U$,}
\\
L(T,U) & \text{if $T\not= U$}\ .
\end{cases}
\end{displaymath}
Then,
\begin{displaymath}
M(T) = \frac{
\frac{1}{2}\left(
\gamma^2\arctan\frac{\sqrt{T}}{\gamma} + T \arctan\frac{\sqrt{T}}{\gamma} + \gamma\sqrt{T}
\right)
}{
\gamma^5+T\gamma^3
}
+\frac{\pi}{4\gamma^3},
\end{displaymath}
and
\begin{displaymath}
L(T,U) = 
\frac{
-\log\left(\gamma^{2} + U\right)+ \log\left(\gamma^{2} + T\right)
}{{\left(4 \, \gamma^{2} + T + 3 \, U\right)} \sqrt{T} - {\left(4 \, \gamma^{2} + 3 \, T + U\right)} \sqrt{U}} +  
\frac{\pi+ \arctan\left(\frac{\sqrt{T}}{\gamma}\right) + \arctan\left(\frac{\sqrt{U}}{\gamma}\right)}
{4 \, \gamma^{3} + T \gamma - 2 \, \sqrt{T} \sqrt{U} \gamma + U \gamma} \ .
\end{displaymath}
To shorten notations, define furthermore
\begin{displaymath}
Z = \frac{2\gamma}{\pi}\quad W=\gamma (N-1) K_{\hat{\mathcal{F}}_N^\uparrow}\quad W' = \frac{W}{N-1}\ . 
\end{displaymath}

With the aforementioned positions, Eq.~\ref{eq:gammahat} becomes
\begin{equation}
\begin{split}
1 &= Z^2 M(0) + W^2 M(N-2) + W^2 M(N) + W'^2 M(2N-2) \\
&\phantom{=}- 2ZW L(0,N-2) - 2ZWL(0,N) -2 ZW' L(0,2N-2) \\
&\phantom{=}+ 2 W^2 L(N-2,N) +2 WW' L(N-2,2N-2) +2WW' L(N,2N-2)\ .
\end{split}
\label{eq:gammahatexp}
\end{equation}

\begin{figure}[!t]
\begin{center}
\begin{tabular}{cc}
\begin{tabular}{rrr}
\multicolumn{1}{c}{N} & \multicolumn{1}{c}{$\overline{\gamma}$} & \multicolumn{1}{c}{$\overline{\gamma}^\uparrow$} \\
\hline
5&  0.4272836 & 0.3866861 \\
10&  0.4517012&  0.4300291\\
50&  0.4752742 & 0.4704579\\
100& 0.4777976  & 0.4753463\\
500 &  0.4787492 &  0.4782538\\
1000 &   0.4785596&    0.4783119\\
10000& 0.4779060& 0.4778813\\
\hline
\end{tabular}&
\raisebox{-2cm}{\includegraphics[width=0.6\textwidth]{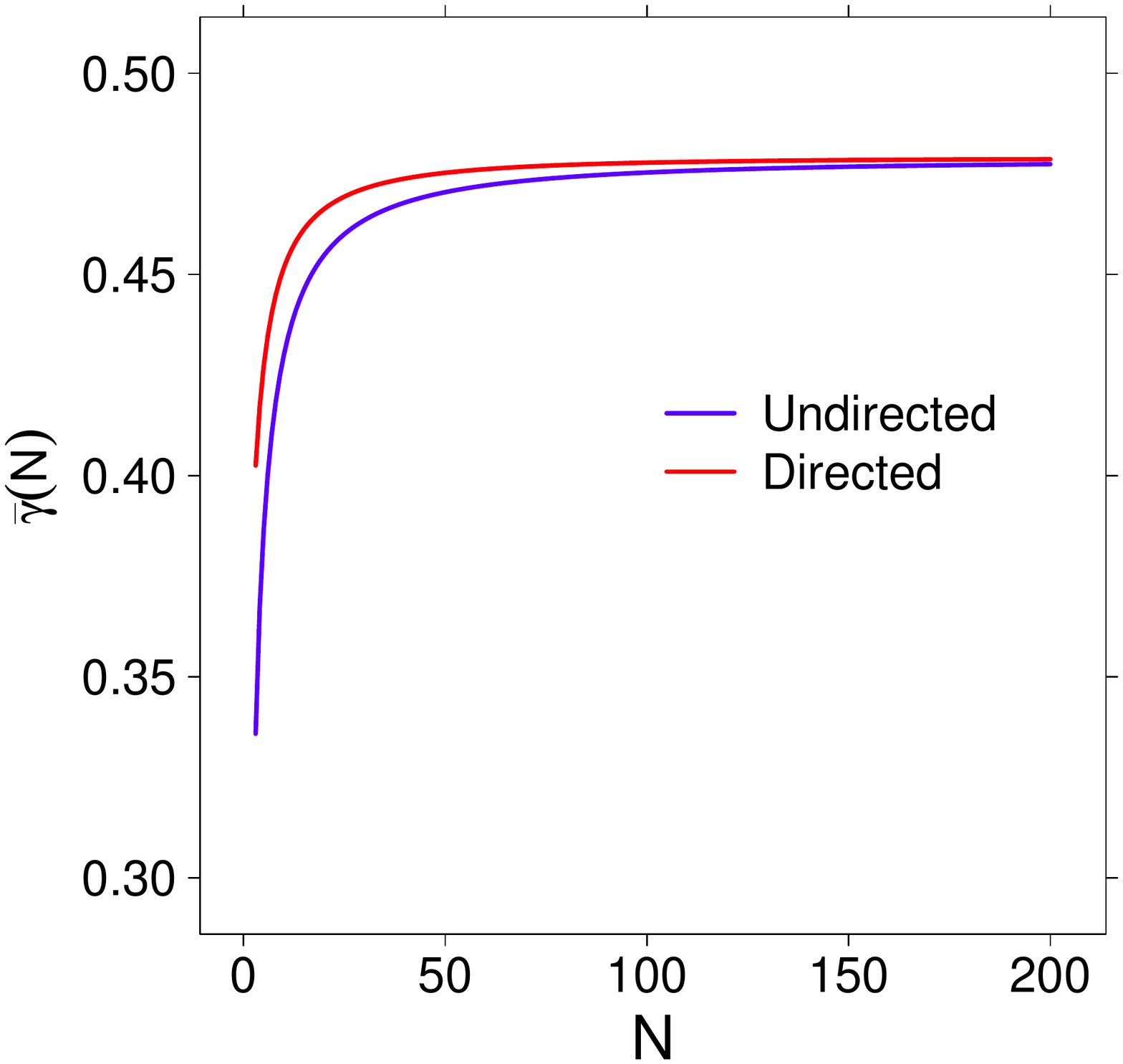}}
\end{tabular}
\end{center}
\caption{Comparison of $\overline{\gamma}$ and $\overline{\gamma}^\uparrow$ for different number of nodes $N$.}
\label{tab:gammas}
\end{figure}

As in the undirected case, for each $N$ Eq.~\ref{eq:gammahatexp} has an unique solution $\overline{\gamma}^\uparrow$, whose value is quite close to $\overline{\gamma}$, as shown in Fig.~\ref{tab:gammas}.

\end{document}